\documentclass{gIPE2e_basti}
\usepackage{color}
\pdfoutput=1

\usepackage{hyperref}
\def\disp{\displaystyle}
\def\DS{\mbox{\boldmath{$\delta\sigma$}}}
\def\DV{\mbox{\boldmath{$\delta V$}}}
\def\eref#1{(\ref{#1})}
\def\f{\frac}
\def\mE{{\mathcal E}}
\def\n{\mathbf{n}}
\def\na{\nabla}
\def\Om{\Omega}

\def\p{\partial}
\def\PHI{\mbox{\boldmath{$\Phi$}}}
\def\q{\quad}
\def\r{\mathbf{r}}
\def\S{\mbox{\boldmath{$S$}}}
\def\v{{\mathbf v}}

\begin{document}
\doi{}
\issn{}
\issnp{}
\jvol{} \jnum{} \jyear{} \jmonth{}

\markboth{}{}


\title{Regularizing a linearized EIT reconstruction method using a sensitivity based factorization method}

\author{
Moon Kyung Choi$\dagger$\footnote{$\dagger$ Computational Science and Engineering, Yonsei University, Korea}, Bastian Harrach$\ddagger$ \footnote{$\ddagger$ Department of Mathematics, University of Stuttgart, Germany}
and Jin Keun Seo$\dagger$
%
}

\maketitle

\vspace*{-6cm}
\raisebox{0cm}{\fbox{\begin{minipage}{13cm}\centering
This is an Author's Accepted Manuscript of an article published in\\
\emph{Inverse Probl.~Sci.~Eng.} \textbf{22}(7), 1029--1044, 2014,\\
published online: 01 Nov 2013, \textcopyright Taylor \& Francis, available online at: \url{http://www.tandfonline.com/doi/abs/10.1080/17415977.2013.850682}
\end{minipage}}}
\vspace*{4cm}

\begin{abstract}
For electrical impedance tomography (EIT), most practical reconstruction methods are based on linearizing the
underlying non-linear inverse problem. Recently, it has been shown that the linearized problem still contains
the exact shape information. However, the stable reconstruction of shape information from measurements of finite accuracy
on a limited number of electrodes remains a challenge.

In this work we propose to regularize the standard linearized reconstruction method (LM) for EIT
using a non-iterative shape reconstruction method (the factorization method). Our main tool is
a discrete sensitivity-based variant of the factorization method (herein called S-FM) which allows us to
formulate and combine both methods in terms of the sensitivity matrix.
We give a heuristic motivation for this new method and show numerical examples
that indicate its good performance in the localization of anomalies and the alleviation of ringing artifacts.
\end{abstract}

\begin{classcode}
35R30; 65F22
\end{classcode}

\section{Introduction}

In electrical impedance tomography (EIT), we try to reconstruct the
electrical conductivity distribution of an imaging subject from boundary measurements that are collected by placing surface electrodes around the subject, injecting different linearly independent currents using chosen pairs of electrodes, and measuring the required voltages \cite{Webster1978,Cheney1999,Holder2004}.
EIT seems to be a unique technique capable of low-cost measuring continuous real-time tomographic impedance of human body non-invasively. EIT has been used for regional lung function monitoring \cite{Barber1984,Wolf2007}.
Although numerous EIT reconstruction algorithms  have been developed last three decades,
it has not yet reached a satisfactory level of clinical competency 
\cite{Brown1985,Brown1987,Yorkey1987,Wexler1985,Gisser1988,Newell1988,Cheney1990,Nachman1996,Novikov1988,Borcea2003,Adler2006,Adler2009,Novikov2009,Bikowski2011}.

In EIT, the boundary data-set describing the boundary current-potential relation along the attached electrodes is determined mainly by the boundary geometry, electrode positions and the conductivity distribution. The nonlinear and ill-posed relation between the data-set and the conductivity distribution are entangled with the modeling errors caused by inaccurate boundary geometry, uncertainty in electrode positions, measurement noises, and unknown contact impedance \cite{Holder2004}. The ill-posedness and the forward modeling errors can be alleviated by using
the difference between a reference data-set and a received data-set to produce a {\it difference image} of the conductivity distribution.
The most commonly used algorithms in practical EIT systems are linearized reconstruction methods (LM) \cite{Isaacson1991,Cheney1999,Holder2004,Adler2009}. They are based on the linearized approximation of the relation between the conductivity difference $\DS$ and the difference data-set
{\boldmath $\delta V$}
\begin{equation}\label{eq:Ssigma_eq_V}
  \S \;\DS \q\approx\q {\mbox{\boldmath $\delta V$}},
\end{equation}
where $\S$ is the so-called sensitivity matrix (see section \ref{sect:method}).

The linearized relation (\ref{eq:Ssigma_eq_V}) is ill-conditioned. Small error in data produces large error in the reconstructed image.
A natural and frequently used regularization strategy is to only reconstruct those features of the image that are least affected by noise
which leads to the truncated singular value decomposition (tSVD), see section \ref{subsect:hybrid}. However, this approach is known to lead to strong ringing artifacts around boundaries in the reconstructed images.

Recently, two of authors showed rigorously that the linearized EIT equation contains full information about the shape and position of conductivity anomalies \cite{Harrach2010a}.
Hence, linearization should not be the cause for shape and position artifacts.
However, this theoretical result requires a continuous sampling on the boundary. In a realistic EIT setting, we have to deal with a limited number of electrodes and with measurement errors. Standard regularization methods as tSVD are not adapted to
shape recovery and do not stably recover shape and position information.

The aim of the present work is to enhance the conventional LM by regularizing it with
the factorization method that was specifically designed for shape recovery, cf.\ \cite{Kirsch1998,Bruhl2000,Bruhl2001,Hanke2003} for the origins of this method, and \cite{Kirsch2008,hanke2011sampling,harrach2013recent} for recent overviews.
Our main tool is a {\it sensitivity matrix based} variant of the factorization method (S-FM),
that uses the columns of the sensitivity matrix to form a pixel-wise index indicating the possibility of the presence or absence of an anomaly in this pixel. This allows us to formulate and combine the LM and the S-FM in the same
framework. We perform various numerical simulations to show that our new combination effectively alleviates ringing artifacts.

\section{Method}\label{sect:method}

Let $\Omega$ be a two- or three-dimensional domain with $C^\infty$-boundary $\p\Om$.
Assume that  $\sigma_0=1$ is a reference conductivity and $\sigma(\r)=\sigma_0(\r)+\delta\sigma(\r)\chi_D\in L^\infty_+(\Omega)$ 
is a conductivity distribution with a perturbation $\delta\sigma(\r)$ on the anomaly $D$,
where $\overline D\subseteq \Omega$ and $\r$ denotes the position in $\Omega$.

We begin with the standard linearized method (LM) in $n_E$-channel EIT system. In order to inject $n_E$ currents into
the imaging subject $\Omega$, we attach $n_E$ point electrodes $\mathcal E_j$, $j=1,\cdots,n_E$, on the boundary
$\partial \Omega$ and inject a current of $I$ mA between each pair of adjacent electrodes
$(\mE_1,\mE_2),(\mE_2,\mE_3),\cdots, (\mE_{n_E},\mE_1)$.
For simplicity, we assume $I=1$ and ignore the effects of the unknown contact impedances between electrodes with the boundary.
The electrical potential due to the $j$-th current is denoted by $u_j$ and satisfies the following Neumann boundary value problem:
\begin{alignat}{2}
\label{govern1} \nabla\cdot\left( \sigma(\r) \nabla u_j(\r) \right) &=0 && \q \mbox{ in }\Omega\\
\label{govern2} \sigma(\r)\f{\p}{\p \n} u_j(\r) & = \delta(\r-\mE_j)-\delta(\r-\mE_{j+1}) && \q \mbox{ on } \p\Om
\end{alignat}
where ${\bf n}$ is the outward unit normal vector on $\p\Om$. For the 
justification and well-posedness of the point electrode model we refer to
\cite{Hanke2011}.

We  measure a boundary voltage between an adjacent pair of electrodes, $\mE_k$ and $\mE_{k+1}$ for $k=1,\cdots,n_E$. The $k$-th boundary 
voltage difference subject to the $j$-th injection current is denoted as
\begin{equation}\label{EIT-f-data}
V_{j,k}[\sigma]= u_j(\mE_k)-u_j(\mE_{k+1}) \q\q\mbox{ for } j,k=1,\cdots,n_E,
\end{equation}
see our remark on the end of this subsection about evaluating the potential on current-driven electrodes.

Collecting $(n_E)^2$ number of  boundary voltage data, the measured data-set can be expressed as
\begin{equation}\label{EIT-data}
{\bf V}_1[\sigma]:=
\left(
   \begin{array}{c}
  V_{1,1}[\sigma]  \\
 V_{1,2}[\sigma]   \\
    \vdots  \\
V_{1,n_E}[\sigma]  \\
   \end{array}
 \right),\q\cdots,\q {\bf V}_{n_E}[\sigma]=:
  \left(
   \begin{array}{c}
  V_{n_E,1}[\sigma] \\
 V_{n_E,2}[\sigma]  \\
 \vdots  \\
 V_{n_E,n_E}[\sigma]  \\
   \end{array}
 \right).
\end{equation}
Assume that we also have a data-set ${\bf V}_j[\sigma_0] ~(j=1,\cdots,n_E)$ for the reference conductivity $\sigma_0$. Let $u_j^{0}$ be the solution of (\ref{govern1}), (\ref{govern2}) with $\sigma$ replaced by $\sigma_0$.

The inverse problem is to reconstruct $\delta\sigma:=\sigma- \sigma_{0}$ from the
difference $\DV_j:={\bf V}_j[\sigma]-{\bf V}_j[\sigma_{0}]$.  The difference data-set is a response of the anomaly $D$ such that
\begin{equation}\label{eq:uuV1}
 V_{j,k}[\sigma]- V_{j,k}[\sigma_0] = -\int_{D}(\sigma-\sigma_0) \nabla u_j \cdot \nabla
 u_k^0  ~d\r.
\end{equation}
In LM, we use the following approximation
\begin{equation}\label{linearization}
\int_{D}(\sigma-\sigma_0) \nabla u_j \cdot \nabla
 u_k^{(0)}  ~d\r
 \q\approx\q \int_{D}(\sigma-\sigma_0) \nabla u_j^{(0)} \cdot \nabla
 u_k^{(0)}  ~d\r
\end{equation}
Note that according to the theoretical result in \cite{Harrach2010a},
this linearization does not lead to position and shape errors if a continuous sampling
on the boundary is available. A EIT system with a limited number of electrodes will always produce shape errors; but the shape errors depend mainly on the configuration of electrodes and the linearization is minor cause of the shape and position errors.

Discretizing the domain $\Om$ into $n_p$ elements  as $\Om=\disp\cup_{n=1}^{n_p} q_n$ and assuming that $\sigma$ and $\sigma_0$ are constants on each element $q_n$, we obtain from \eref{eq:uuV1} and \eref{linearization}
\begin{equation}\label{dgTn}
\sum_{n=1}^{n_p} \delta\sigma_n\int_{q_n} \nabla u_j^{(0)} \cdot \nabla u_k^{(0)}~  d\r \approx
V_{j,k}[\sigma_0]-V_{j,k}[\sigma]=:\delta V_{j,k}
\end{equation}
where $\delta\sigma_n:=(\sigma-\sigma_0)|_{q_n}$ is the value 
of $\sigma-\sigma_0$ on the pixel $q_n$. (See the remark at the end of this
section concerning the pixels adjacent to a point electrode.)

\eref{dgTn} can be written as the following linear system
\begin{equation}\label{EIT-Nonlinear-inverse}
  \S \;\DS \q\approx\q {\mbox{\boldmath$\delta V$}}
\end{equation}
where
$$
\S=
\left(
  \begin{array}{cc}
   \S_1\\
   \vdots \\
   \S_{n_E} \\
  \end{array}
\right),\q\q\DS=
\left(
  \begin{array}{cc}
   \delta\sigma_1\\
   \vdots\\
\delta\sigma_{n_p}
  \end{array}
\right),\q\q   \DV=
\left(
  \begin{array}{cc}
   \DV_1 \\
   \vdots \\
  \DV_{n_E} \\
  \end{array}
\right).
$$
Here, $\S_j$ is
the $j$-th block of the linearized sensitivity matrix:
\begin{equation}\label{sensitivity}
{\S}_j= \left(\begin{array}{ccc}
\int_{q_1} \nabla u_j^{(0)} \cdot \nabla u_1^{(0)} ~d\r& \q \cdots\q&\int_{q_{n_p}} \nabla u_j^{(0)} \cdot \nabla u_1^{(0)}~d\r\\
 \vdots& \q\q \vdots\q\q& \vdots\\
\int_{q_1} \nabla u_j^{(0)} \cdot \nabla u_{n_E}^{(0)} ~d\r& \q \cdots\q&\int_{q_{n_p}} \nabla u_j^{(0)} \cdot \nabla
u_{n_E}^{(0)}~d\r\\
\end{array}\right).
\end{equation}
It is well-known that the sensitivity matrix $\S$ is ill-conditioned for realistic numbers of pixels and electrodes,
so that solving (\ref{EIT-Nonlinear-inverse}) requires regularization.

For the sake of mathematical rigorosity we have to note that, in the point electrode model, the electric potential possesses singularities
on current driven electrodes. Hence, the evaluations in (\ref{EIT-f-data})
and (\ref{EIT-data}) may not be well-defined. Moreover, $\nabla u_j^{(0)}$ is not
neccessarily square integrable on pixes adjacent to a point electrodes, so that
entries of $S_j$ belonging to such boundary pixels might also not be well-defined.

A rigorous mathematical way to overcome the first problem 
is to note that we only work with the difference data (\ref{eq:uuV1}).
Since $\overline D\subseteq \Omega$, point evaluations of difference data
are well-defined, even on current-driven electrodes, see \cite{Hanke2011}.
The second problem can be overcome by replacing the partition 
$\Omega=\bigcup_{n=1}^{n_p} q_n$ by a partition $\Omega'=\bigcup_{n=1}^{n_p} q_n$
for a slightly smaller open set $\Omega'$ with $\overline{\Omega'}\subset \Omega$
that is known to contain the inclusions. 

In the scope of this article, we consider both problems to be mathematical 
subtleties, as we work with numerical approximations to the point electrode 
model anyway. Hence, for the sake of readability, we stick to the somewhat sloppy
formulations above

\subsection{Sensitivity matrix based Factorization Method (S-FM)}\label{subsect:S_FM}

Though the linearized equation still contains the correct shape information, standard regularization methods for the linear reconstruction method (LM) are not designed for shape reconstruction and tend to produce ringing artifacts
around conductivity anomalies. To improve this, we will combine it with the factorization method (FM) which is a non-iterative anomaly detection method that provides a criterion for determining the presence or absence of an anomaly at each location in the imaging subject. Hence, we expect the FM to improve the anomaly detection performance of the standard LM.

Standard formulations of the FM utilize special dipole functions to characterize whether a given point belongs to an anomaly or not.
In order to use the FM for regularizing the LM, we develop a {\it sensitivity matrix based factorization method} (S-FM).
The main idea is that the dipole functions can be replaced by the column vectors of the sensitivity matrix $\S_j$.
We will now explain this idea in some detail:

Given a pixel $q_n$ with $\overline{q_n}\subset \Om$ and the $j$-th injection current, let $\phi^n_j$ solve
\begin{equation}\label{phi}
\left\{\begin{array}{ll}
 \na\cdot\na \phi^n_j = \na\cdot( \chi_{q_n} \na  u_j^{(0)})\q\q&\mbox{in }~\Om \\
 \\
\na \phi^n_j\cdot\n = 0 &\mbox{on }\p\Om
 \end{array}\right.
\end{equation}
where $\chi_{q_n}$ is the characteristic function having $1$ on $q_n$ and $0$ otherwise.
$\phi_n^j$ can be regarded as a pixel-based analogue to a dipole function $\phi_{z,d}$ which solves
$\na\cdot\na\phi_{z,d}=d\cdot \nabla\delta_z$.

The boundary voltages of the pixel-based dipole function $\phi_n^j$ agree with the column vectors of the sensitivity matrix $\S_j$.
\begin{lemma}\label{lemma1} Denote $\PHI^n_j:=(\Phi_{j}^n(1)~\cdots~\Phi_{j}^n(n_E))^T$ where \begin{equation}\label{Phi-data}
\Phi_{j}^n(k):= \phi_j^n(\mE_k) -\phi_j^n(\mE_{k+1}) \q\q\mbox{for } k=1,\cdots,n_E.
\end{equation}
Then $\Phi_{j}^n(k)=\int_{q_n} \nabla u_j^{(0)} \cdot \nabla u_k^{(0)} ~ d\r$ and hence
\begin{equation}\label{Phi3-0}
\PHI^n_j\q=\q\mbox{ $n$-th column of } {\S}_j
\end{equation}
\end{lemma}
{\it Proof.} From the Neumann boundary condition of  $u_k^{(0)}$ in \eref{govern1},\eref{govern2}, we have
\begin{equation}\label{Phi2}
\int_{\Om } \nabla \phi_j^n \cdot \nabla u_k^{(0)} ~ d\r = \phi_j^n(\mE_k) -\phi_j^n(\mE_{k+1})
\end{equation}
for $k=1,\cdots,n_E$. On the other hand,  it follows from the definition of $\phi^n_j$ in \eref{phi} that
\begin{equation}\label{Phi2-1}
\int_{\Om } \nabla \phi_j^n \cdot \nabla u_k^{(0)} ~ d\r = \;\int_{q_n} \nabla u_j^{(0)} \cdot \nabla u_k^{(0)} ~ d\r.
\end{equation}
Hence, \eref{Phi3-0} follows from the definition of ${\S}_j$ in \eref{sensitivity}.\hfill $\Box$

Now, we will derive the pixel-based anomaly detection criterion of S-FM.
First, we interpret the $n_E^2$ values in the data set as a matrix, i.e.,
we replace the $n_E^2$-dimensional vector $\DV$ by the $n_E\times n_E$-matrix
\begin{equation}\label{BbbV}
\delta\Bbb{V}:=(\DV_1\cdots\DV_{n_E}).
\end{equation}

To simplify the presentation let us assume in the following
that the conductivity difference $\DS$ satisfies the linear approximation $\S \DS=\mbox{\boldmath $\delta V$}$,
and that $\delta\Bbb{V}$ is invertible, cf.\ \cite{Harrach2010b} for a justification of the FM without these assumptions.
Moreover, we assume that the conductivity contrast is equal to one inside the inclusion, i.e.
$\delta \sigma=\chi_D$, where $D$ is compactly supported in $\Om$.

The following theorem is the basis for the S-FM:
\begin{theorem}\label{thm:S-FM}
\begin{enumerate}
\item[(a)] For all pixel $q_n$ and $j$-th injection current
\begin{equation}\label{key7-0}
\PHI^n_j\cdot (\delta\mathbb V)^{-1}\PHI^n_j \q\ge\q
\max_{{\bf h}\in \mathbb{R}^{{n_E}}}\f{ \left| \int_{q_n} \nabla u_j^{(0)} \cdot \nabla u_{\bf h}^{(0)} ~ d\r \right|^2}{
   \int_D |\nabla u^{(0)}_{\bf h}|^2 d\r}.
\end{equation}
where $u^{(0)}_{\bf h}(\r):=\sum_{k=1}^{n_E}\; h_k\; u_{k}^{(0)}(\r)$.
\item[(b)] If the pixel $q_n$ lies inside the inclusion $D$ then for all injection currents
  \begin{equation}\label{key4}
\left|\PHI^n_j\cdot (\delta\mathbb V)^{-1}\PHI^n_j\right| \q \le\q  \int_{q_n}\left|\nabla
 u_j^{(0)}\right|^2 d\r.
\end{equation}
\end{enumerate}
\end{theorem}
\emph{Proof.}
By lemma \ref{lemma1} and superposition we have for all ${\bf h}\in \mathbb{R}^{{n_E}}$
\begin{equation}\label{Phi3}
\;{\bf h}\cdot\PHI_{j}^n\; =\; \sum_{k=1}^{n_E}\; h_k\; \int_{q_n} \nabla u_j^{(0)} \cdot \nabla u_k^{(0)} ~ d\r
\; =\; \int_{q_n} \nabla u_j^{(0)} \cdot \nabla u_{\bf h}^{(0)} ~ d\r
\end{equation}

On the other hand, we can define
\begin{equation}\label{hjn}
{\bf h}_{j}^n=(h_{j1}^n,\cdots, h_{jn_E}^n)^T:=(\delta\mathbb V)^{-1}  \PHI_{j}^n.
\end{equation}
and use that
\begin{equation}\label{dV-j}
\mbox{ $k$-th column of } \delta\Bbb{V}\q=\q \DV_k \q=\q \S_k\DS
\end{equation}
to obtain
\begin{eqnarray}
\nonumber \;{\bf h}\cdot\PHI_{j}^n & = &
{\bf h}\cdot ({\delta\mathbb V}~{\bf h}_{j}^n) ={\bf h}\cdot \left(\sum_{k=1}^{n_E}\; h_{jk}^n \S_k\DS\right)=
\sum_{k=1}^{n_E}\; h_{jk}^n \left( {\bf h}\cdot  \S_k\DS\right)\\
\nonumber &=& \sum_{k=1}^{n_E}\; h_{jk}^n \left(\sum_{\ell=1}^{n_E} h_\ell\int_{\Om} \delta\sigma\nabla u_k^{(0)} \cdot \nabla u_{\ell}^{(0)} ~ d\r\right)\\
\nonumber &=&\sum_{k=1}^{n_E}\; h_{jk}^n \left(\int_{\Om} \delta\sigma\nabla u_k^{(0)} \cdot \na u_{\bf h}^{(0)} ~ d\r\right)
 =\int_\Om  \delta\sigma  \nabla u_{{\bf h}_{j}^n}^{(0)}\cdot\na u_{\bf h}^{(0)} d\r\\
\label{Phi_other} &=& \int_D  \nabla u_{{\bf h}_{j}^n}^{(0)}\cdot\na u_{\bf h}^{(0)} d\r.
\end{eqnarray}
Now we are ready to prove the two assertions (a) and (b).
\begin{enumerate}
\item[(a)] Using (\ref{Phi3}), (\ref{Phi_other}), and the Cauchy-Schwartz inequality, we have that
\begin{eqnarray*}
\left| \int_{q_n} \nabla u_j^{(0)} \cdot \nabla u_{\bf h}^{(0)} ~ d\r \right|^2
&=& \left| {\bf h}\cdot\PHI_{j}^n \right|^2 =
\left| \int_D  \nabla u_{{\bf h}_{j}^n}^{(0)}\cdot\na u_{\bf h}^{(0)} d\r \right|^2\\
&\leq & \int_D  \left| \na u_{{\bf h}_{j}^n}^{(0)} \right|^2 d\r \; \int_D  \left| \na u_{\bf h}^{(0)} \right|^2 d\r \\
& = &  \left({{\bf h}_{j}^n}\cdot\PHI_{j}^n \right)\; \int_D  \left| \na u_{\bf h}^{(0)} \right|^2 d\r \\
& = &  \left( \PHI^n_j\cdot (\delta\mathbb V)^{-1}\PHI^n_j \right)\; \int_D  \left| \na u_{\bf h}^{(0)} \right|^2 d\r,
\end{eqnarray*}
which shows (a).
\item[(b)] Again, we use (\ref{Phi3}), (\ref{Phi_other}), and the Cauchy-Schwartz inequality, to obtain that
\begin{eqnarray*}
\left| \PHI^n_j\cdot (\delta\mathbb V)^{-1}\PHI^n_j  \right| &= & \left| {{\bf h}_{j}^n}\cdot\PHI_{j}^n \right|
= \left| \int_{q_n} \nabla u_j^{(0)} \cdot \nabla u_{{\bf h}_{j}^n}^{(0)} ~ d\r \right|\\
&\leq & \left(\int_{q_n} |\nabla u_j^{(0)}|^2 ~ d\r \right)^{1/2}\; \left(\int_{q_n} |\nabla u_{{\bf h}_{j}^n}^{(0)}|^2 ~ d\r \right)^{1/2}\\
&\leq & \left(\int_{q_n} |\nabla u_j^{(0)}|^2 ~ d\r \right)^{1/2}\; \left(\int_{D} |\nabla u_{{\bf h}_{j}^n}^{(0)}|^2 ~ d\r \right)^{1/2}\\
&=& \left(\int_{q_n} |\nabla u_j^{(0)}|^2 ~ d\r \right)^{1/2}\; \left( \PHI^n_j\cdot (\delta\mathbb V)^{-1}\PHI^n_j  \right) ^{1/2},
\end{eqnarray*}
which shows (b). \hfill $\Box$
\end{enumerate}

The estimates in theorem~\ref{thm:S-FM}(a) and (b) provide the following observations that can be viewed as a sensitivity matrix based version of the previous work in \cite{harrach2009detecting,Harrach2010b}:
\begin{itemize}
\item If the pixel $q_n$ lies outside the inclusion $D$, then we can expect that
there exists an excitation pattern $h$ so that the resulting current will be large
in the pixel $q_n$ but small inside the inclusion $D$. See \cite{Gebauer2008} for the rigorous mathematical
basis of this argument.

Hence, we can expect that for some $j$ the expression
\begin{equation*}
\PHI^n_j\cdot (\delta\mathbb V)^{-1}\PHI^n_j \q\ge\q
\max_{{\bf h}\in \mathbb{R}^{{n_E}}}\f{ \left| \int_{q_n} \nabla u_j^{(0)} \cdot \nabla u_{\bf h}^{(0)} ~ d\r \right|^2}{
   \int_D |\nabla u^{(0)}_{\bf h}|^2 d\r}.
\end{equation*}
will be very large.
\item If the pixel $q_n$ lies inside the inclusion $D$, then
\begin{equation*}
\left|\PHI^n_j\cdot (\delta\mathbb V)^{-1}\PHI^n_j\right| \q \le\q  \int_{q_n}\left|\nabla
 u_j^{(0)}\right|^2 d\r.
\end{equation*}
will be not so large.
\end{itemize}
The above observations  enable us to predict the presence and absence of anomaly at each location $q_n$ by evaluating
\begin{equation}\label{Diag}
\zeta_j^n~:=\PHI^n_j\cdot \left( \delta\mathbb V\right)^{-1}\PHI^n_j.
\end{equation}

\subsection{Hybrid method combining S-FM with LM}\label{subsect:hybrid}
Now, we are ready to explain a hybrid method in which S-FM and LM cooperate with each other.

\subsubsection{Drawbacks of the standard linearized method}
\begin{figure}
\centering
\subfigure[]{\label{0thhybasis}
\begin{tabular}{|c|c|c|c|}
\hline
$\v_{\small{1}}$ & $\v_{\small{17}}$ & $\v_{\small{33}}$ & $\v_{\small{49}}$\\
\hline
\includegraphics[keepaspectratio=true,width=2.46cm]{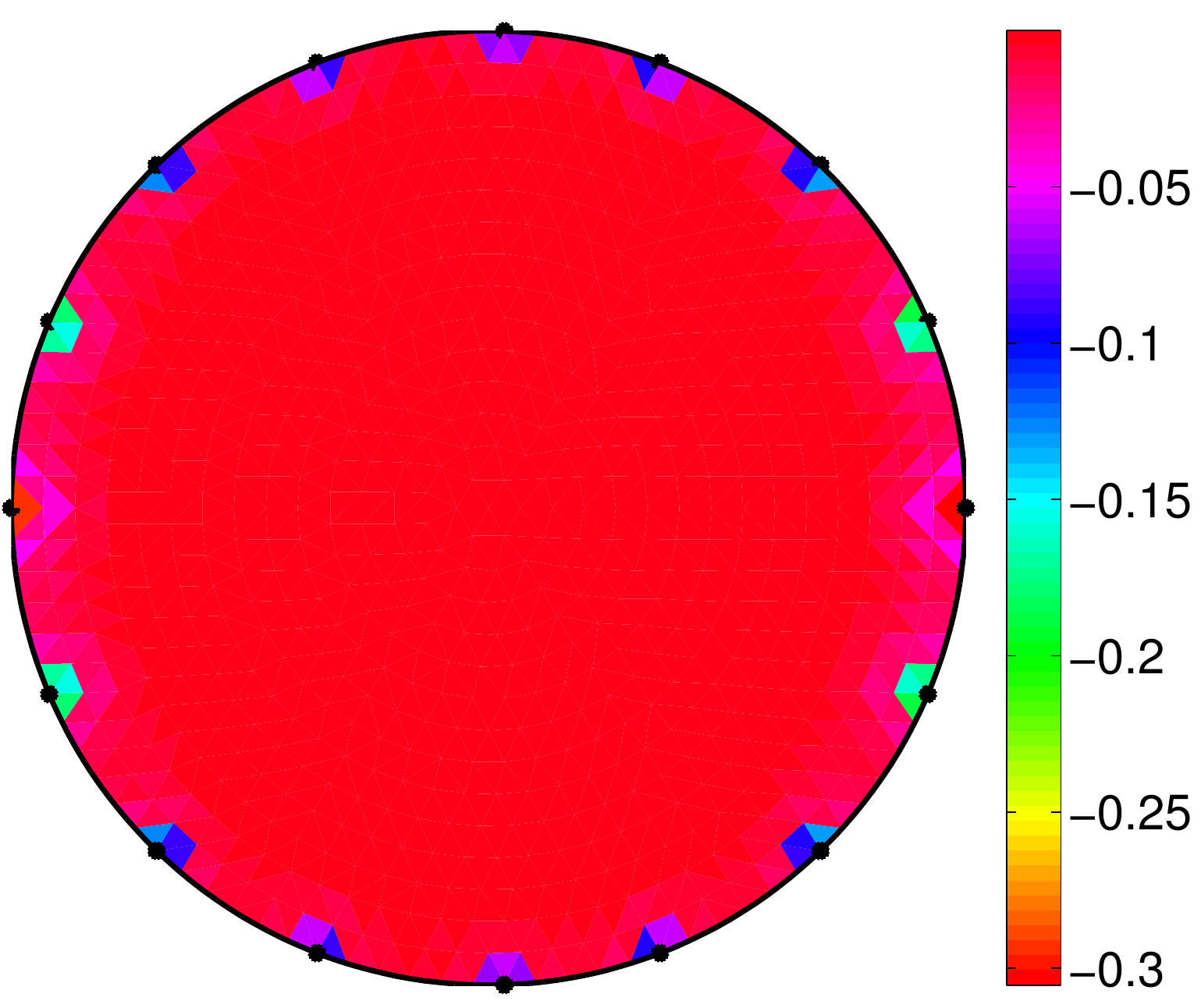}&
\includegraphics[keepaspectratio=true,width=2.46cm]{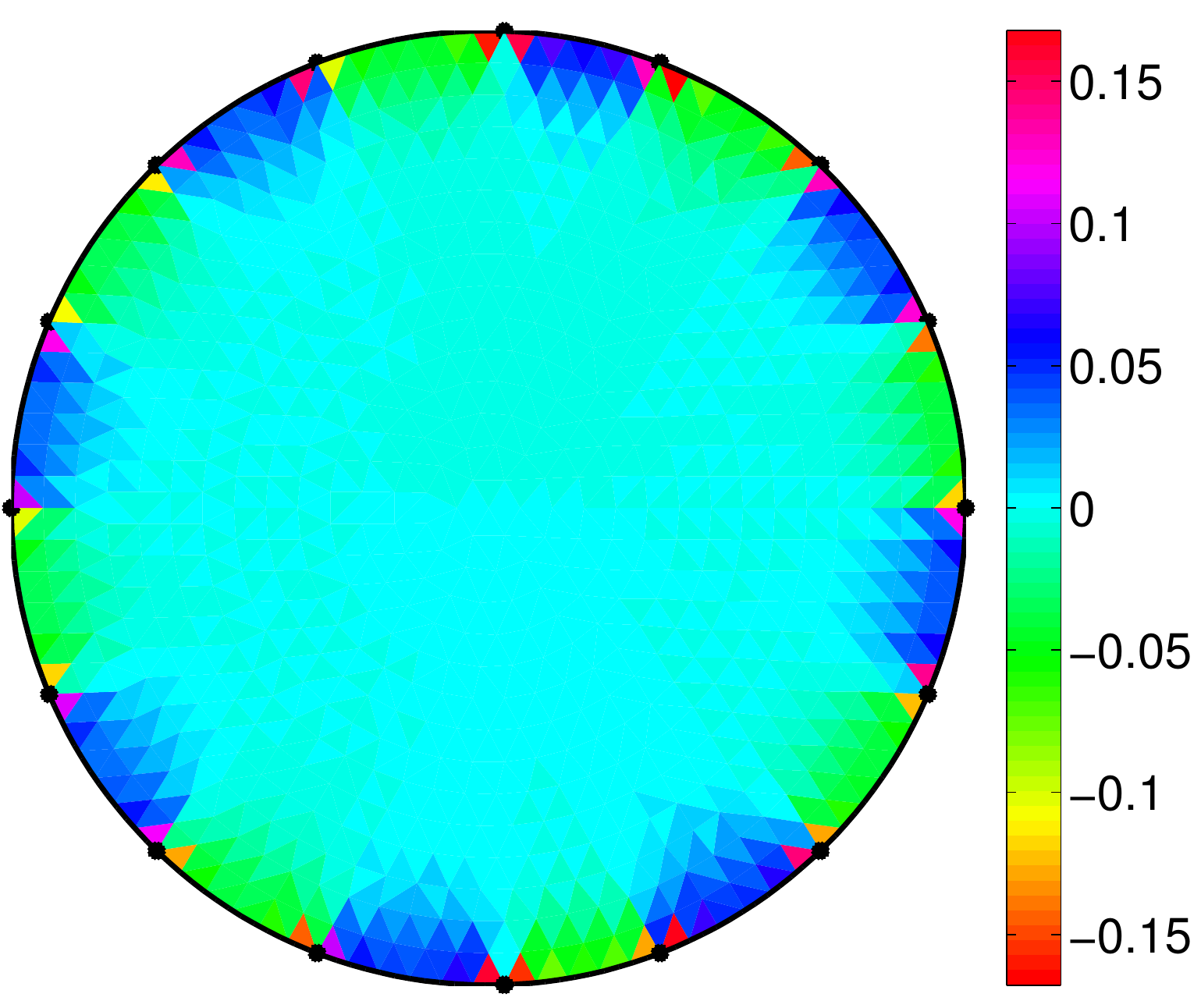}&
\includegraphics[keepaspectratio=true,width=2.46cm]{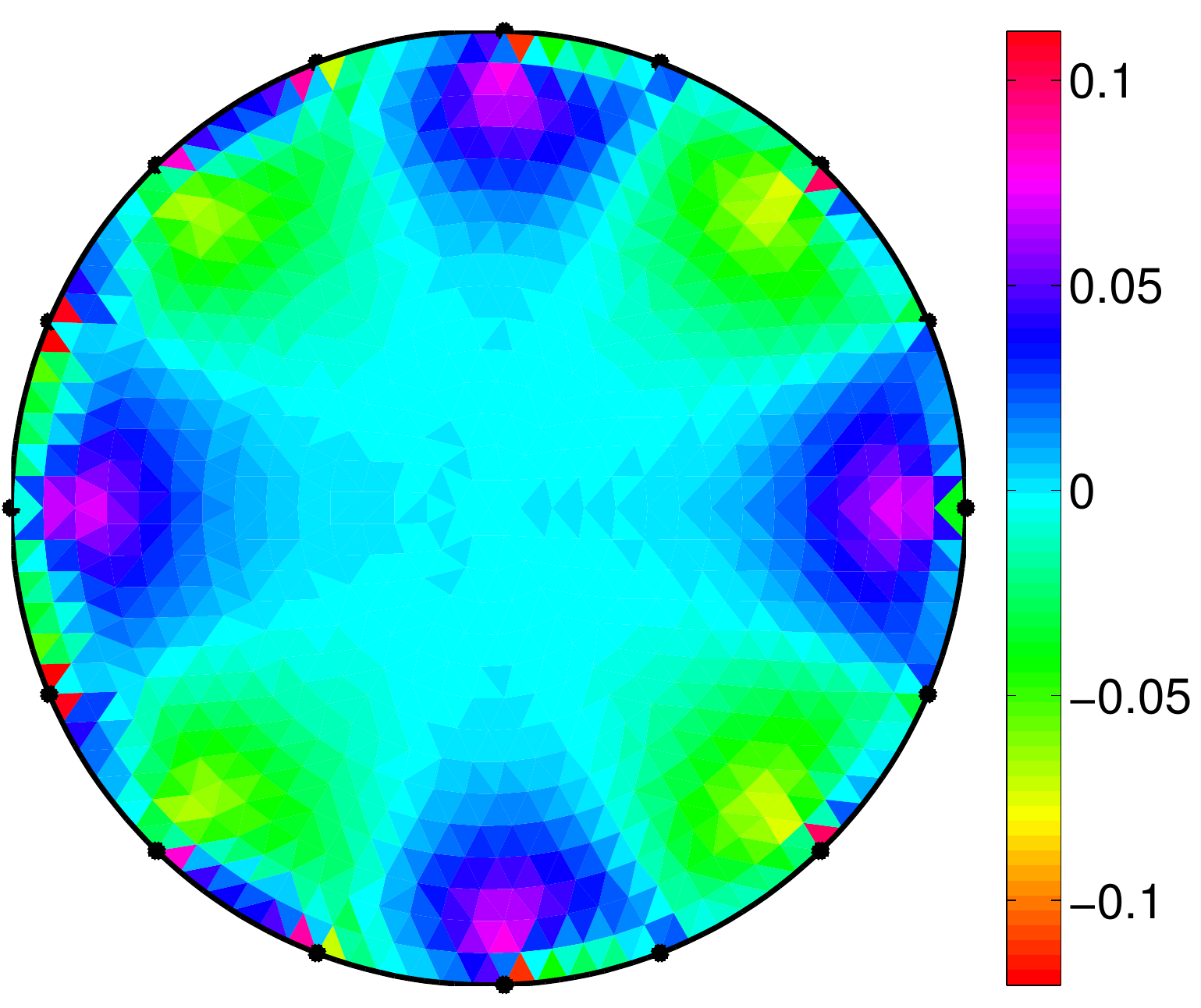}&
\includegraphics[keepaspectratio=true,width=2.46cm]{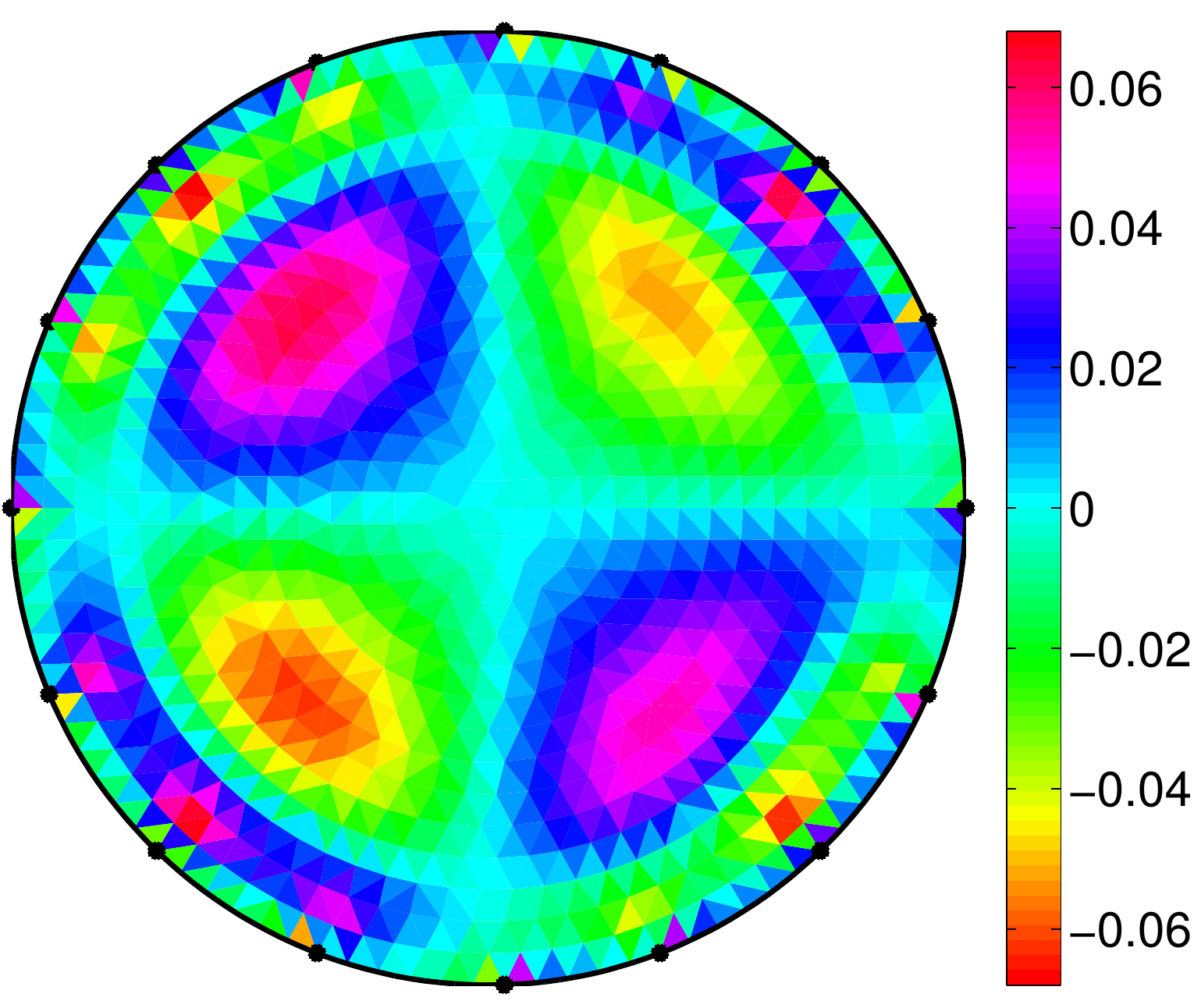}\\
\hline
$\v_{\small{16}}$& $\v_{\small{32}}$ & $\v_{\small{48}}$ & $\v_{\small{64}}$\\
\hline
\includegraphics[keepaspectratio=true,width=2.46cm]{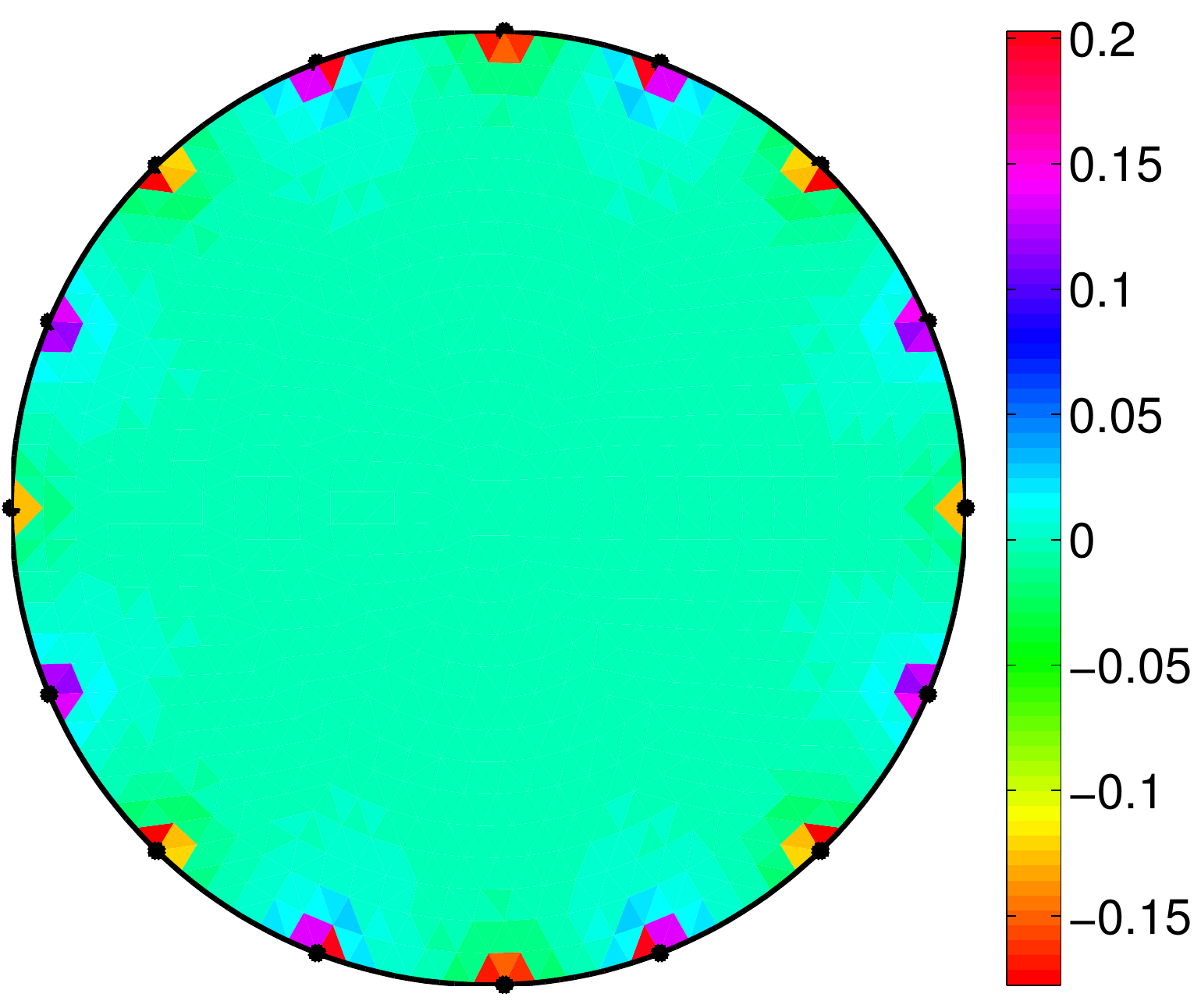}&
\includegraphics[keepaspectratio=true,width=2.46cm]{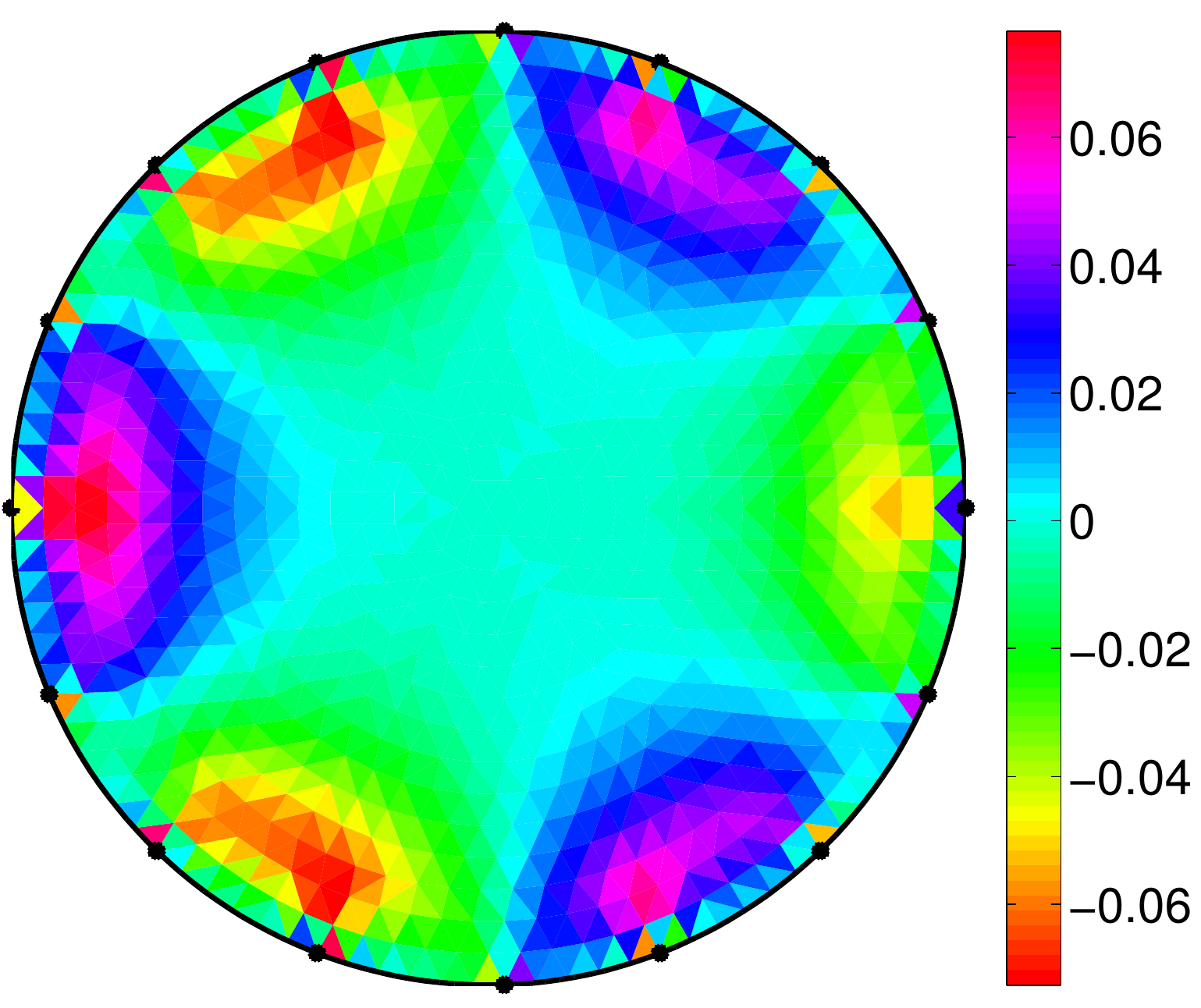}&
\includegraphics[keepaspectratio=true,width=2.46cm]{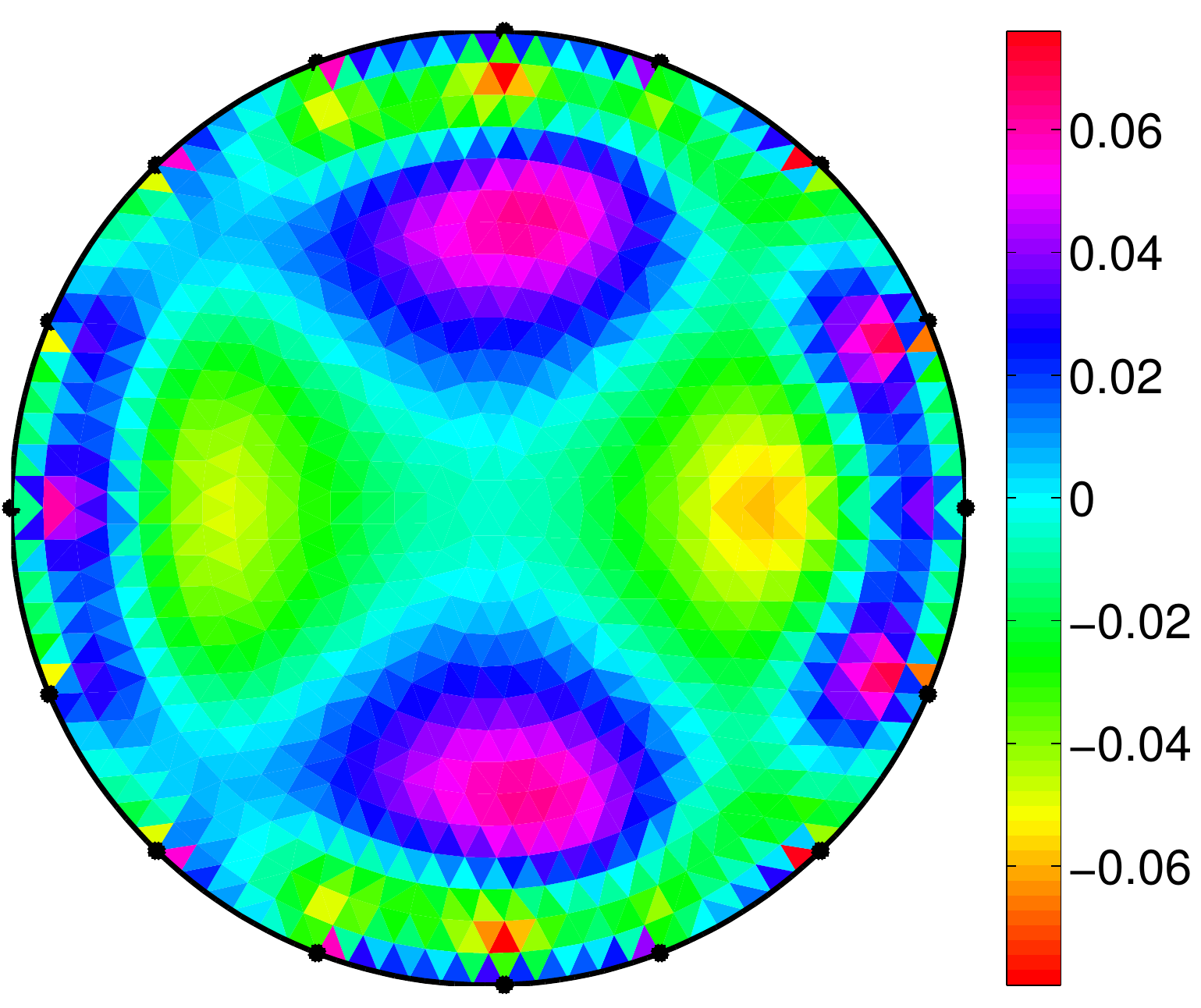}&
\includegraphics[keepaspectratio=true,width=2.46cm]{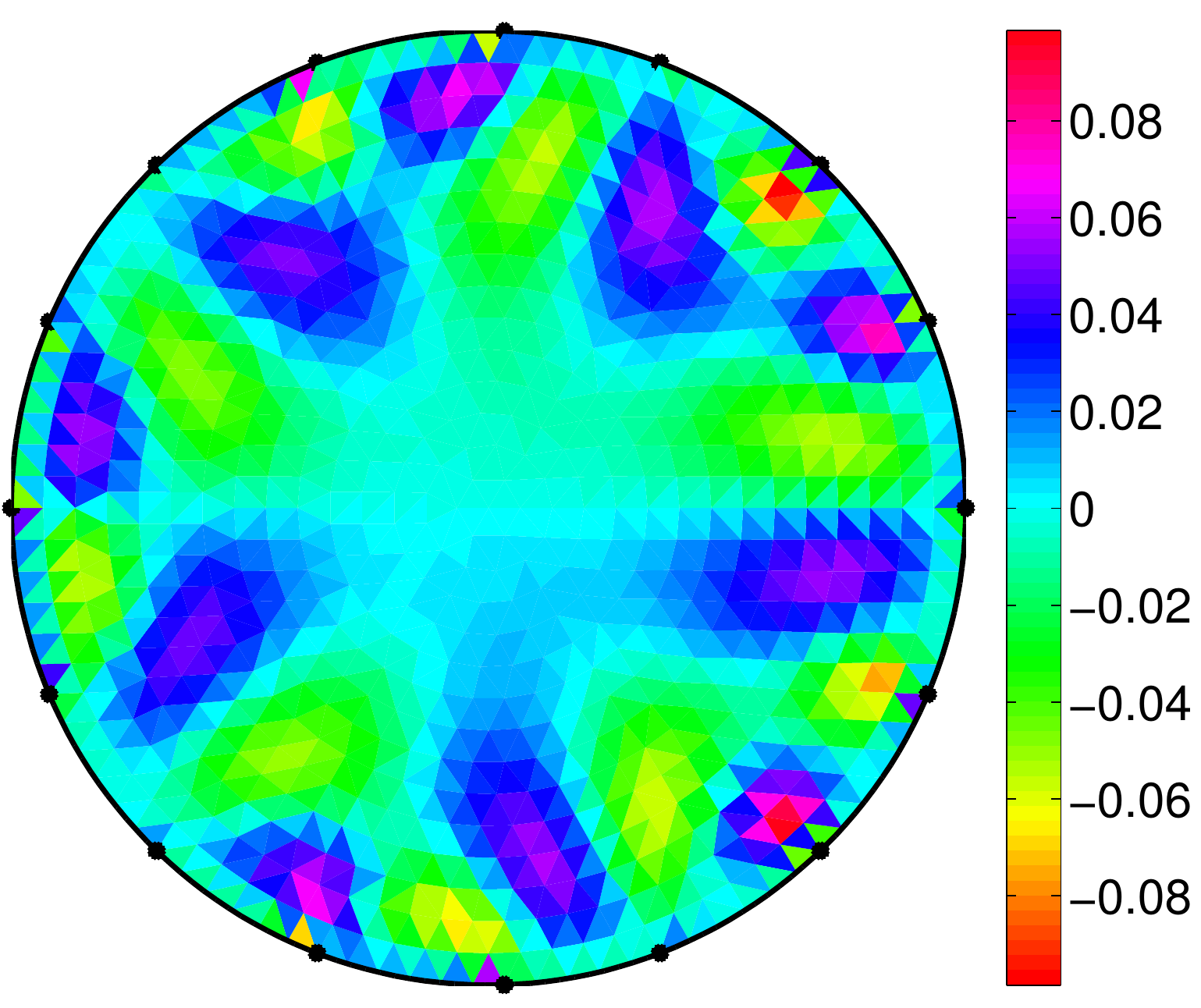}\\
\hline
\end{tabular}}
\phantom{.}\\
\subfigure[]{\label{fig:projected_image}
\begin{tabular}{cccc}
\includegraphics[keepaspectratio=true,width=2.5cm]{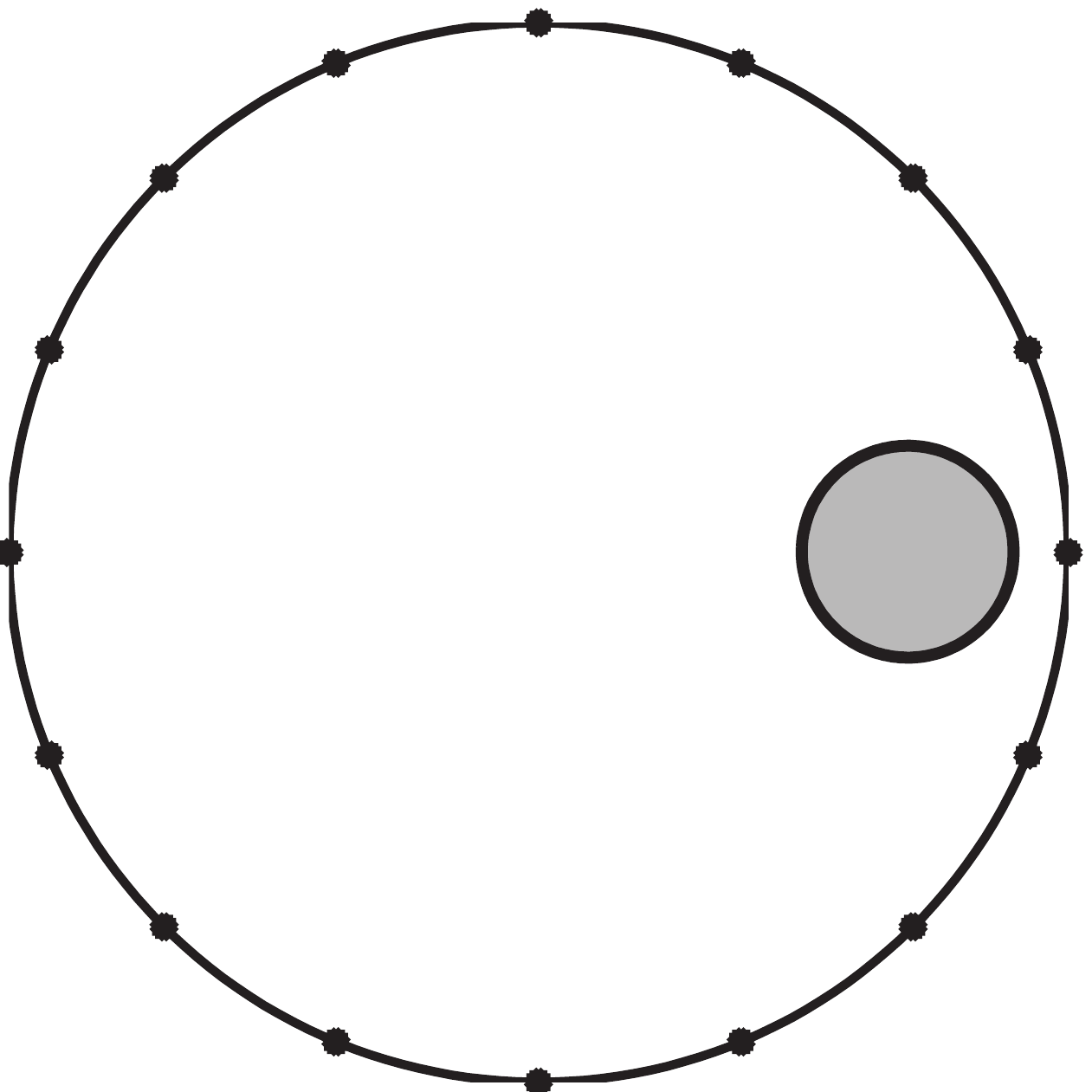}&
\includegraphics[keepaspectratio=true,width=3.1cm]{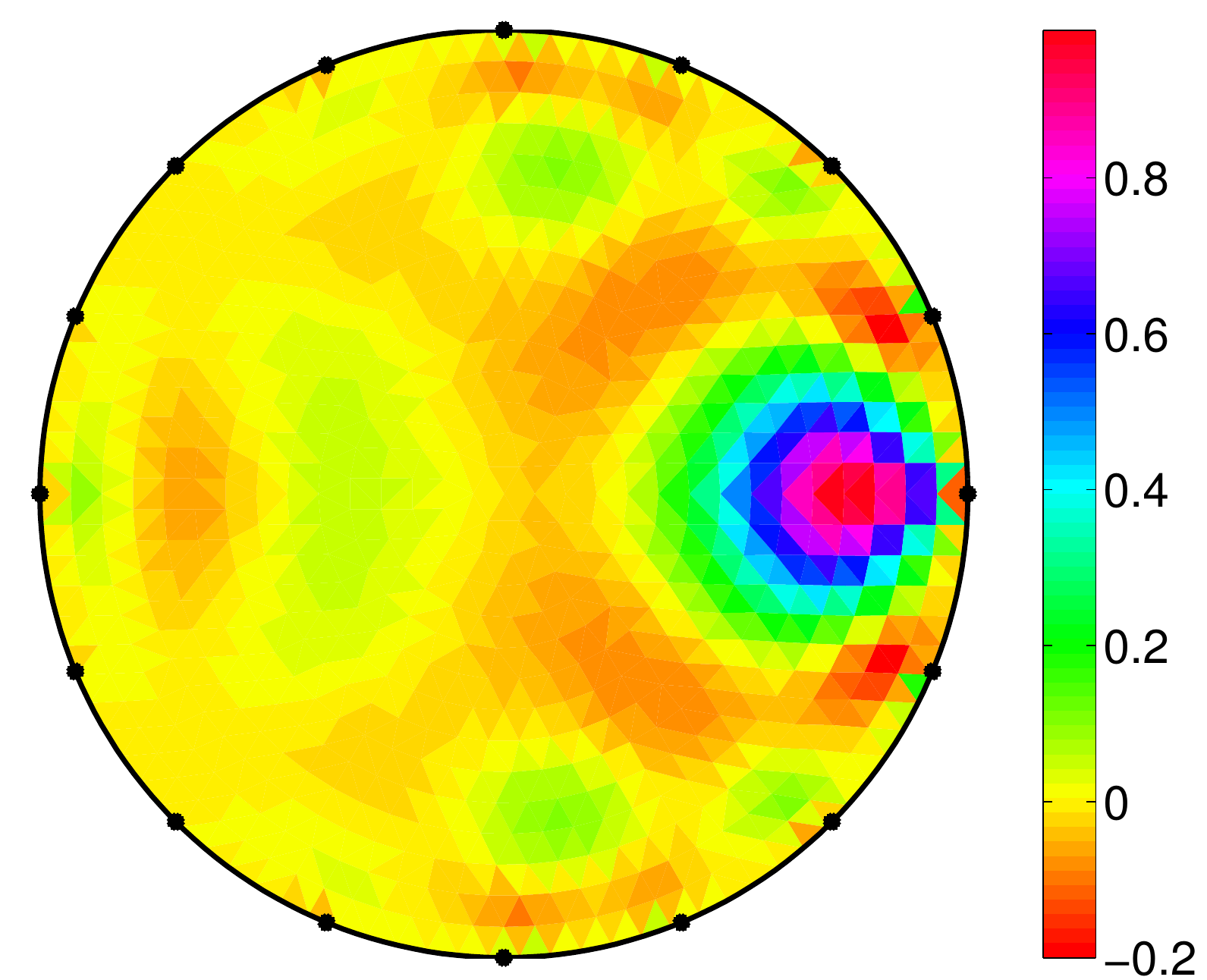}$\q$&$\q$
\includegraphics[keepaspectratio=true,width=2.5cm]{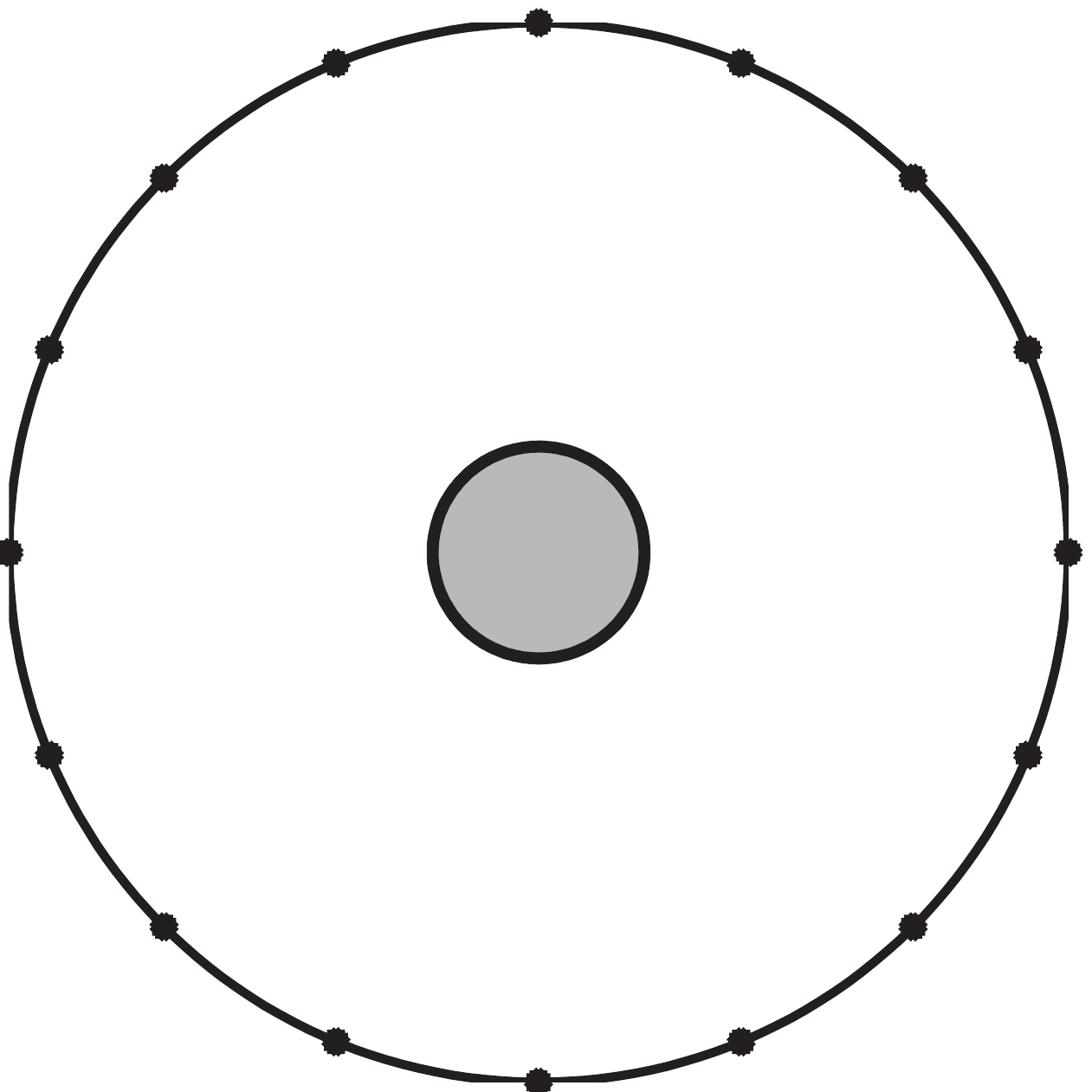}&
\includegraphics[keepaspectratio=true,width=3.1cm]{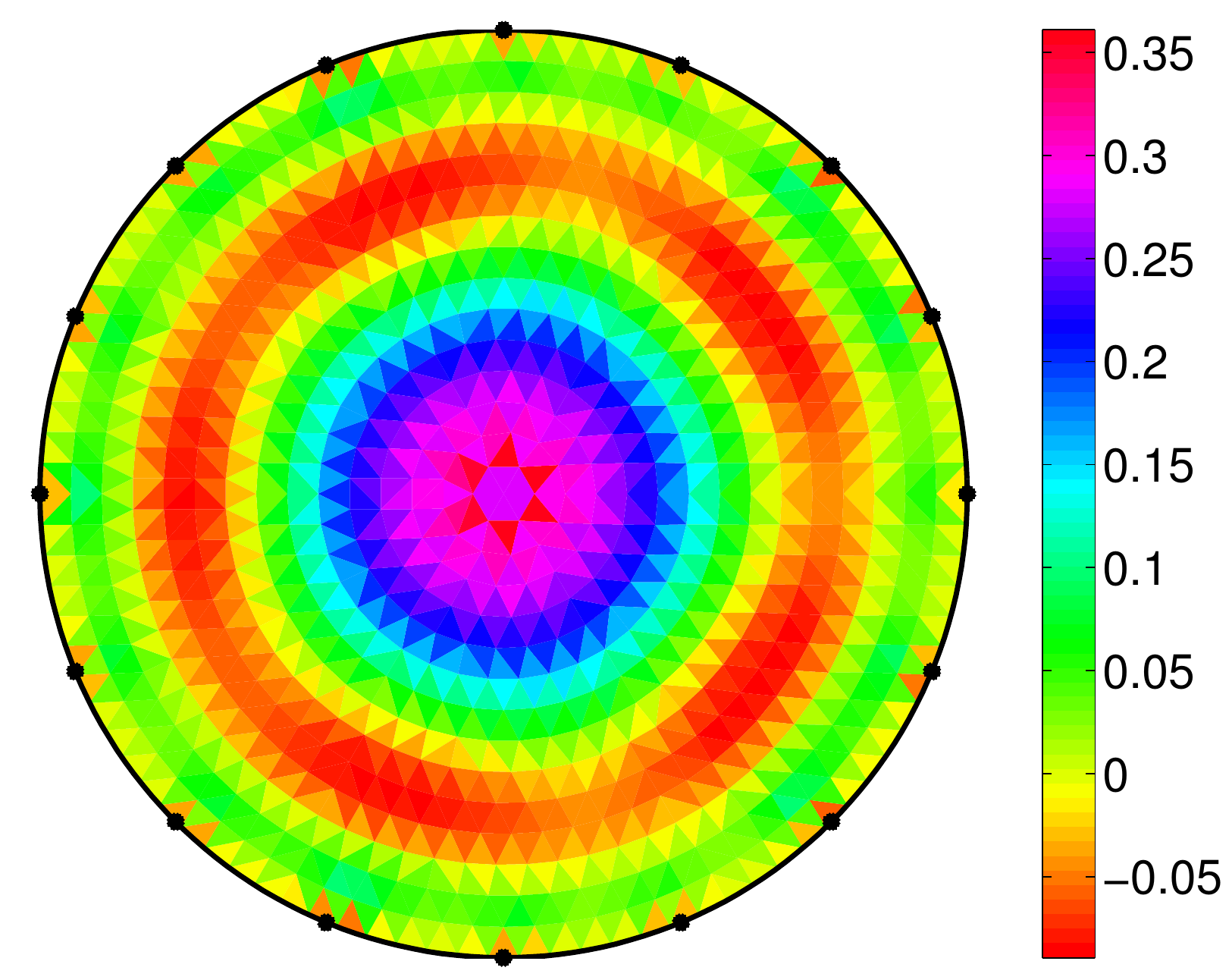}
\end{tabular}}
\caption{(a) Images of eigenvectors $\v_1,\cdots\v_{64}$ of $\S^*\S$
(b) The color images (right images) are projected images of the true  $\DS$ (left images) onto  the vector space spanned by $\{\v_1,\cdots\v_{64}\}$.  } \label{basis_proj}
\end{figure}

For the image reconstruction using LM, assume that  we use tSVD of the form $\S\approx\mathbf{U}_{t_0}\mathbf{\Lambda}_{t_0} \mathbf{V}_{t_0}^*$  where $\mathbf{\Lambda}_{t_0}$ is a diagonal matrix with the $t_0$ largest singular values $\lambda_1\ge \cdots\ge  \lambda_{t_0}>0$, $\mathbf{U}_{t_0}$ and $\mathbf{V}_{t_0}$ is a matrix with the first $t_0$-th columns of the left singular matrix $\mathbf{U}$ and the right singular matrix $\mathbf{V}$, respectively. Here,  $t_0$ can be determined by the noise level in the measurements. The smallest singular value $\lambda_{t_0}$ must be sufficiently large in such a way that the following reconstruction method is reasonably robust against noise in the data-set $\DV$:
\begin{equation}\label{LM2}
\DS_{S}~:=~\S^\dag~\DV=~ \sum_{t=1}^{t_0}\frac{1}{\lambda_t}\langle \DV,\mathbf{u}_t\rangle\v_t
\end{equation}
where $\S^\dag=\mathbf{V}_{t_0}\mathbf{\Lambda}_{t_0}^{-1} \mathbf{U}_{t_0}^*$ is the tSVD pseudoinverse of $\S$.
To explain the underlying idea  intuitively,
we focus on the standard 16-channel EIT system. In this EIT system, we usually choose $t_0=64$ which is the  number corresponding to $\f{\lambda_{64}}{\lambda_{1}}\approx 10^{-3}$. Figure \ref{0thhybasis} shows a part of images of $\{\v_1,\cdots, \v_{64}\}$. Note that the best achievable image by LM would be the projected image of the true
$\DS$  into the vector space spanned by $\{\v_1,\cdots, \v_{t_0}\}$ as shown in Figure \ref{fig:projected_image}. Clearly, the basis $\{\v_1,\cdots, \v_{t_0}\}$ is insufficient to achieve localization of conductivity anomalies and  causes   Gibbs ringing artifacts.

\subsubsection{A naive combination of LM and S-FM}
To deal with the drawbacks in LM, it seems natural to use the additional information obtained from S-FM as a regularization term, i.e., to
\begin{equation}\label{regularization}
\text{minimize}~~||\S\DS-\DV||_2^2+\beta^2||\mathbf{W}\DS||_2^2
\end{equation}
where $\beta$ is a regularization parameter and $\mathbf{W}$ a positive diagonal matrix having $n$-th entry \begin{equation}\label{SFMimage}
w_n:=~\ln\left(1+\sum_{j=1}^{n_E}\left|\frac{\zeta_j^n}{\S_j^n\cdot\S_j^n}\right|\right)=
~\ln\left(1+\sum_{j=1}^{n_E}\left|\frac{\S_j^n\cdot\delta\Bbb{V}^{-1}\S_j^n}{\S_j^n\cdot\S_j^n}\right|\right).
\end{equation}
Recalling the property of $\zeta_j^n$ from subsection \ref{subsect:S_FM}, a large value of  $w_n$ corresponds to a higher chance of getting $\delta\sigma_n=0$.
 The least squares formulation of \eref{regularization}  takes the form
$$\min\left|\left|
\underbrace{\left(\begin{array}{ccc}
&\S&\\
&&\\
&\beta\mathbf{W}&\\
\end{array}\right)}_{:=\mathbf{B}}\DS-\underbrace{\left(\begin{array}{c}
\DV\\
\\
\mathbf{0}\\
\end{array}\right)}_{:=\mathbf{b}}\right|\right|_2$$
where $\mathbf{0}$ is a zero vector. Note that this method tries to make zeros in which a subregion having large value of $w_n$ has a higher likelihood of $\delta\sigma=0$. 

We demonstrate in several numerical examples in the next section that this naive combination of LM and S-FM methods performs better than the standard linearized 
method, but still produces some ringing artifacts. We believe that this is due to the singular value structure of the matrix ${\bf B}$. 
Let the tSVD pseudoinverse of ${\bf B}$ be
${\bf B}^\dag=\hat{\mathbf{V}}_{t_1} \hat{\mathbf{\Lambda}}^{-1}_{t_1} \hat{\mathbf{U}}^*_{t_1}$ where $\hat{\mathbf{V}}_{t_1}$ and $\hat{\mathbf{U}}_{t_1}$, respectively, is the truncated right and left singular matrix of ${\bf B}$. There exists a significant difference between the right singular vectors  $\{\hat\v_1,\cdots, \hat\v_{t_1}\}$ and $\{\v_1,\cdots, \v_{t_0}\}$.
Roughly speaking, the right singular vectors of ${\bf B}$ corresponding to high singular values can be viewed as a basis to reconstruct images in the region of those pixels having high values of $w_n$. Hence, according to S-FM, the right singular vectors associated with high singular values are used to reconstruct images in the region of $\delta\sigma=0$. On the other hand, the right singular vectors associated with low singular values are used to reconstruct images of anomaly region ( $\delta\sigma\neq 0$), that is, the expression of $\DS_{B}:={\bf B}^\dag\mathbf{b}$  provides image of anomaly region ( $\delta\sigma\neq 0$) by a linear combination of right singular vectors associated with low singular values. Thus, $\DS_{B}={\bf B}^\dag\mathbf{b}$ may still contain undesirable artifacts as shown in Figure \ref{recon_circle}.

\subsubsection{The proposed combination of LM and S-FM}
The very intuitve above approach of using the information from the S-FM as a regularization term did not produce satisfactory numerical results. Summing up the above arguments, we believe that this is due to the fact that the right singular vectors 
(corresponding to the higher singular values) of the resulting linear system tend to have their support outside the anomaly. Based on this idea we try to design a linear system where these vectors tend to have their support
inside the anomaly. We suggest to use the following new matrix
\begin{equation}\label{AA}
\mathbf{A}=\left[\begin{array}{ccc}
&\S&\\
&&\\
&\alpha\mathbf{W}^{-1}&\\
\end{array}\right]
\end{equation}
where $\alpha$ is a positive constant. Based on the above heuristic arguments, we can expect that using the auxiliary matrix $\alpha\mathbf{W}^{-1}$ should change the singular values in such a way that the right singular vectors of ${\bf A}$ associated with high singular values tend to be useful for reconstructing images in the region of anomaly (i.e. the region with $\delta\sigma\neq 0$). 
Thus, using ${\bf A}$ instead of the sensitivity matrix $\S$ adds some location information of conductivity anomalies to the 
matrix used in the reconstruction process.

To extend the linear system $\S \DS = \DV$ to one with ${\bf A}\DS$ on the left hand side, 
we also have to extend the right hand side. The extension is only done to change the singular vector structure.
In the noise-free, unregularized case the new system should still be equivalent to solving $\S \DS = \DV$. 
Hence, we choose
\begin{equation}\label{recon}
\mathbf{A}\DS=\left[\begin{array}{ccc}
&\S&\\
&&\\
&\alpha\mathbf{W}^{-1}&\\
\end{array}\right]\DS=\left[\begin{array}{c}
\DV\\
\\
\alpha\mathbf{W}^{-1}\S^{\dag}\DV\\
\end{array}\right].
\end{equation}

Note that the second line in (\ref{recon}) is not to be interpreted as a penalization term
for $\DS-\S^{\dag}\DV$. For all $\alpha$, (\ref{recon}) obviously possesses the minimum norm solution $\DS=\S^{\dag}\DV$.
In that sense, without any further regularization, (\ref{recon}) is equivalent to $\S \DS = \DV$. However, a TSVD regularization of
(\ref{recon}) can be expected to yield
superior results since $\mathbf{A}$ has been designed so that its right singular vectors 
are useful for reconstructing images in the region of the anomaly.
 
We numerically demonstrate in the next section that this new reconstruction method indeed shows a better performance in localization with alleviating  Gibbs ringing artifacts. 
Let us however stress again, that the motivation for the proposed method comes from purely heuristic arguments and numerical evidence. 
At the moment, we have no further theoretical explanation for the promising performance of this new method.

\section{Numerical results}

\begin{figure}
\centering
\begin{tabular}{|c|cccc|c|}
\hline
Case &  $\DS$ & $\DS_{S}$ & $\DS_{B}$ &  $\DS_{A}$ &$\mathbf{W1}$ \\
\hline
\raisebox{4ex}{\footnotesize\begin{tabular}{c}
(a)
\end{tabular}}
&
\includegraphics[keepaspectratio=true,height=1.5cm]{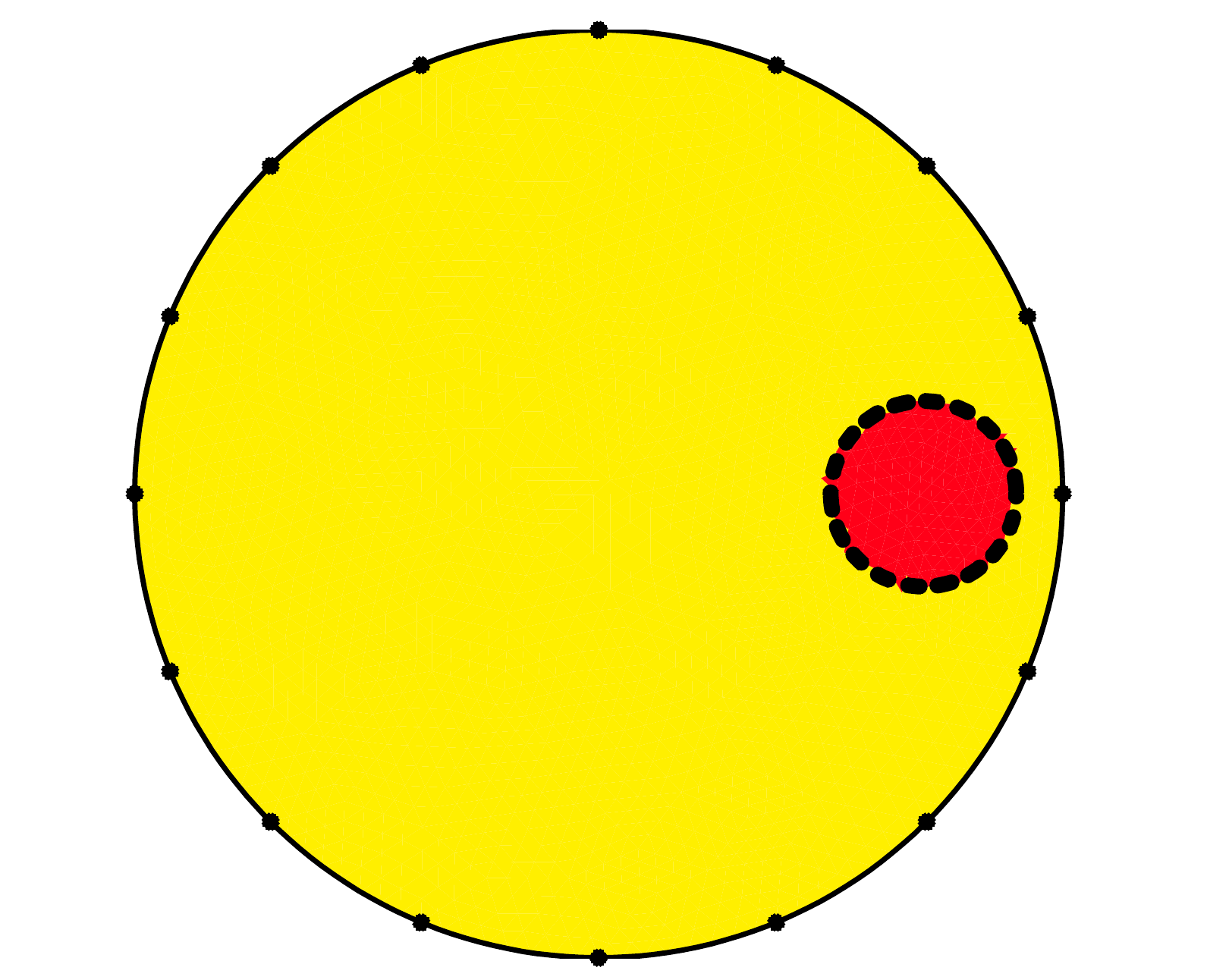}&
\includegraphics[keepaspectratio=true,height=1.5cm]{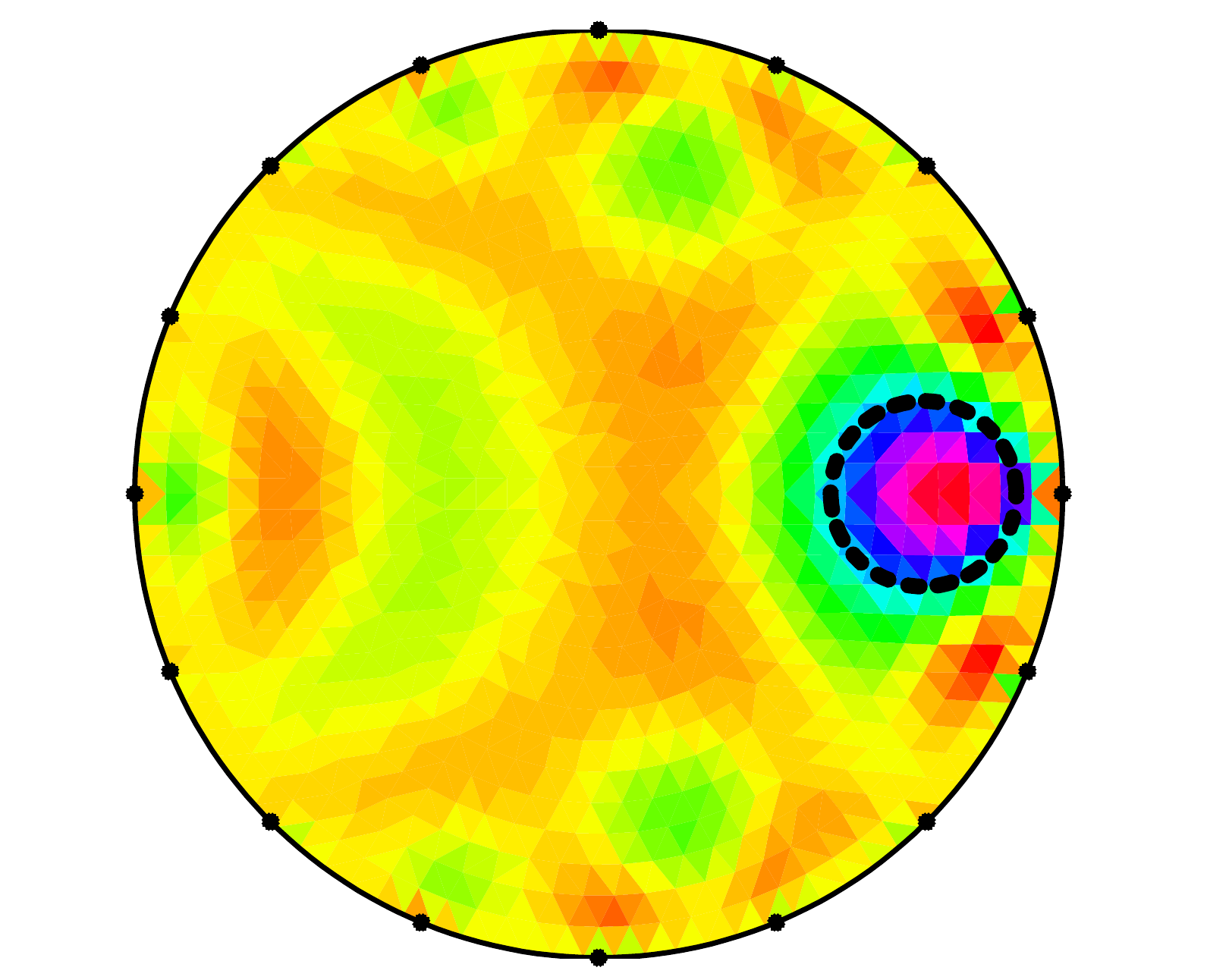}&
\includegraphics[keepaspectratio=true,height=1.5cm]{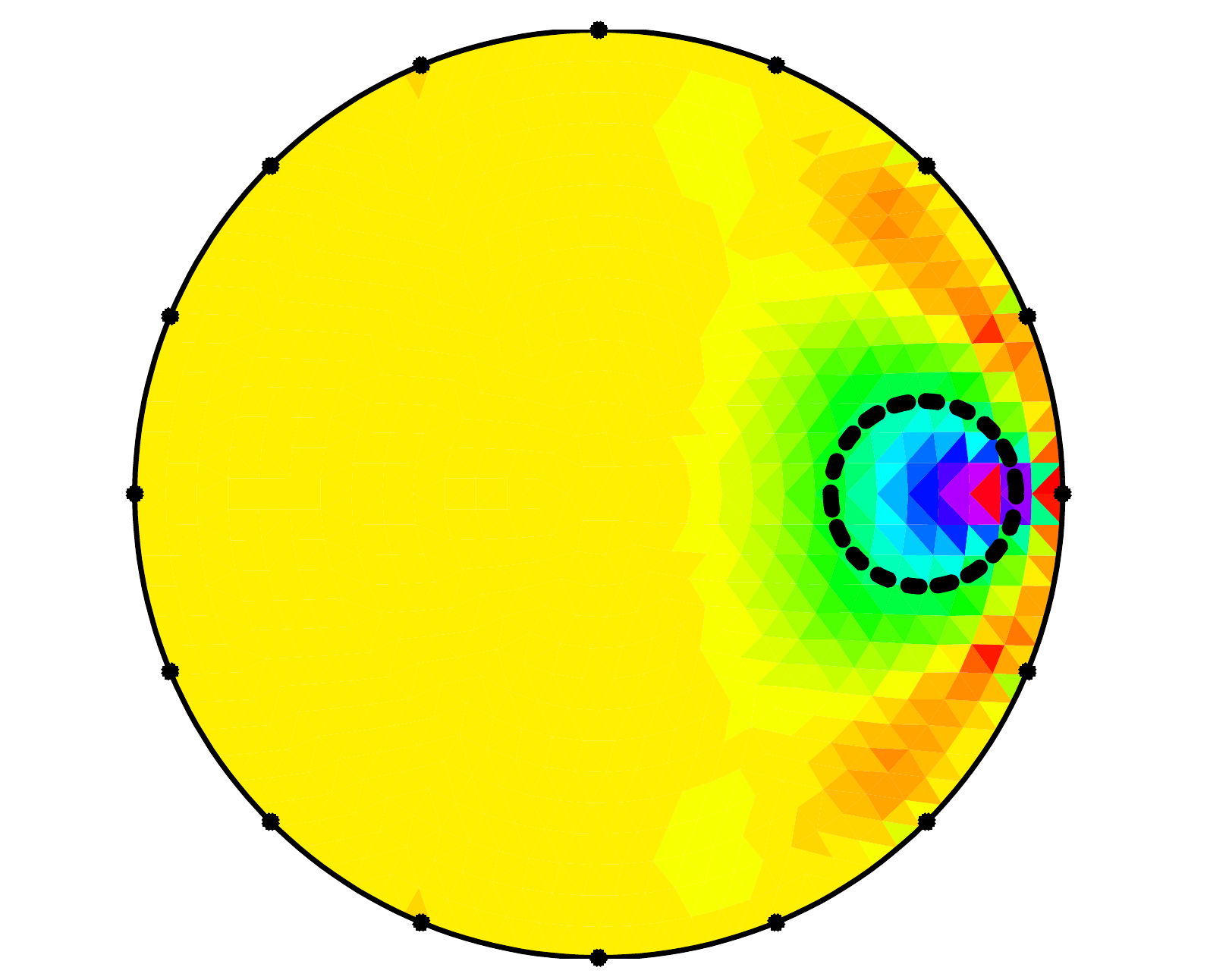}&
\includegraphics[keepaspectratio=true,height=1.5cm]{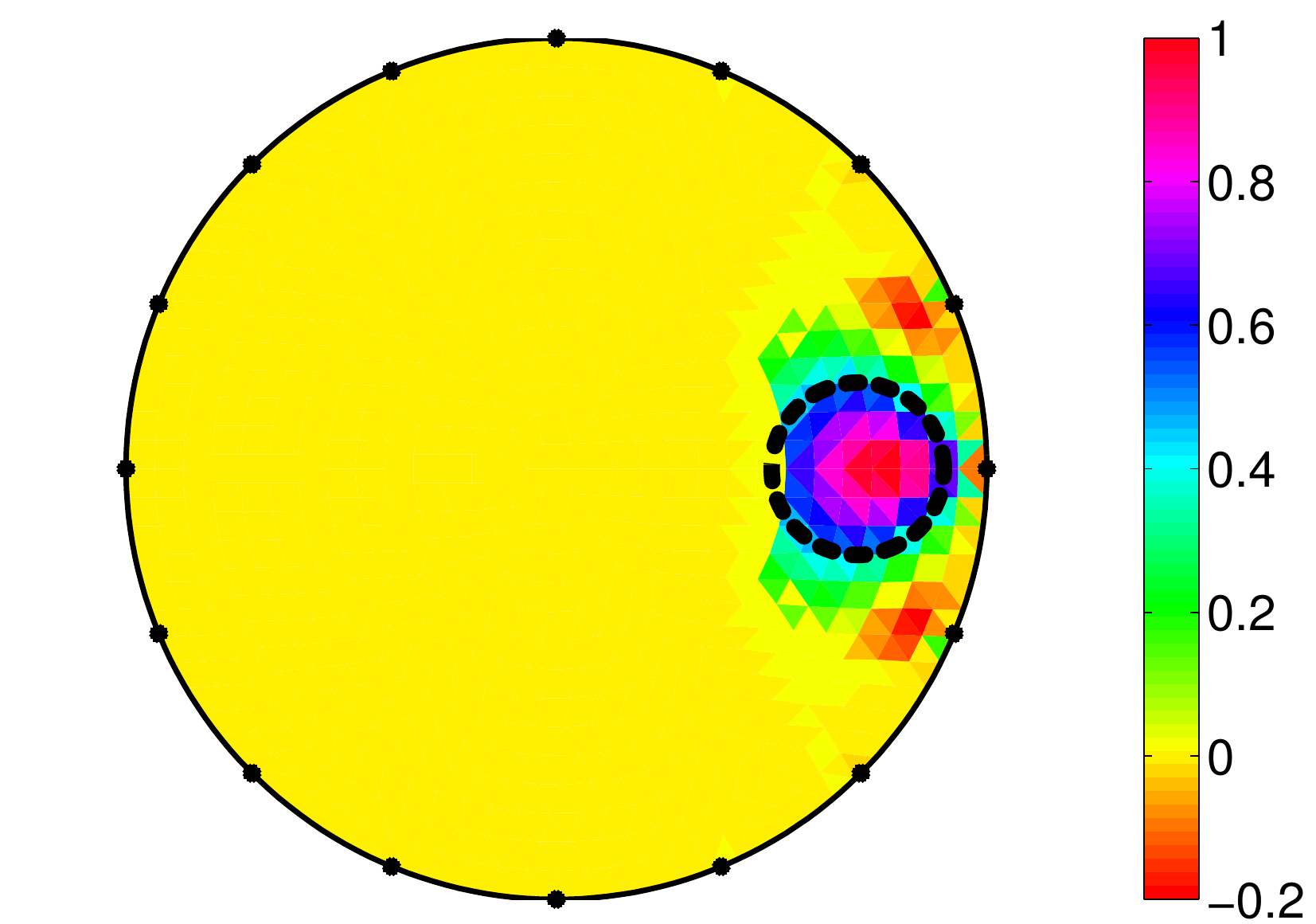}&
\includegraphics[keepaspectratio=true,height=1.5cm]{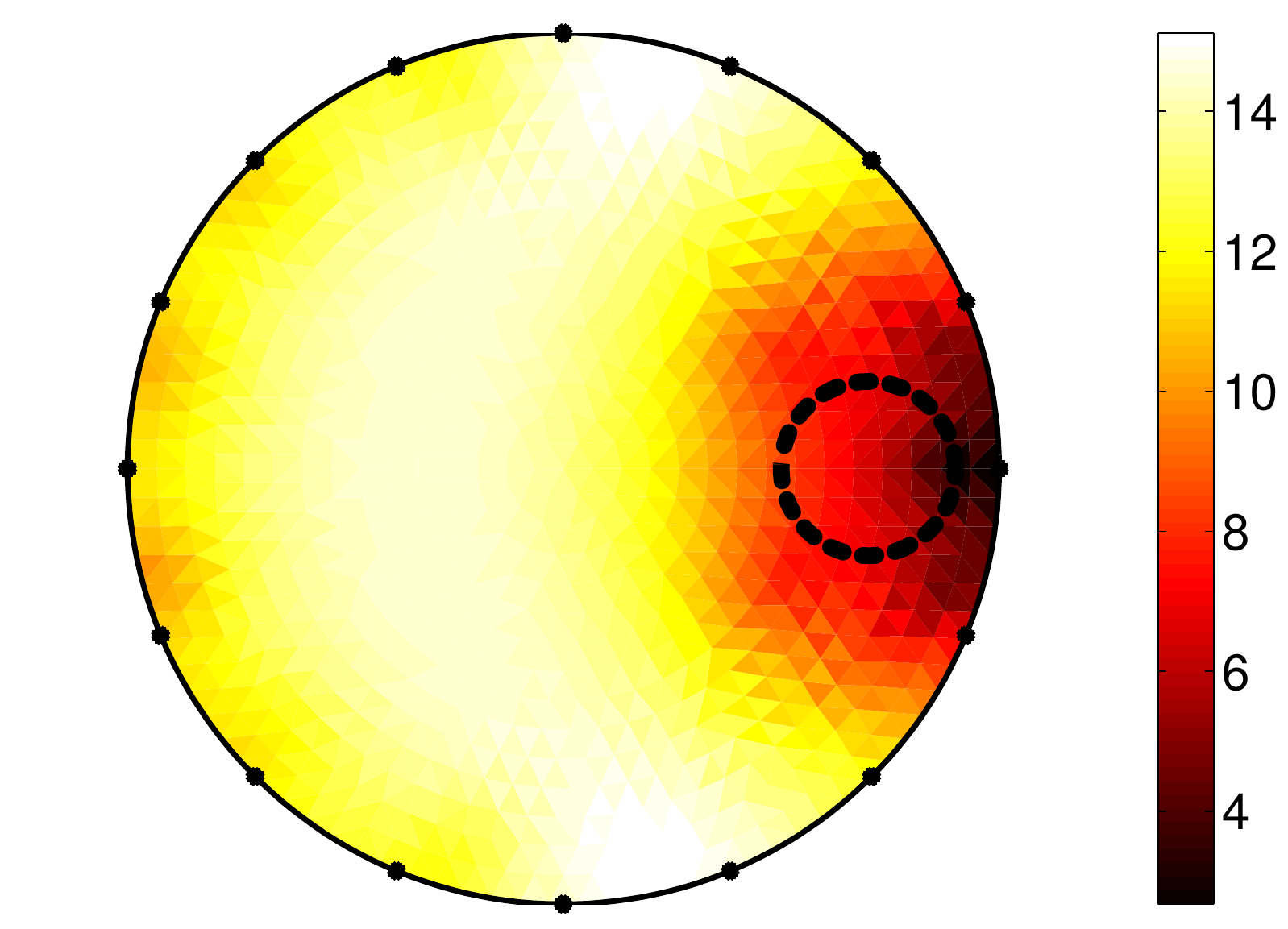}\\
\hline
\raisebox{4ex}{\footnotesize\begin{tabular}{c}
(b)
\end{tabular}}
&
\includegraphics[keepaspectratio=true,height=1.5cm]{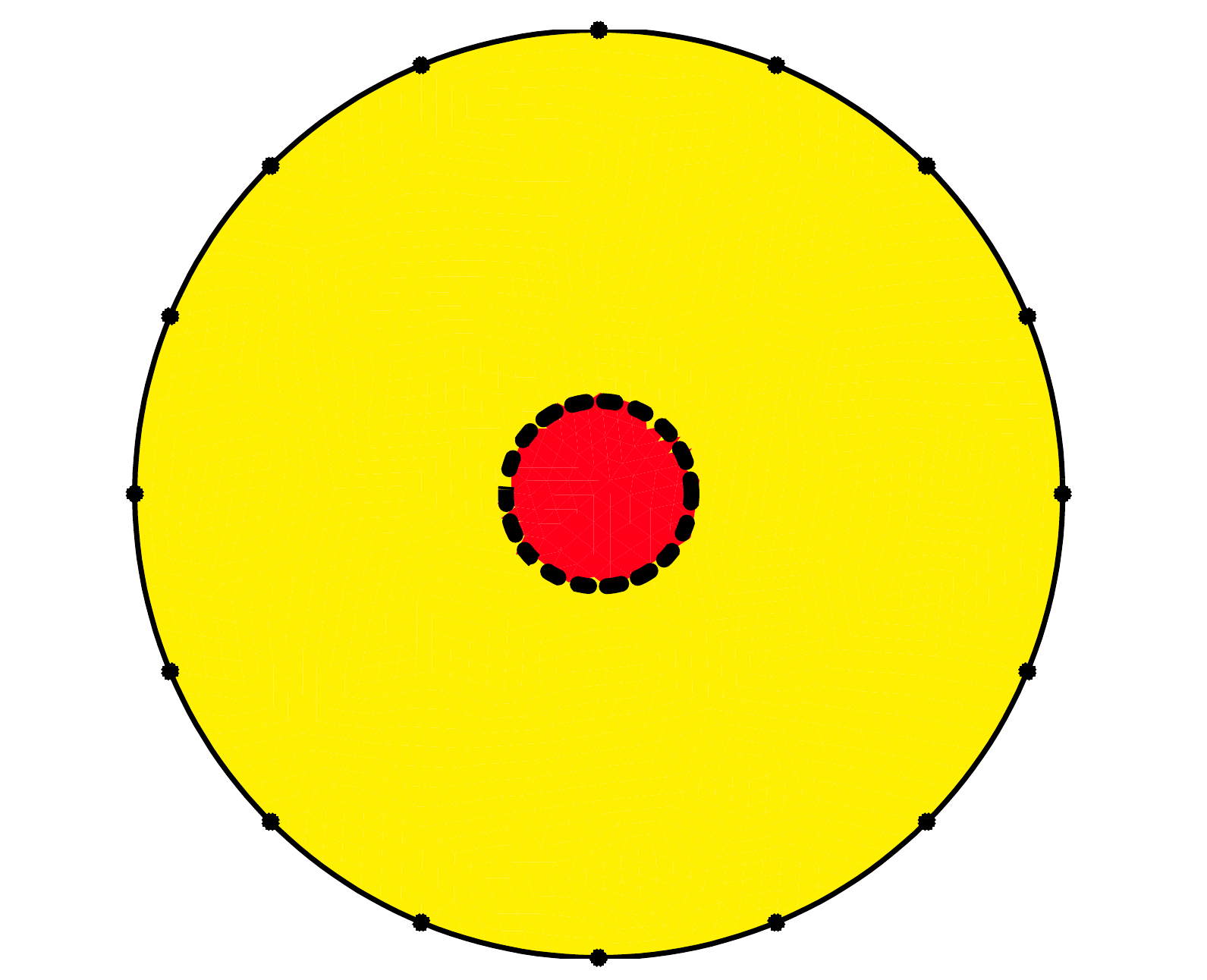}&
\includegraphics[keepaspectratio=true,height=1.5cm]{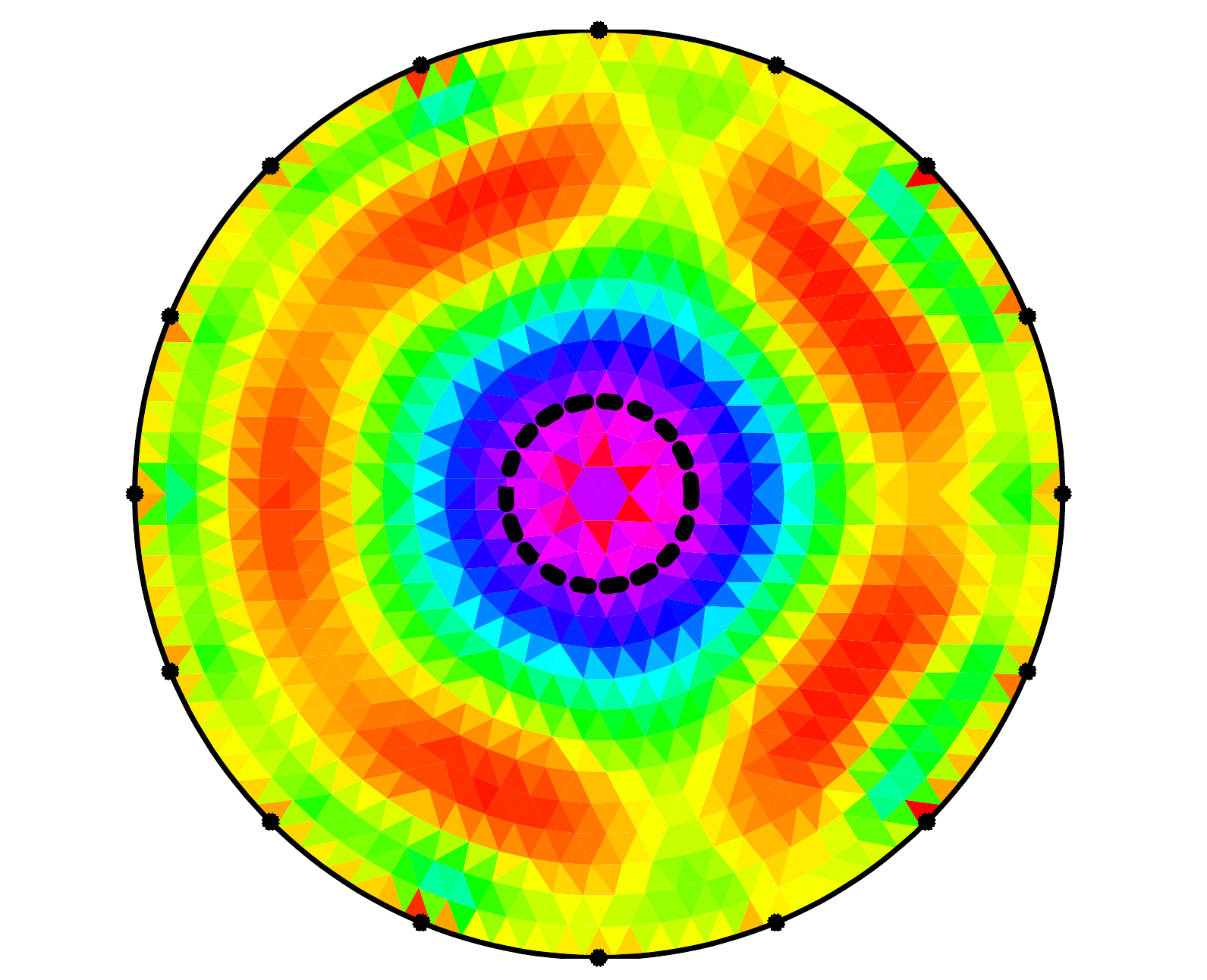}&
\includegraphics[keepaspectratio=true,height=1.5cm]{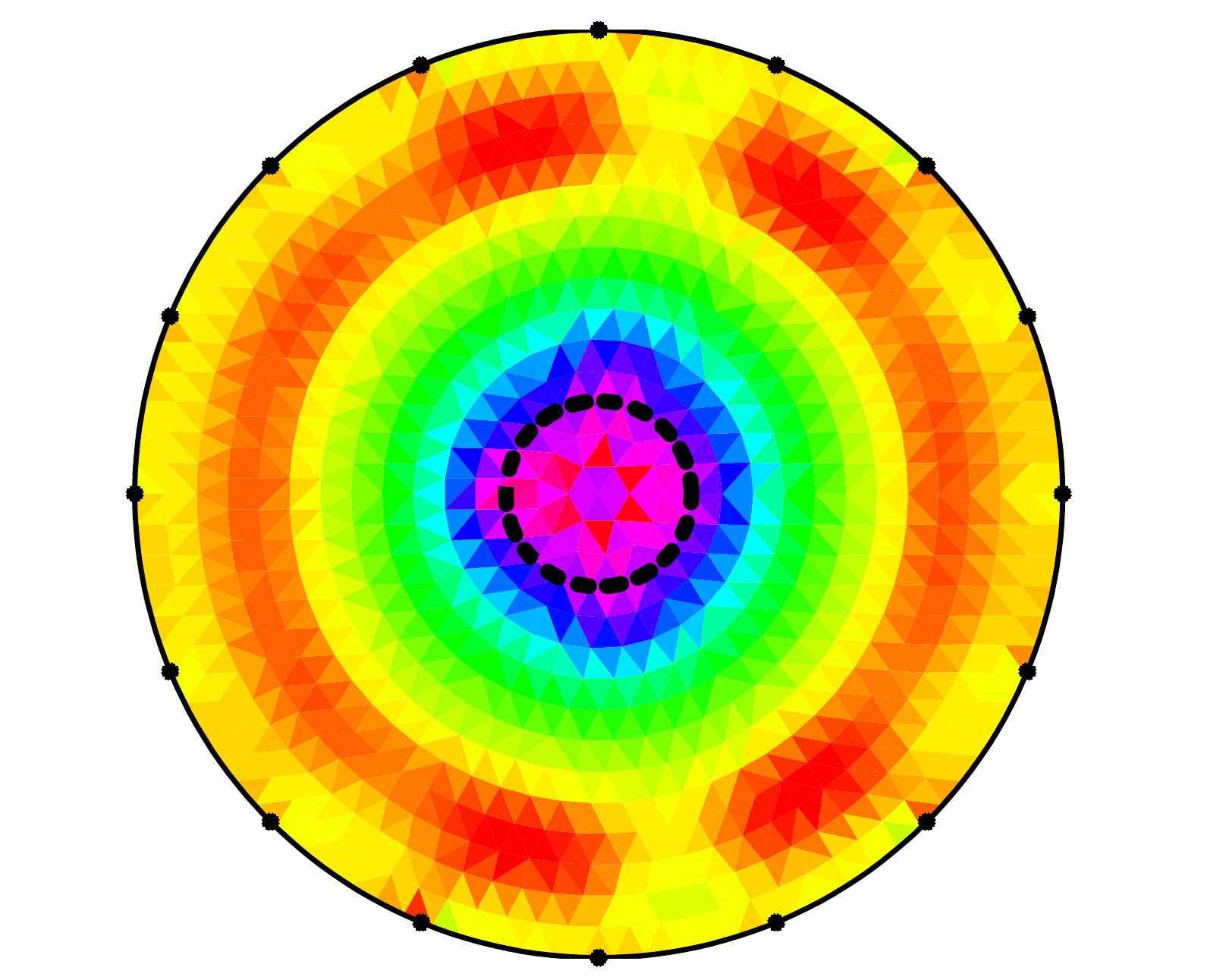}&
\includegraphics[keepaspectratio=true,height=1.5cm]{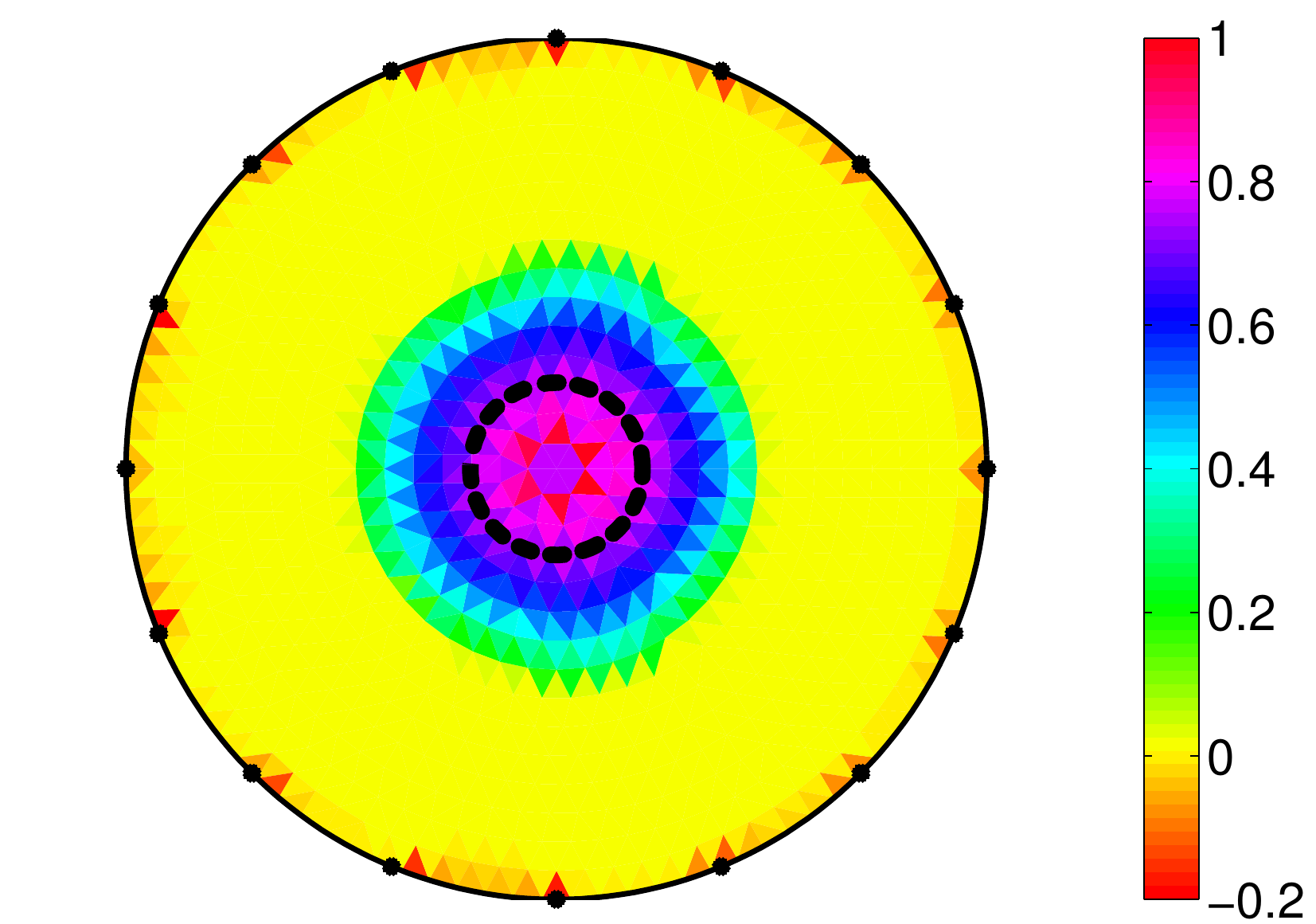}&
\includegraphics[keepaspectratio=true,height=1.5cm]{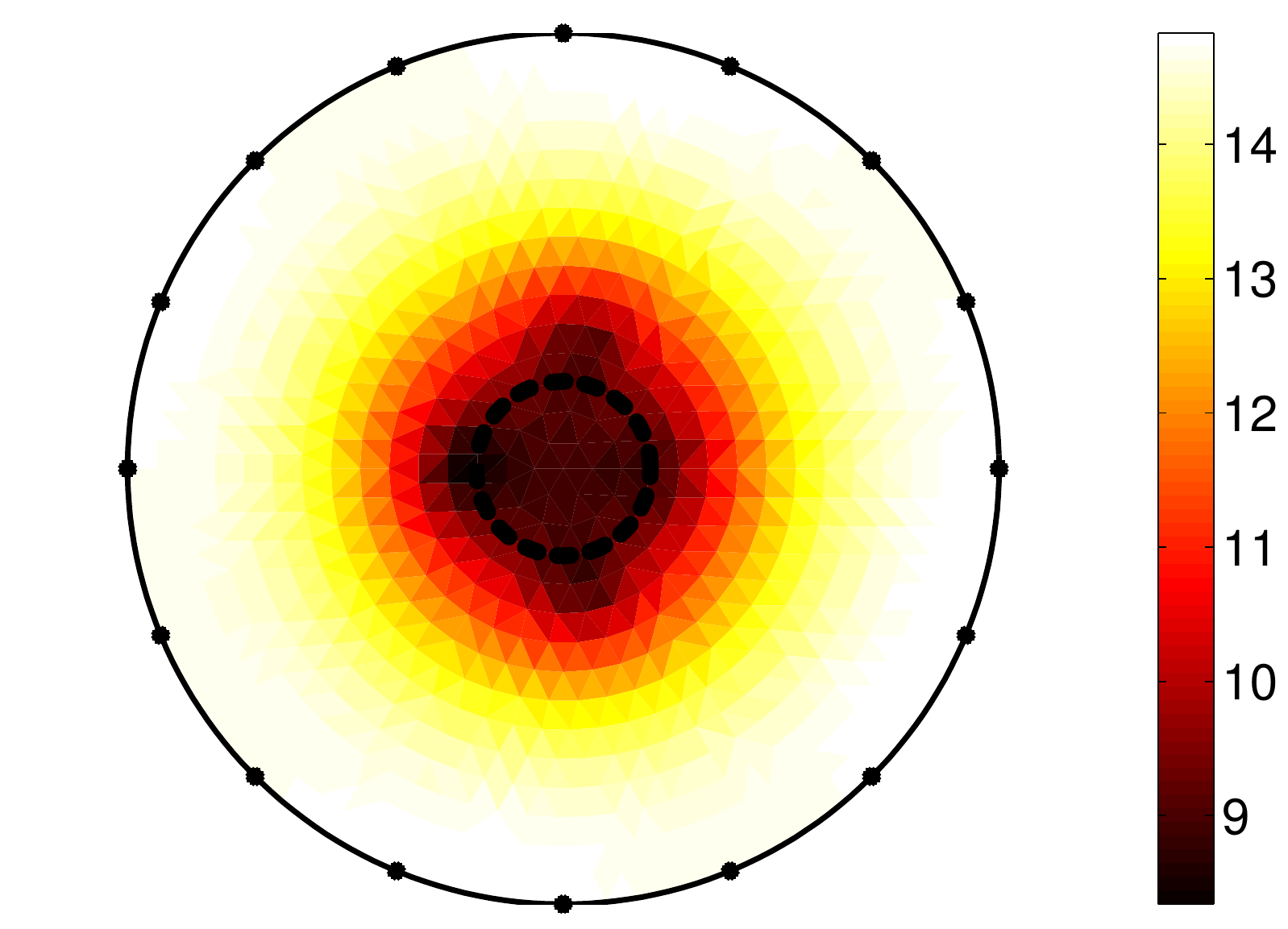}\\
\hline
\raisebox{4ex}{\footnotesize\begin{tabular}{c}
(c)
\end{tabular}}
&
\includegraphics[keepaspectratio=true,height=1.5cm]{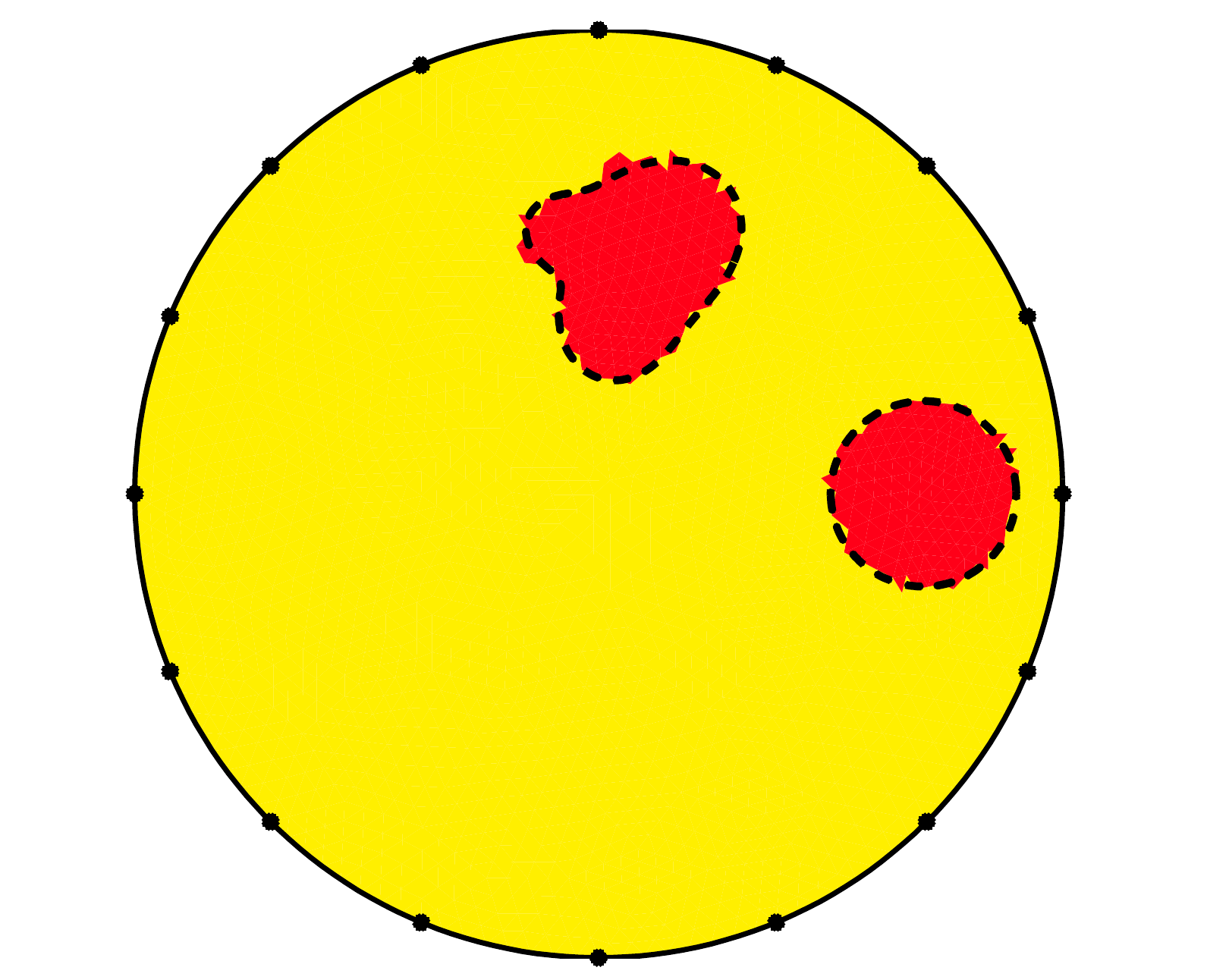}&
\includegraphics[keepaspectratio=true,height=1.5cm]{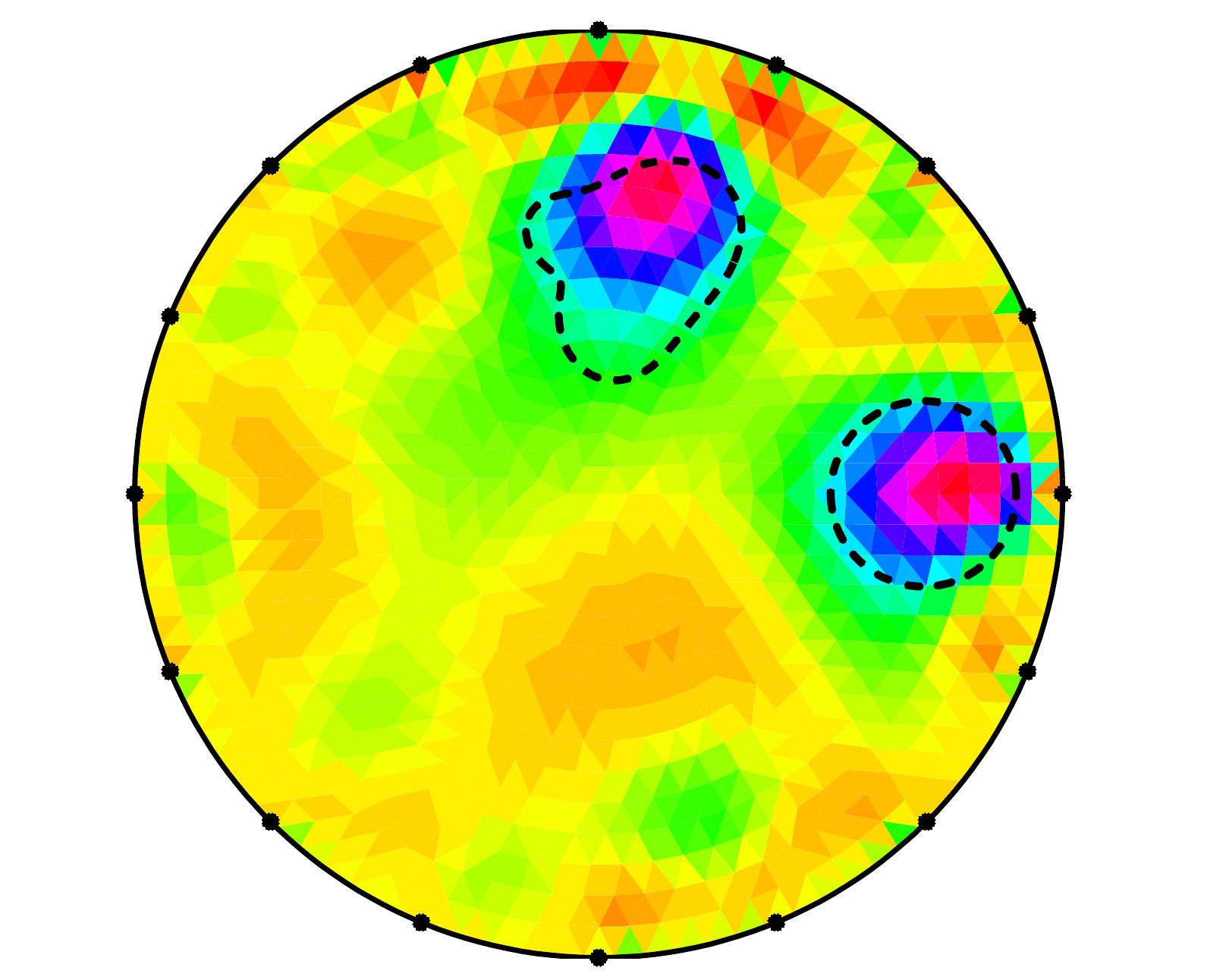}&
\includegraphics[keepaspectratio=true,height=1.5cm]{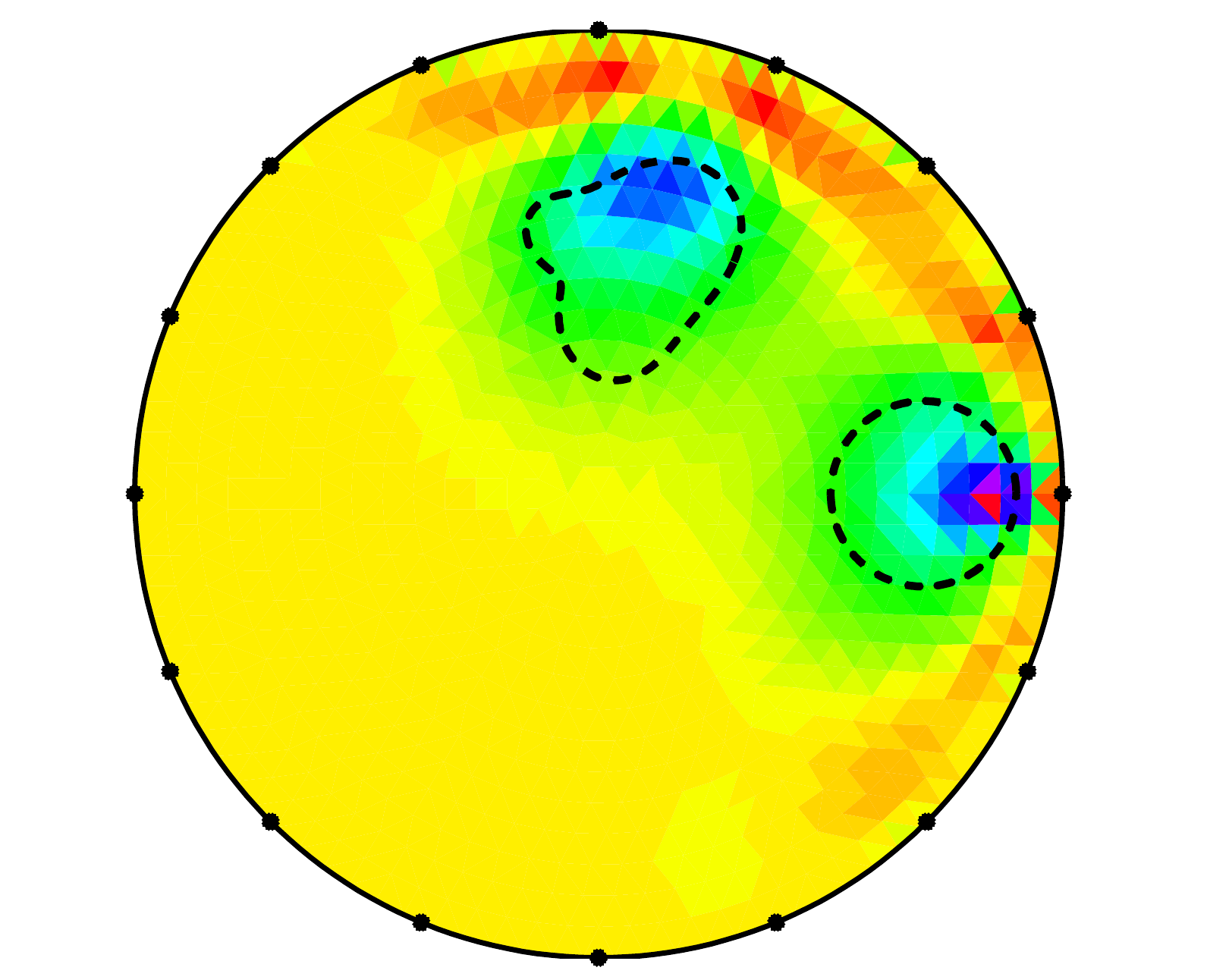}&
\includegraphics[keepaspectratio=true,height=1.5cm]{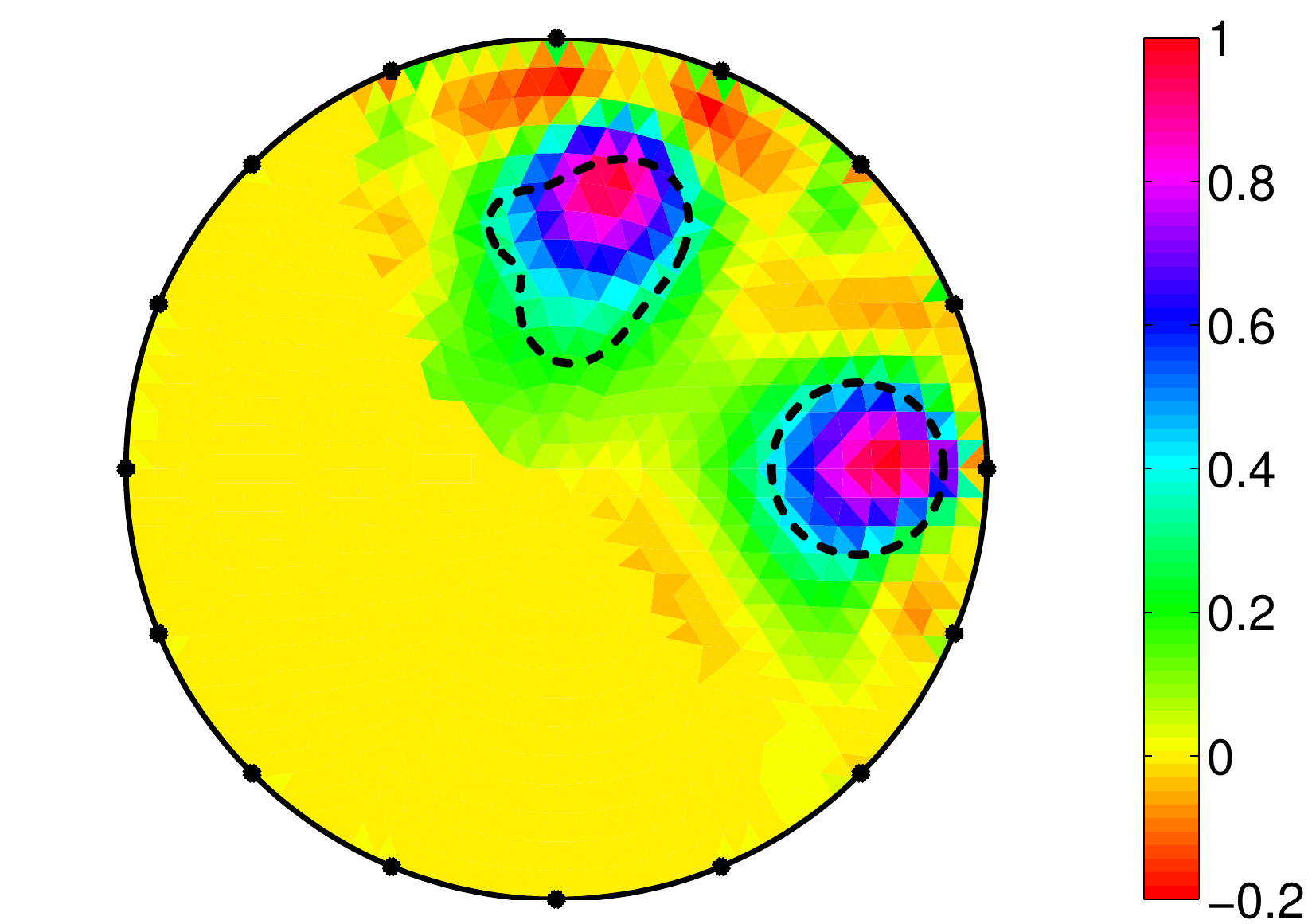}&
\includegraphics[keepaspectratio=true,height=1.5cm]{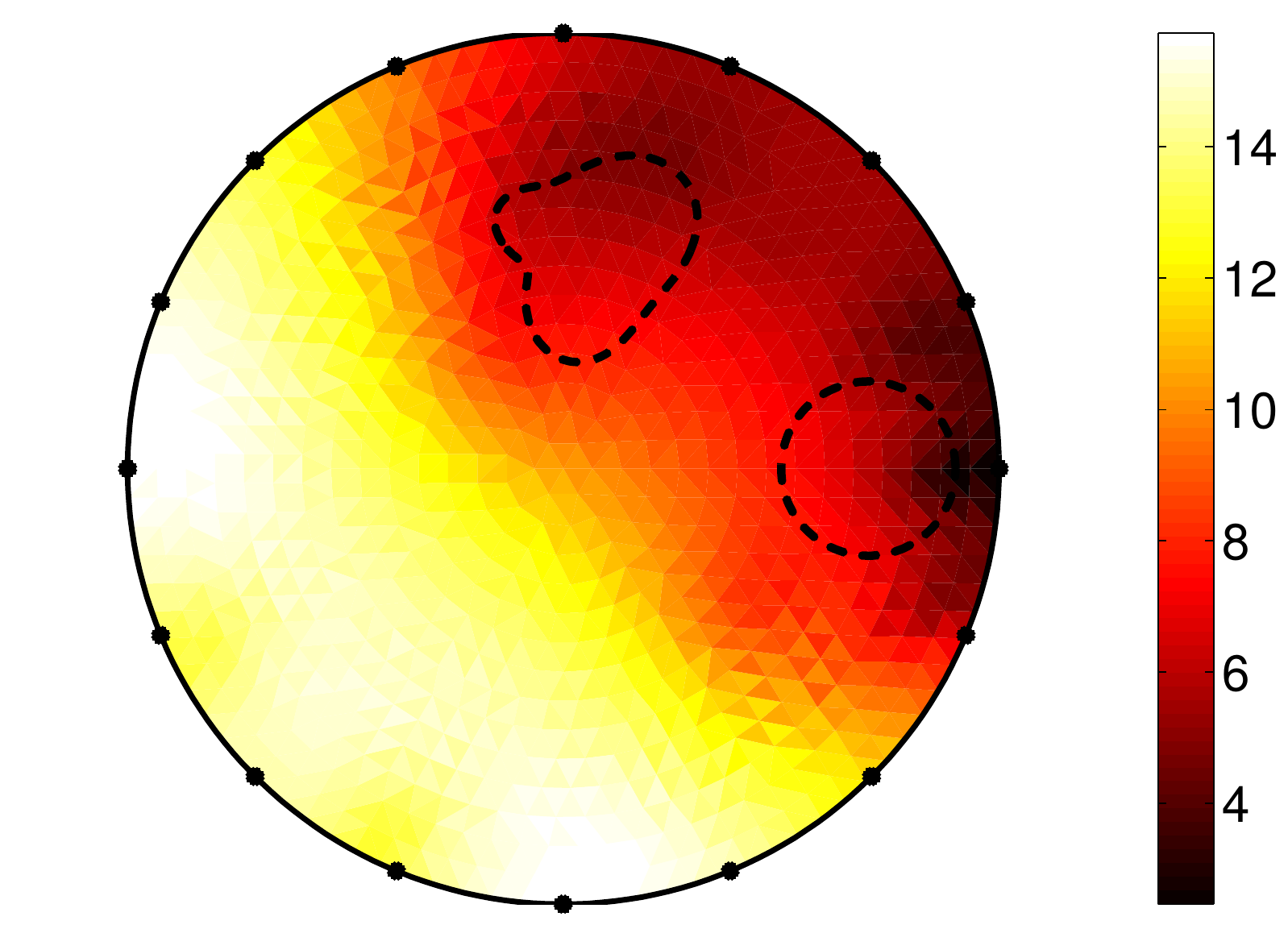}\\
\hline
\raisebox{4ex}{\footnotesize\begin{tabular}{c}
(d)
\end{tabular}}
&
\includegraphics[keepaspectratio=true,height=1.5cm]{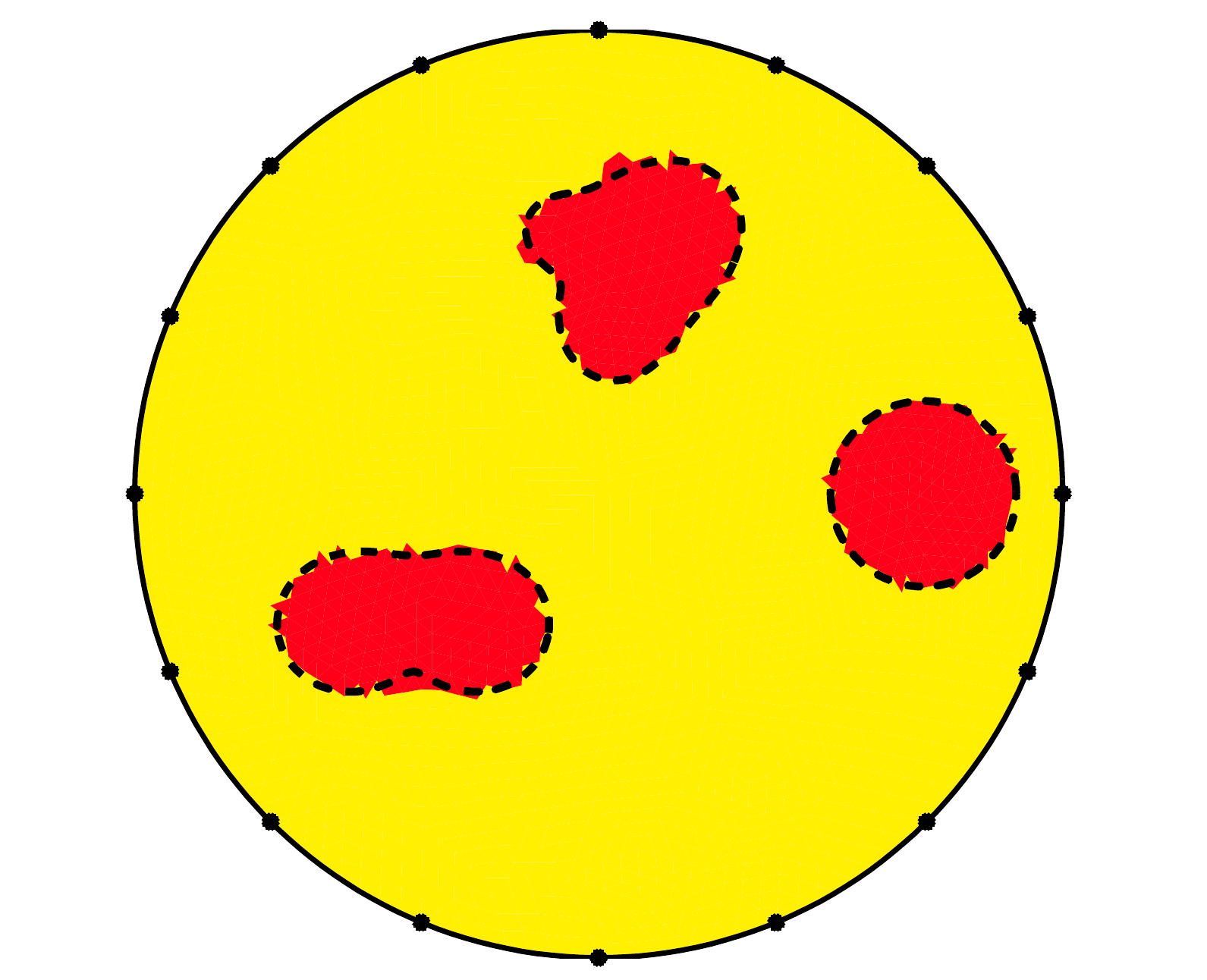}&
\includegraphics[keepaspectratio=true,height=1.5cm]{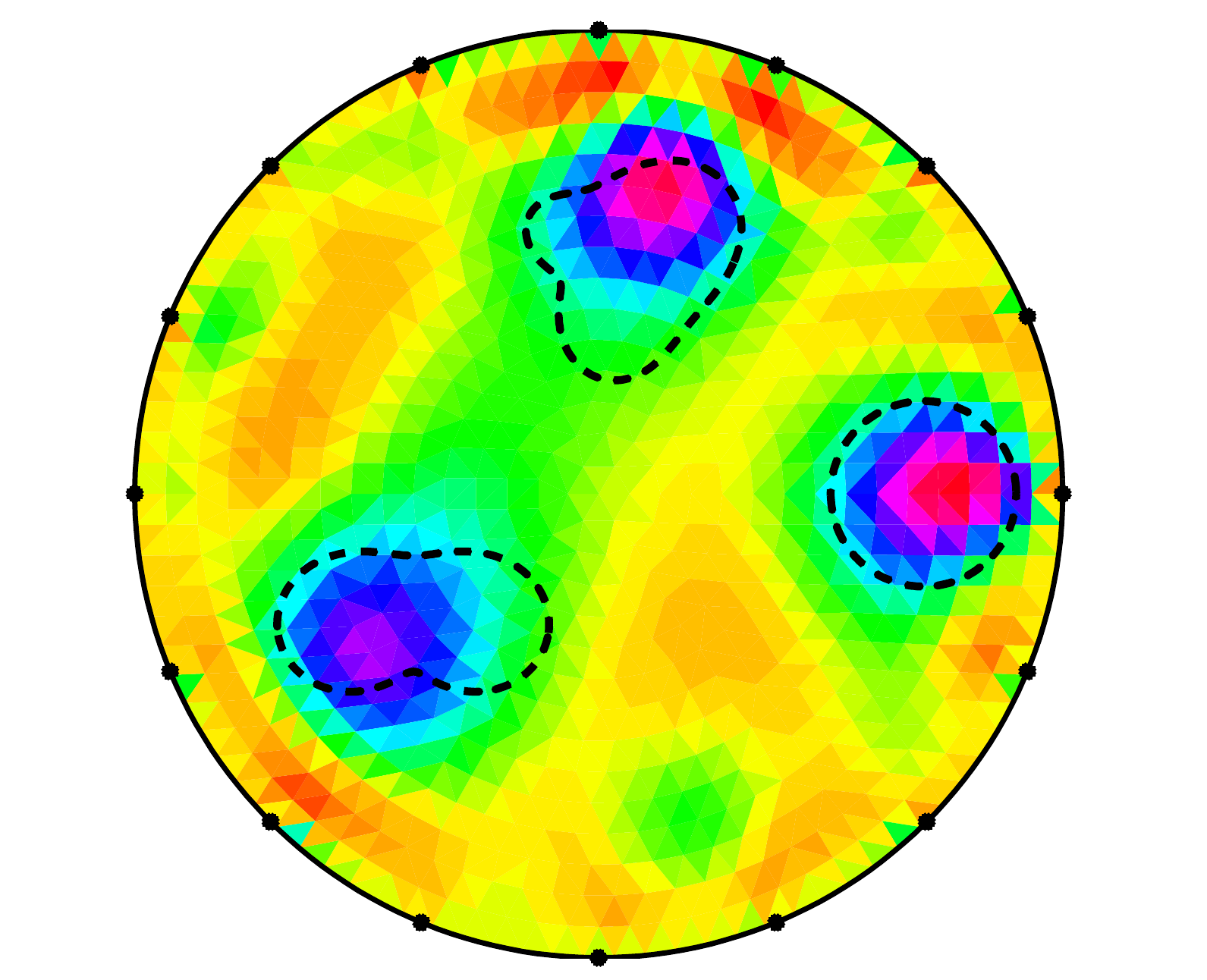}&
\includegraphics[keepaspectratio=true,height=1.5cm]{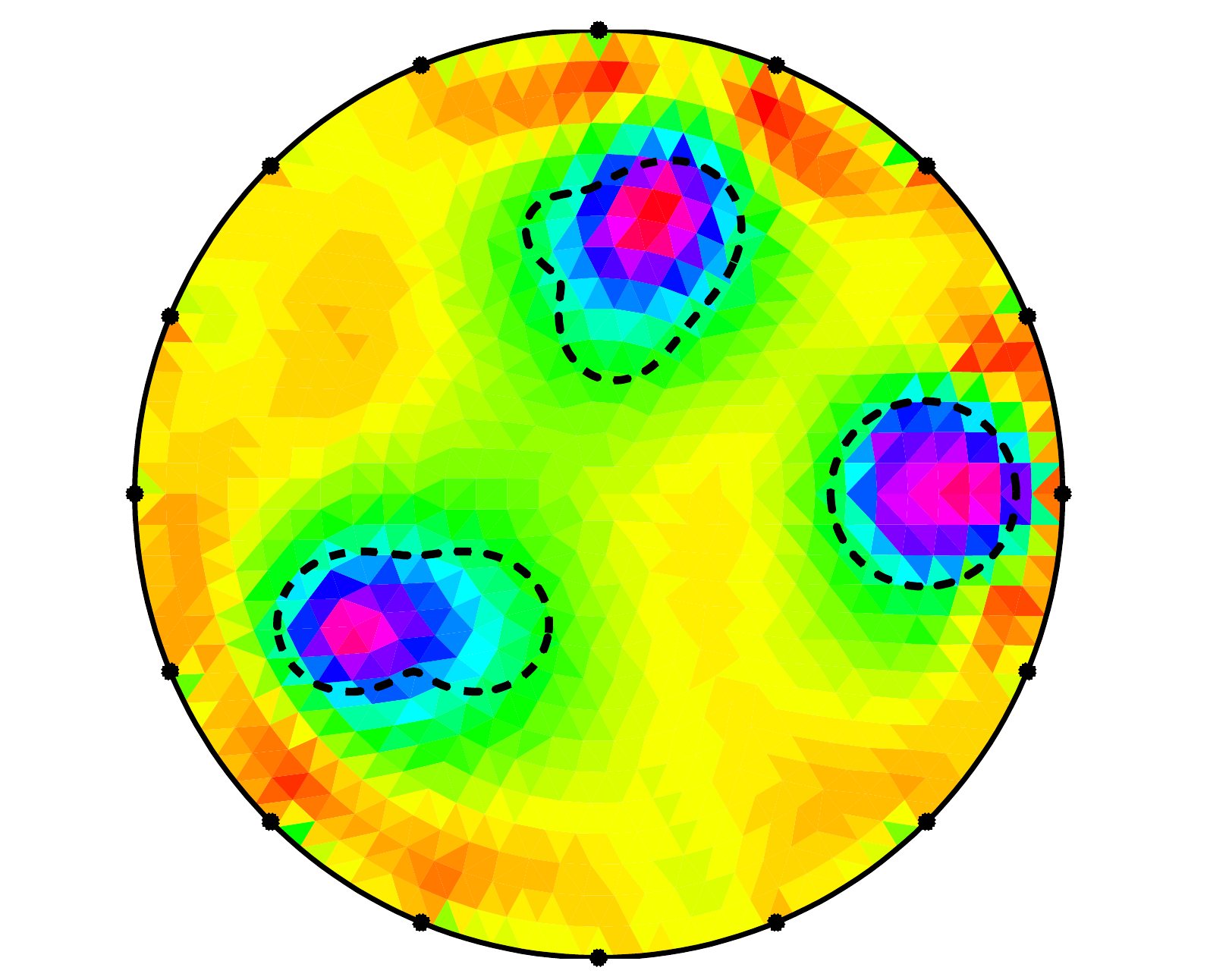}&
\includegraphics[keepaspectratio=true,height=1.5cm]{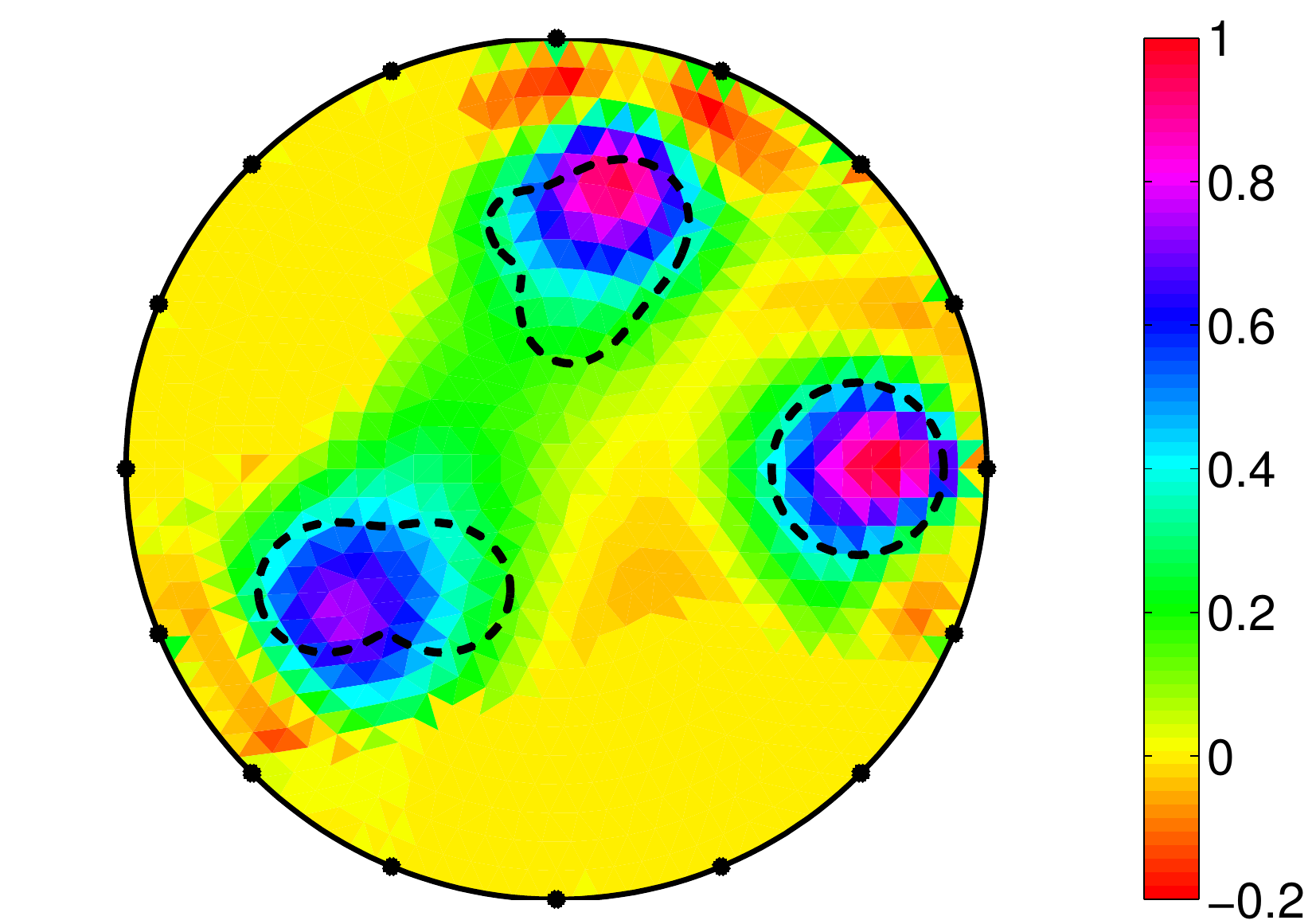}&
\includegraphics[keepaspectratio=true,height=1.5cm]{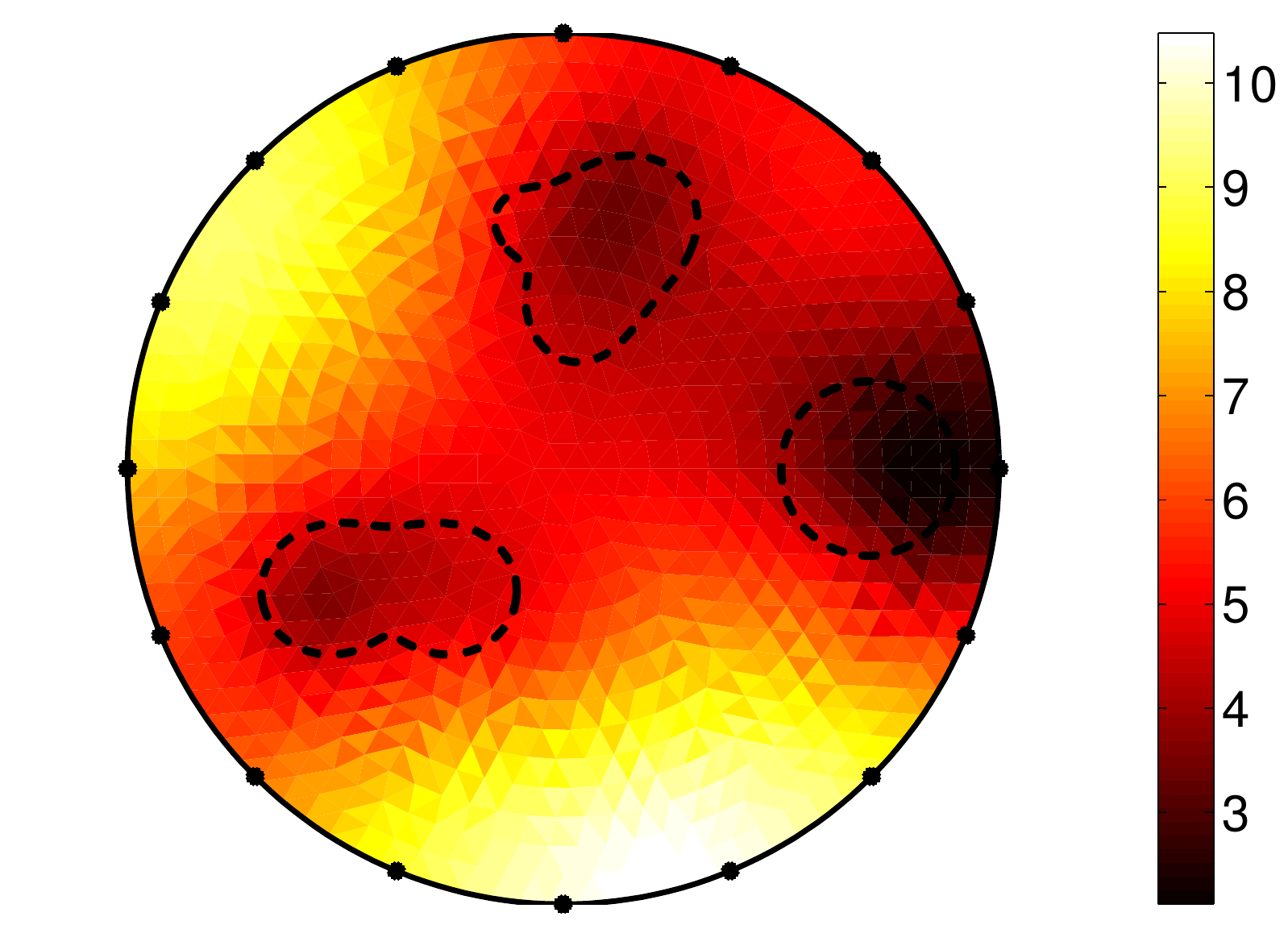}\\
\hline
\end{tabular}
\caption{ Reconstructed difference EIT images in circular domain. $\DS$: true difference image,
$\DS_{S}$: standard linearized method, cf.\ \eref{LM2},
$\DS_{B}$: naive combination of LM and S-FM, cf.\ \eref{regularization},
$\DS_{A}$: proposed combination of LM and S-FM, cf.\ \eref{recon},
$\mathbf{W1}$: S-FM alone, cf.\ \eref{SFMimage}.}
\label{recon_circle}
\end{figure}

\begin{figure}
\centering
\begin{tabular}{|c|cccc|c|}
\hline
Case &  $\DS$ & $\DS_{S}$ & $\DS_{B}$ &  $\DS_{A}$ &$\mathbf{W1}$ \\
\hline
\raisebox{4ex}{\footnotesize\begin{tabular}{c}
(e)
\end{tabular}}
&
\includegraphics[keepaspectratio=true,height=1.5cm]{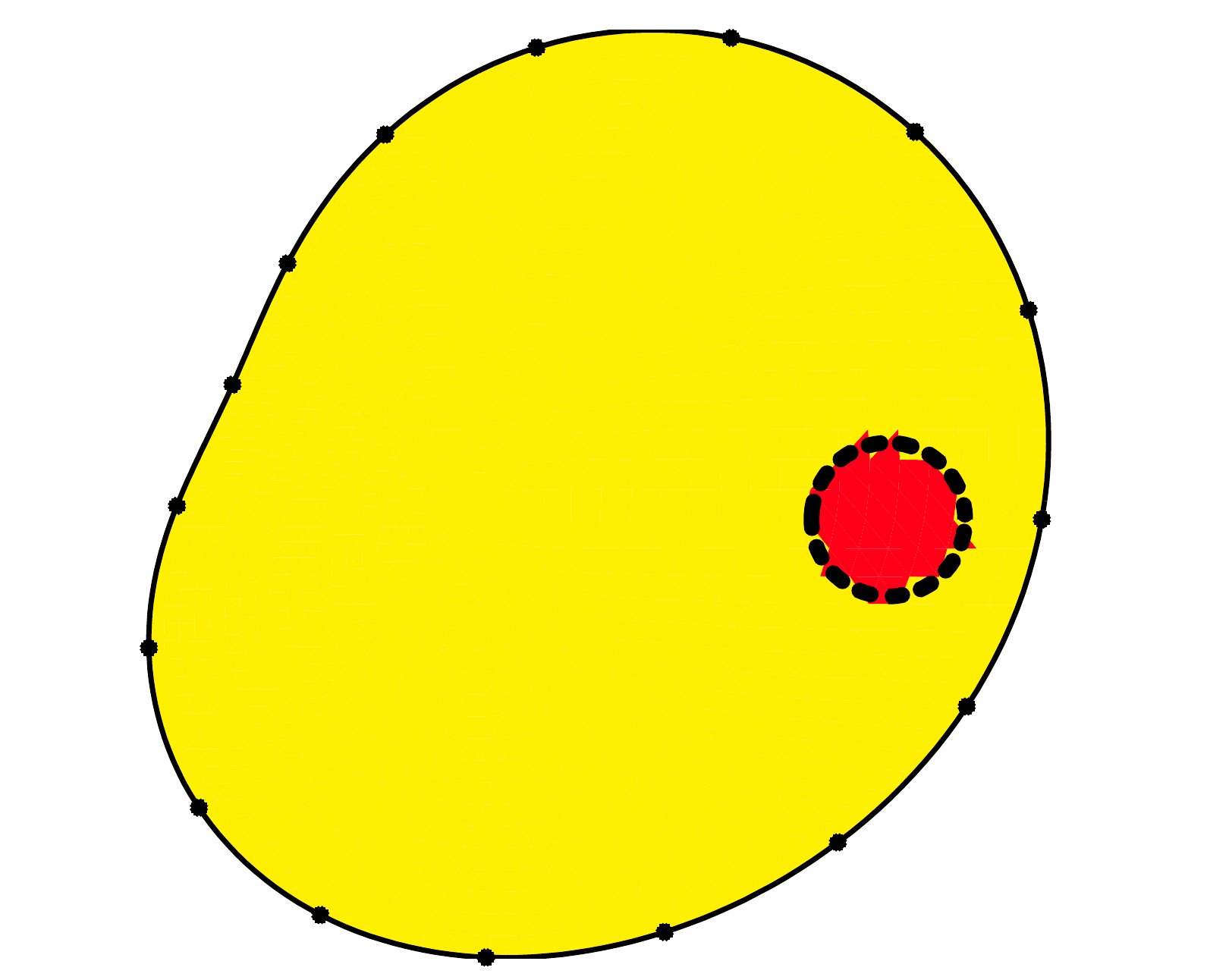}&
\includegraphics[keepaspectratio=true,height=1.5cm]{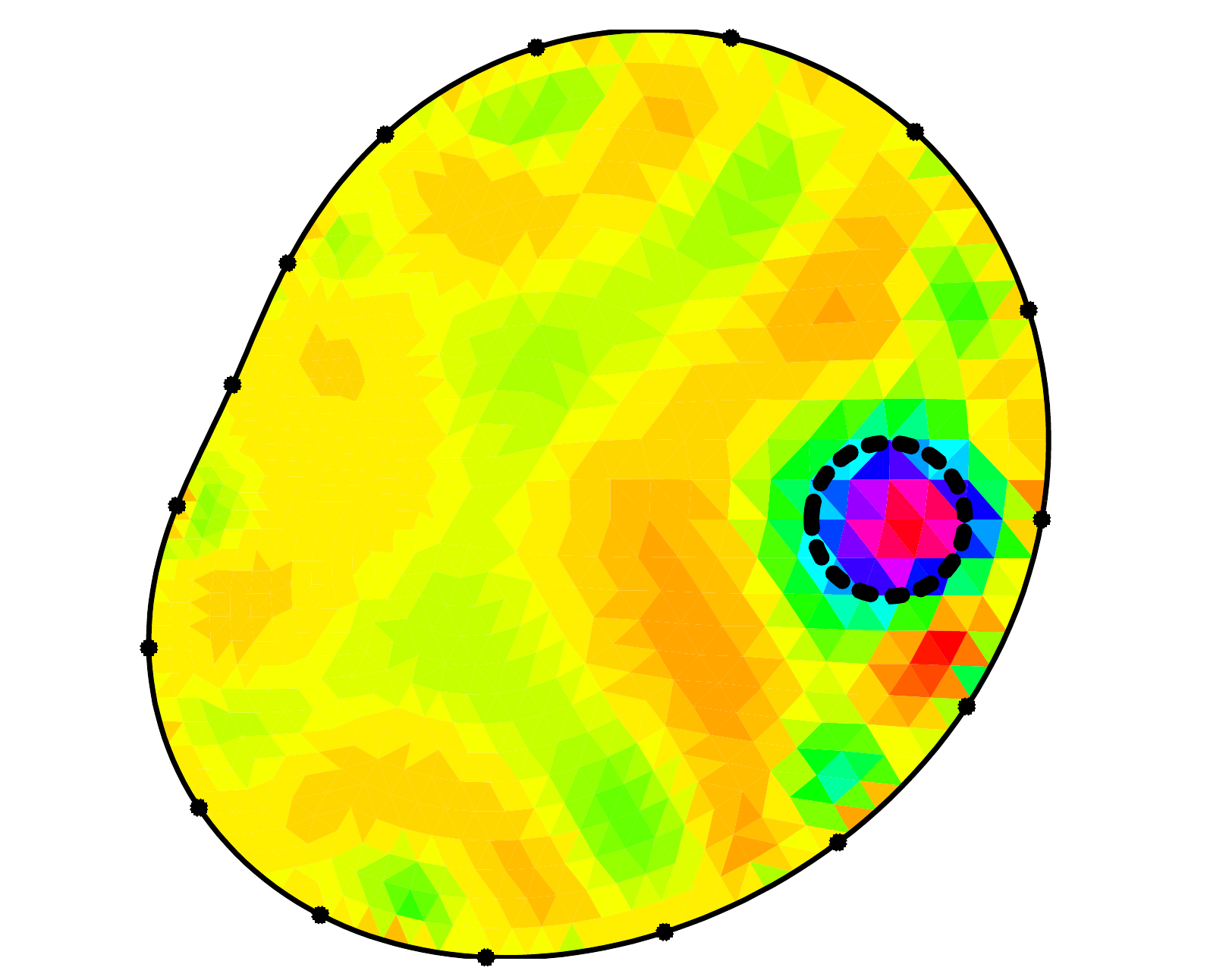}&
\includegraphics[keepaspectratio=true,height=1.5cm]{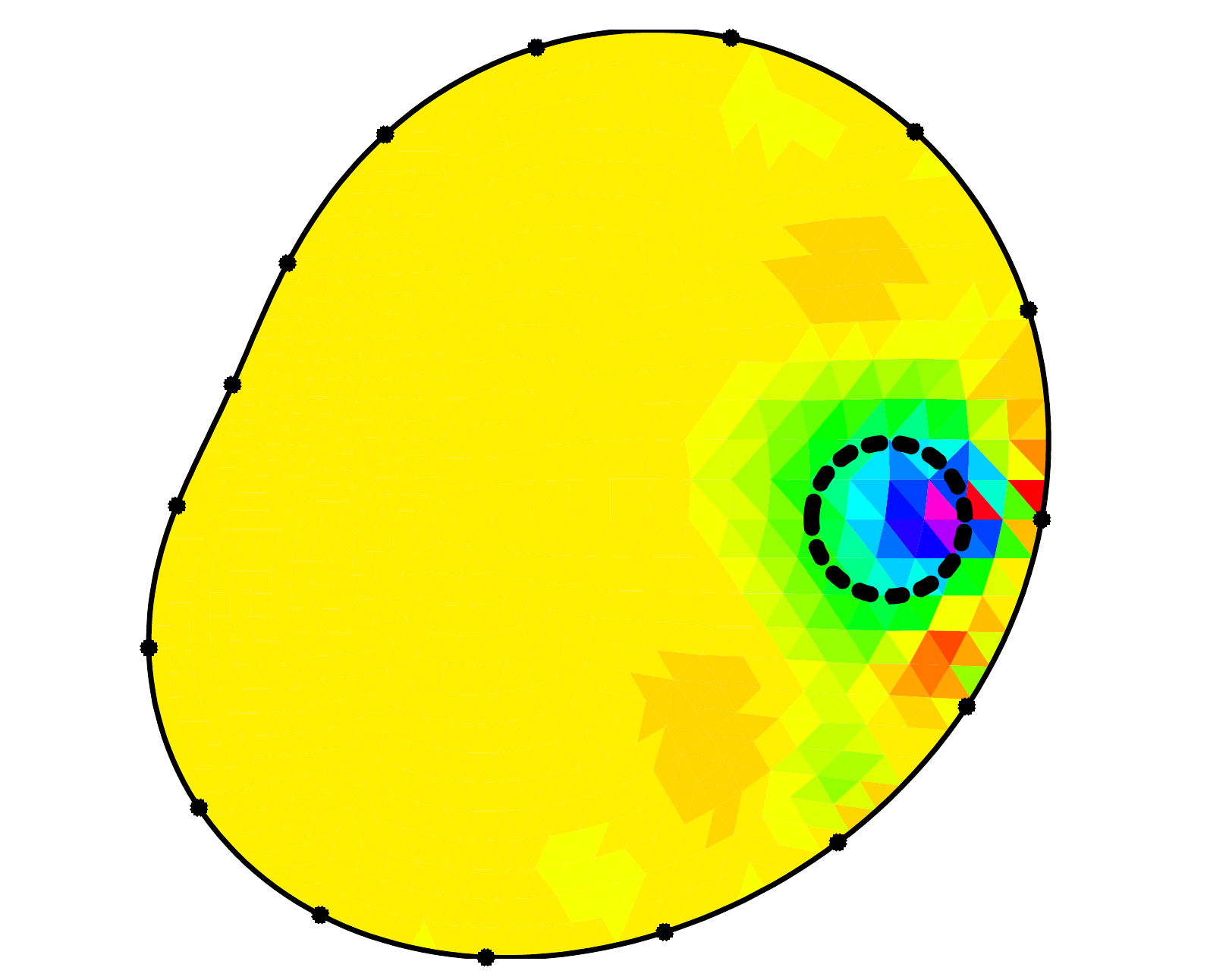}&
\includegraphics[keepaspectratio=true,height=1.5cm]{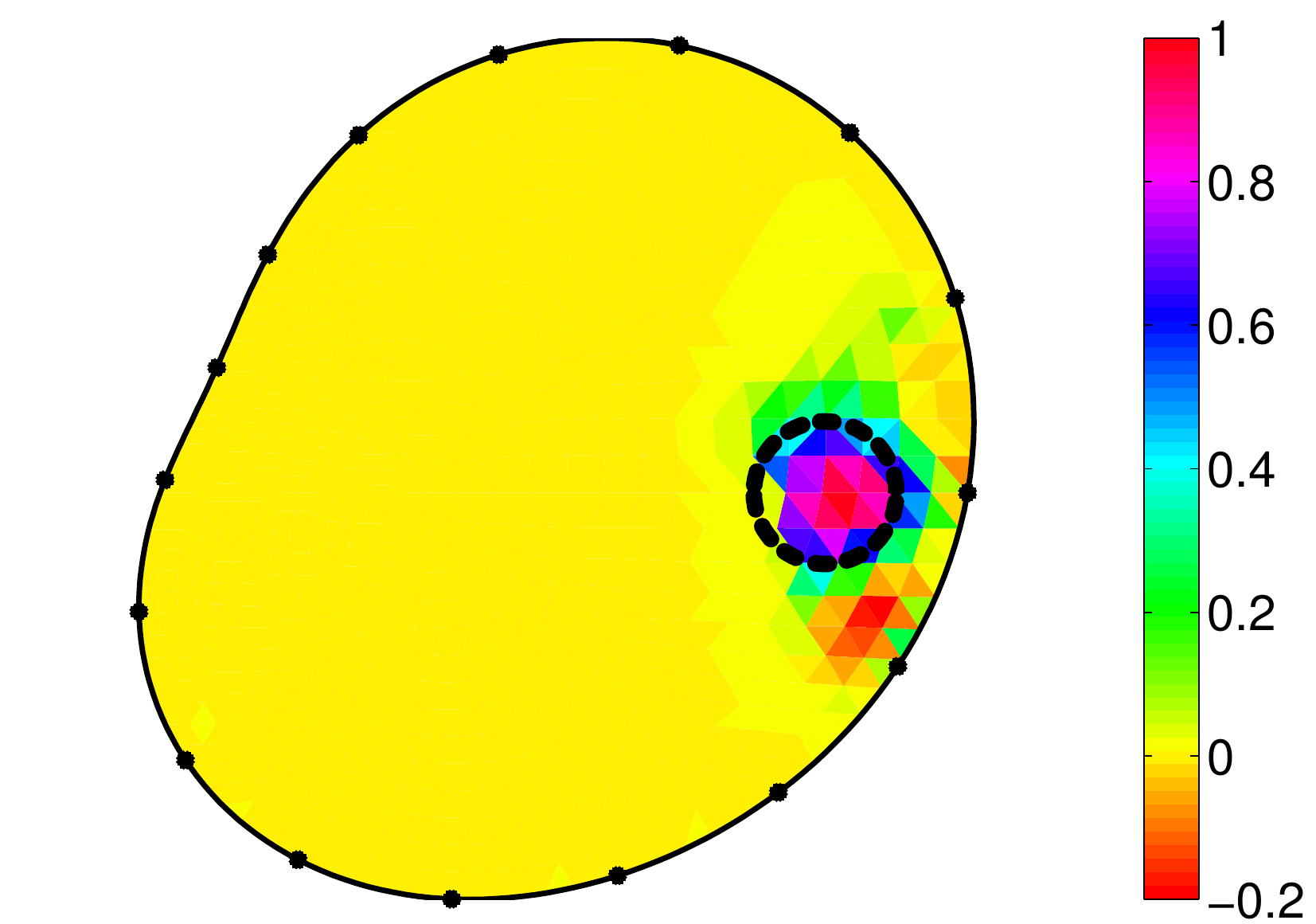}&
\includegraphics[keepaspectratio=true,height=1.5cm]{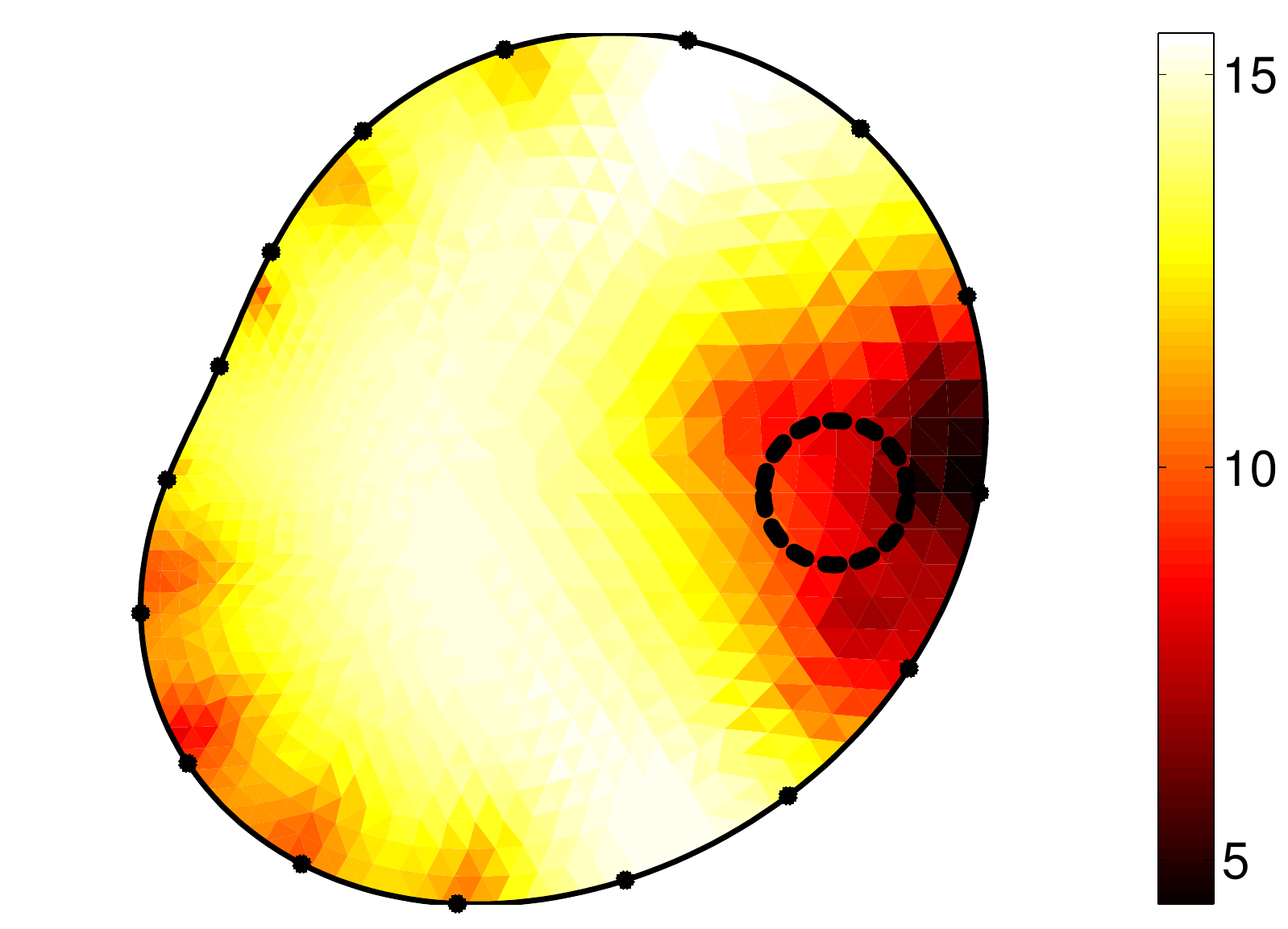}\\
\hline
\raisebox{4ex}{\footnotesize\begin{tabular}{c}
(f)
\end{tabular}}
&
\includegraphics[keepaspectratio=true,height=1.5cm]{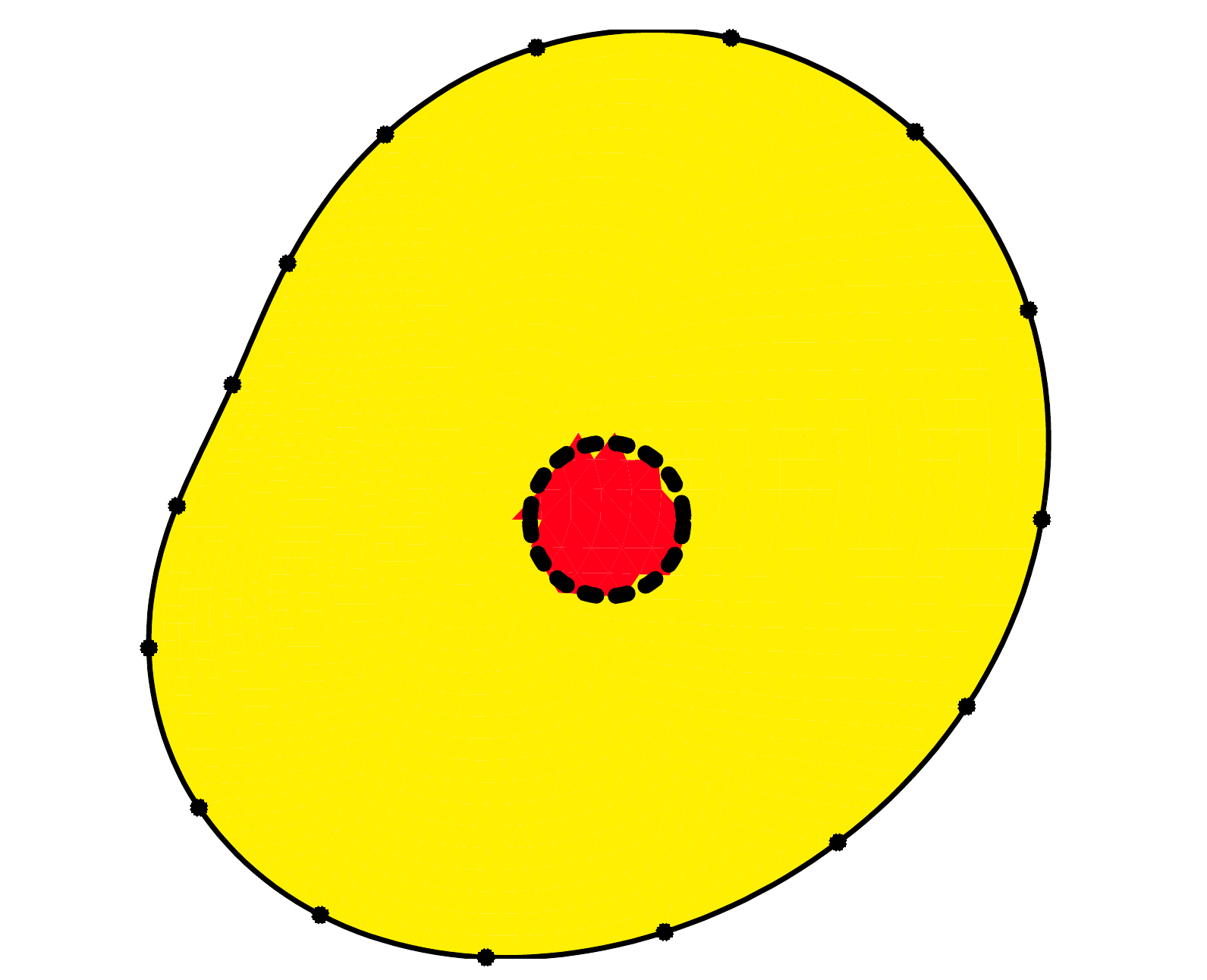}&
\includegraphics[keepaspectratio=true,height=1.5cm]{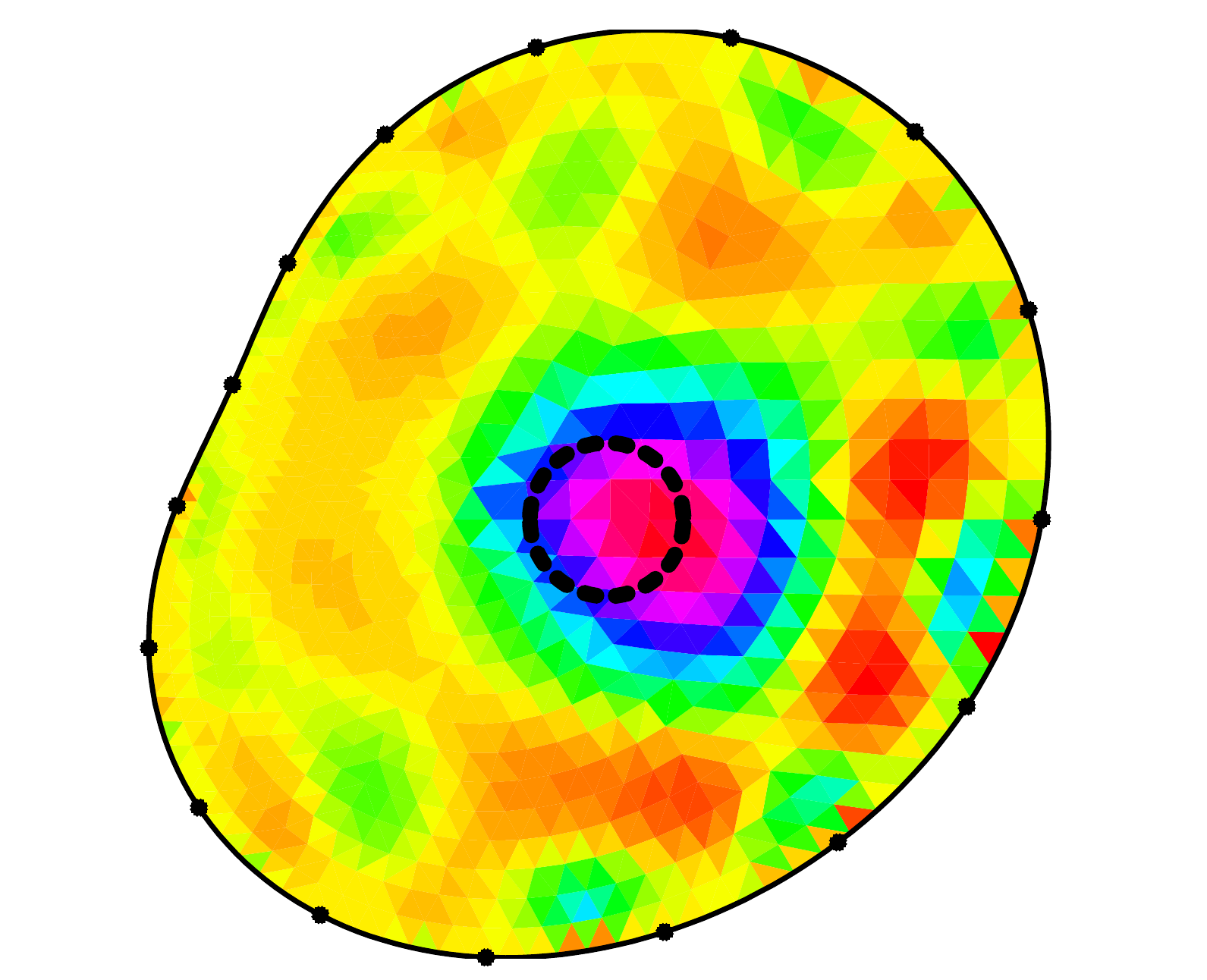}&
\includegraphics[keepaspectratio=true,height=1.5cm]{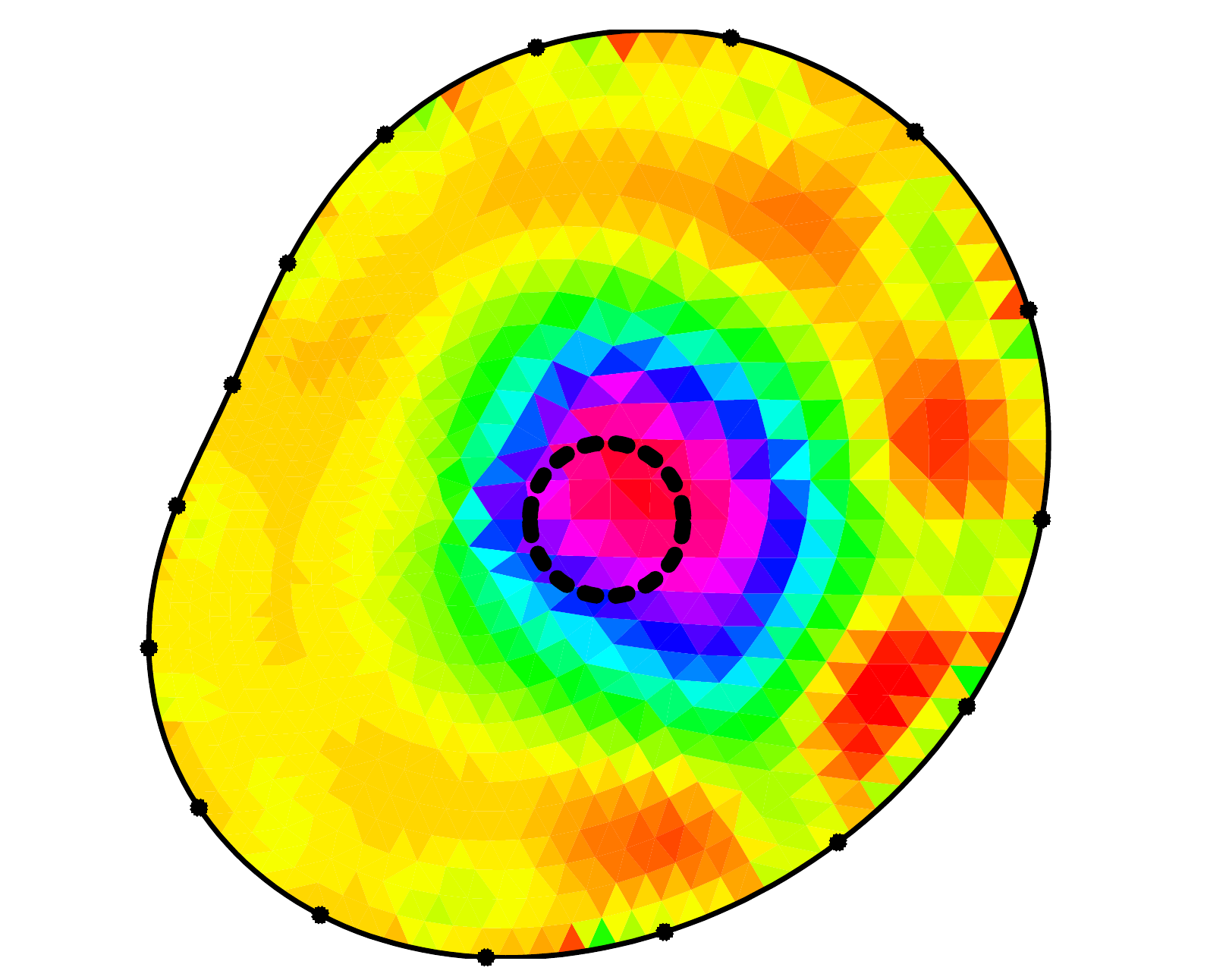}&
\includegraphics[keepaspectratio=true,height=1.5cm]{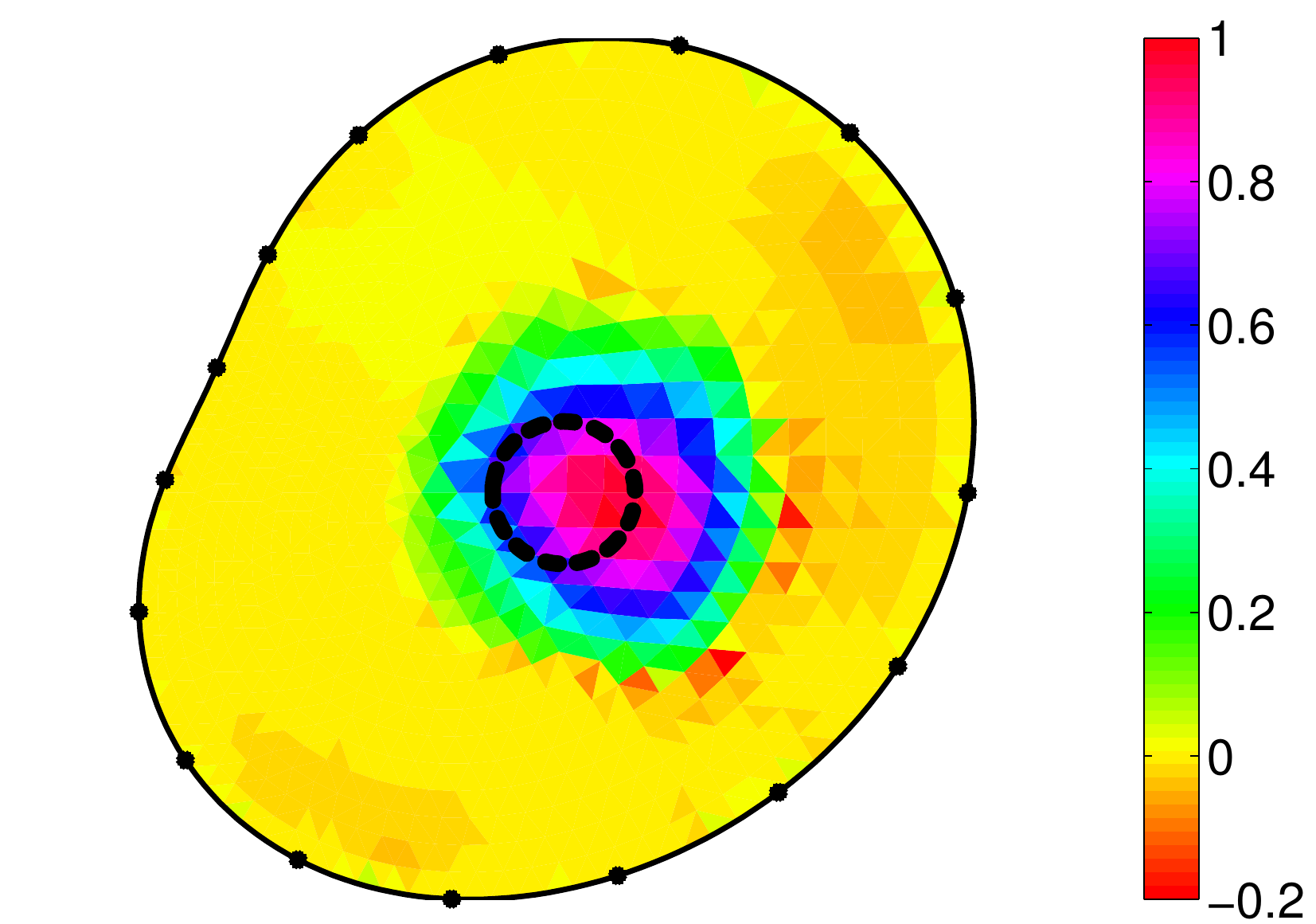}&
\includegraphics[keepaspectratio=true,height=1.5cm]{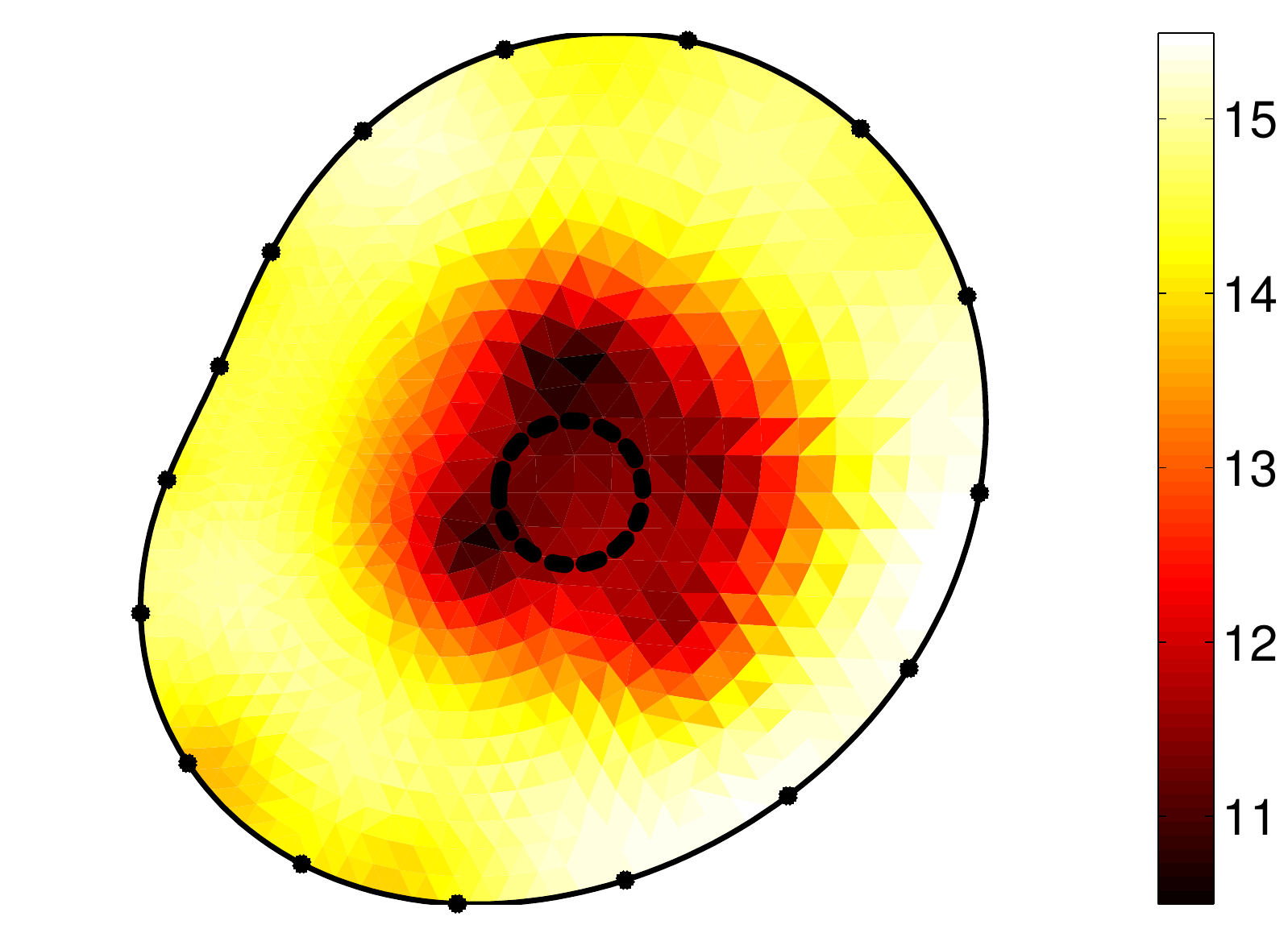}\\
\hline
\raisebox{4ex}{\footnotesize\begin{tabular}{c}
(g)
\end{tabular}}
&
\includegraphics[keepaspectratio=true,height=1.5cm]{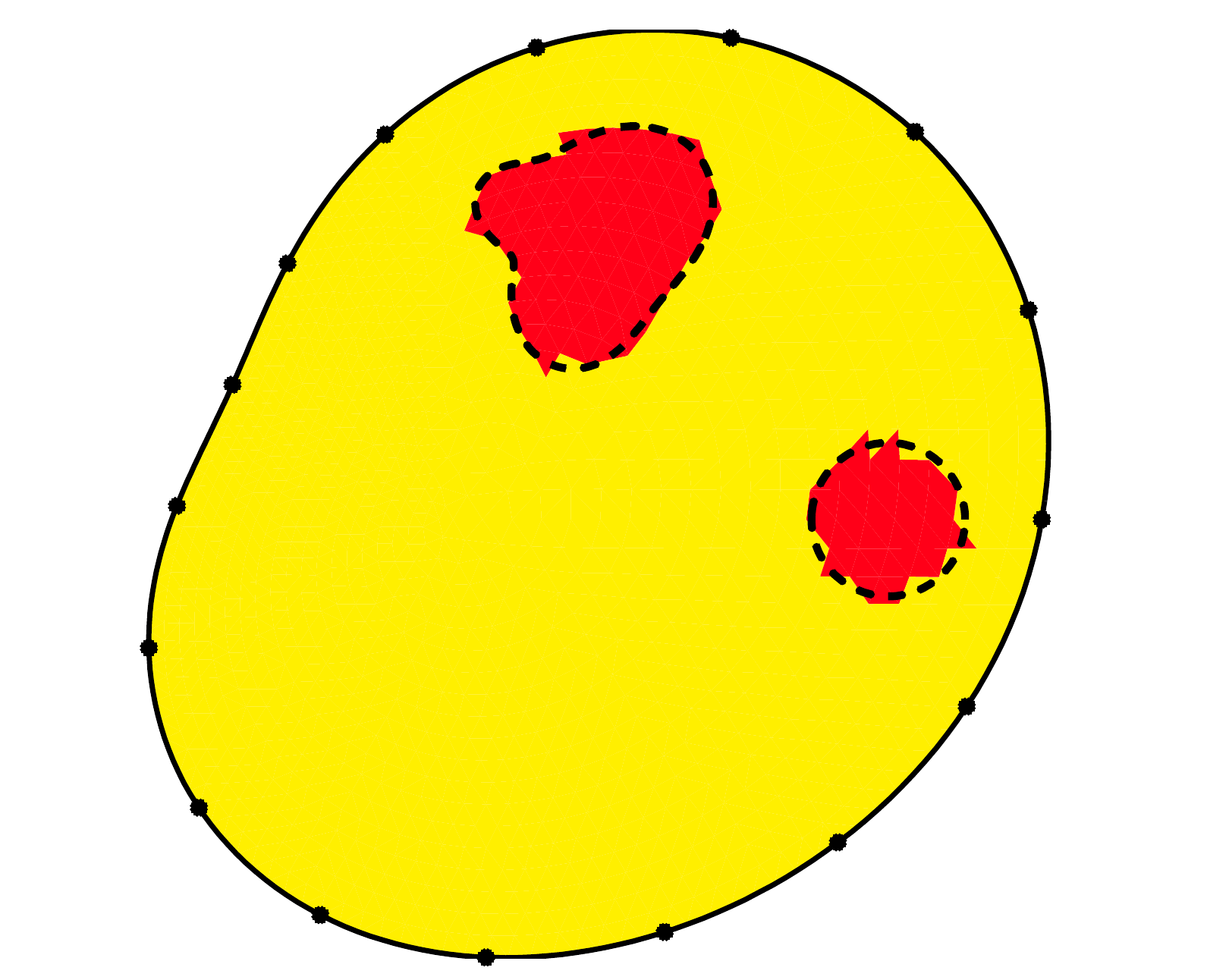}&
\includegraphics[keepaspectratio=true,height=1.5cm]{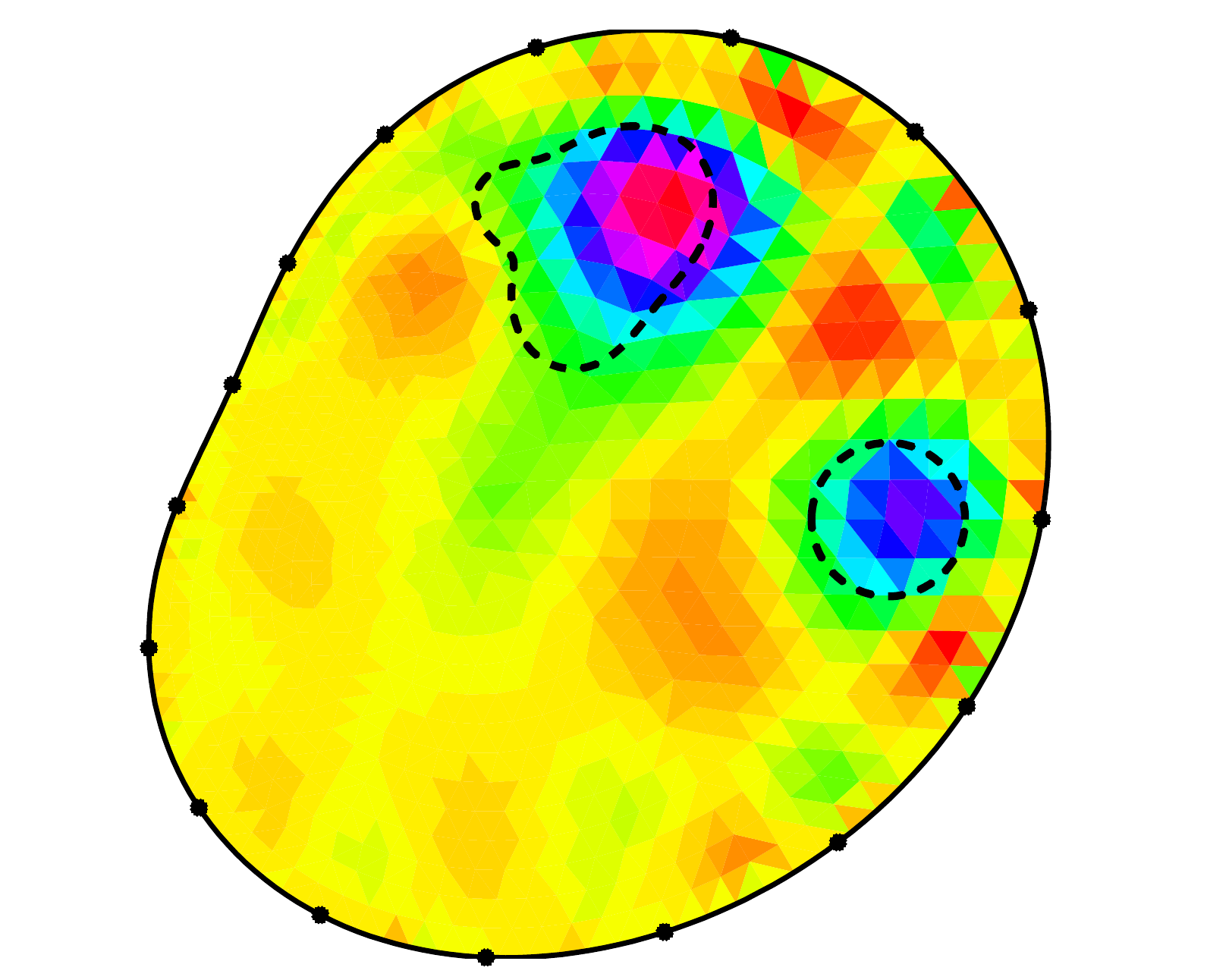}&
\includegraphics[keepaspectratio=true,height=1.5cm]{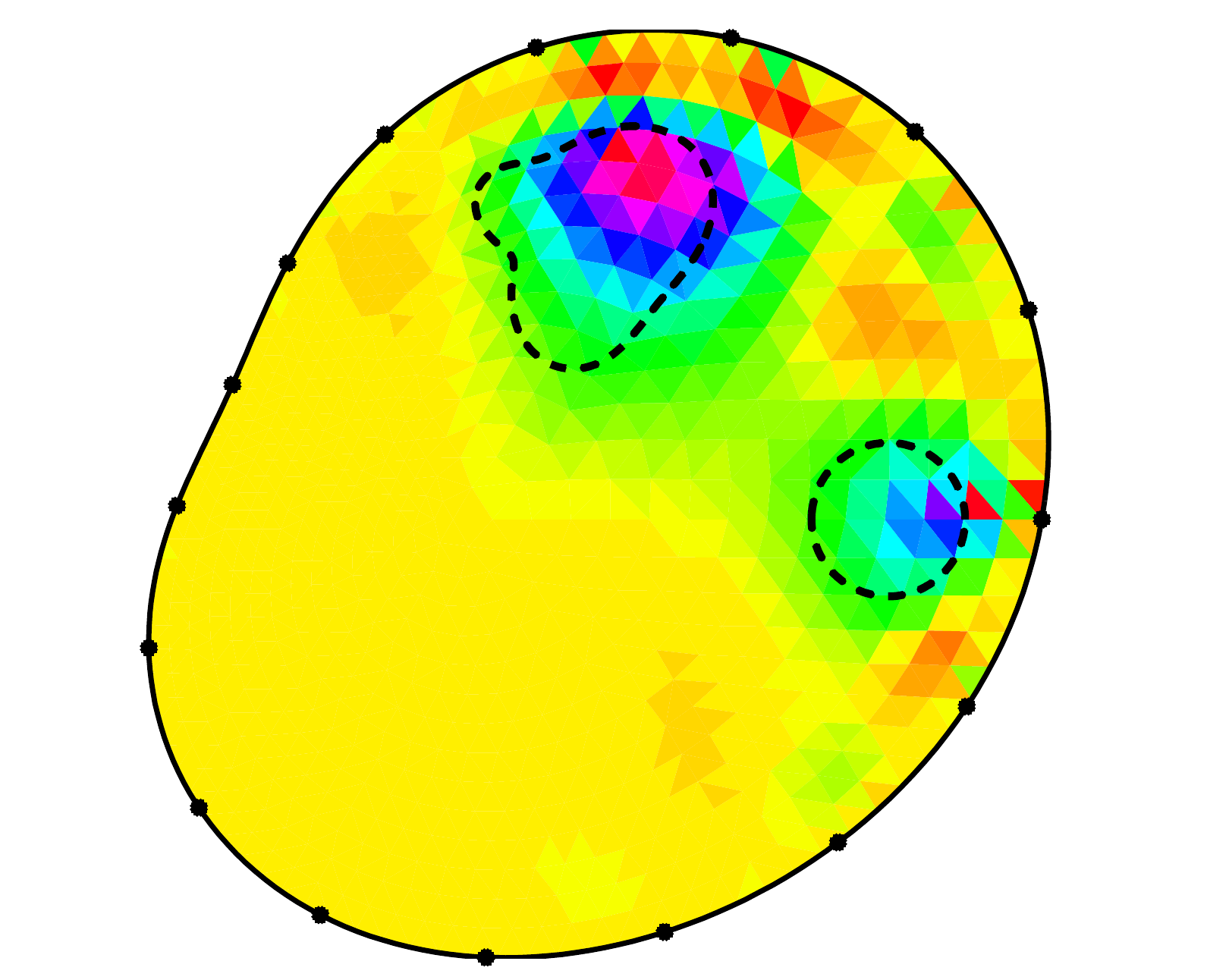}&
\includegraphics[keepaspectratio=true,height=1.5cm]{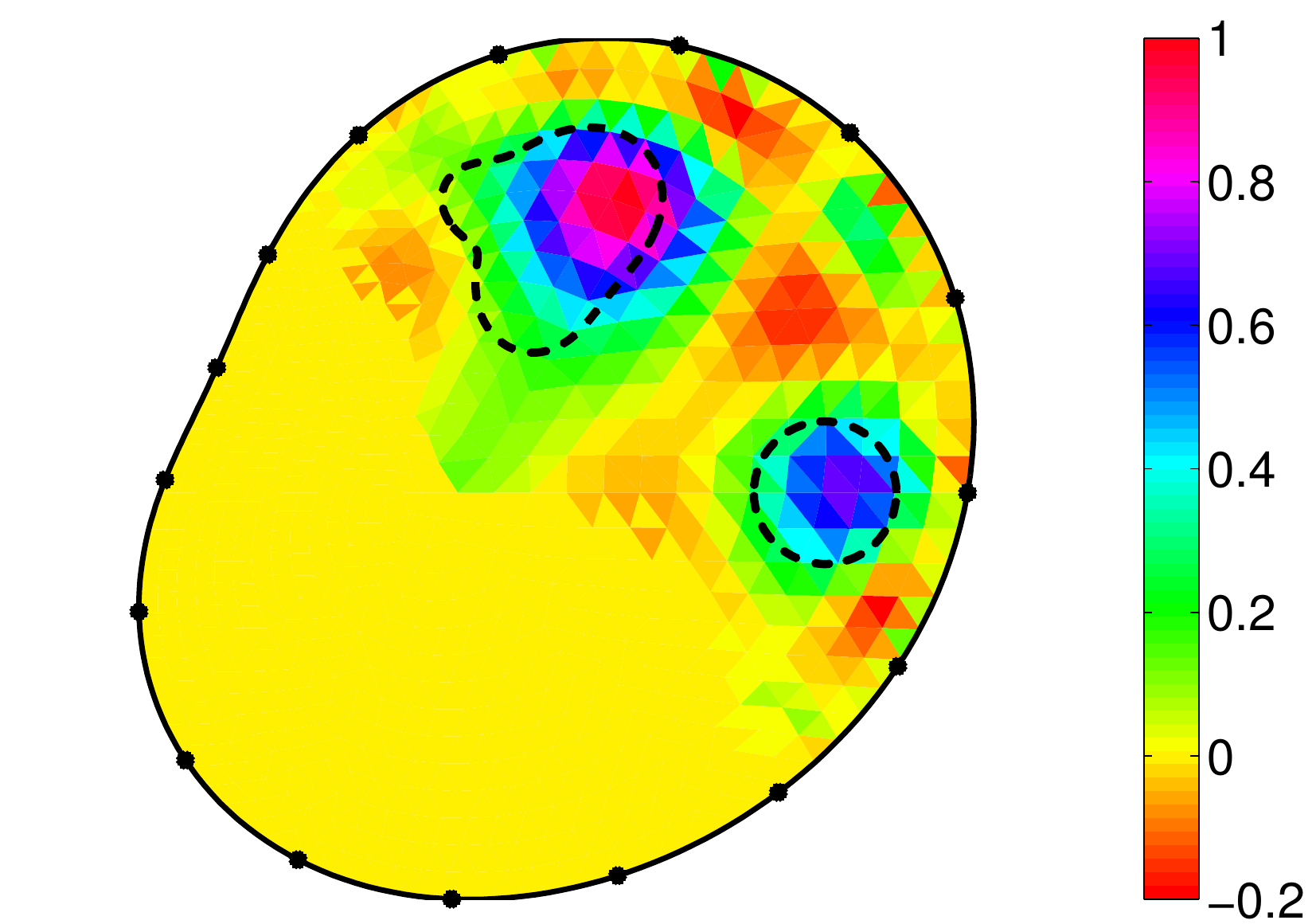}&
\includegraphics[keepaspectratio=true,height=1.5cm]{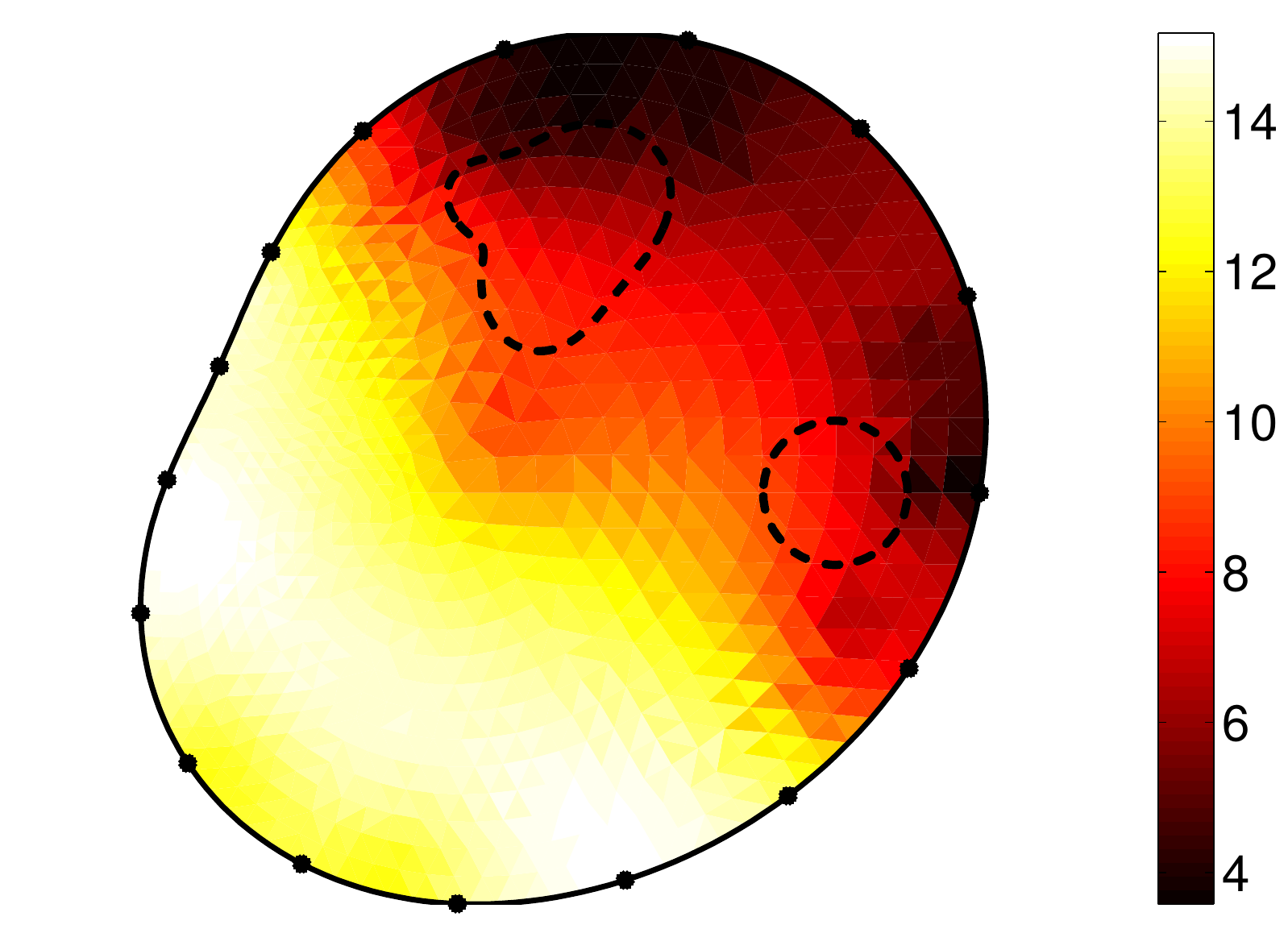}\\
\hline
\raisebox{4ex}{\footnotesize\begin{tabular}{c}
(h)
\end{tabular}}
&
\includegraphics[keepaspectratio=true,height=1.5cm]{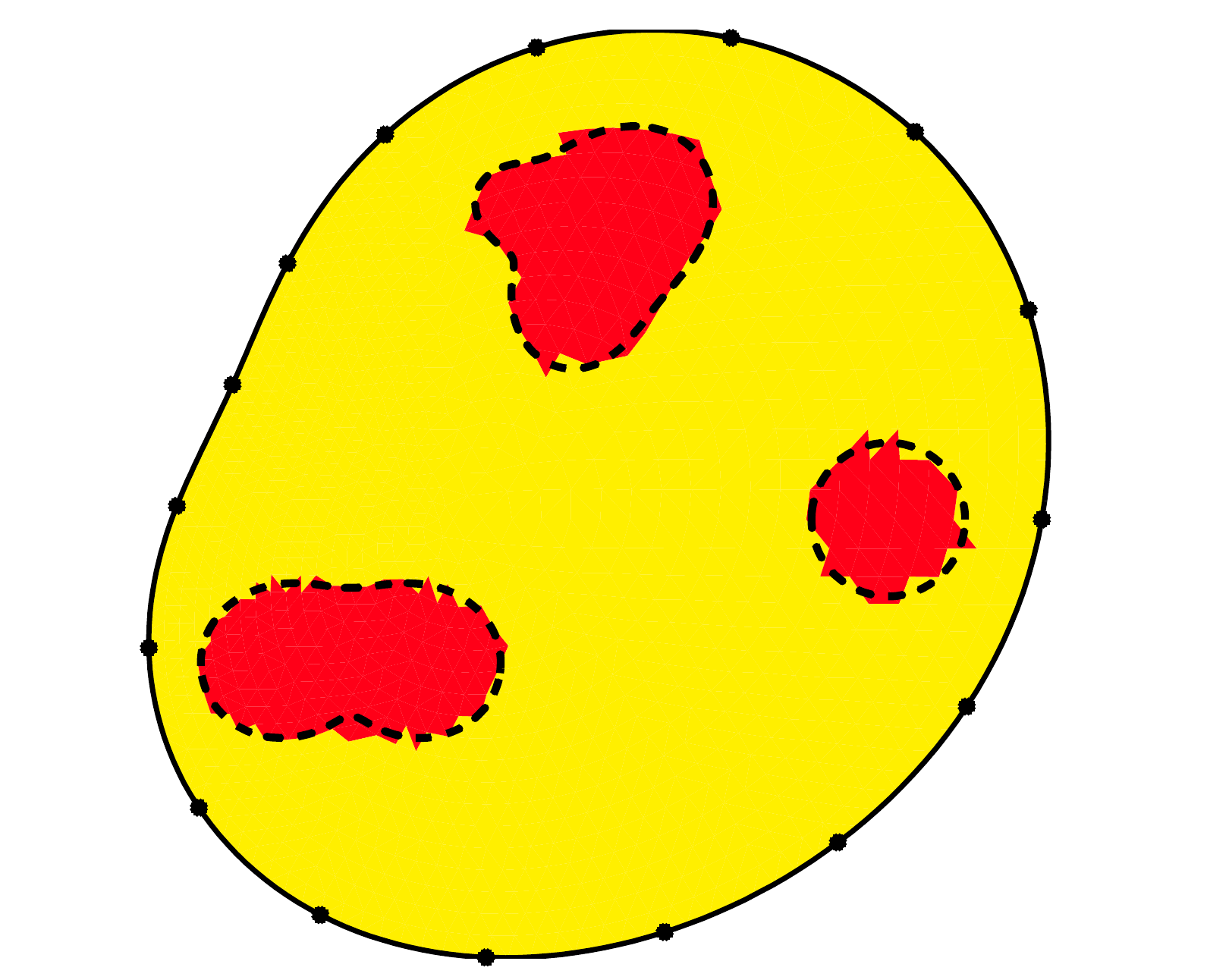}&
\includegraphics[keepaspectratio=true,height=1.5cm]{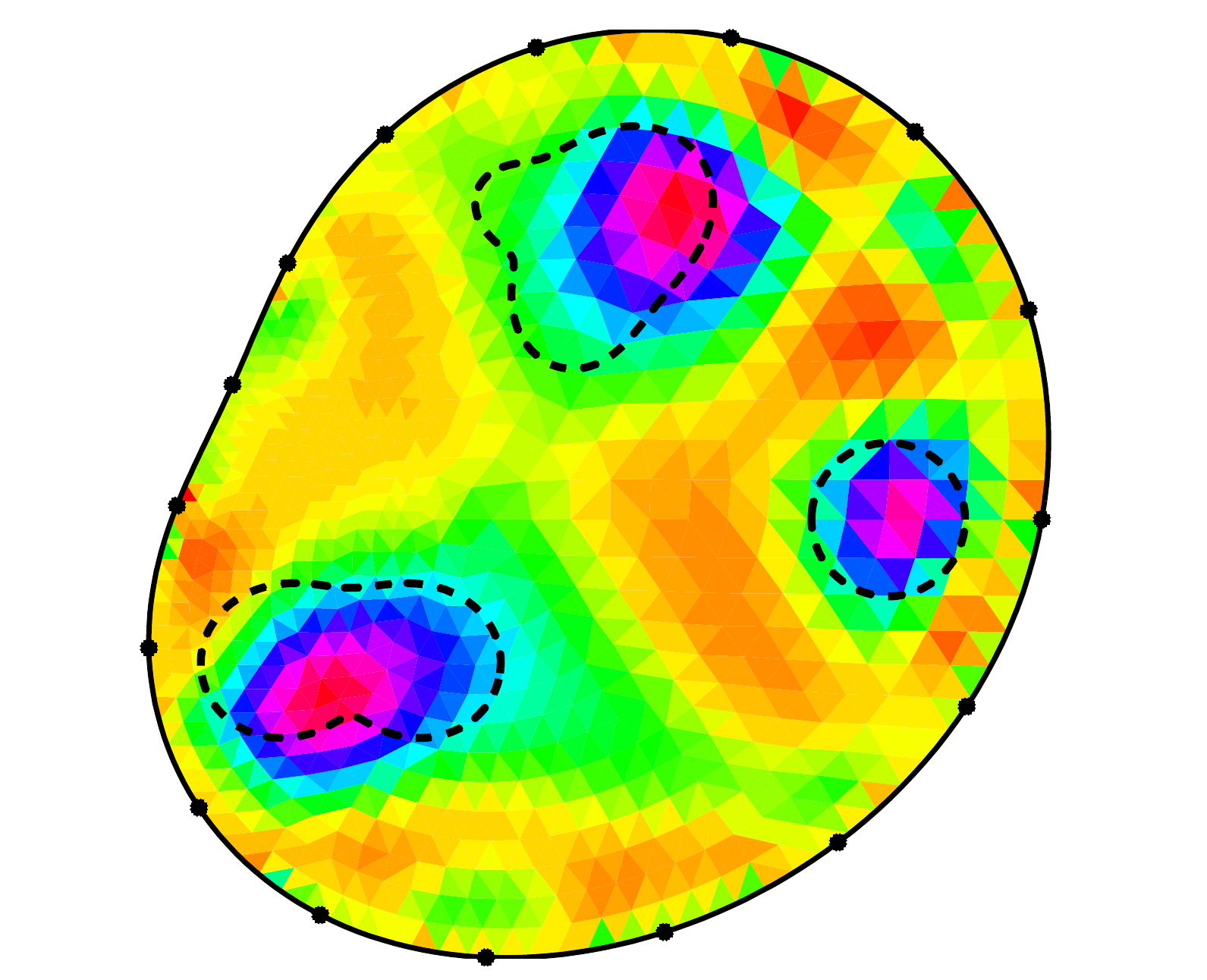}&
\includegraphics[keepaspectratio=true,height=1.5cm]{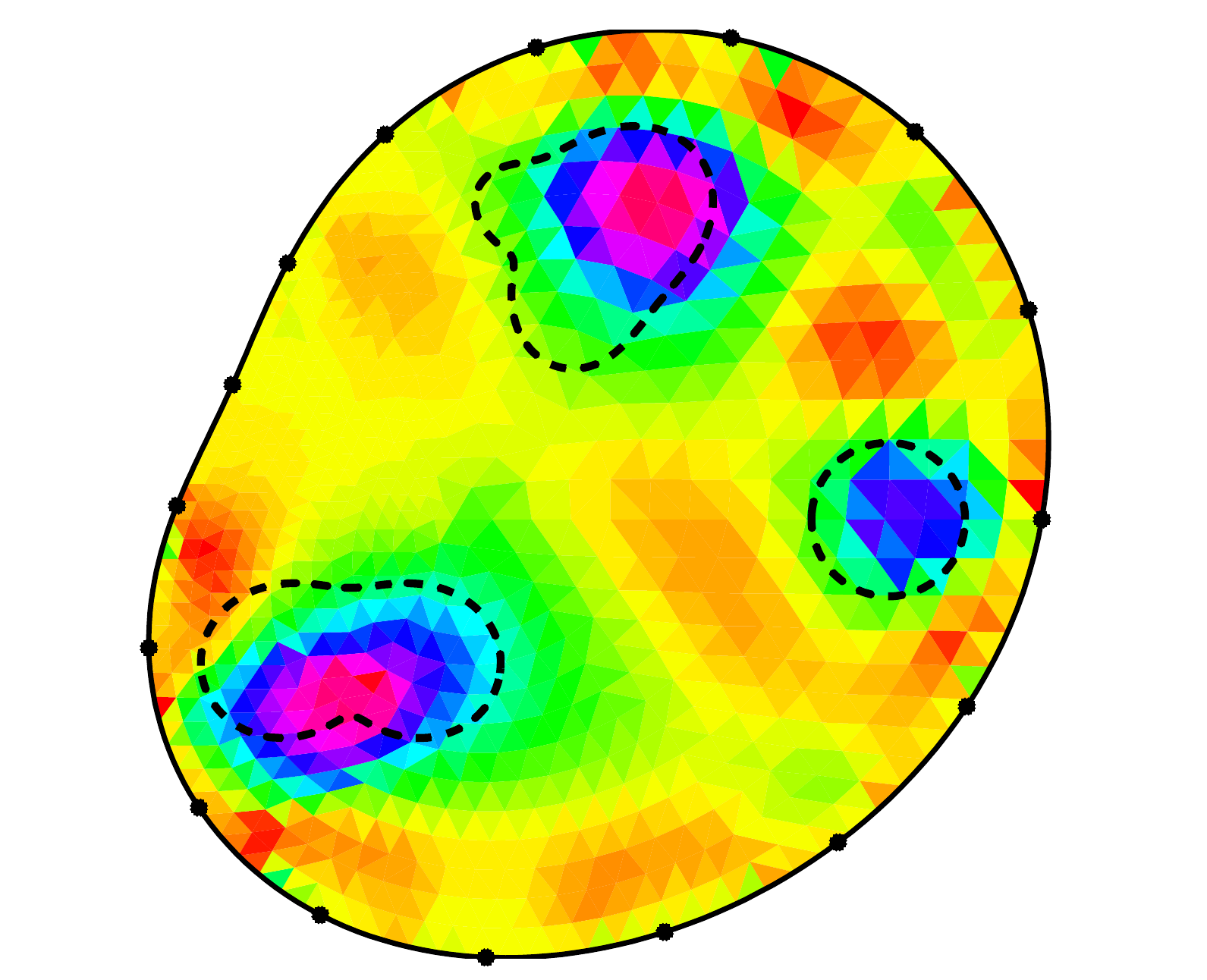}&
\includegraphics[keepaspectratio=true,height=1.5cm]{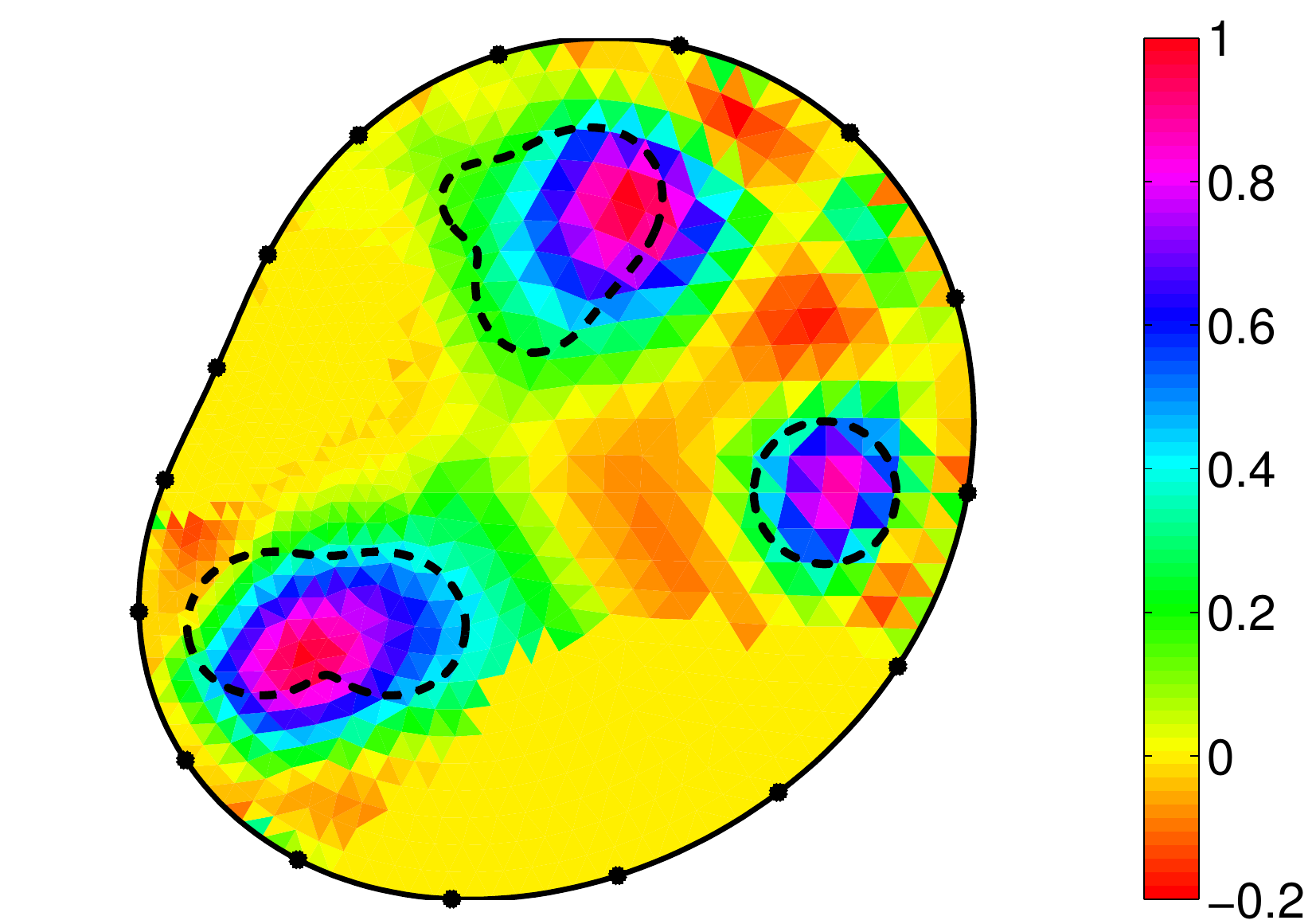}&
\includegraphics[keepaspectratio=true,height=1.5cm]{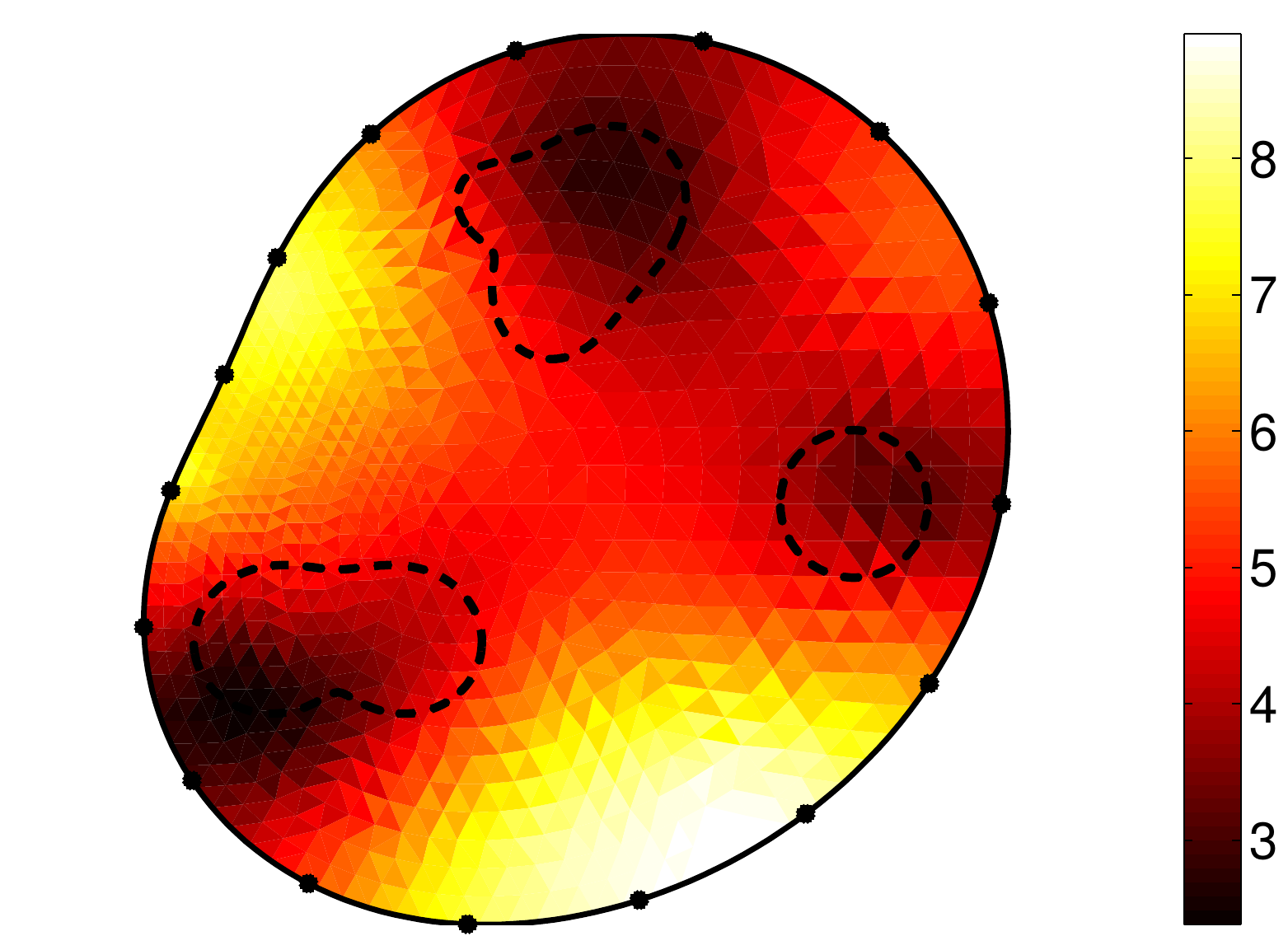}\\
\hline
\end{tabular}
\caption{ Reconstructed difference EIT images in non-circular domain. $\DS$: true difference image,
$\DS_{S}$: standard linearized method, cf.\ \eref{LM2},
$\DS_{B}$: naive combination of LM and S-FM, cf.\ \eref{regularization},
$\DS_{A}$: proposed combination of LM and S-FM, cf.\ \eref{recon},
$\mathbf{W1}$: S-FM alone, cf.\ \eref{SFMimage}.}
\label{recon_noncircle}
\end{figure}

We tested the performance of the proposed hybrid method through numerical simulations with 16-channel EIT system with two different domains. The background conductivity is  $\sigma_0=1$.
Inside the domain $\Om$,  we placed  various anomalies  as shown in  the second column of Figure \ref{recon_circle}-\ref{recon_noncircle}. We generate a mesh of $\Om$ using 4128 and 4432 triangular elements and 2129  and 2291 nodes for the unit circle and deformed circle, respectively. We numerically compute the forward problem in \eref{govern1} and \eref{govern2} to generate the data-set $\DV$. The sensitivity matrix $\S$ is computed by solving \eref{govern1} and \eref{govern2} with $\sigma$ replaced by $1$. The number of pixel for reconstructed images is $n_p=1414$ and $1424$ for the unit circle and deformed circle, respectively, and  it should be different from meshes in the forward model.

The last column in Figure \ref{recon_circle}-\ref{recon_noncircle} shows images of $\mathbf {W1}$:
$$
n\mbox{-th element of } \mathbf{W1}=\ln\left( 1+
\sum_{j=1}^{n_E}\left|\frac{\S_j^n\cdot \delta\Bbb{V}^{-1}\S_j^n}{
\S_j^n\cdot\S_j^n}\right|\right).$$
Here, we use a fixed regularization parameter to implement $\delta\Bbb{V}^{-1}\S_j^n$.

The third column in Figure \ref{recon_circle}-\ref{recon_noncircle} shows the reconstructed images using tSVD pseudoinverse $\S^\dag=\mathbf{V}_{t_0}\mathbf{\Lambda}_{t_0}^{-1} \mathbf{U}_{t_0}^*$  with $t_0=64$ :
$$\DS_{S}=\sum_{t=1}^{t_0}\frac{1}{\lambda_t}\langle \DV,\mathbf{u}_t\rangle\v_t \q\q \Big(\f{\lambda_{t_0}}{\lambda_1}\approx\frac{0.001}{0.8}= 0.0012\Big).
$$

\begin{figure}
\centering
\begin{tabular}{|c|cccc|c|}
\hline
Case &  $\DS$ & $\DS_{S}$ & $\DS_{B}$ &  $\DS_{A}$ &$\mathbf{W1}$ \\
\hline
\raisebox{4ex}{\footnotesize\begin{tabular}{c}
(a)
\end{tabular}}
&
\includegraphics[keepaspectratio=true,width=1.8cm]{Fig/plot_final_20120602/deltasigma_1_woColorbar-eps-converted-to.pdf}&
\includegraphics[keepaspectratio=true,width=1.8cm]{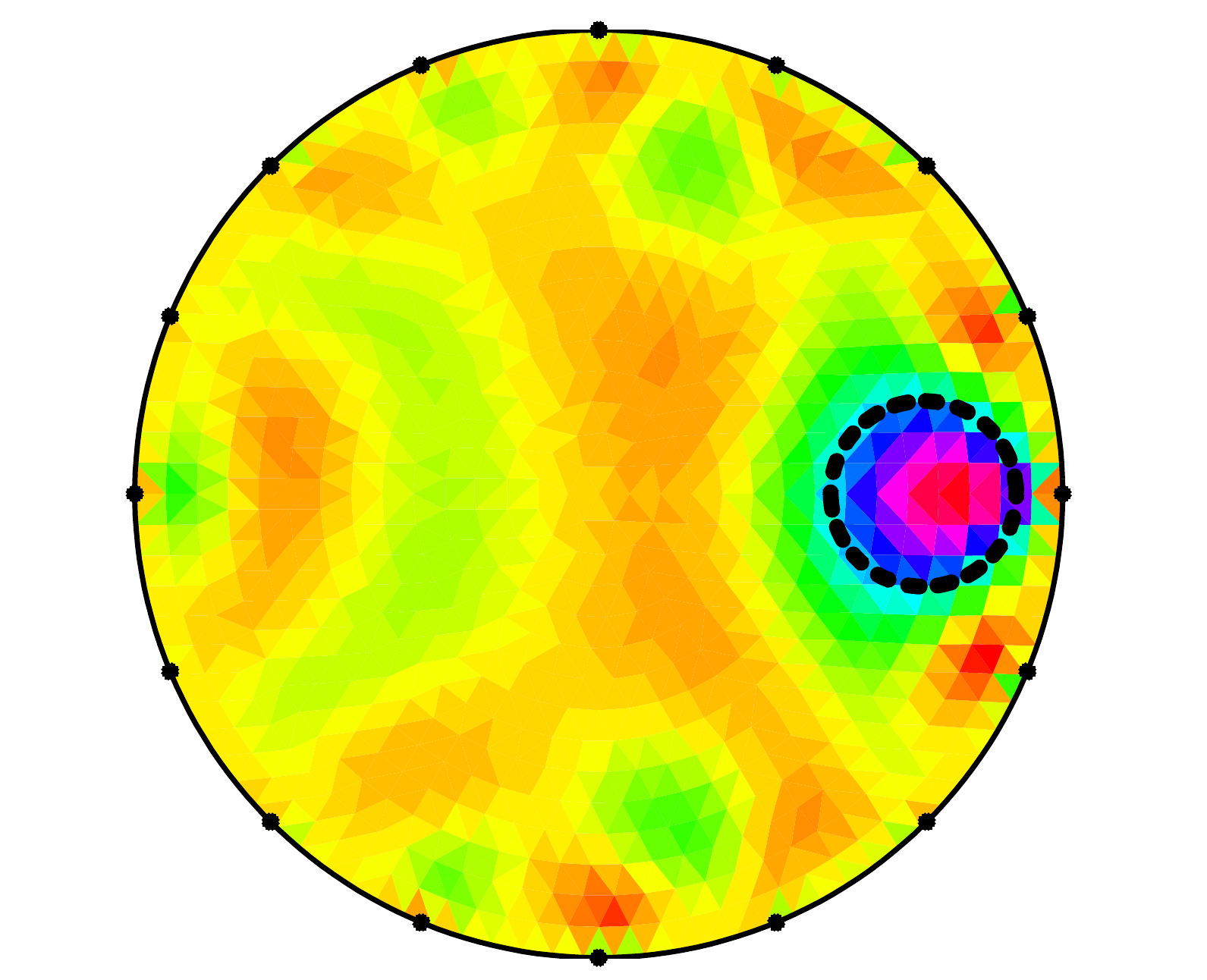}&
\includegraphics[keepaspectratio=true,width=1.8cm]{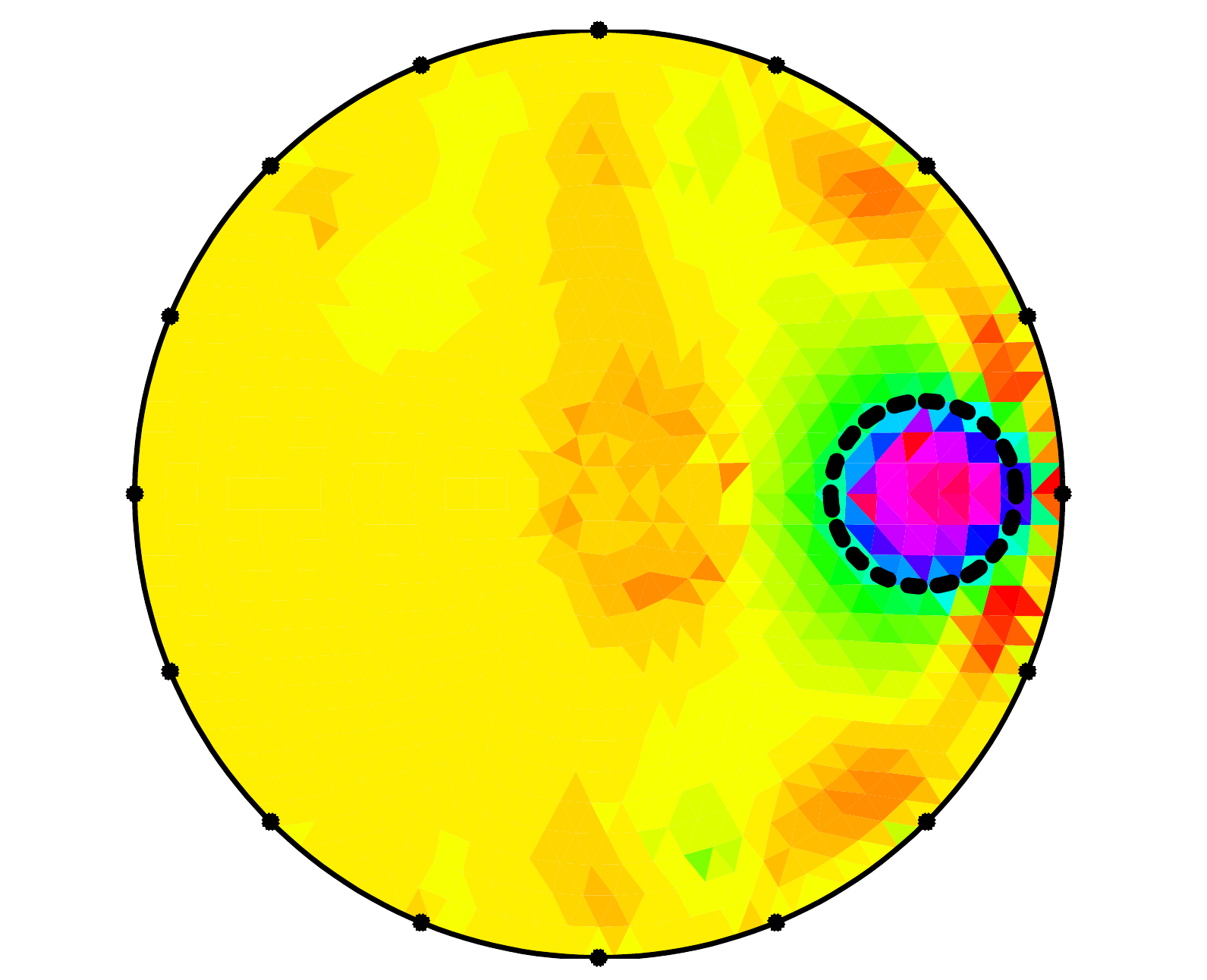}&
\includegraphics[keepaspectratio=true,width=2.1cm]{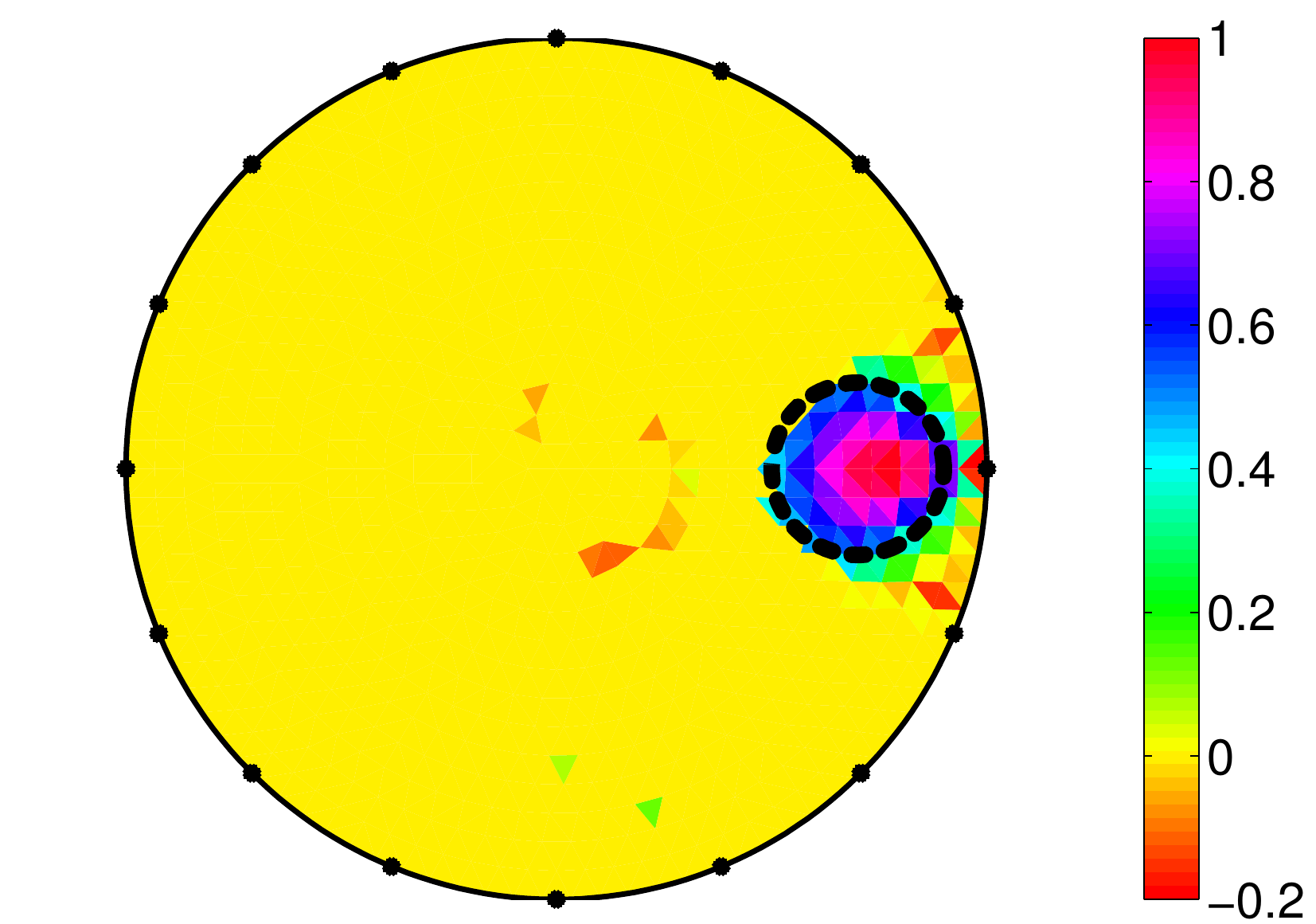}&
\includegraphics[keepaspectratio=true,width=2.05cm]{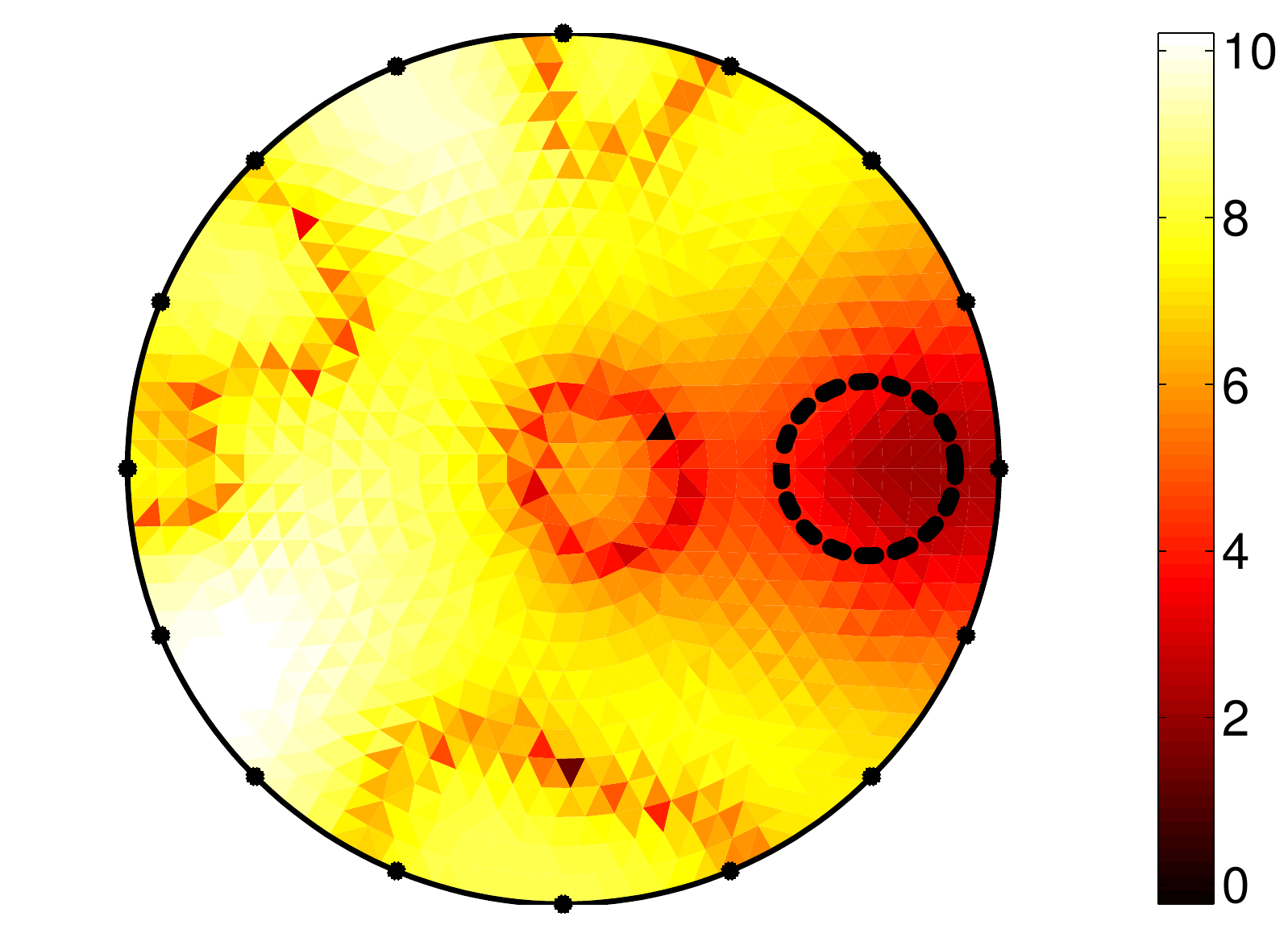}\\
\hline
\raisebox{4ex}{\footnotesize\begin{tabular}{c}
(b)
\end{tabular}}
&
\includegraphics[keepaspectratio=true,height=1.5cm]{Fig/plot_final_20120602/deltasigma_2_woColorbar-eps-converted-to.pdf}&
\includegraphics[keepaspectratio=true,height=1.5cm]{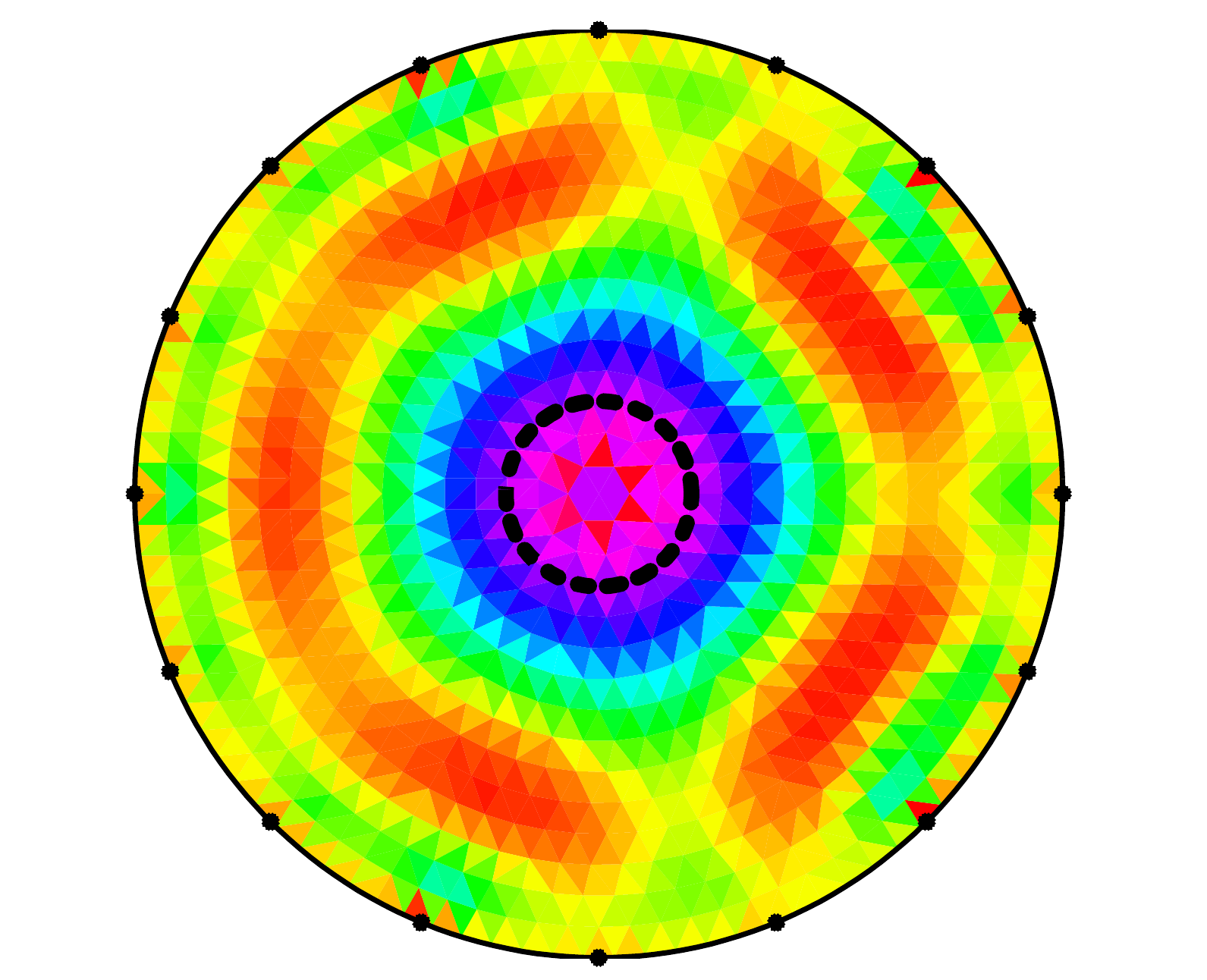}&
\includegraphics[keepaspectratio=true,height=1.5cm]{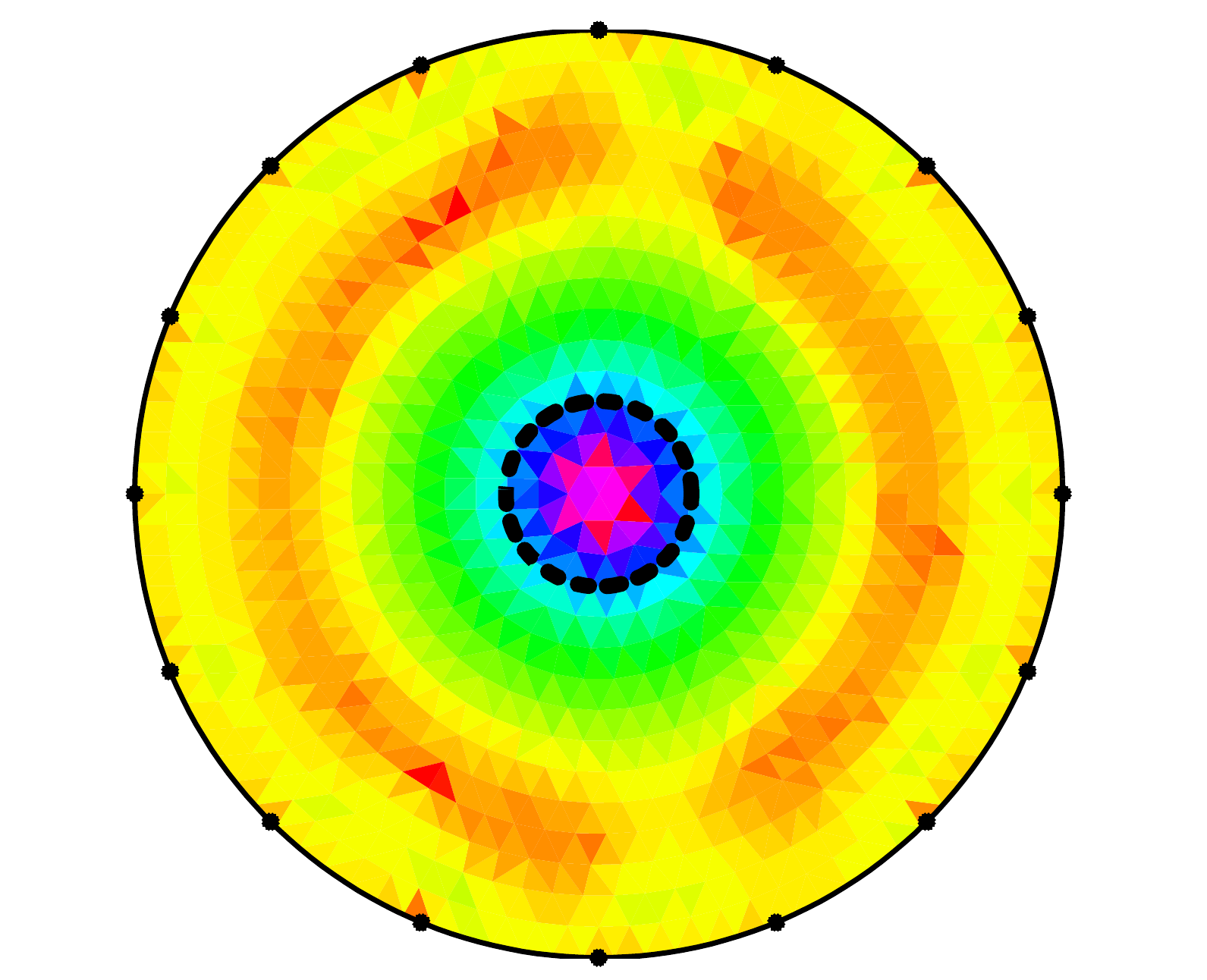}&
\includegraphics[keepaspectratio=true,height=1.5cm]{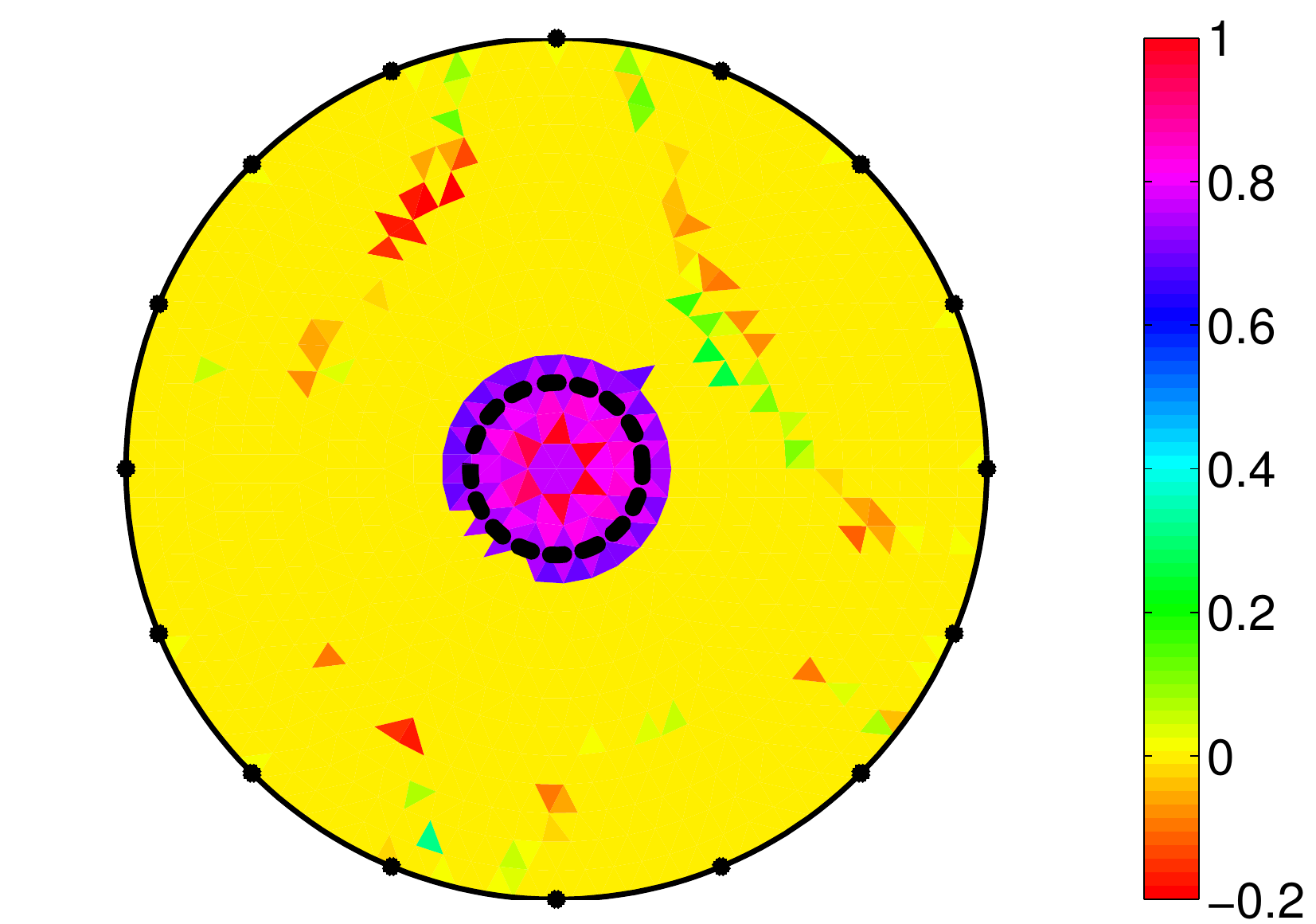}&
\includegraphics[keepaspectratio=true,height=1.5cm]{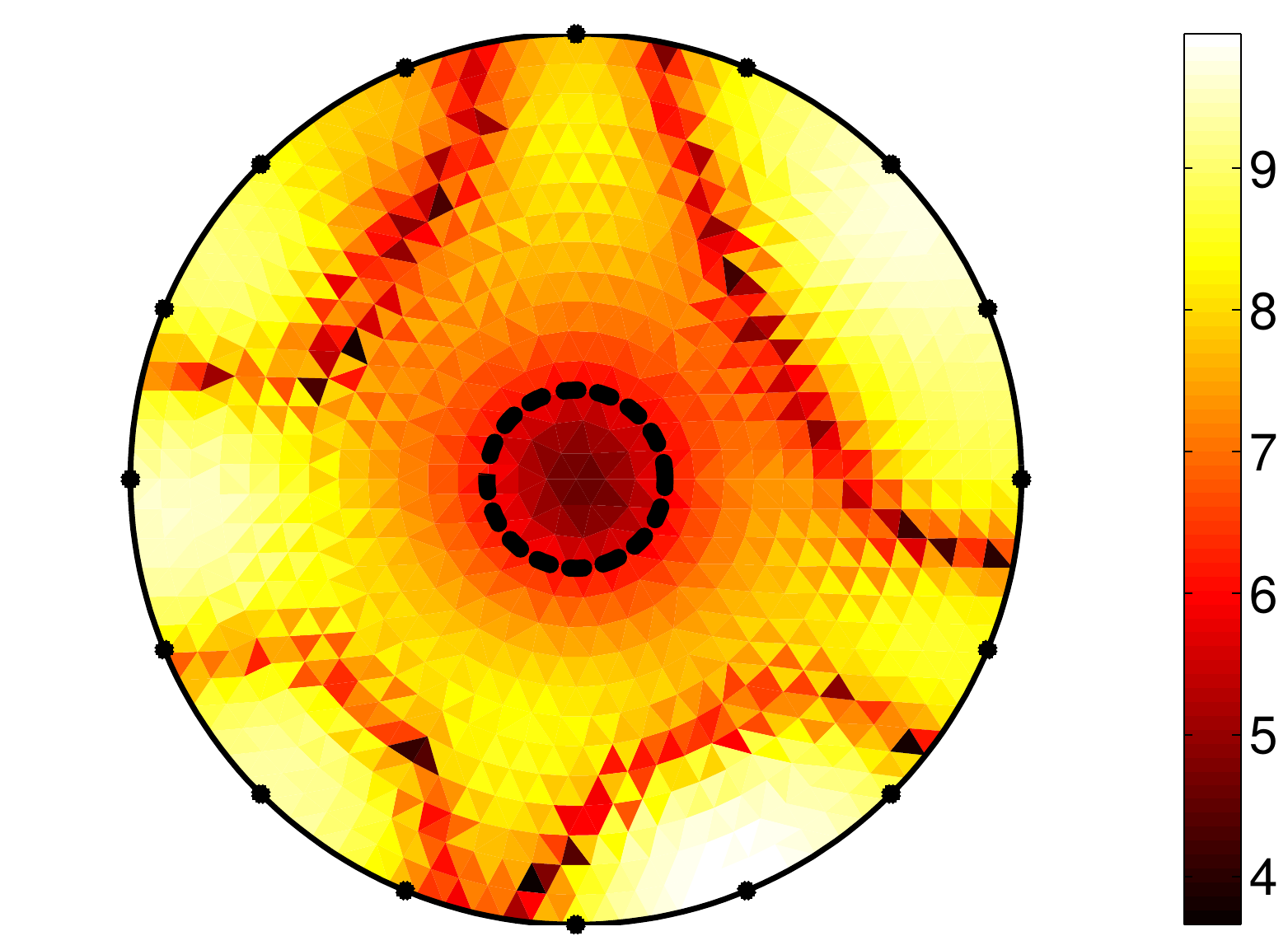}\\
\hline
\raisebox{4ex}{\footnotesize\begin{tabular}{c}
(c)
\end{tabular}}
&
\includegraphics[keepaspectratio=true,height=1.5cm]{Fig/plot_final_20120602/deltasigma_3_woColorbar-eps-converted-to.pdf}&
\includegraphics[keepaspectratio=true,height=1.5cm]{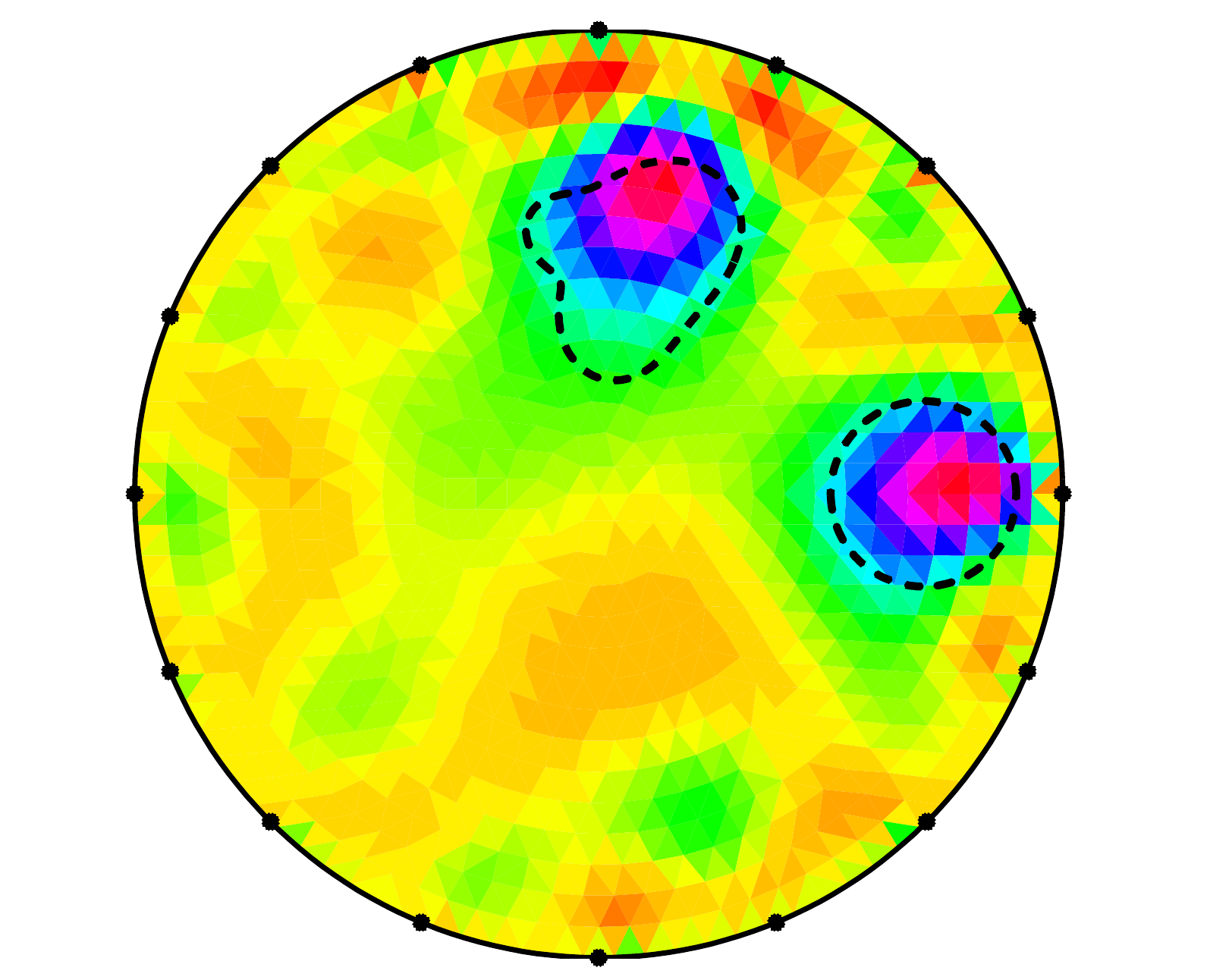}&
\includegraphics[keepaspectratio=true,height=1.5cm]{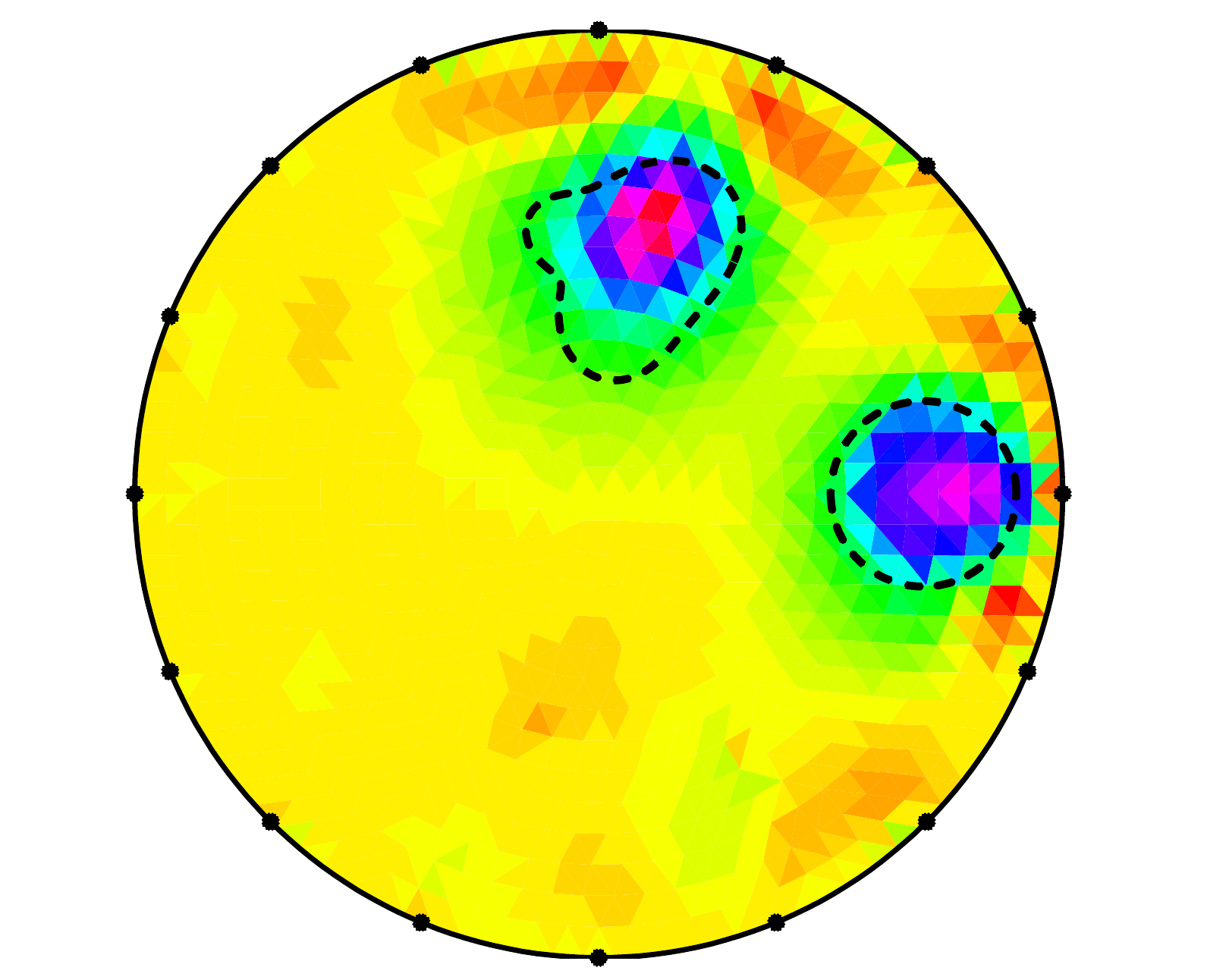}&
\includegraphics[keepaspectratio=true,height=1.5cm]{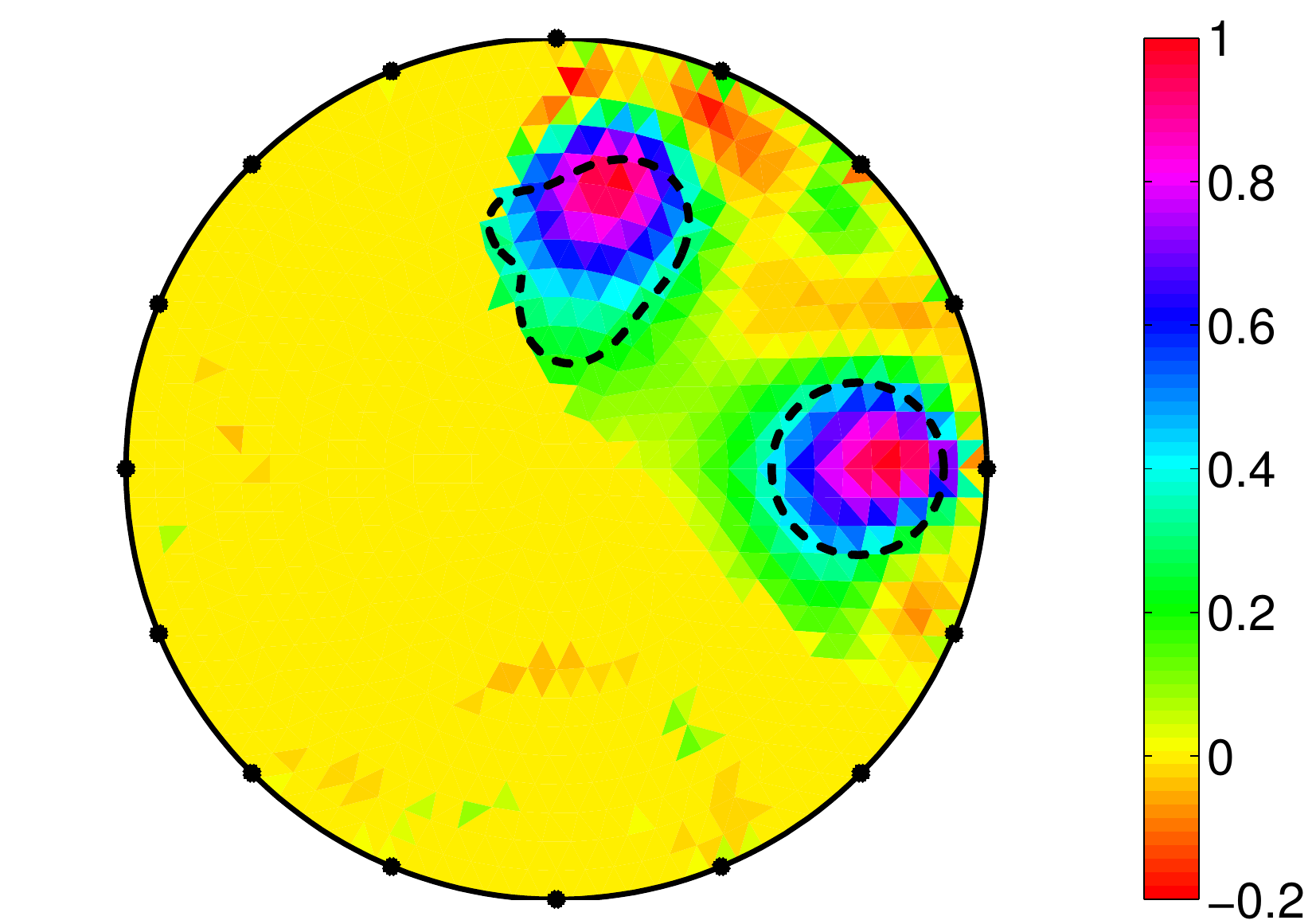}&
\includegraphics[keepaspectratio=true,height=1.5cm]{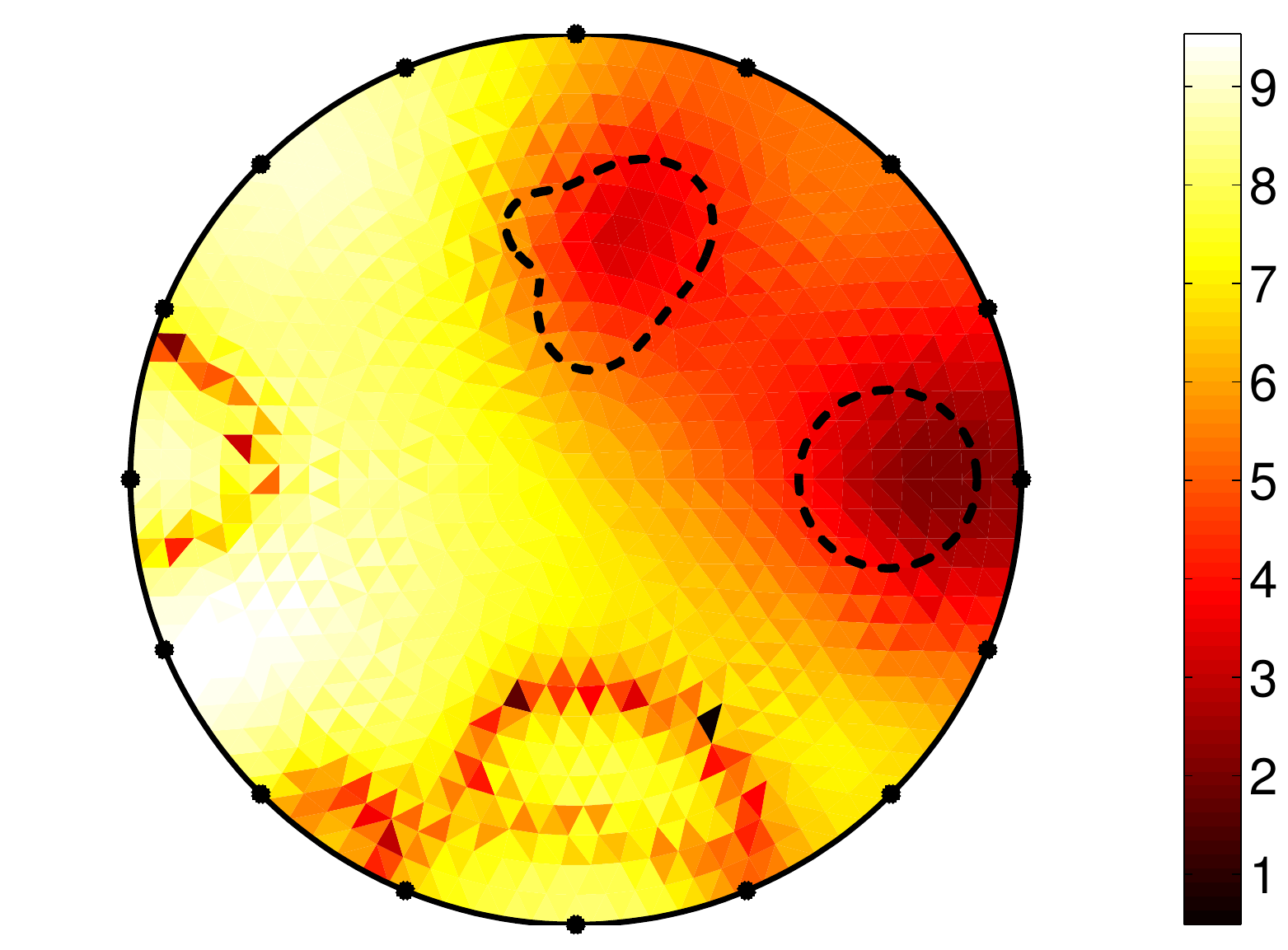}\\
\hline
\raisebox{4ex}{\footnotesize\begin{tabular}{c}
(d)
\end{tabular}}
&
\includegraphics[keepaspectratio=true,height=1.5cm]{Fig/plot_final_20120602/deltasigma_4_woColorbar-eps-converted-to.pdf}&
\includegraphics[keepaspectratio=true,height=1.5cm]{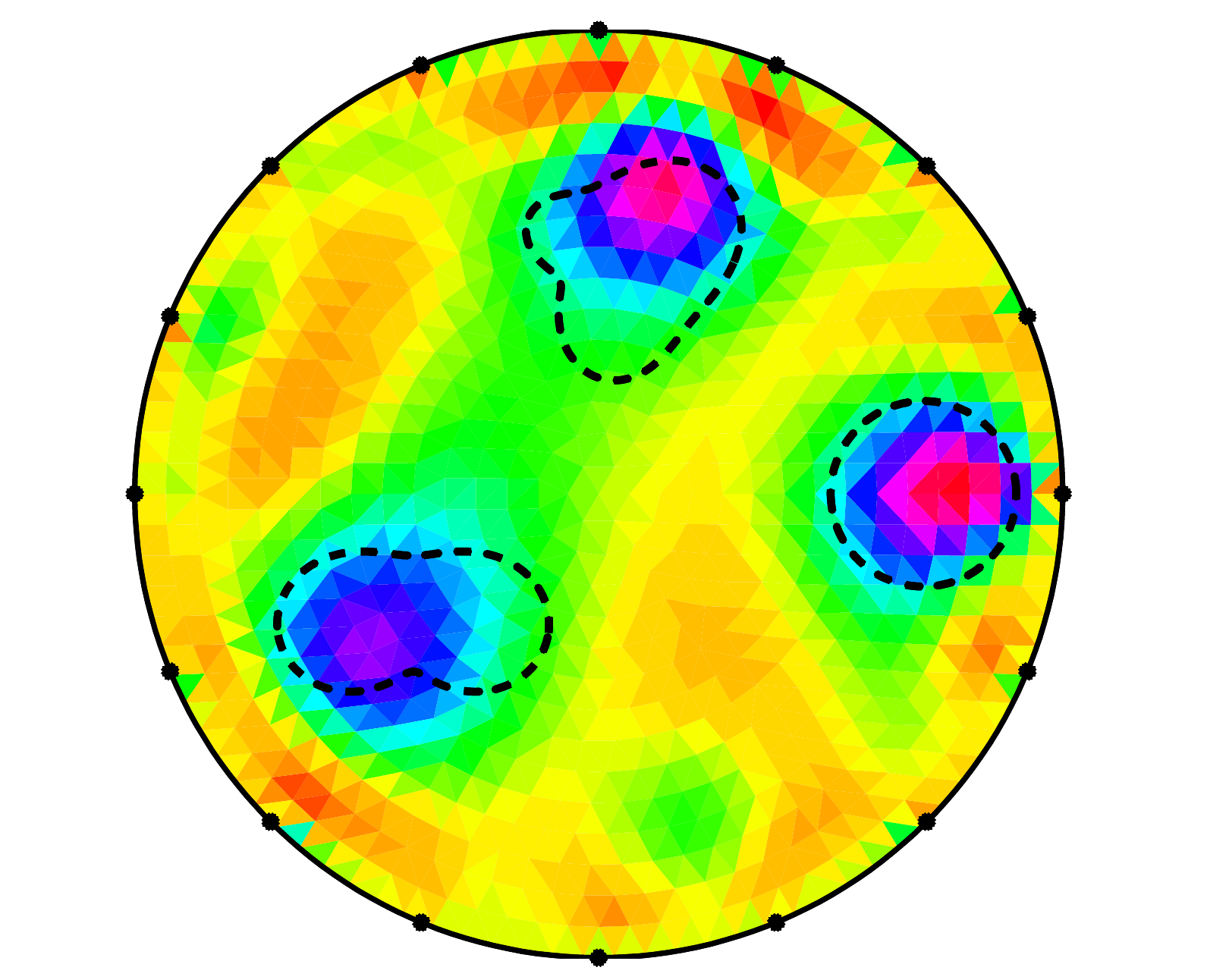}&
\includegraphics[keepaspectratio=true,height=1.5cm]{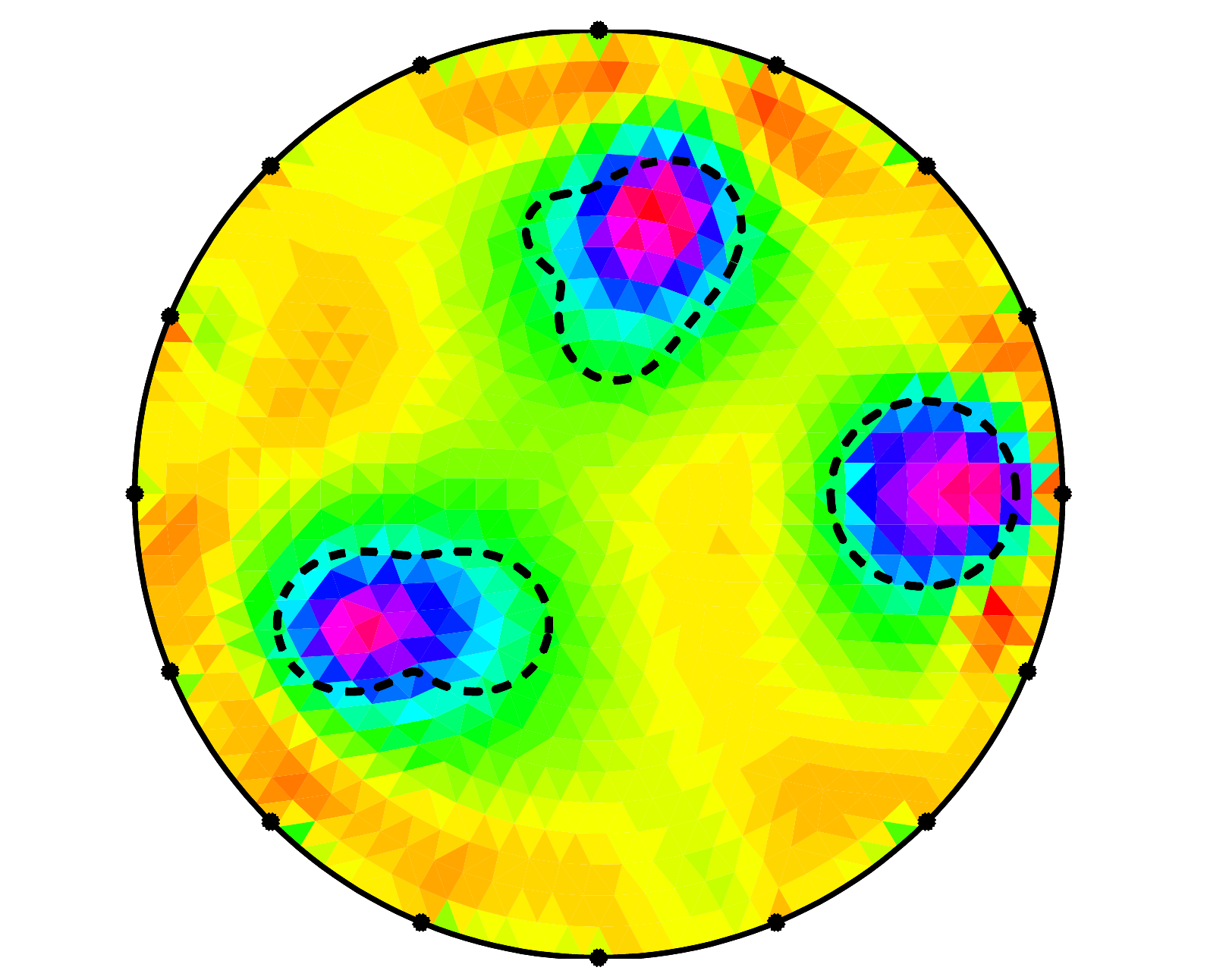}&
\includegraphics[keepaspectratio=true,height=1.5cm]{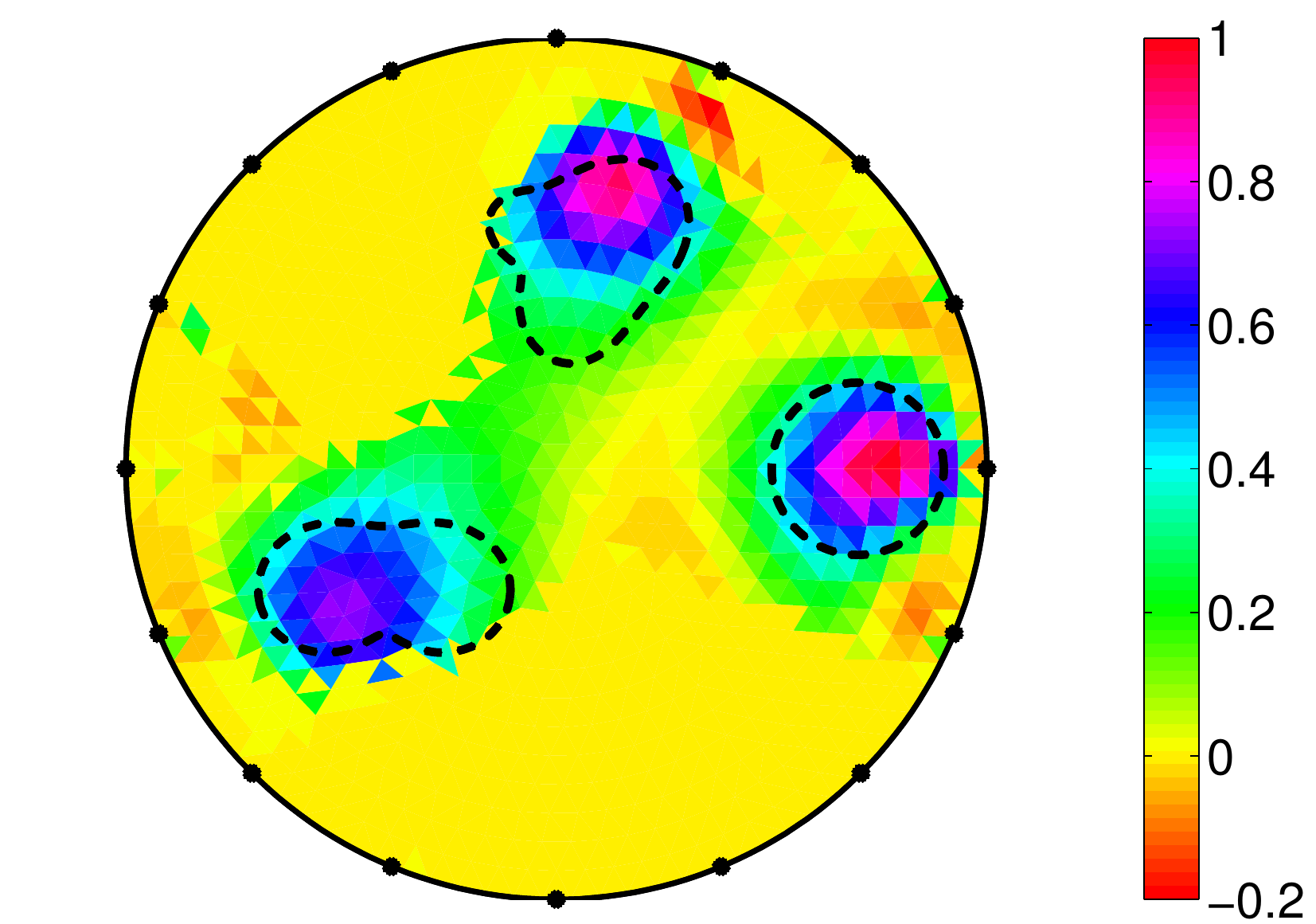}&
\includegraphics[keepaspectratio=true,height=1.5cm]{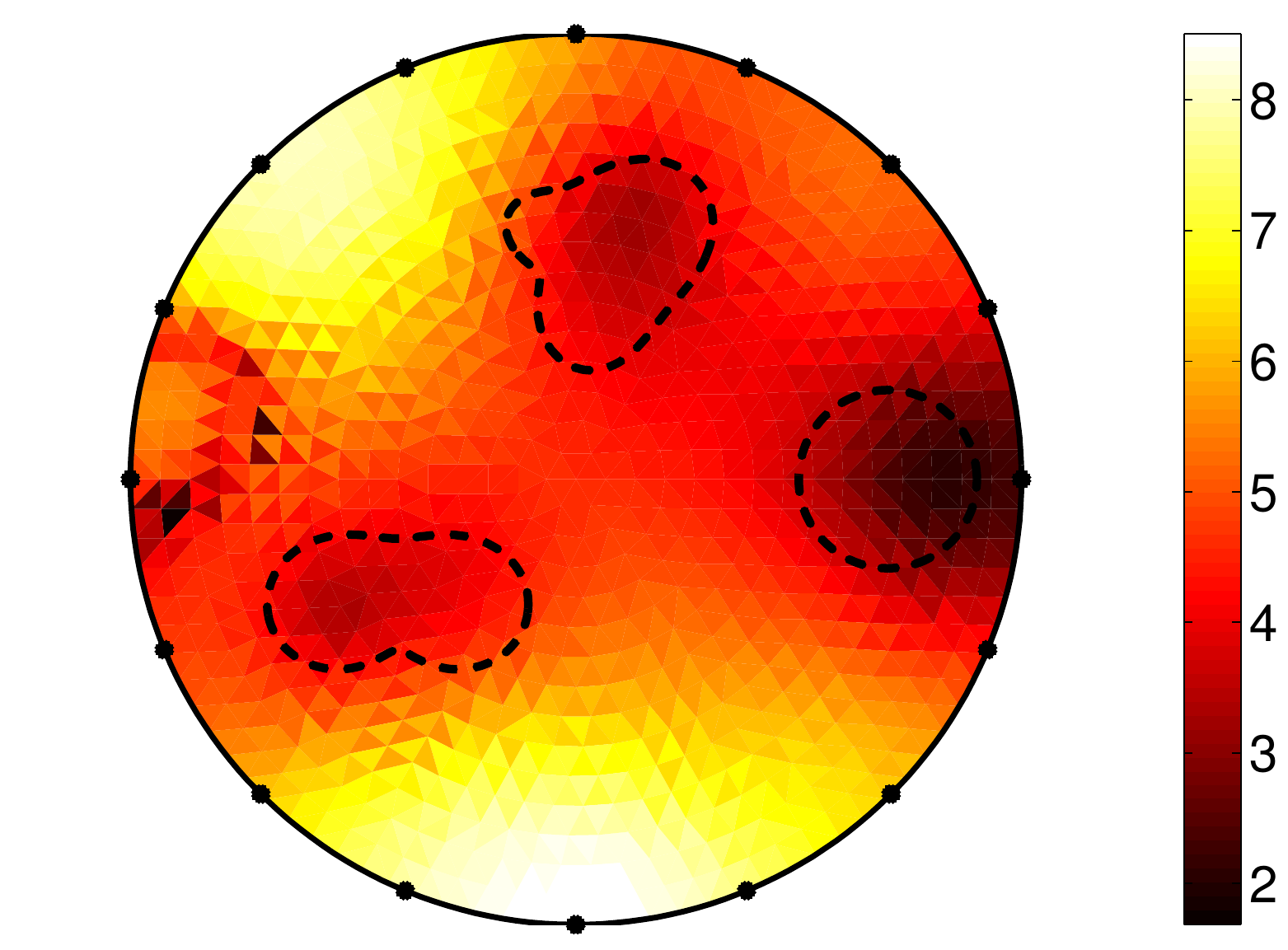}\\
\hline
\end{tabular}
\caption{ Reconstructed difference EIT images in circular domain with a data which adds 1$\%$ random noise. $\DS$: true difference image,
$\DS_{S}$: standard linearized method, cf.\ \eref{LM2},
$\DS_{B}$: naive combination of LM and S-FM, cf.\ \eref{regularization},
$\DS_{A}$: proposed combination of LM and S-FM, cf.\ \eref{recon},
$\mathbf{W1}$: S-FM alone, cf.\ \eref{SFMimage}.}
\label{recon_circle_1p}
\end{figure}

\begin{figure}
\centering
\begin{tabular}{|c|cccc|c|}
\hline
Case &  $\DS$ & $\DS_{S}$ & $\DS_{B}$ &  $\DS_{A}$ &$\mathbf{W1}$ \\
\hline
\raisebox{4ex}{\footnotesize\begin{tabular}{c}
(e)
\end{tabular}}
&
\includegraphics[keepaspectratio=true,height=1.5cm]{Fig/plot_final_20120602/deltasigma_5_woColorbar-eps-converted-to.pdf}&
\includegraphics[keepaspectratio=true,height=1.5cm]{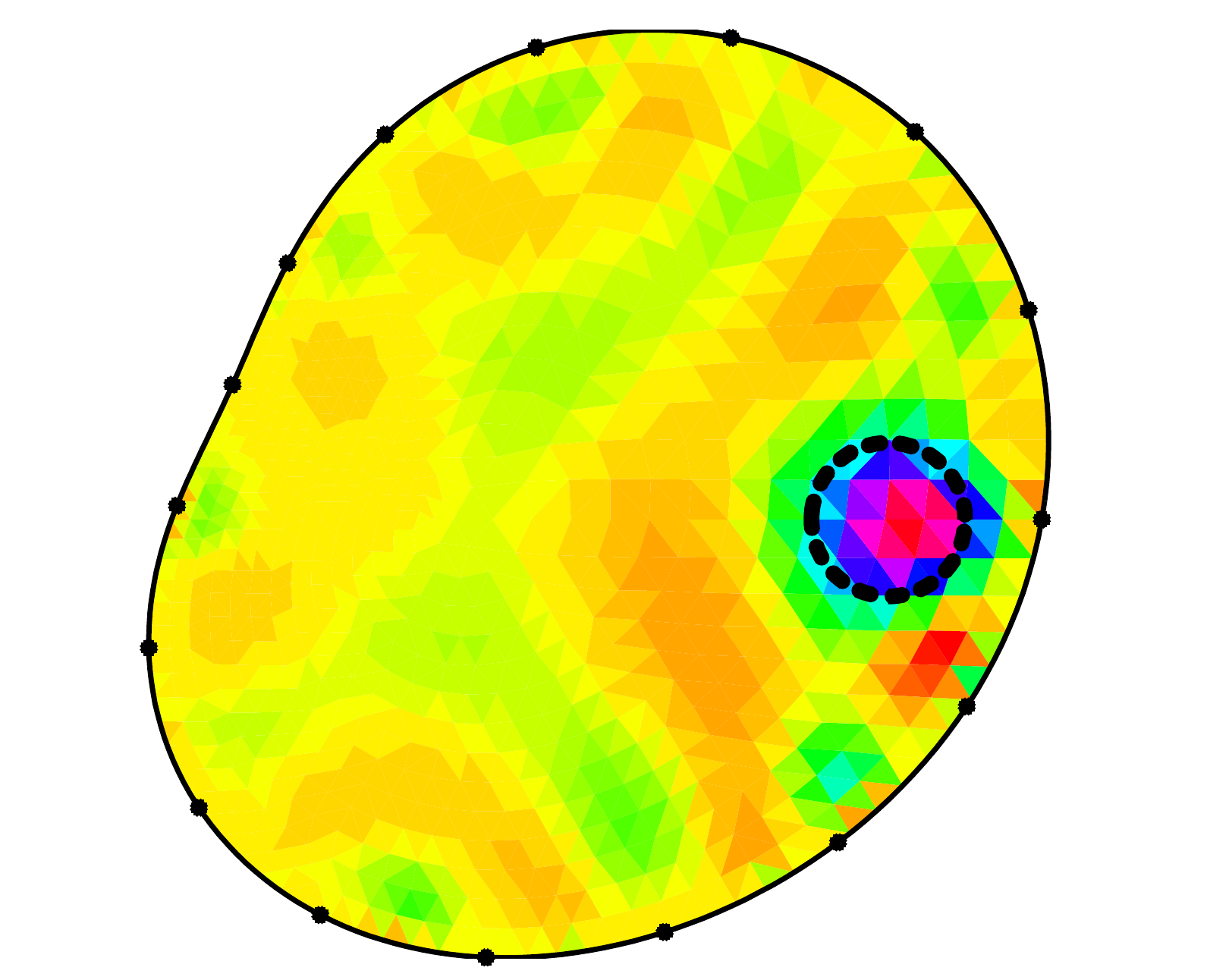}&
\includegraphics[keepaspectratio=true,height=1.5cm]{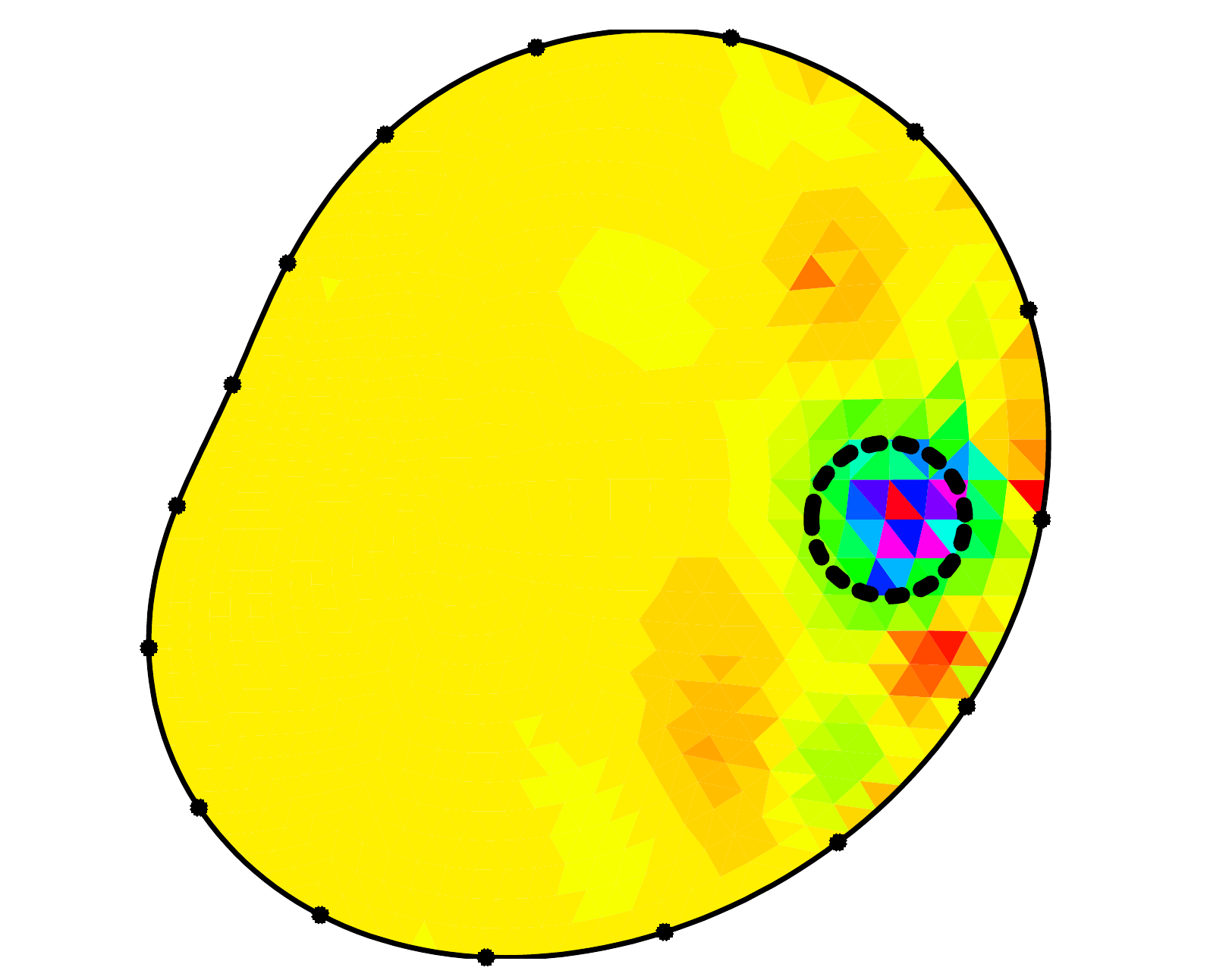}&
\includegraphics[keepaspectratio=true,height=1.5cm]{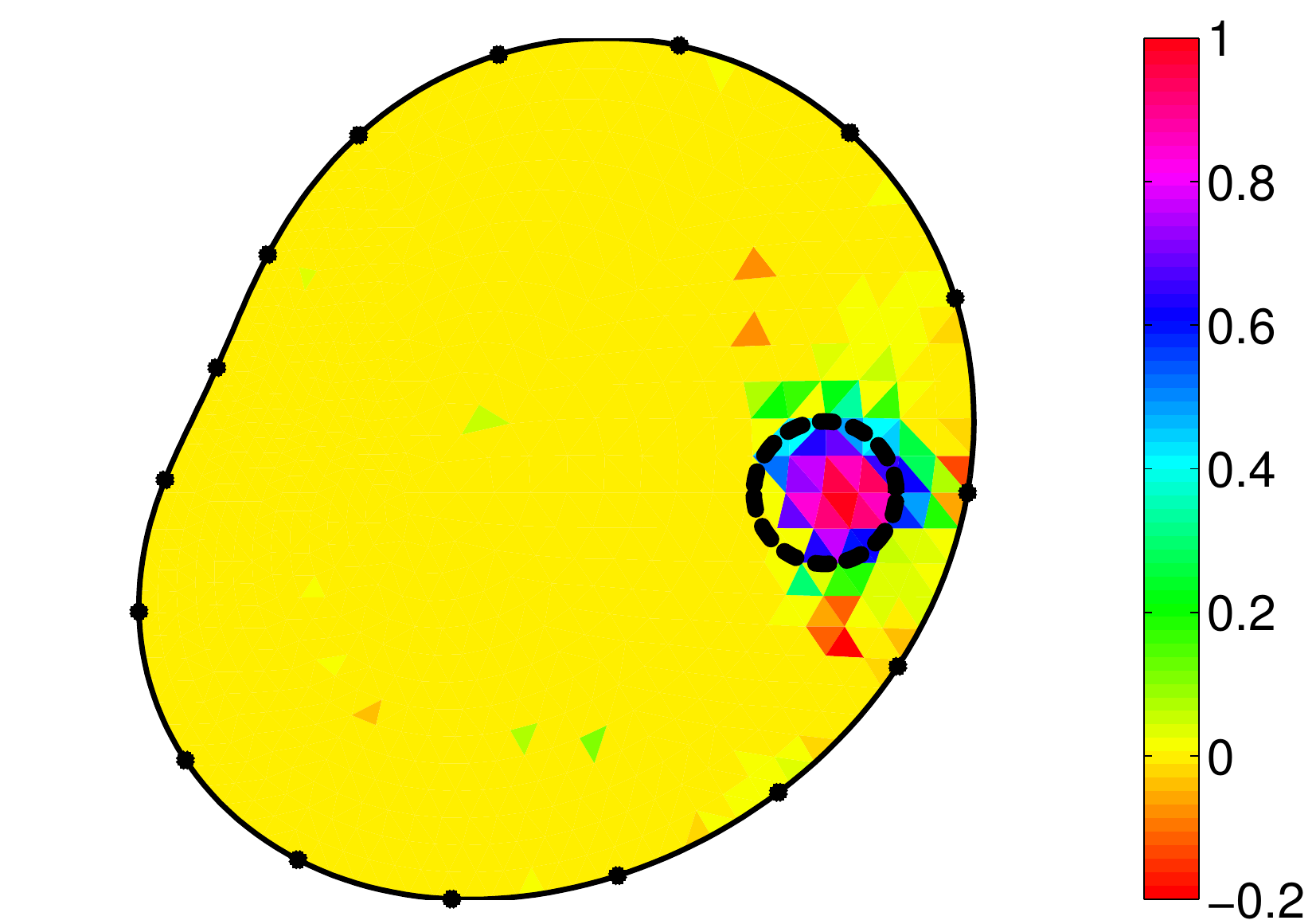}&
\includegraphics[keepaspectratio=true,height=1.5cm]{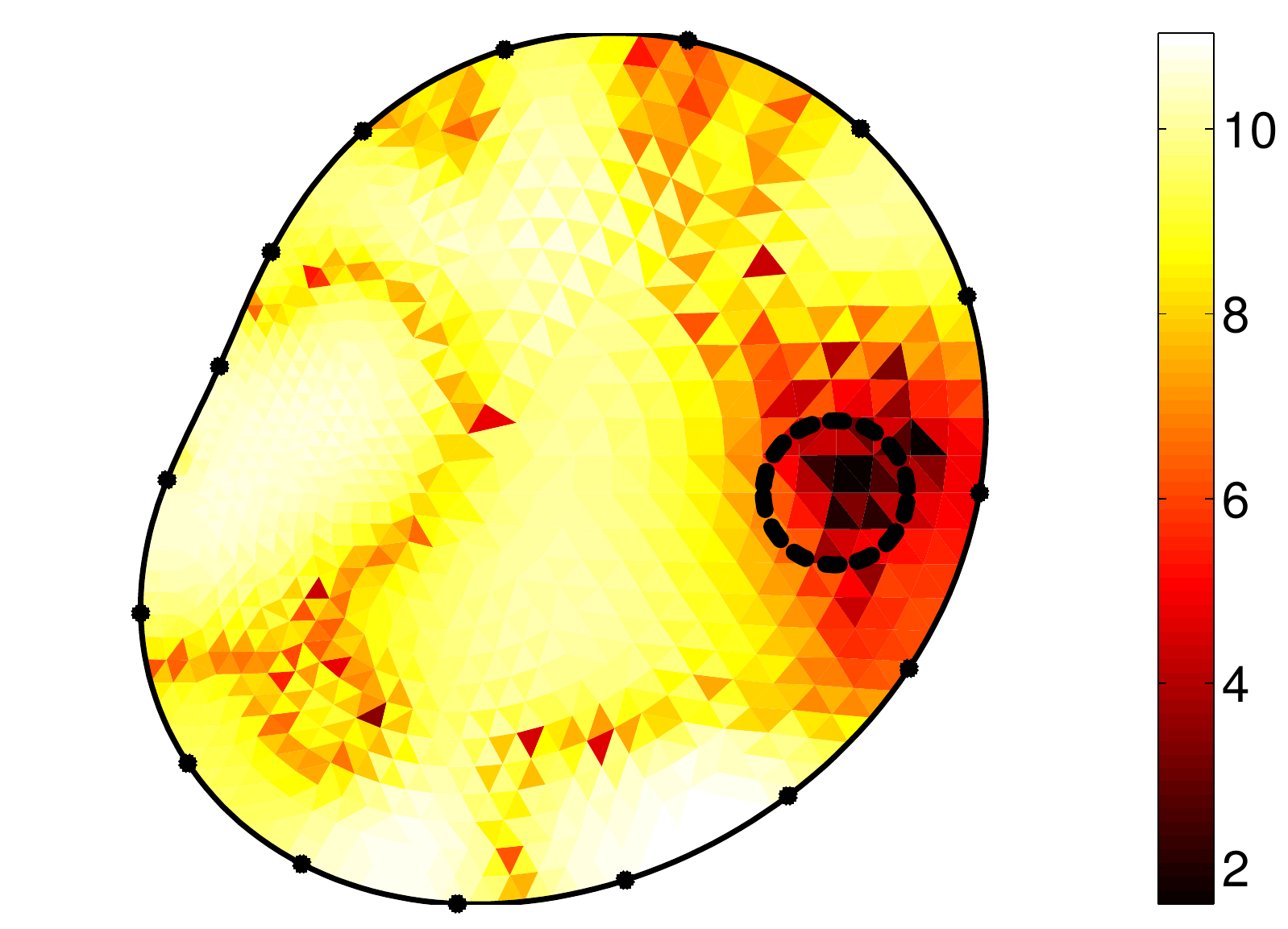}\\
\hline
\raisebox{4ex}{\footnotesize\begin{tabular}{c}
(f)
\end{tabular}}
&
\includegraphics[keepaspectratio=true,height=1.5cm]{Fig/plot_final_20120602/deltasigma_6_woColorbar-eps-converted-to.pdf}&
\includegraphics[keepaspectratio=true,height=1.5cm]{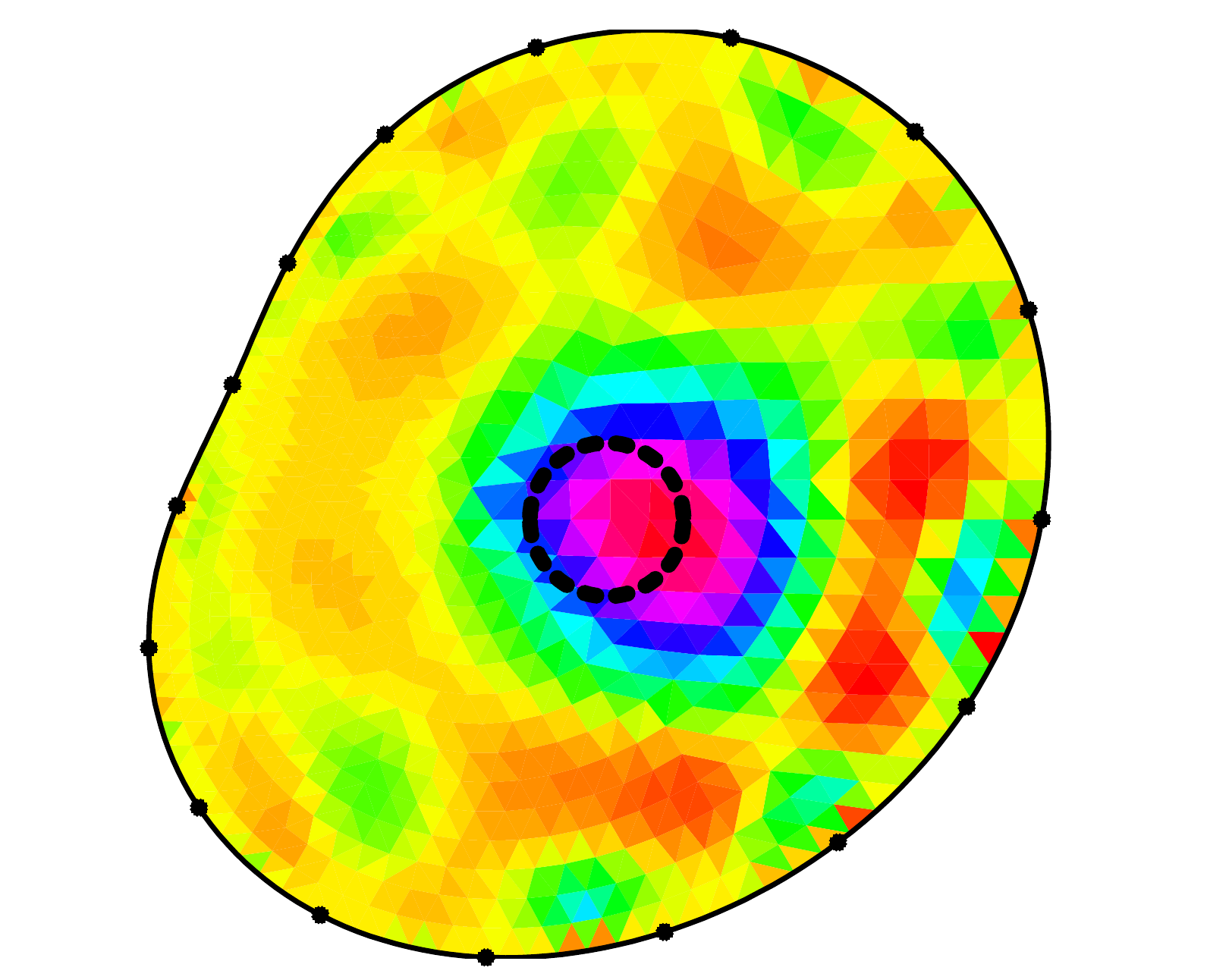}&
\includegraphics[keepaspectratio=true,height=1.5cm]{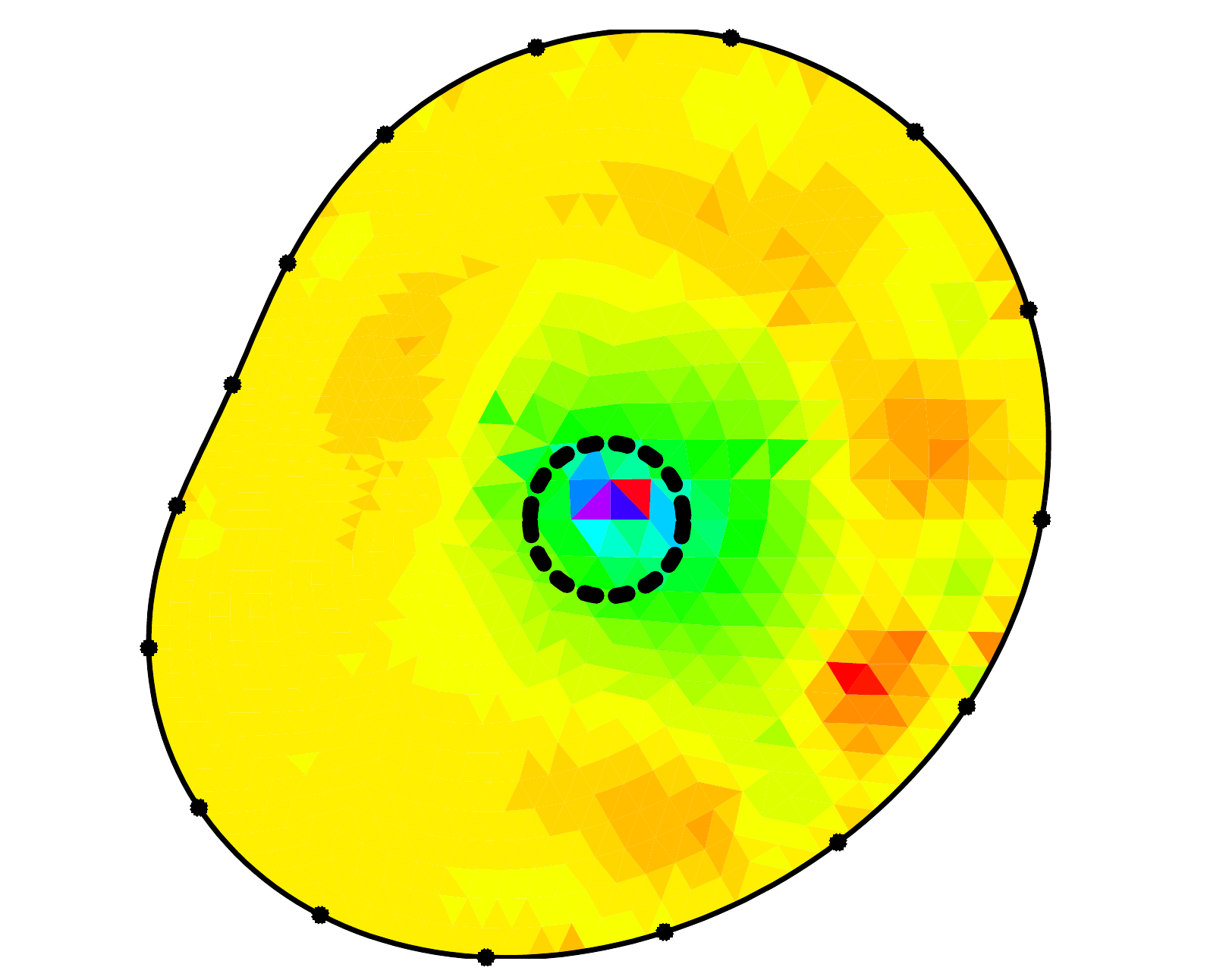}&
\includegraphics[keepaspectratio=true,height=1.5cm]{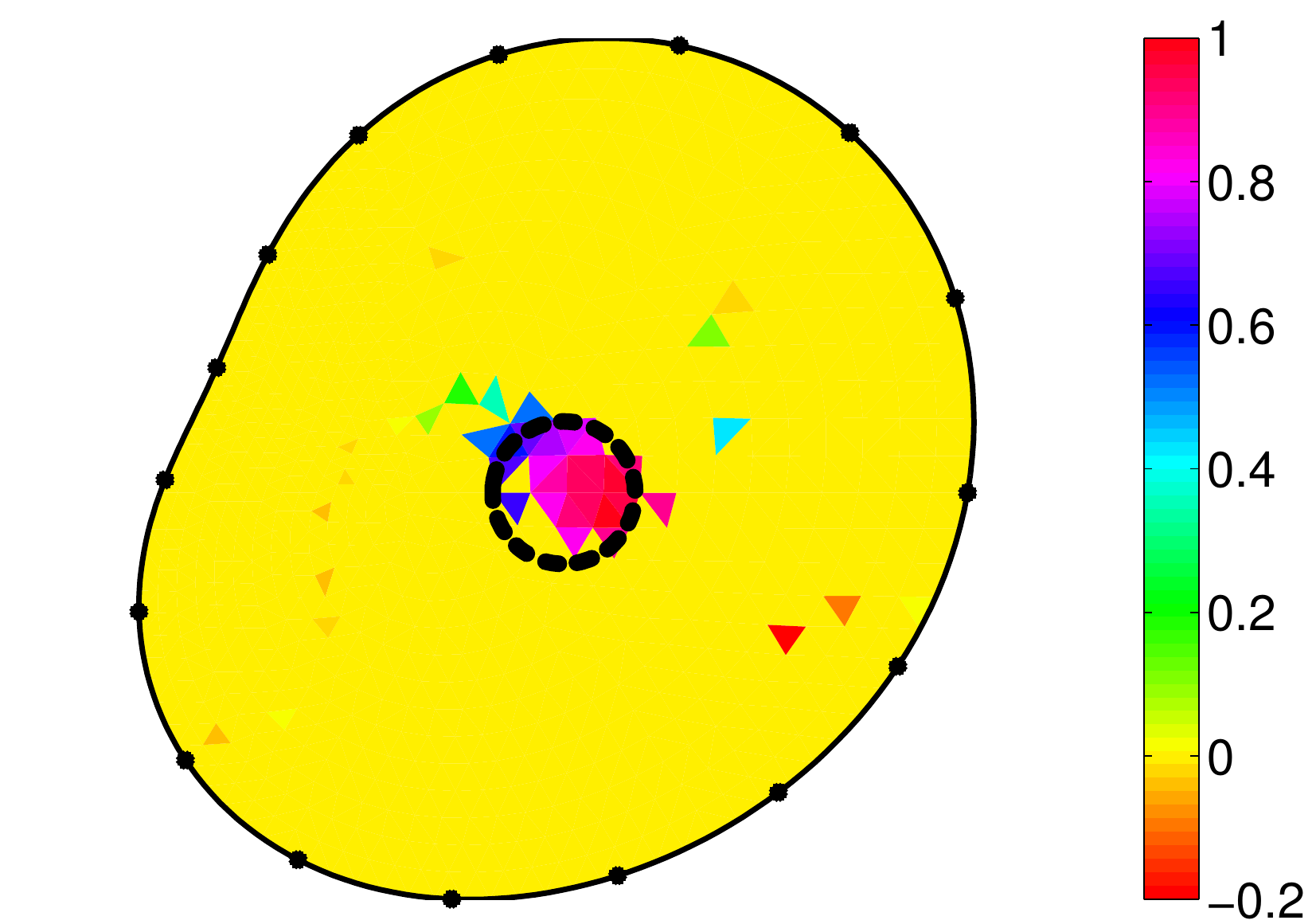}&
\includegraphics[keepaspectratio=true,height=1.5cm]{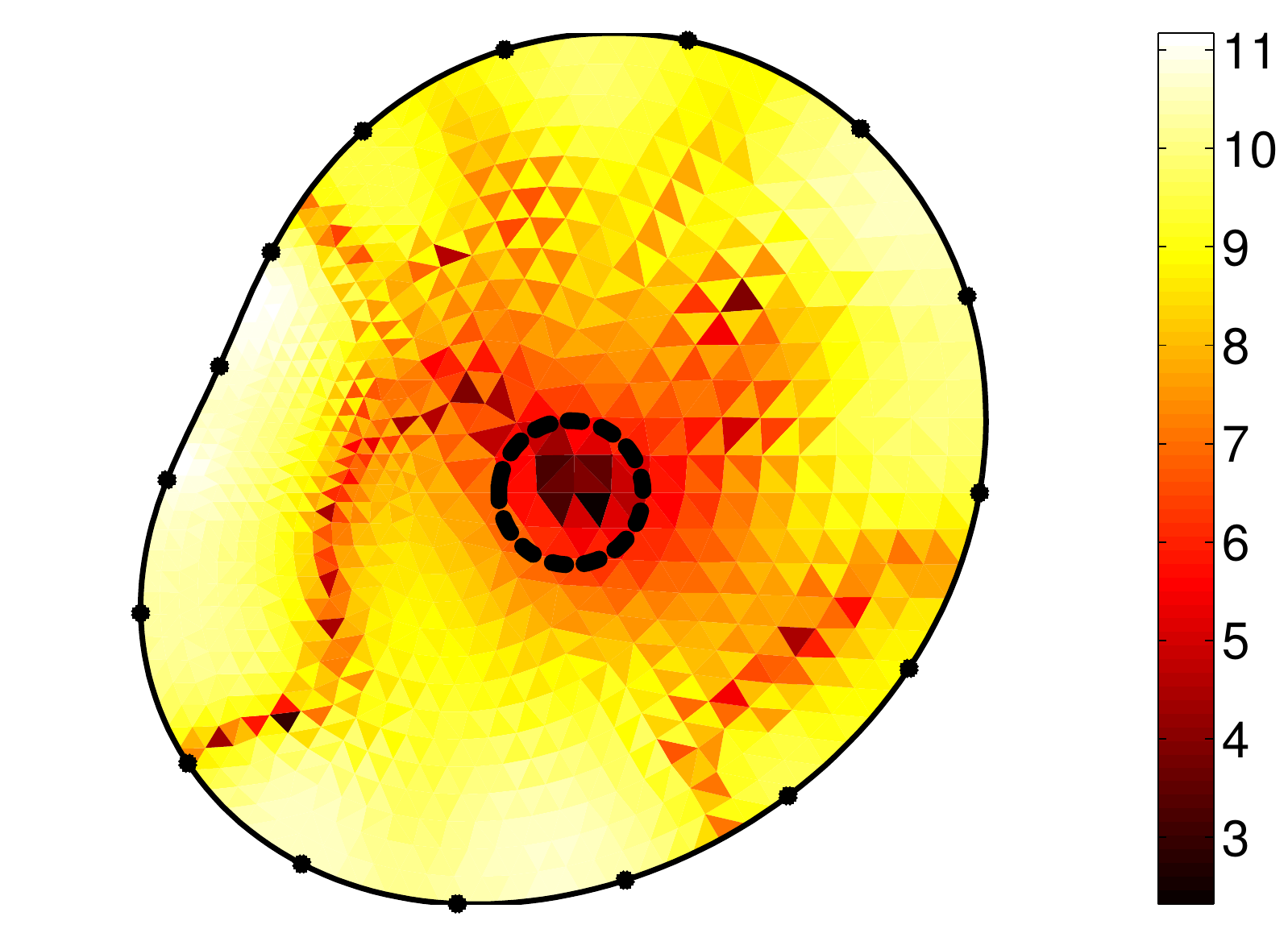}\\
\hline
\raisebox{4ex}{\footnotesize\begin{tabular}{c}
(g)
\end{tabular}}
&
\includegraphics[keepaspectratio=true,height=1.5cm]{Fig/plot_final_20120602/deltasigma_7_woColorbar-eps-converted-to.pdf}&
\includegraphics[keepaspectratio=true,height=1.5cm]{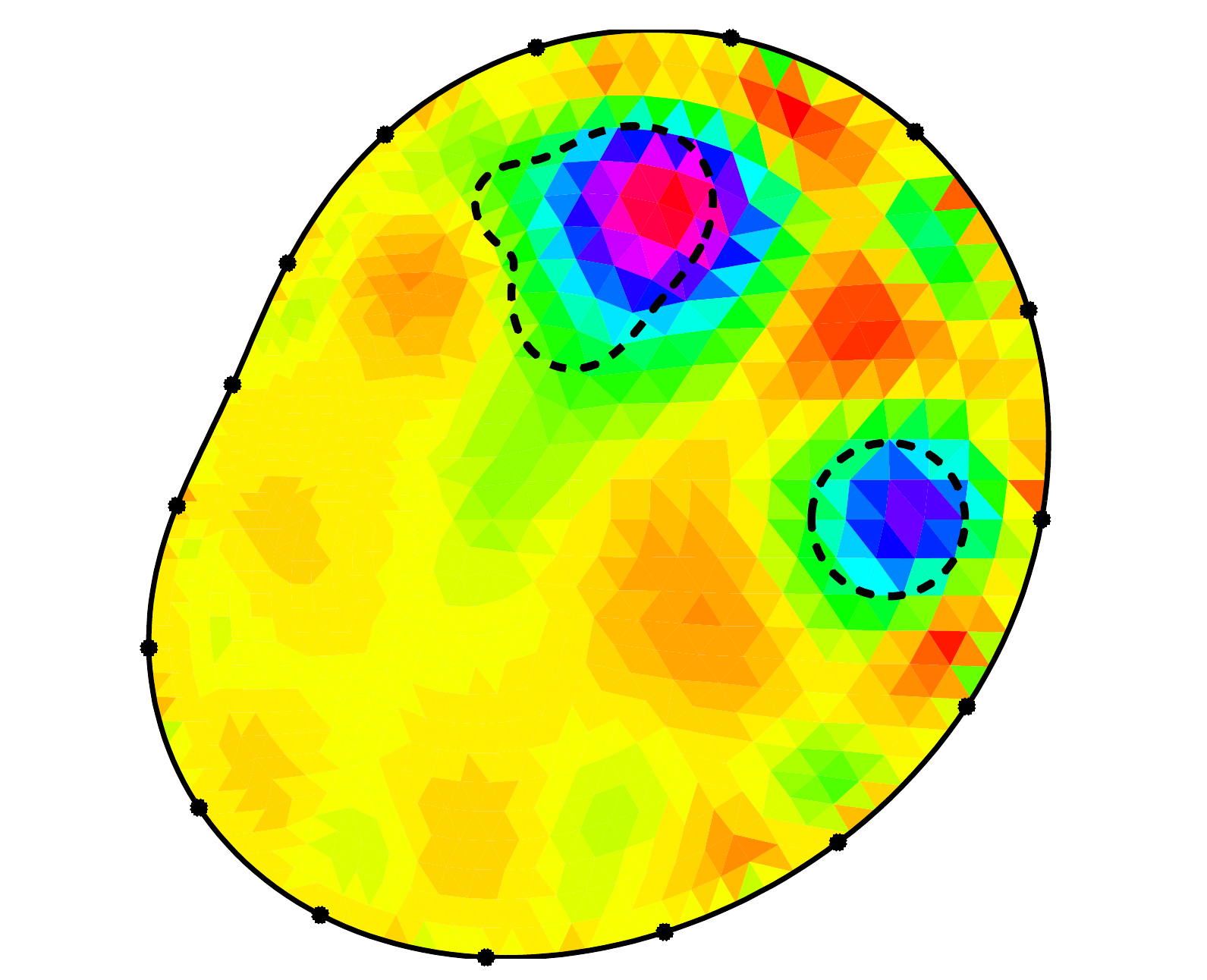}&
\includegraphics[keepaspectratio=true,height=1.5cm]{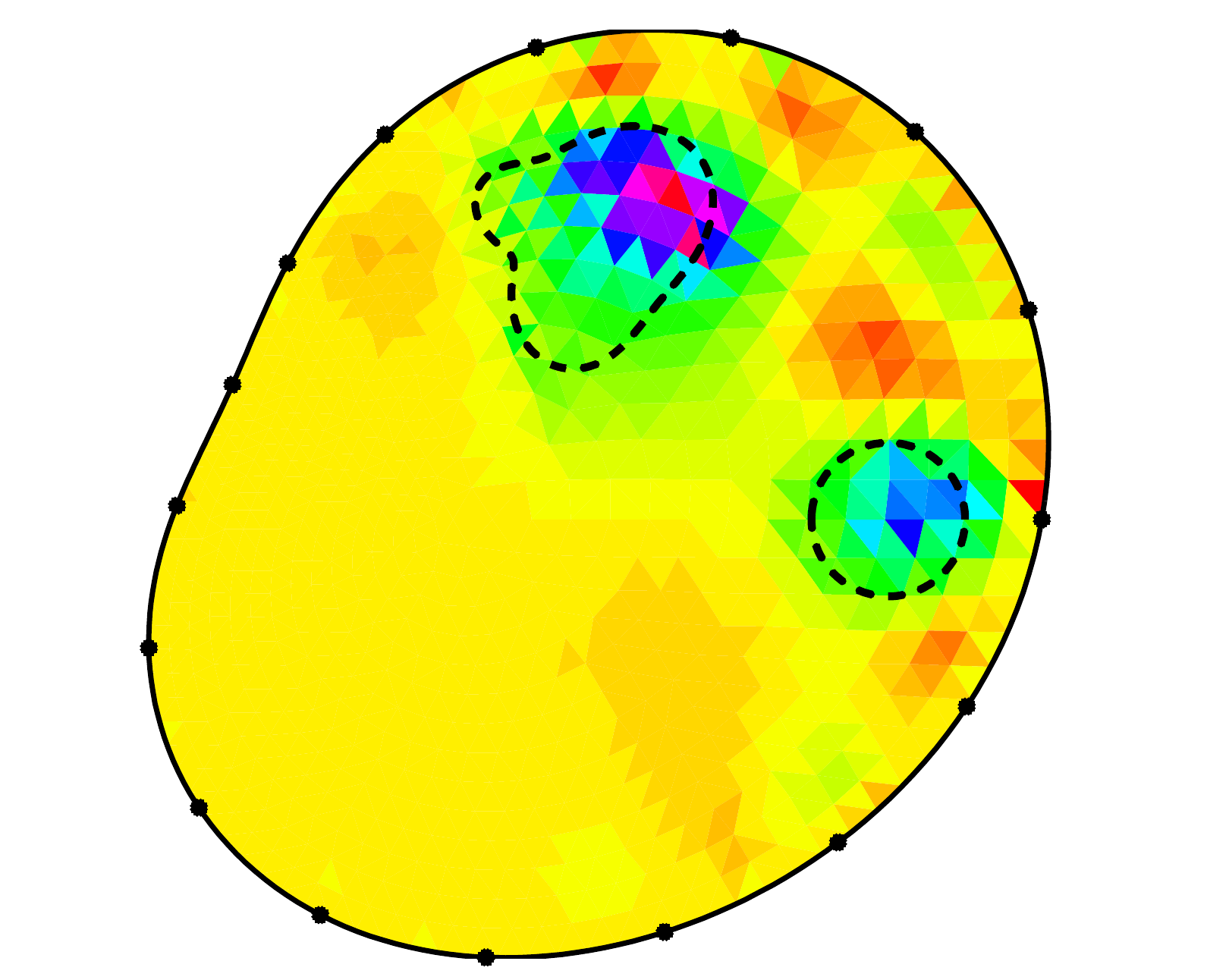}&
\includegraphics[keepaspectratio=true,height=1.5cm]{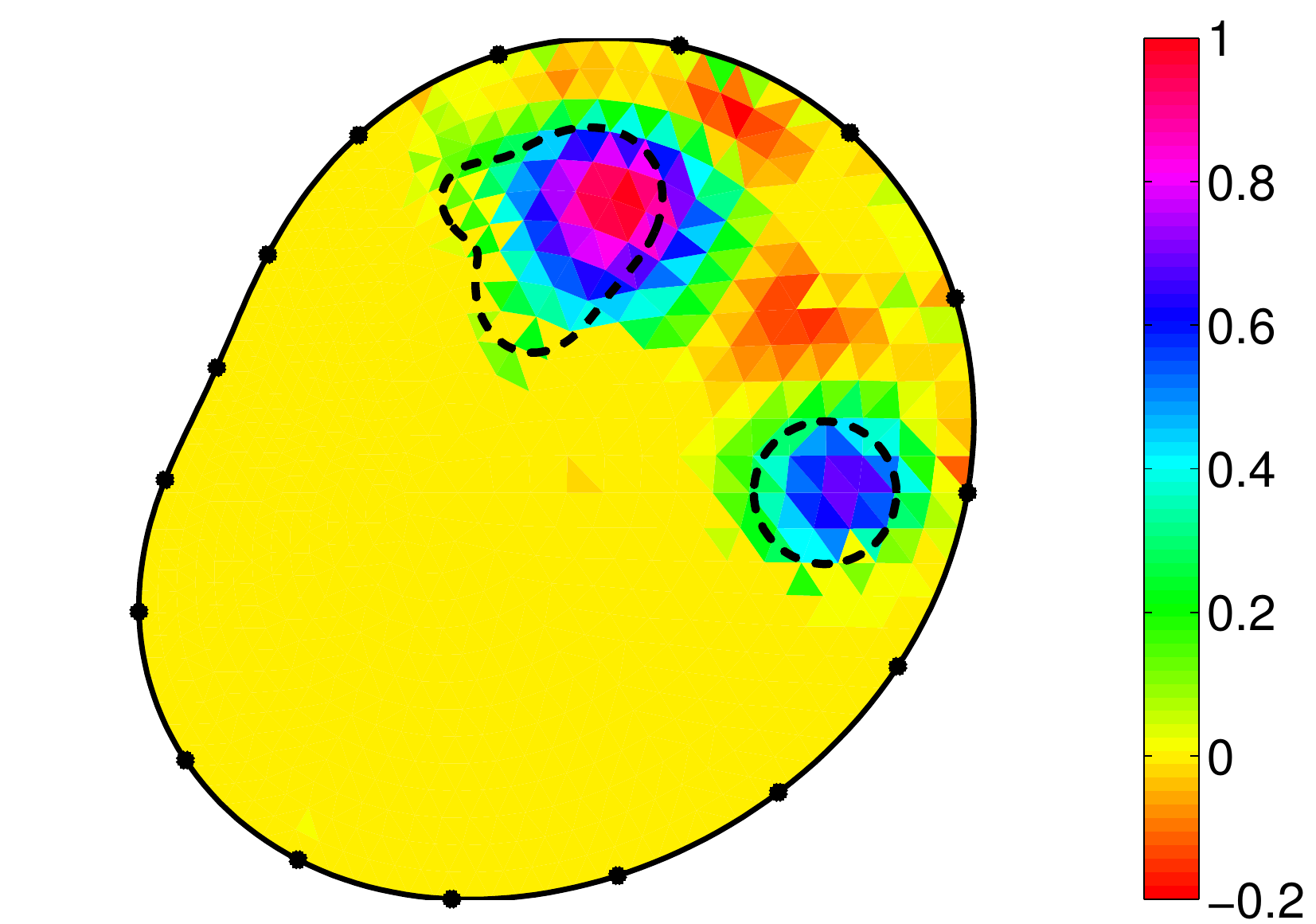}&
\includegraphics[keepaspectratio=true,height=1.5cm]{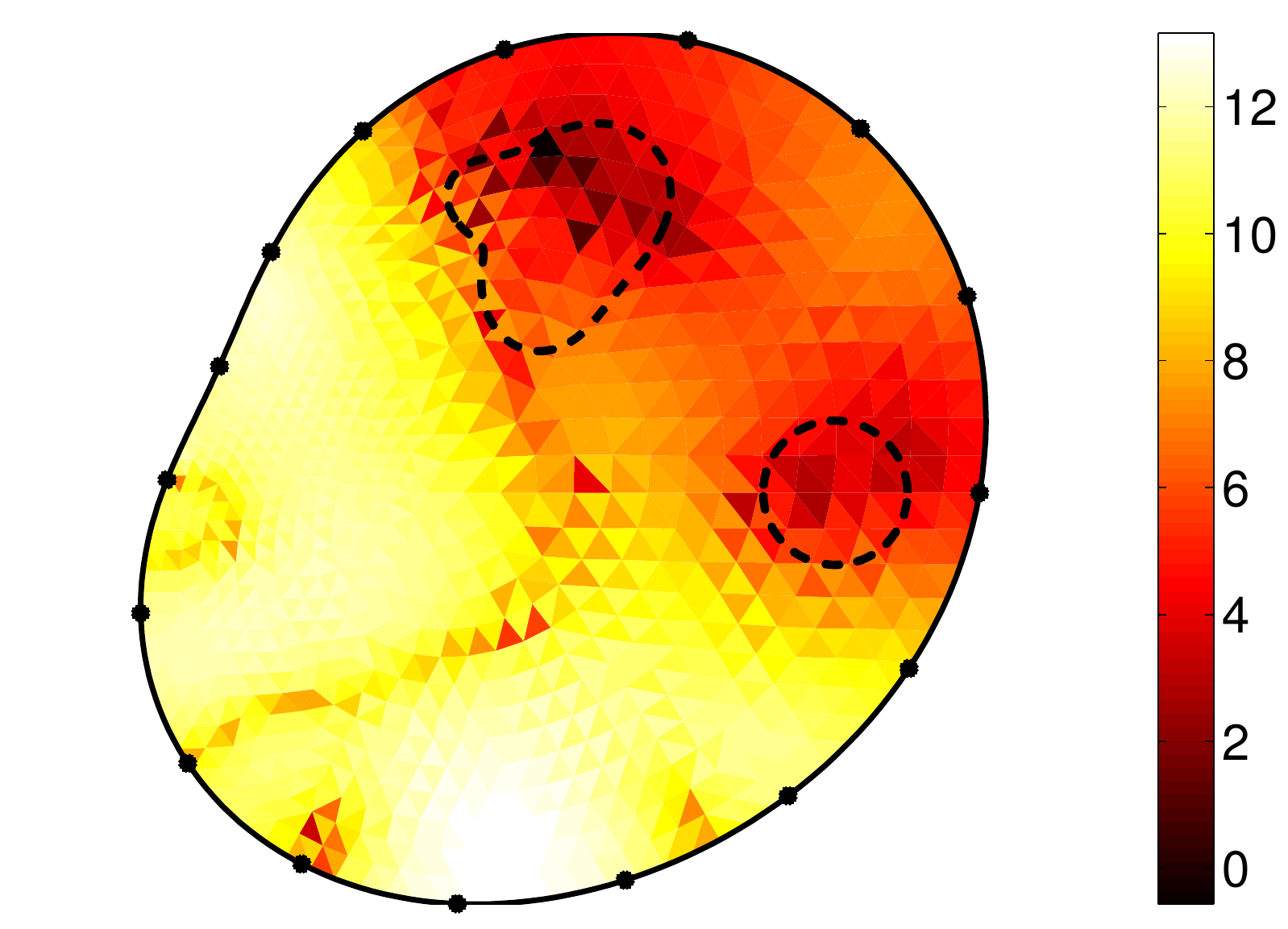}\\
\hline
\raisebox{4ex}{\footnotesize\begin{tabular}{c}
(h)
\end{tabular}}
&
\includegraphics[keepaspectratio=true,height=1.5cm]{Fig/plot_final_20120602/deltasigma_8_woColorbar-eps-converted-to.pdf}&
\includegraphics[keepaspectratio=true,height=1.5cm]{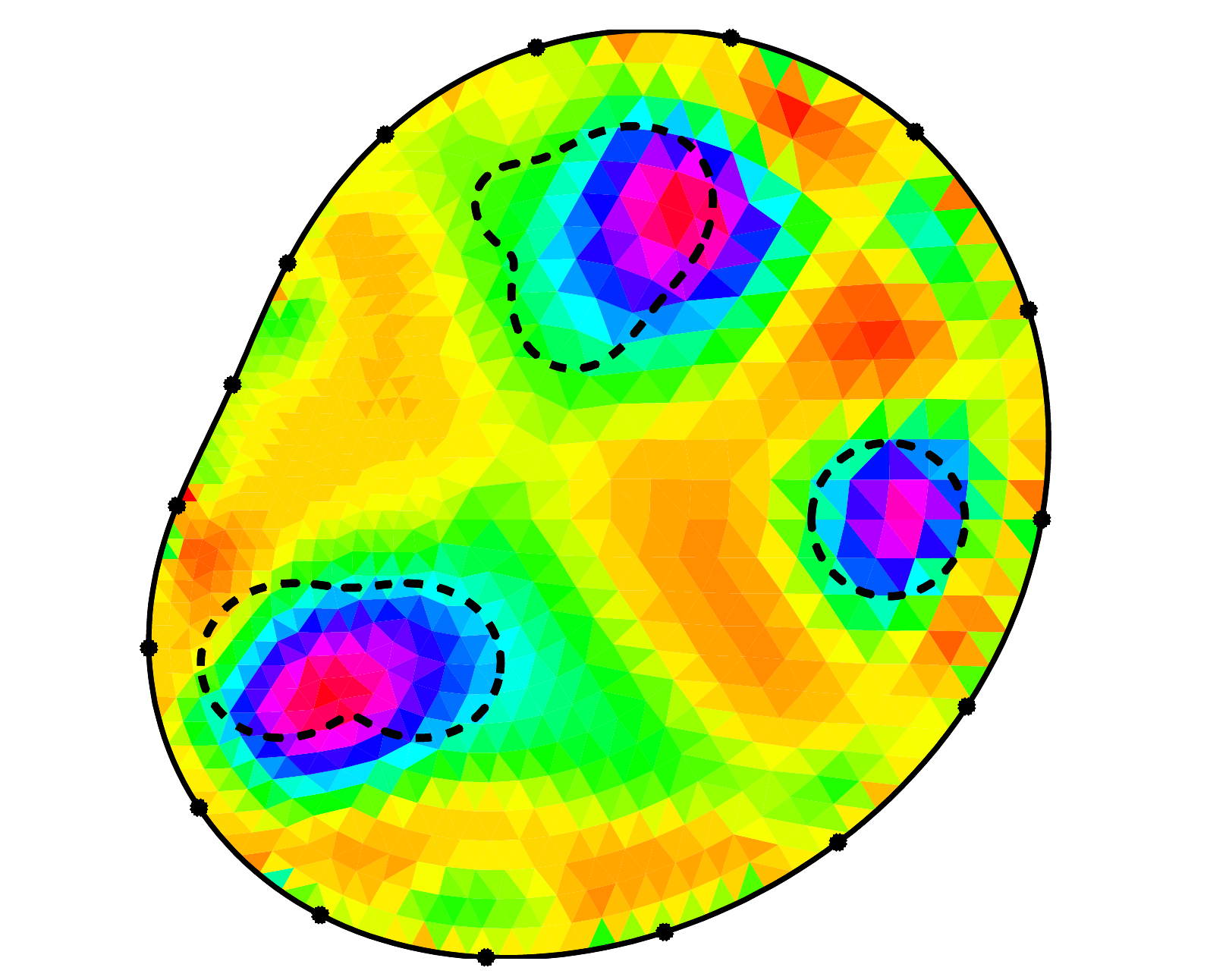}&
\includegraphics[keepaspectratio=true,height=1.5cm]{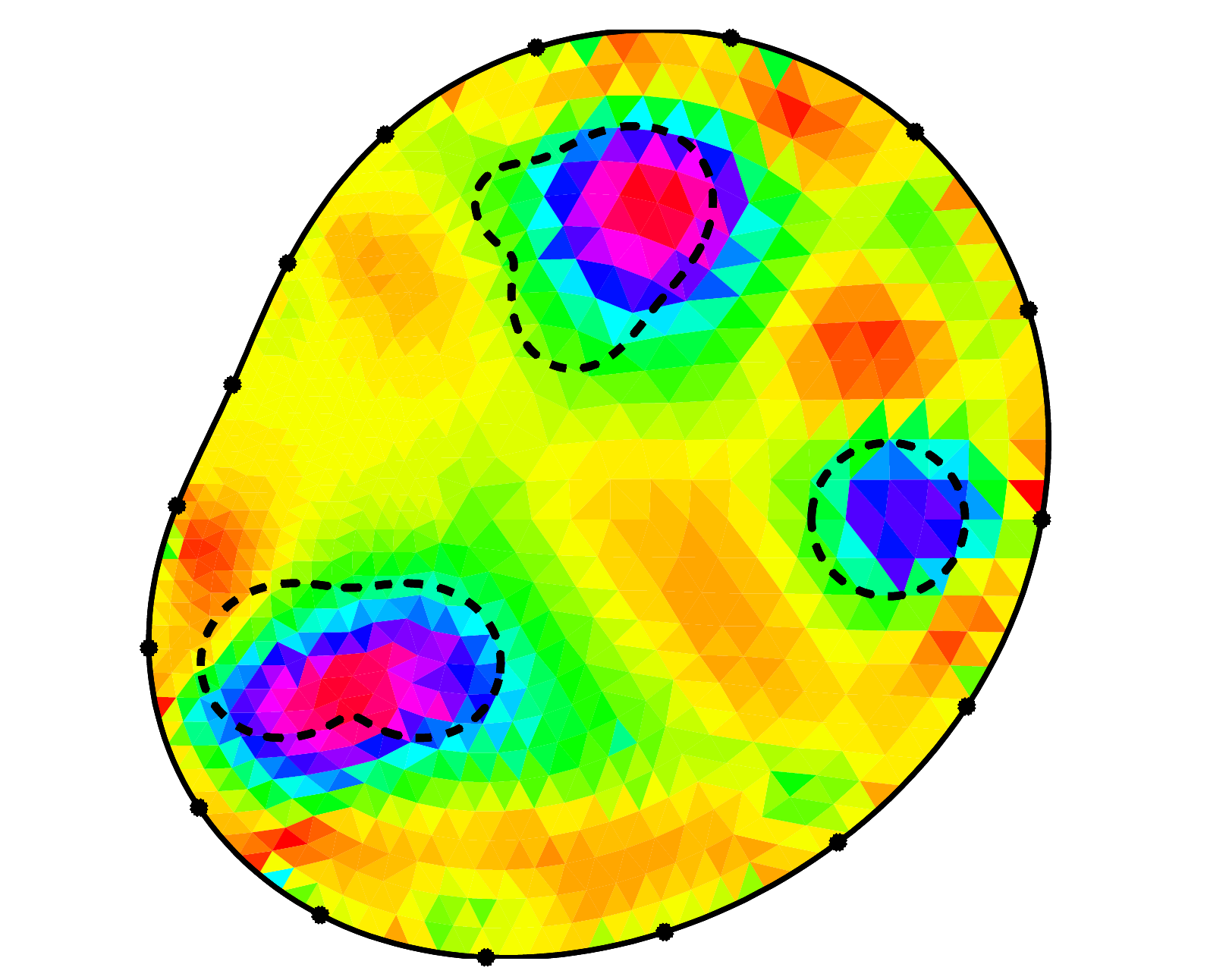}&
\includegraphics[keepaspectratio=true,height=1.5cm]{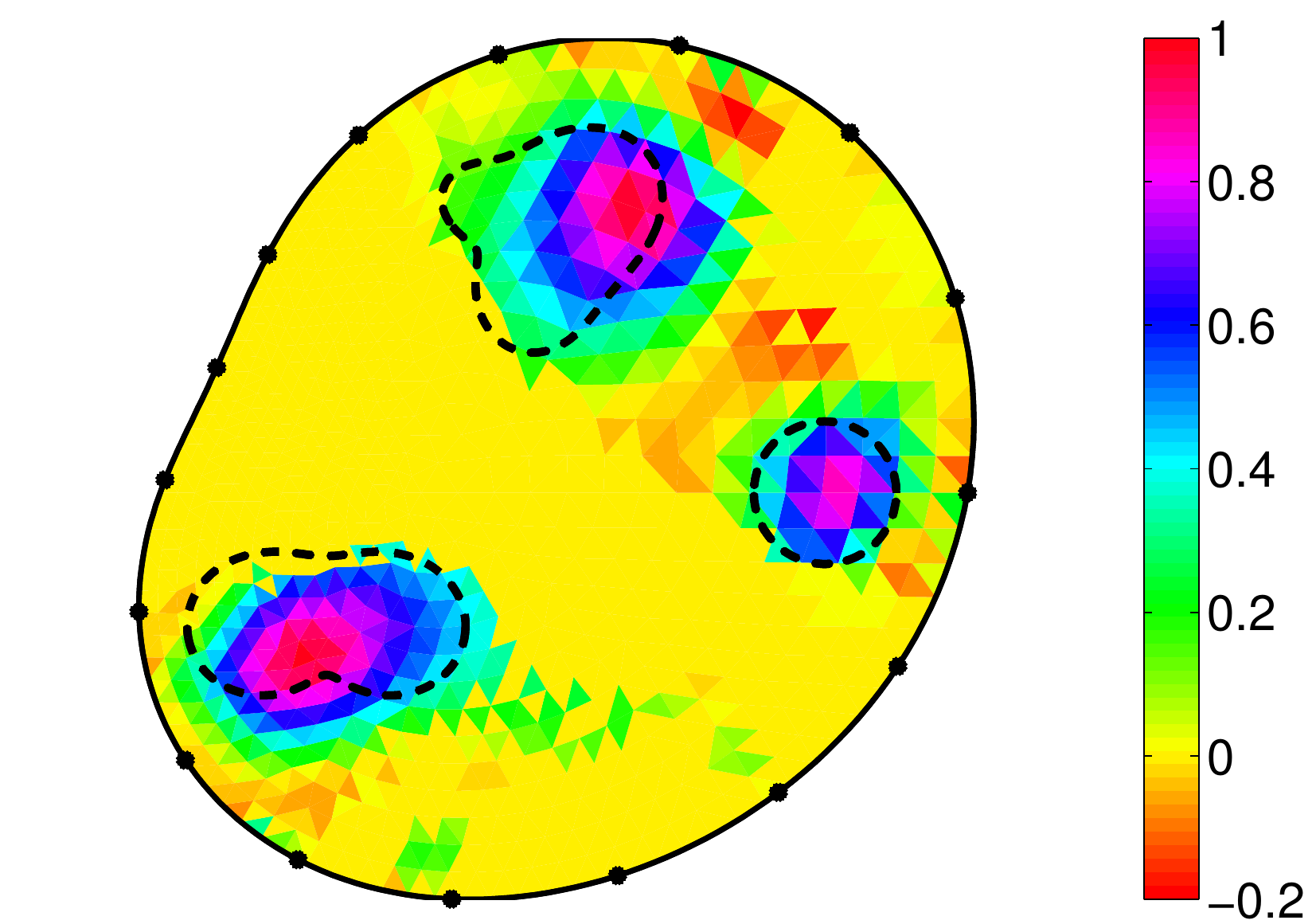}&
\includegraphics[keepaspectratio=true,height=1.5cm]{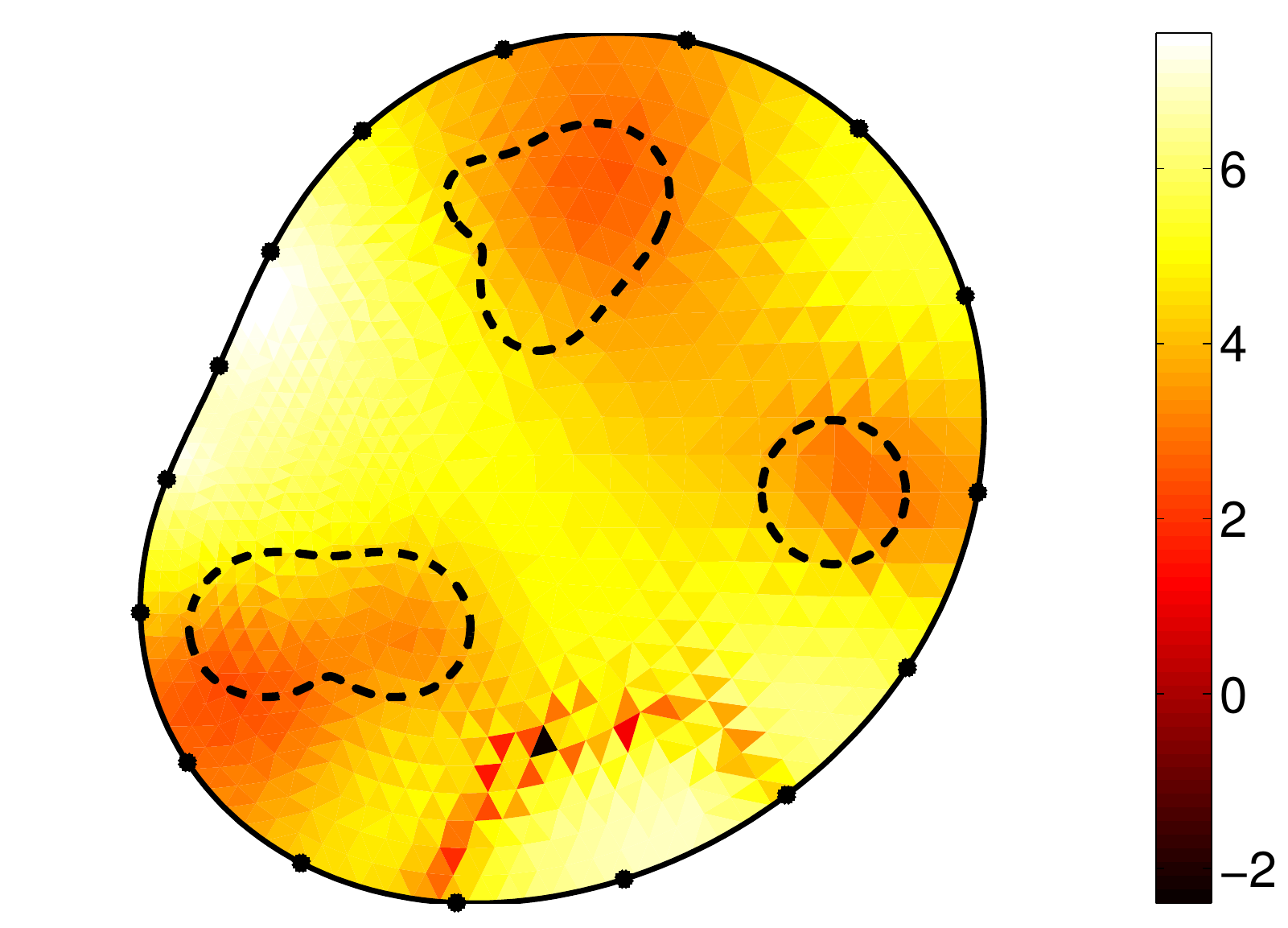}\\
\hline
\end{tabular}
\caption{ Reconstructed difference EIT images in non-circular domain with a data which adds 1$\%$ random noise. $\DS$: true difference image,
$\DS_{S}$: standard linearized method, cf.\ \eref{LM2},
$\DS_{B}$: naive combination of LM and S-FM, cf.\ \eref{regularization},
$\DS_{A}$: proposed combination of LM and S-FM, cf.\ \eref{recon},
$\mathbf{W1}$: S-FM alone, cf.\ \eref{SFMimage}.}
\label{recon_noncircle_1p}
\end{figure}

\begin{figure}
\centering
\begin{tabular}{|c|cccc|c|}
\hline
Case &  $\DS$ & $\DS_{S}$ & $\DS_{B}$ &  $\DS_{A}$ &$\mathbf{W1}$ \\
\hline
\raisebox{4ex}{\footnotesize\begin{tabular}{c}
(a)
\end{tabular}}
&
\includegraphics[keepaspectratio=true,width=1.8cm]{Fig/plot_final_20120602/deltasigma_1_woColorbar-eps-converted-to.pdf}&
\includegraphics[keepaspectratio=true,width=1.8cm]{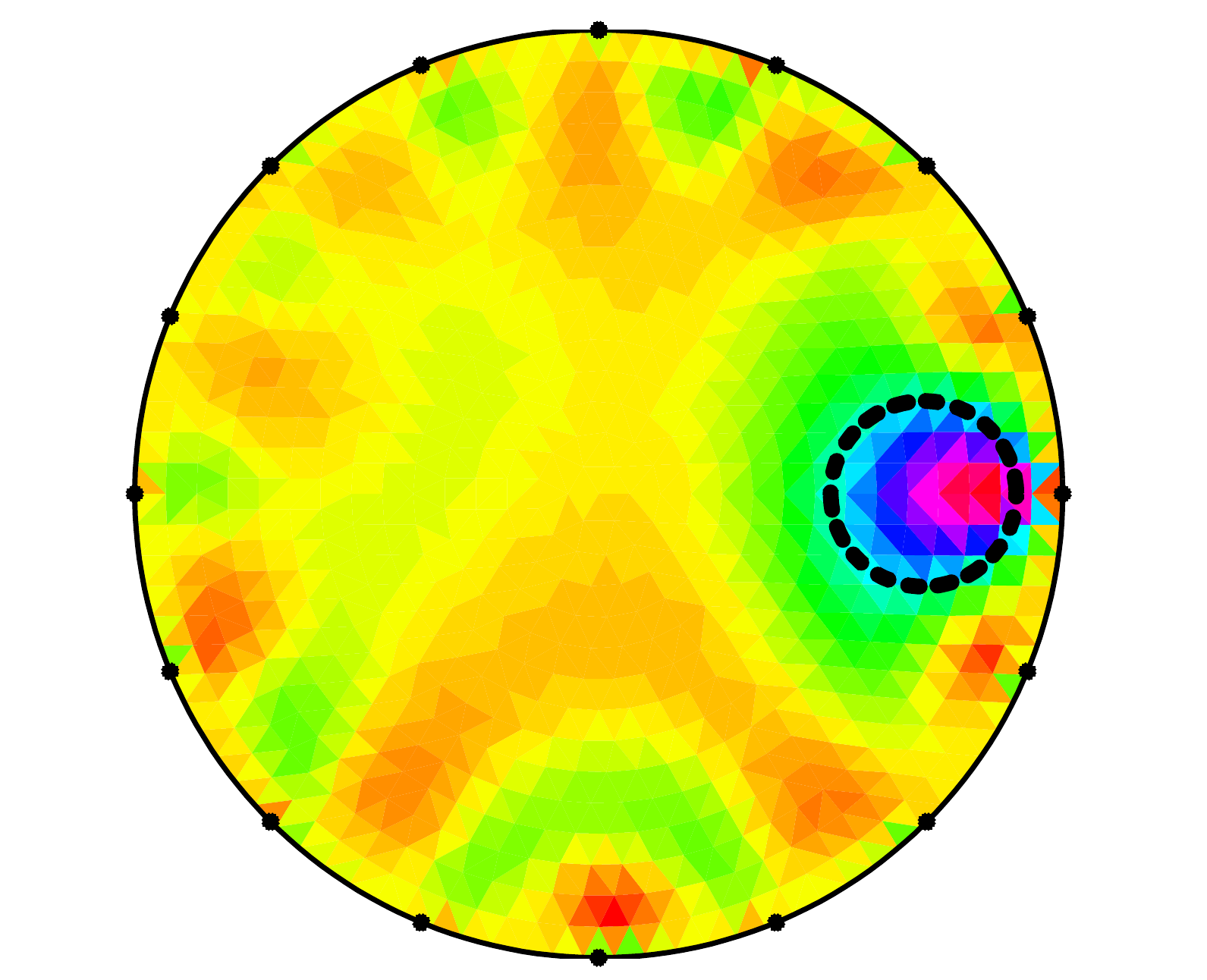}&
\includegraphics[keepaspectratio=true,width=1.8cm]{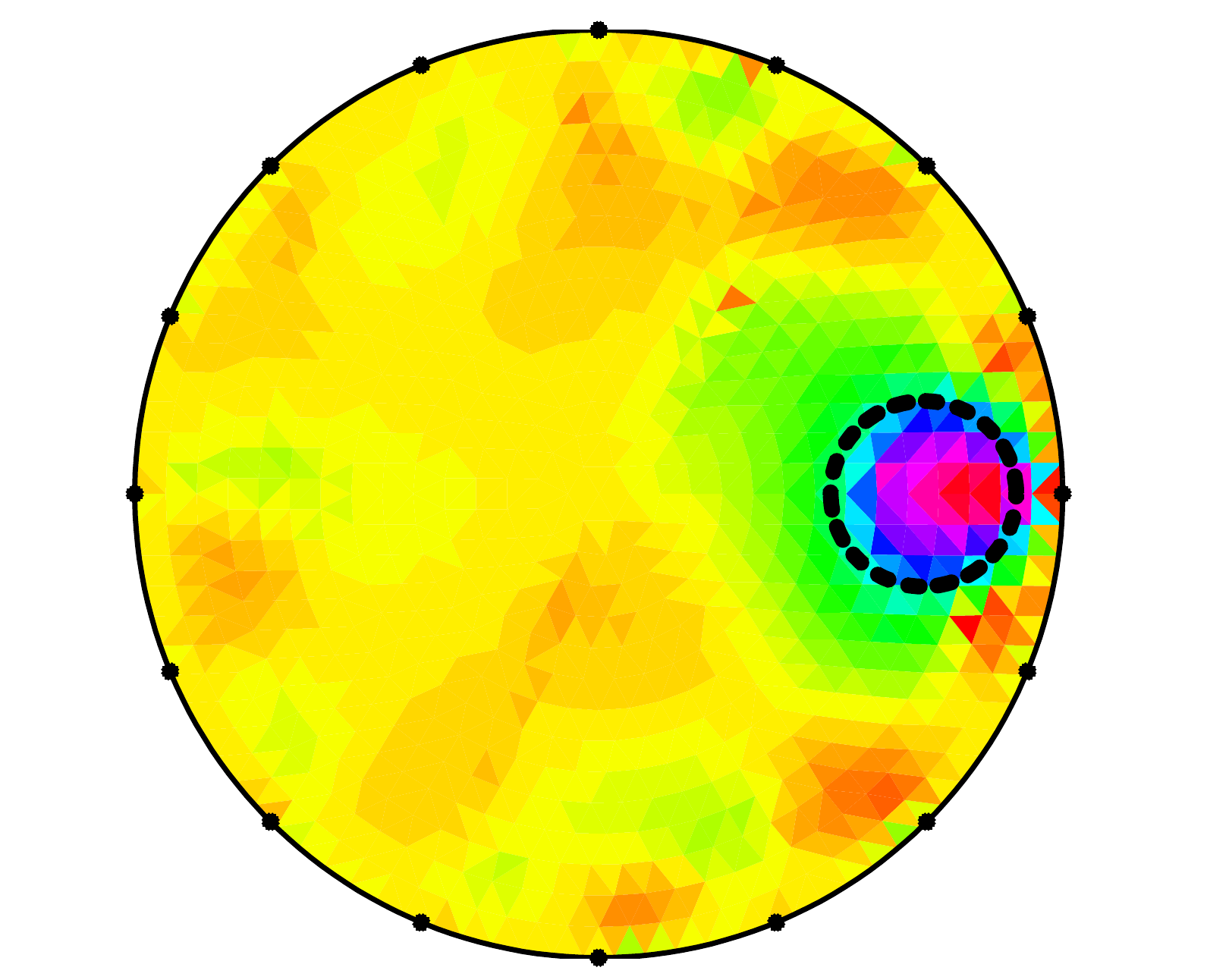}&
\includegraphics[keepaspectratio=true,width=2.1cm]{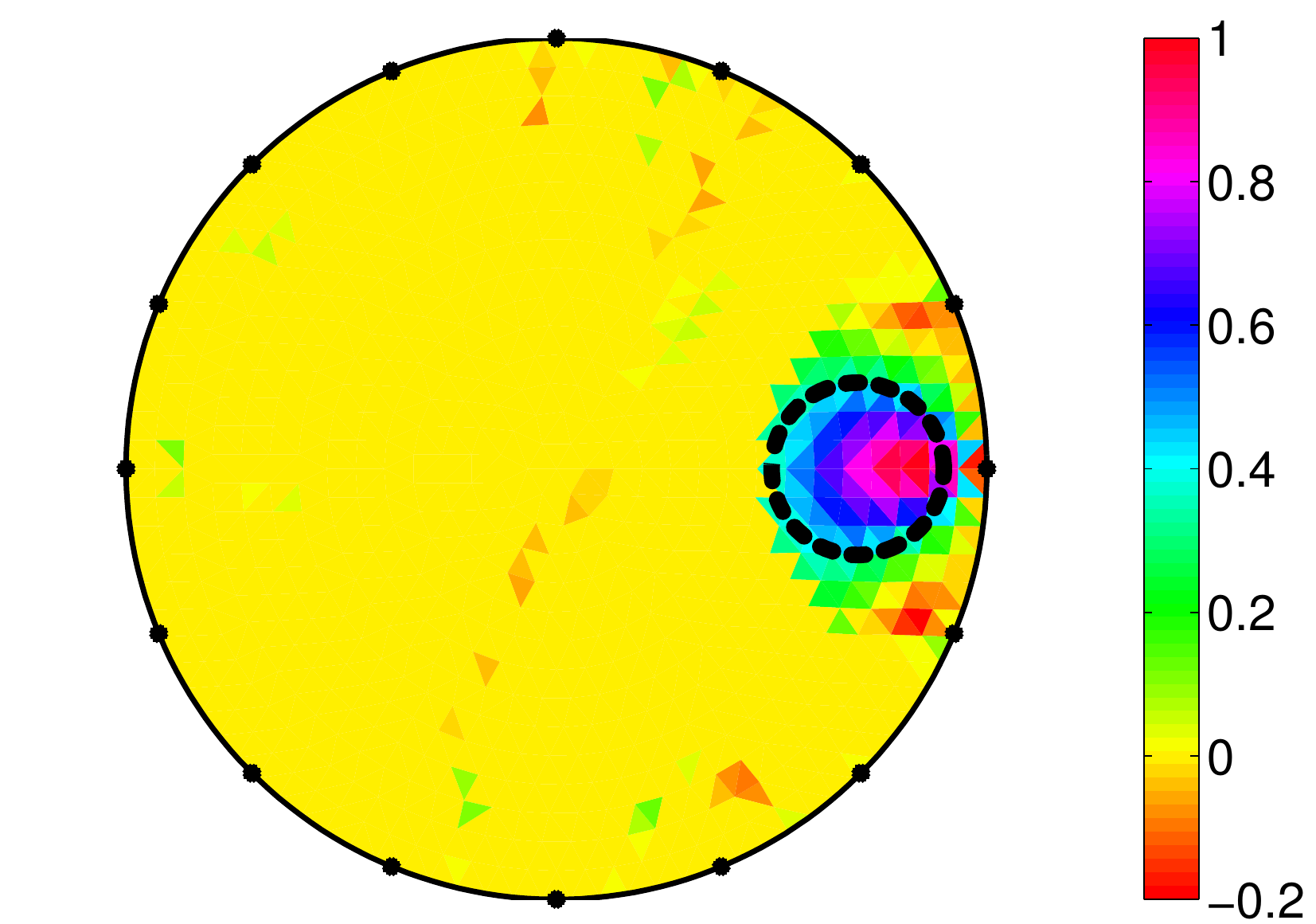}&
\includegraphics[keepaspectratio=true,width=2.cm]{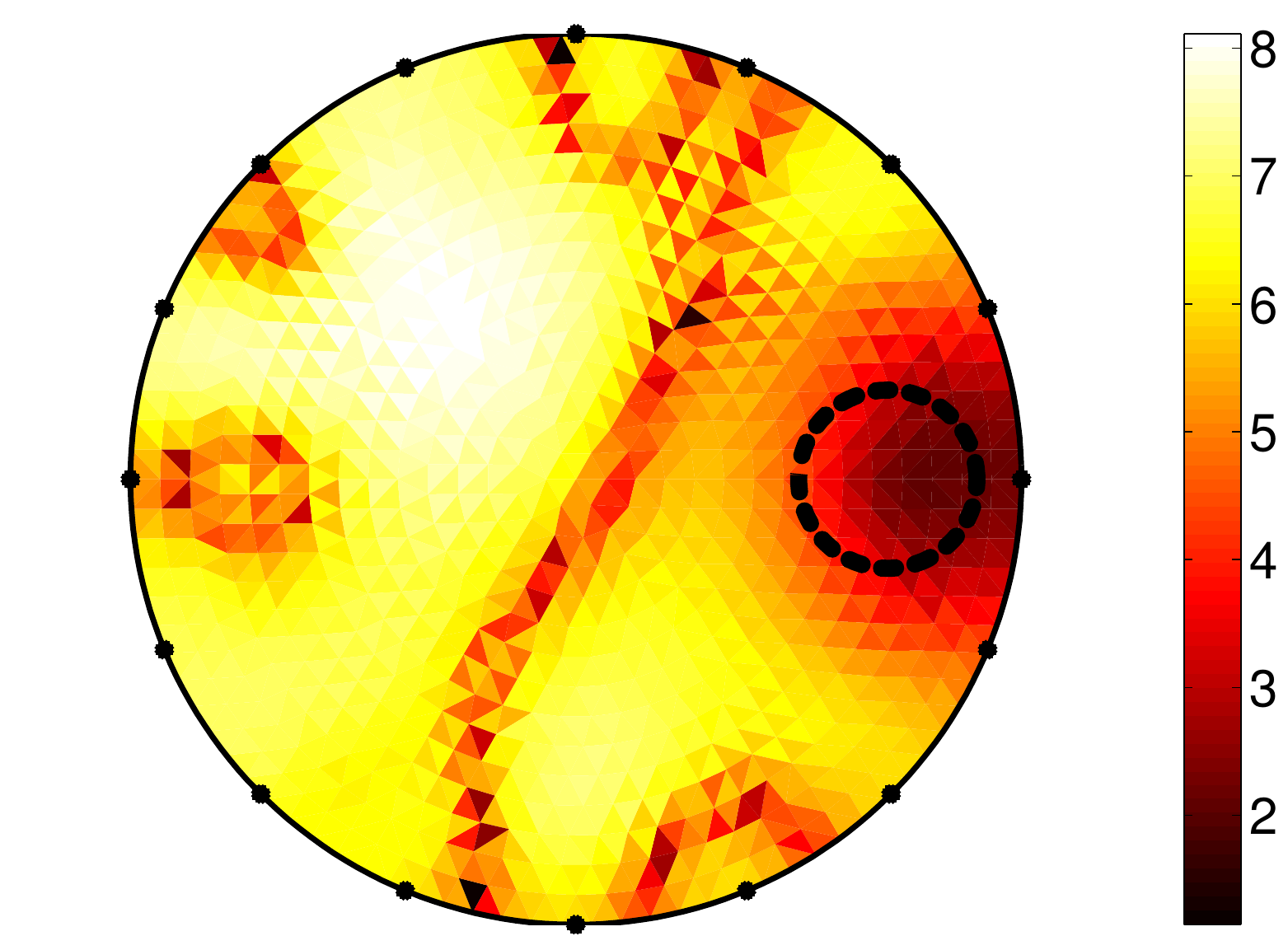}\\
\hline
\raisebox{4ex}{\footnotesize\begin{tabular}{c}
(b)
\end{tabular}}
&
\includegraphics[keepaspectratio=true,height=1.5cm]{Fig/plot_final_20120602/deltasigma_2_woColorbar-eps-converted-to.pdf}&
\includegraphics[keepaspectratio=true,height=1.5cm]{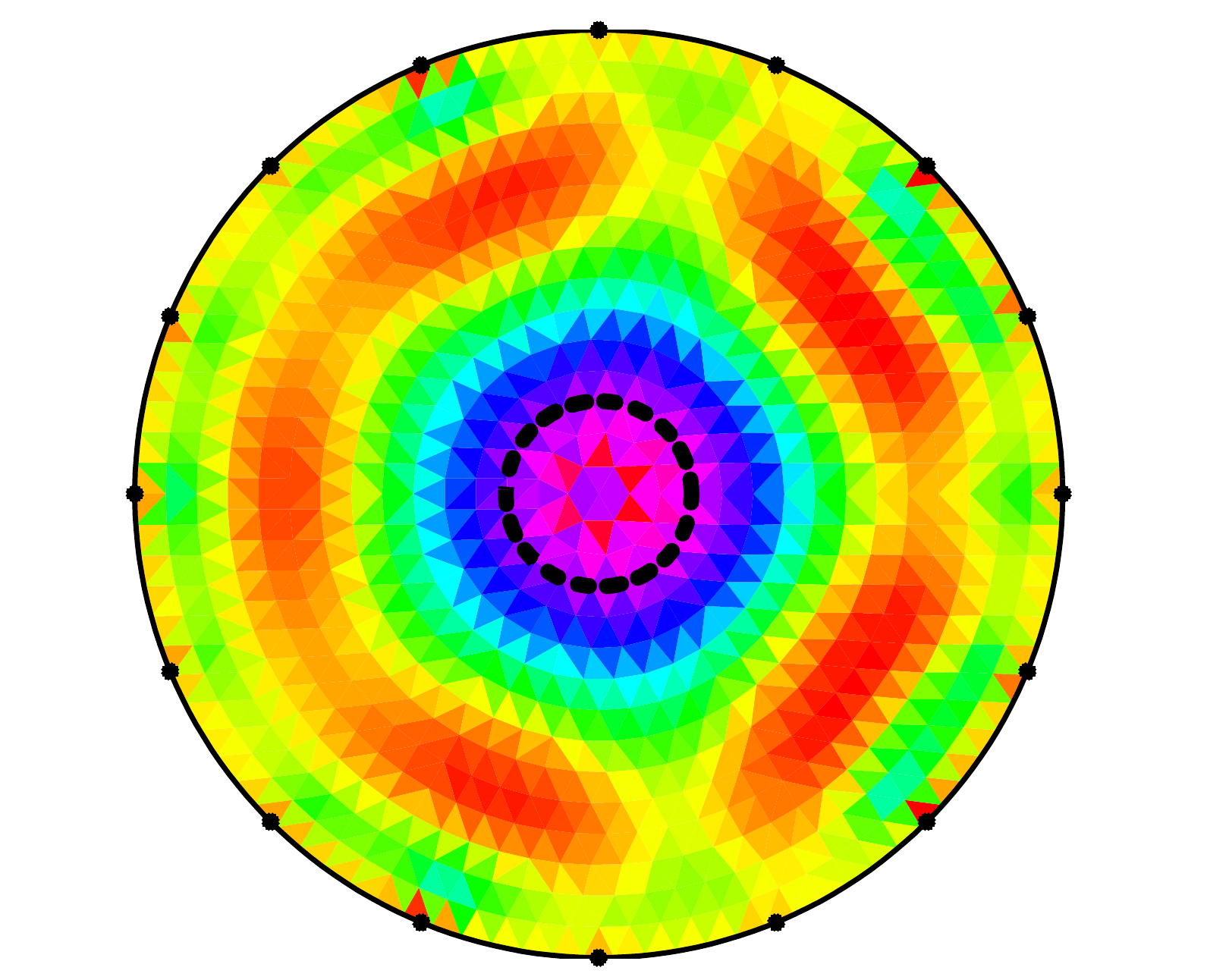}&
\includegraphics[keepaspectratio=true,height=1.5cm]{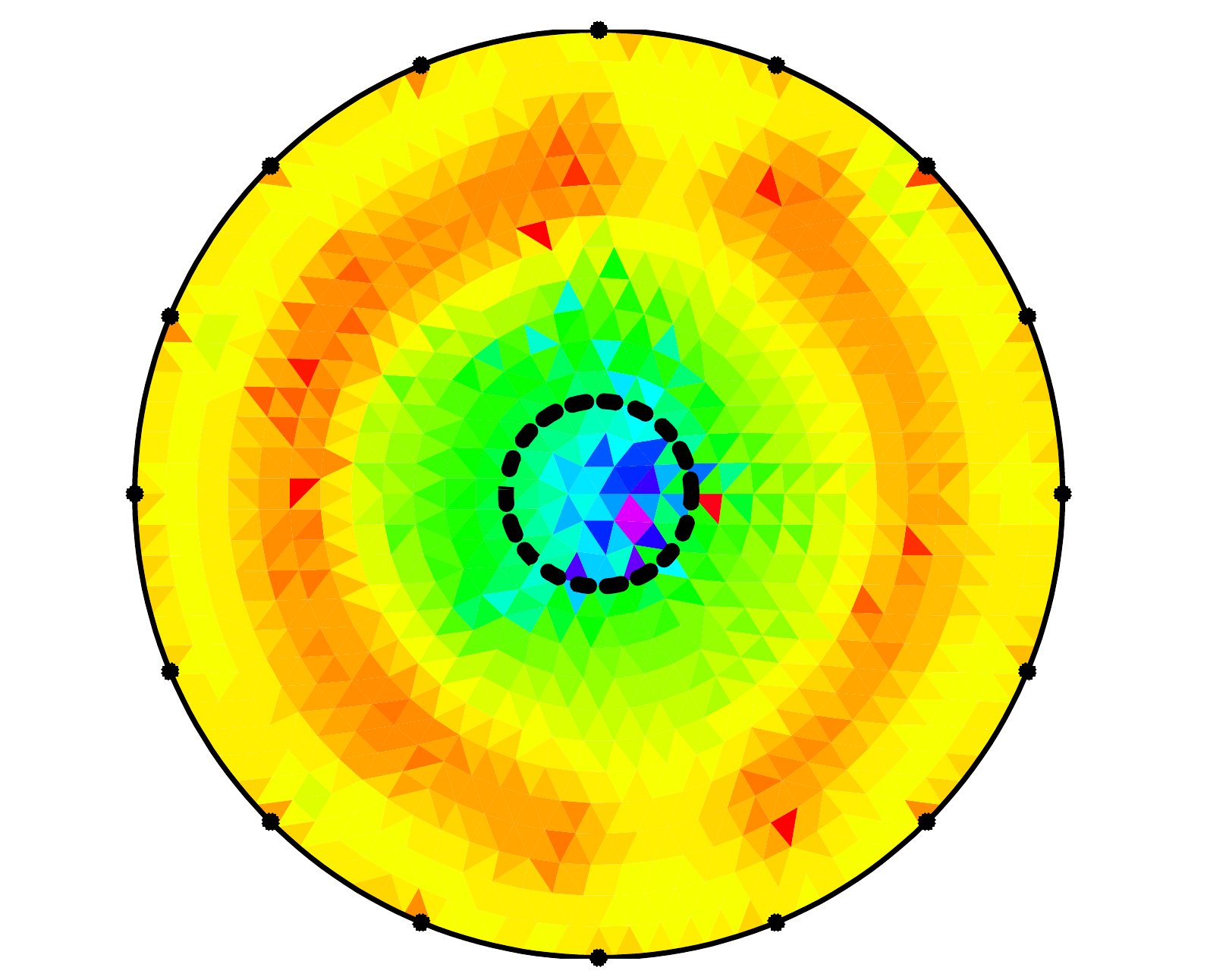}&
\includegraphics[keepaspectratio=true,height=1.5cm]{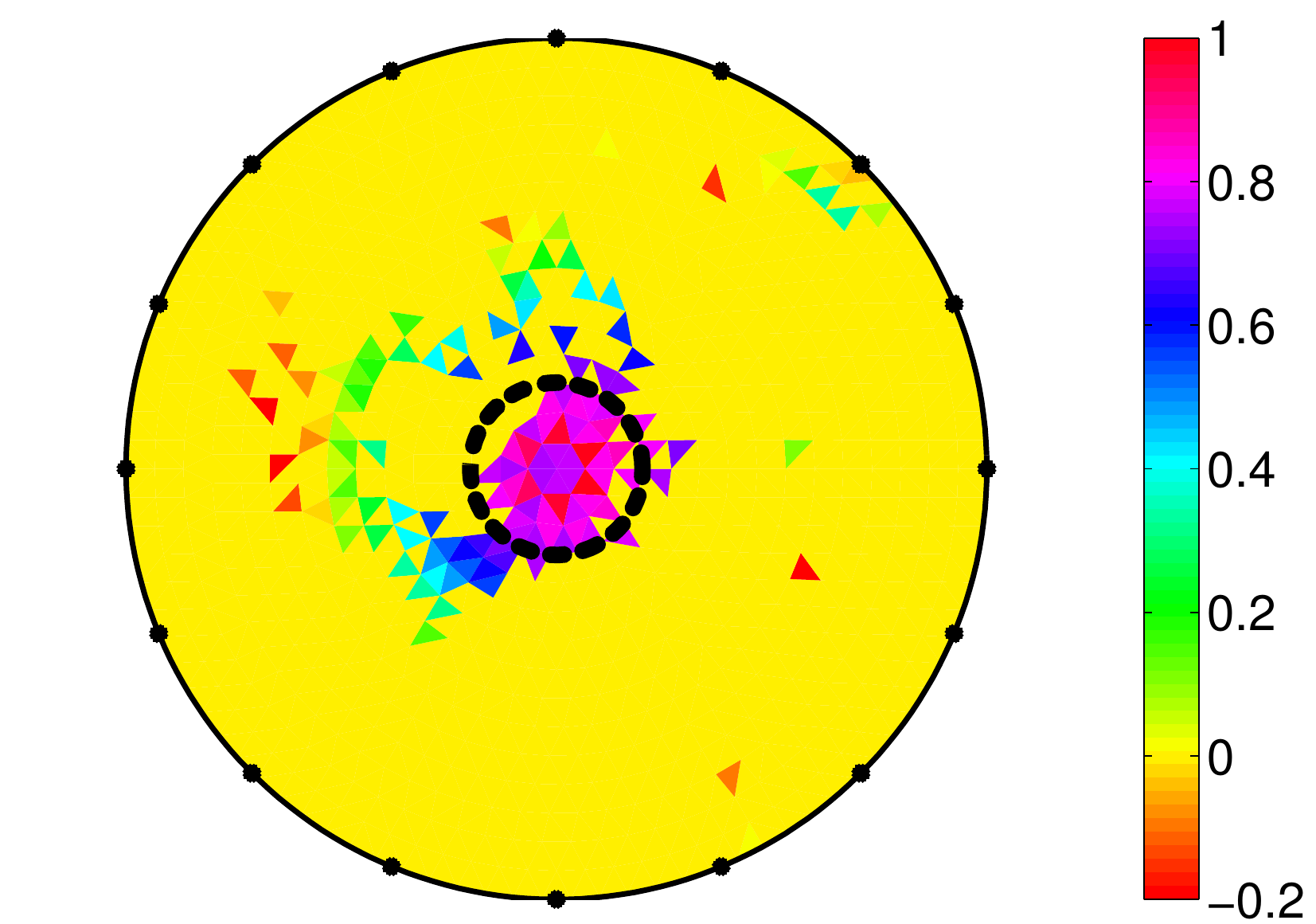}&
\includegraphics[keepaspectratio=true,height=1.5cm]{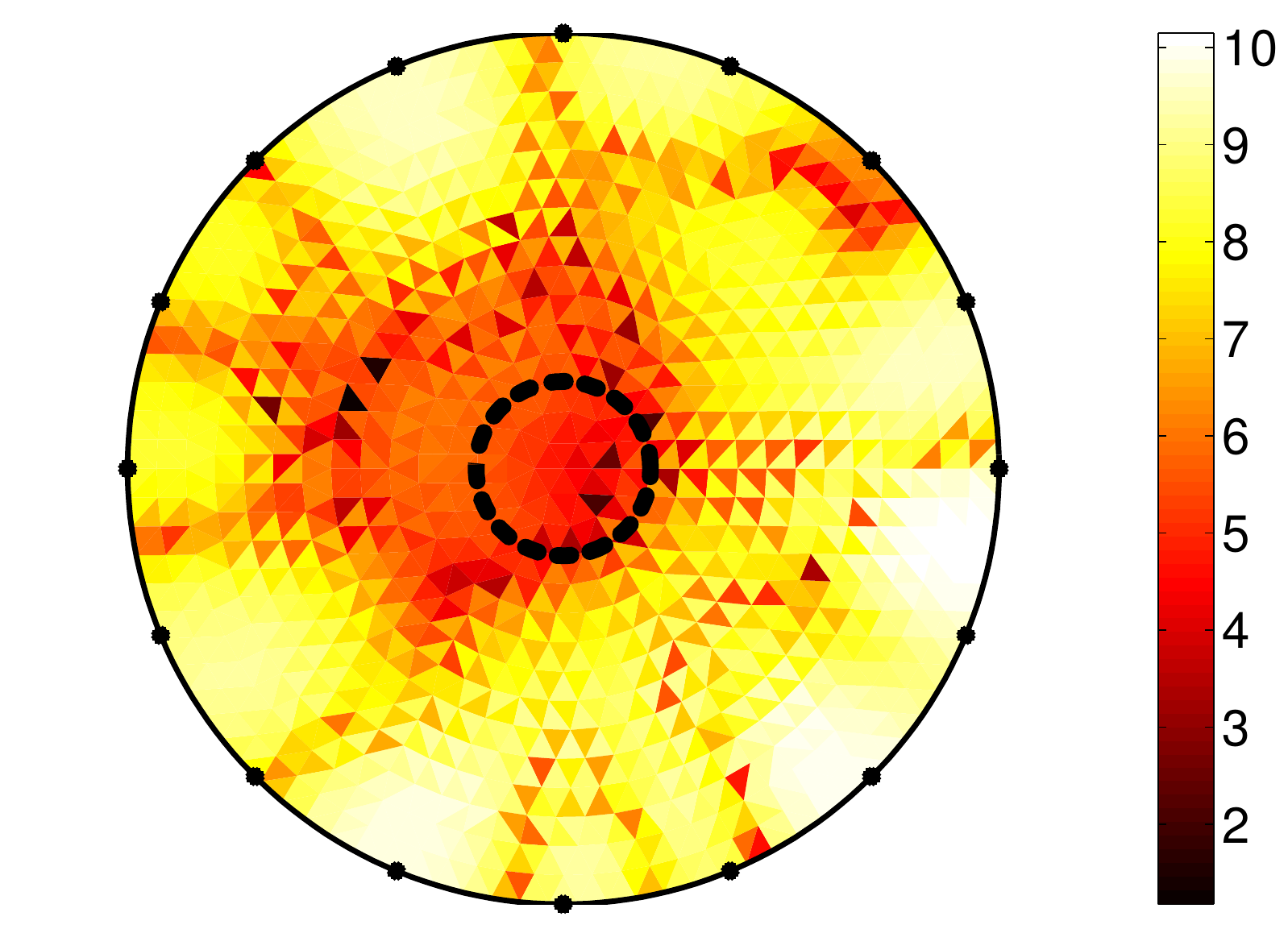}\\
\hline
\raisebox{4ex}{\footnotesize\begin{tabular}{c}
(c)
\end{tabular}}
&
\includegraphics[keepaspectratio=true,height=1.5cm]{Fig/plot_final_20120602/deltasigma_3_woColorbar-eps-converted-to.pdf}&
\includegraphics[keepaspectratio=true,height=1.5cm]{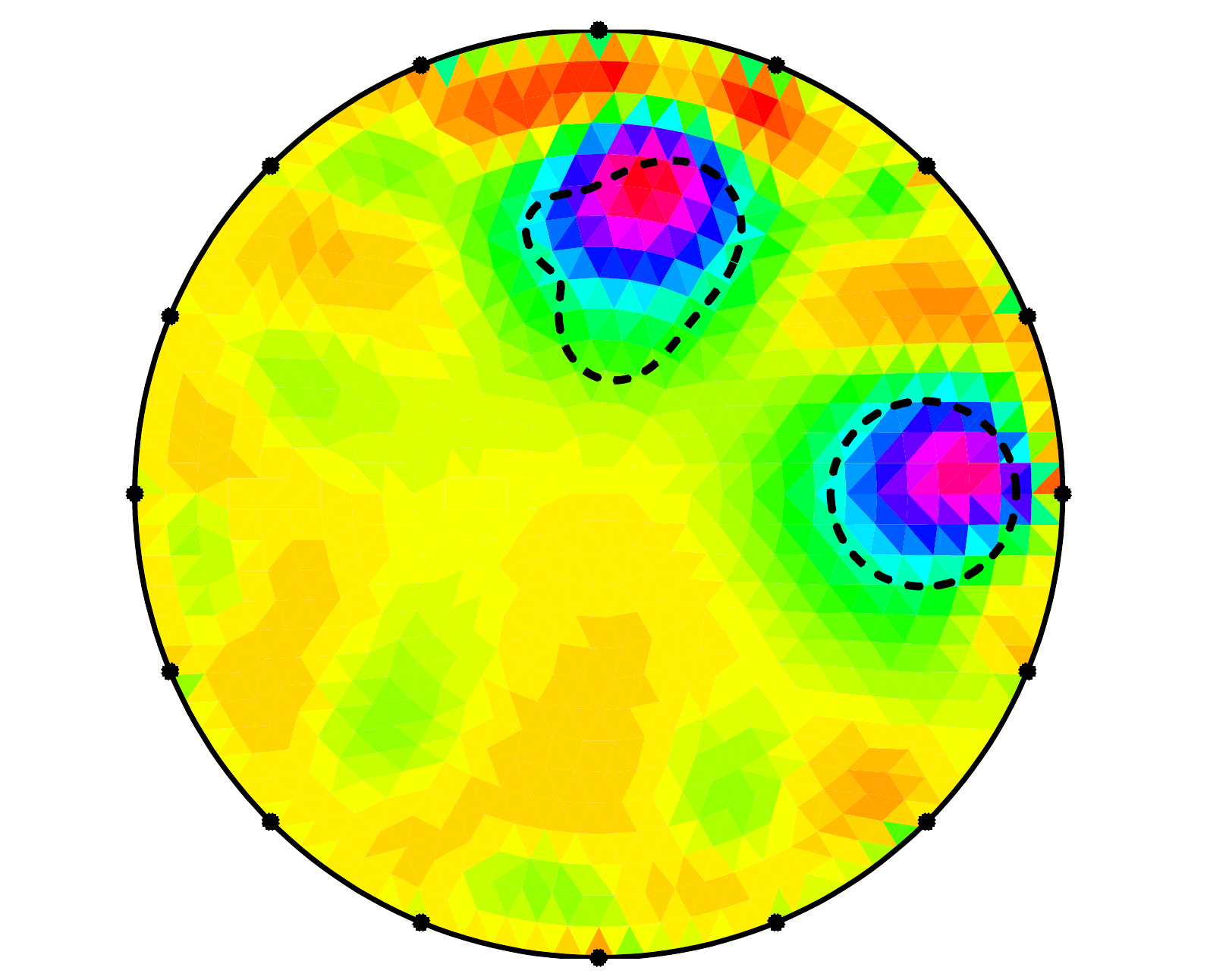}&
\includegraphics[keepaspectratio=true,height=1.5cm]{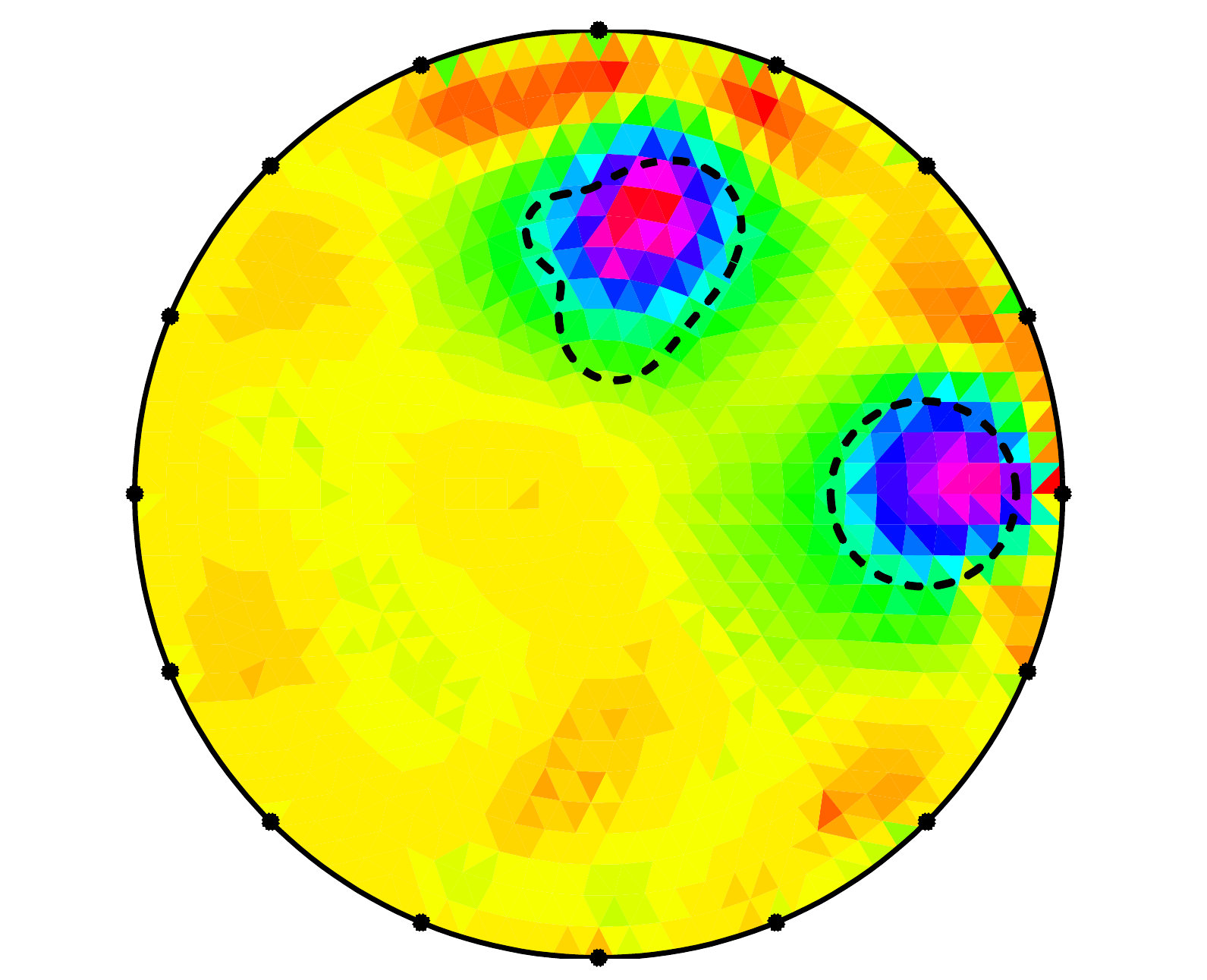}&
\includegraphics[keepaspectratio=true,height=1.5cm]{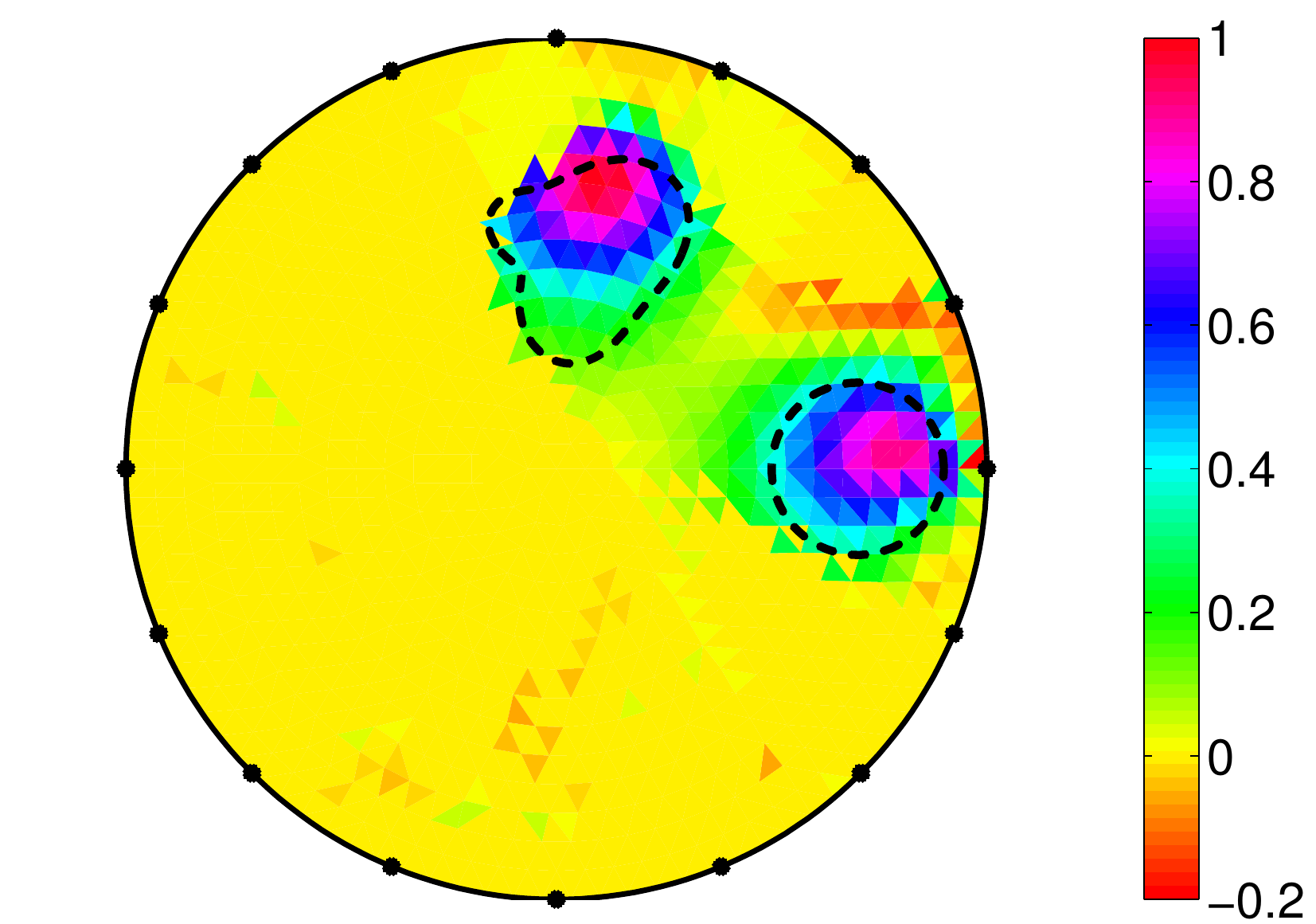}&
\includegraphics[keepaspectratio=true,height=1.5cm]{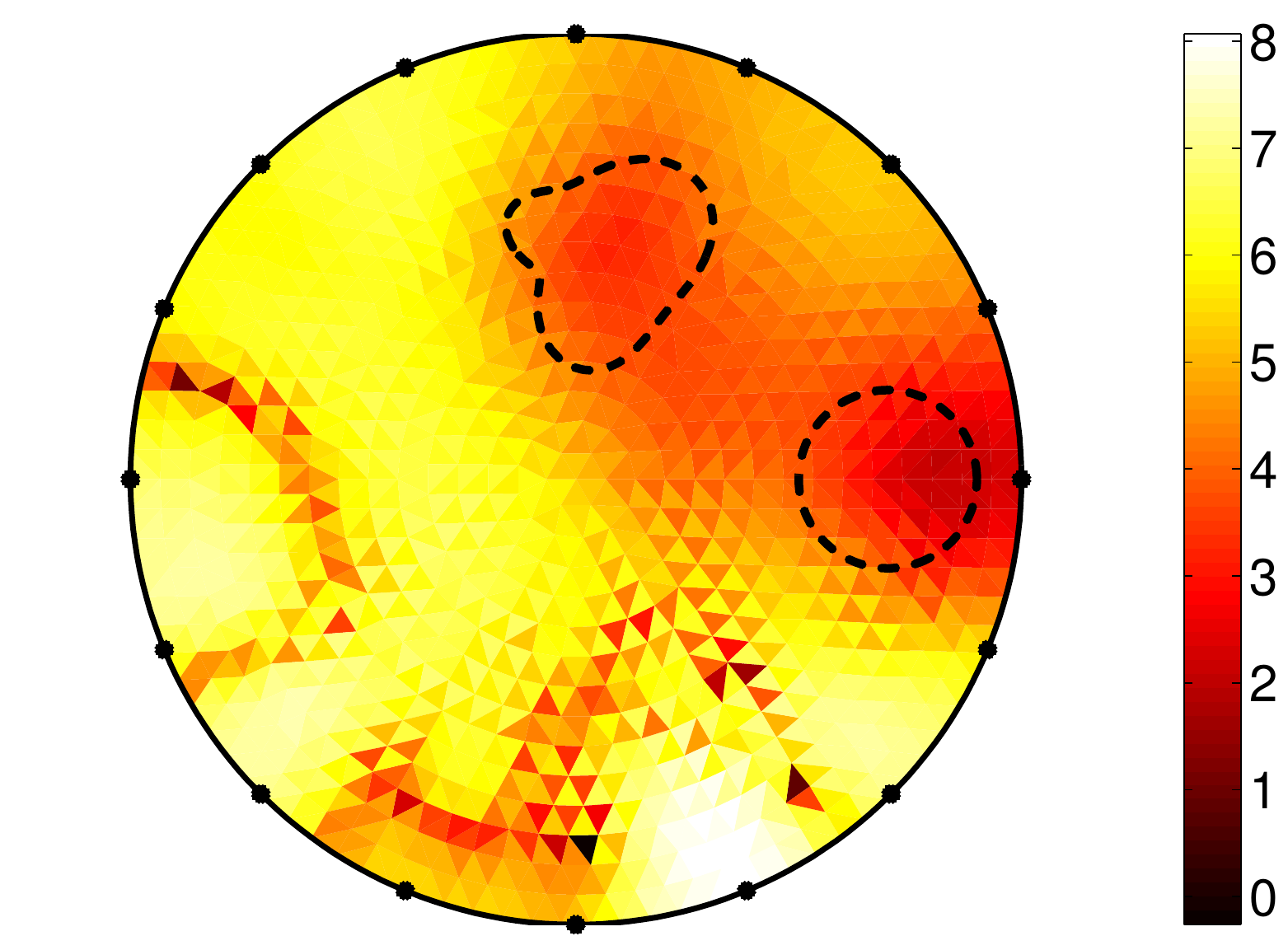}\\
\hline
\raisebox{4ex}{\footnotesize\begin{tabular}{c}
(d)
\end{tabular}}
&
\includegraphics[keepaspectratio=true,height=1.5cm]{Fig/plot_final_20120602/deltasigma_4_woColorbar-eps-converted-to.pdf}&
\includegraphics[keepaspectratio=true,height=1.5cm]{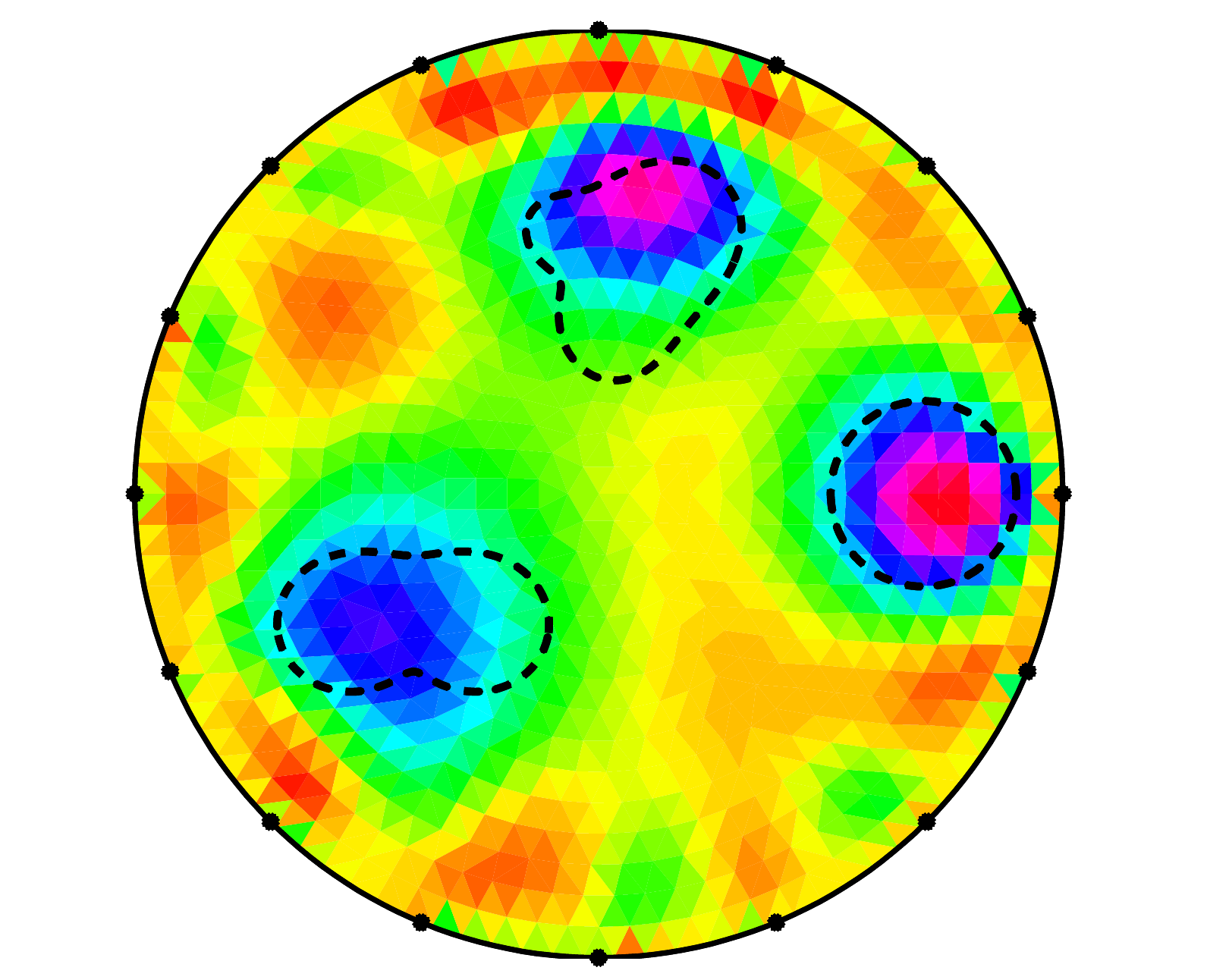}&
\includegraphics[keepaspectratio=true,height=1.5cm]{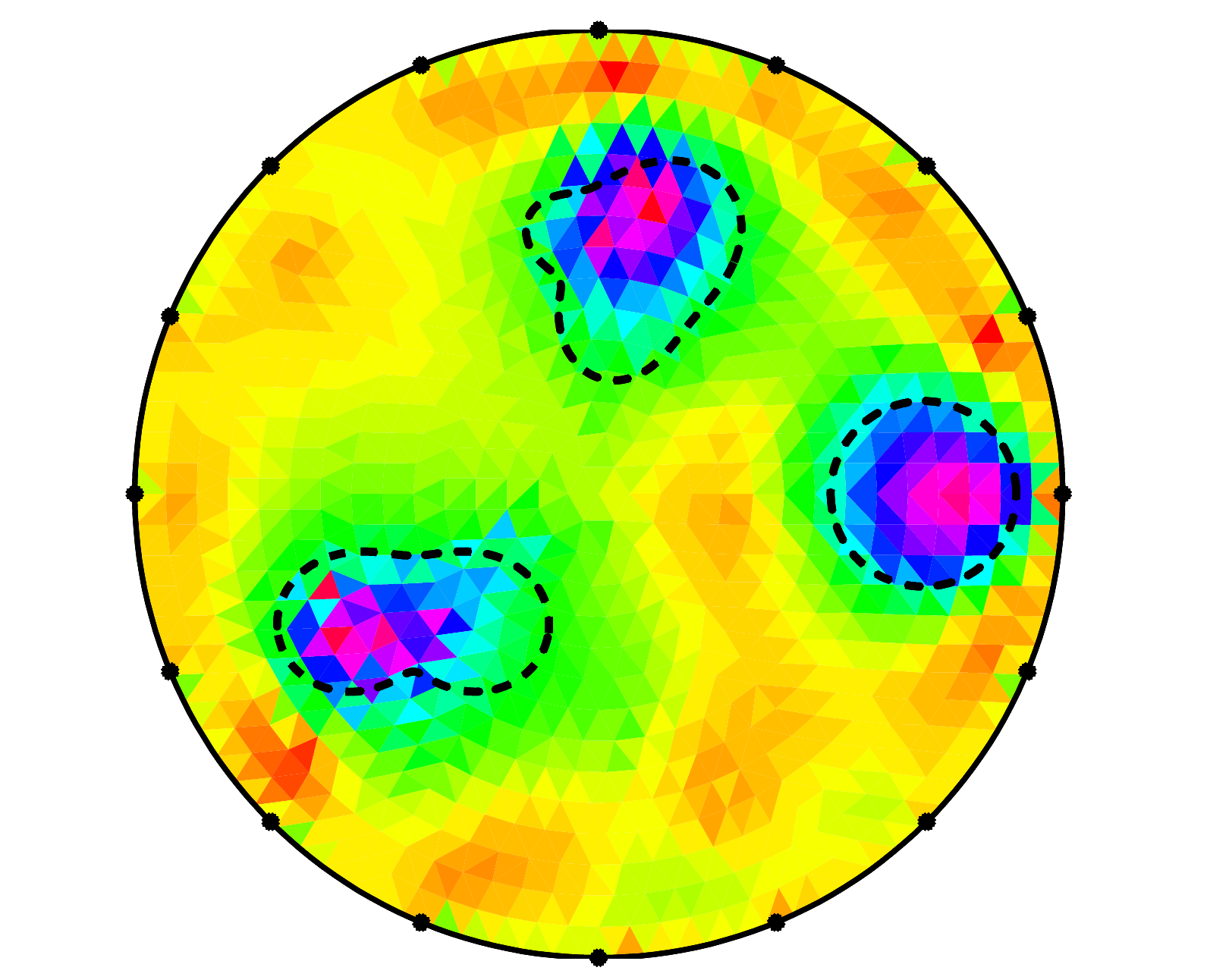}&
\includegraphics[keepaspectratio=true,height=1.5cm]{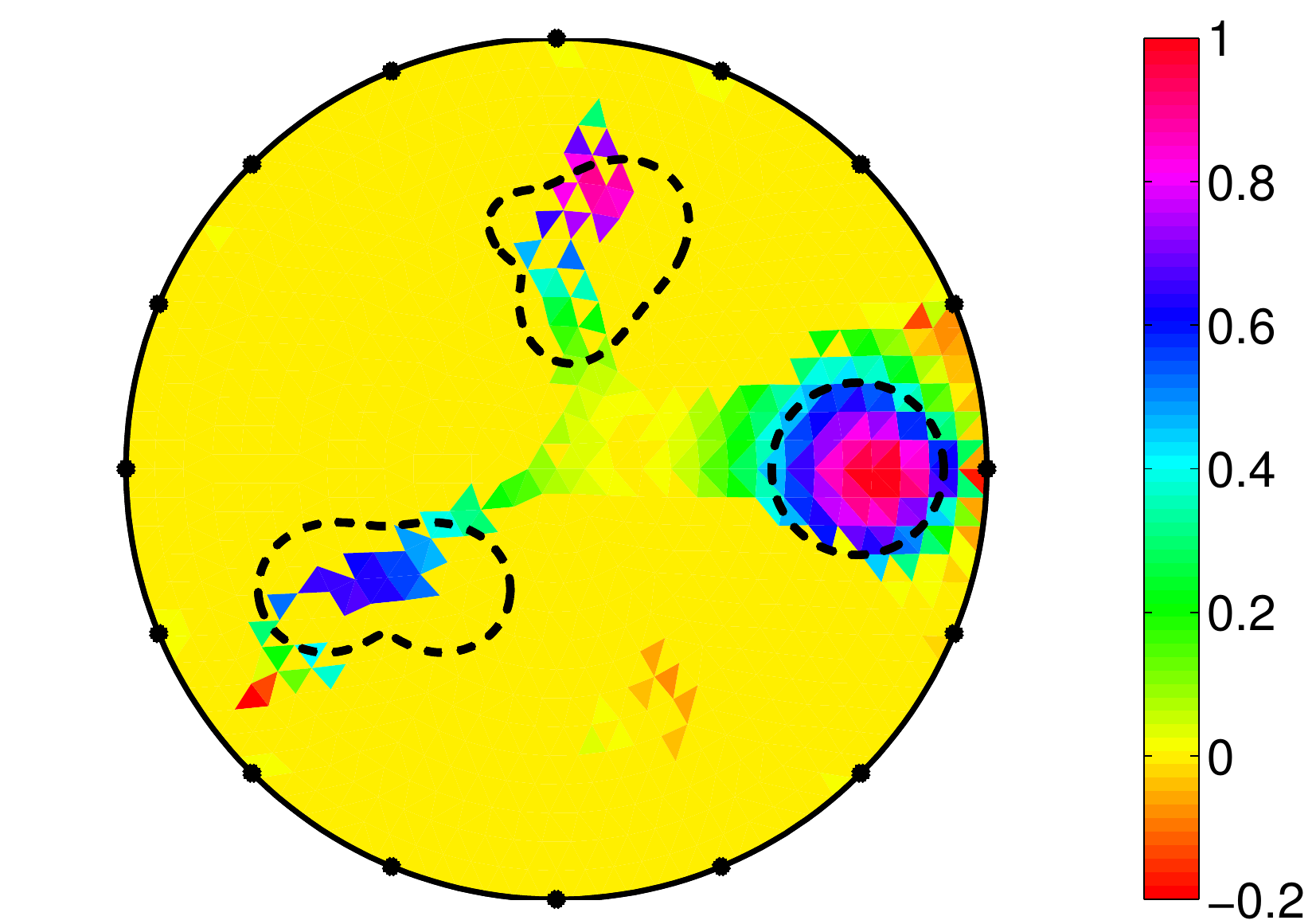}&
\includegraphics[keepaspectratio=true,height=1.5cm]{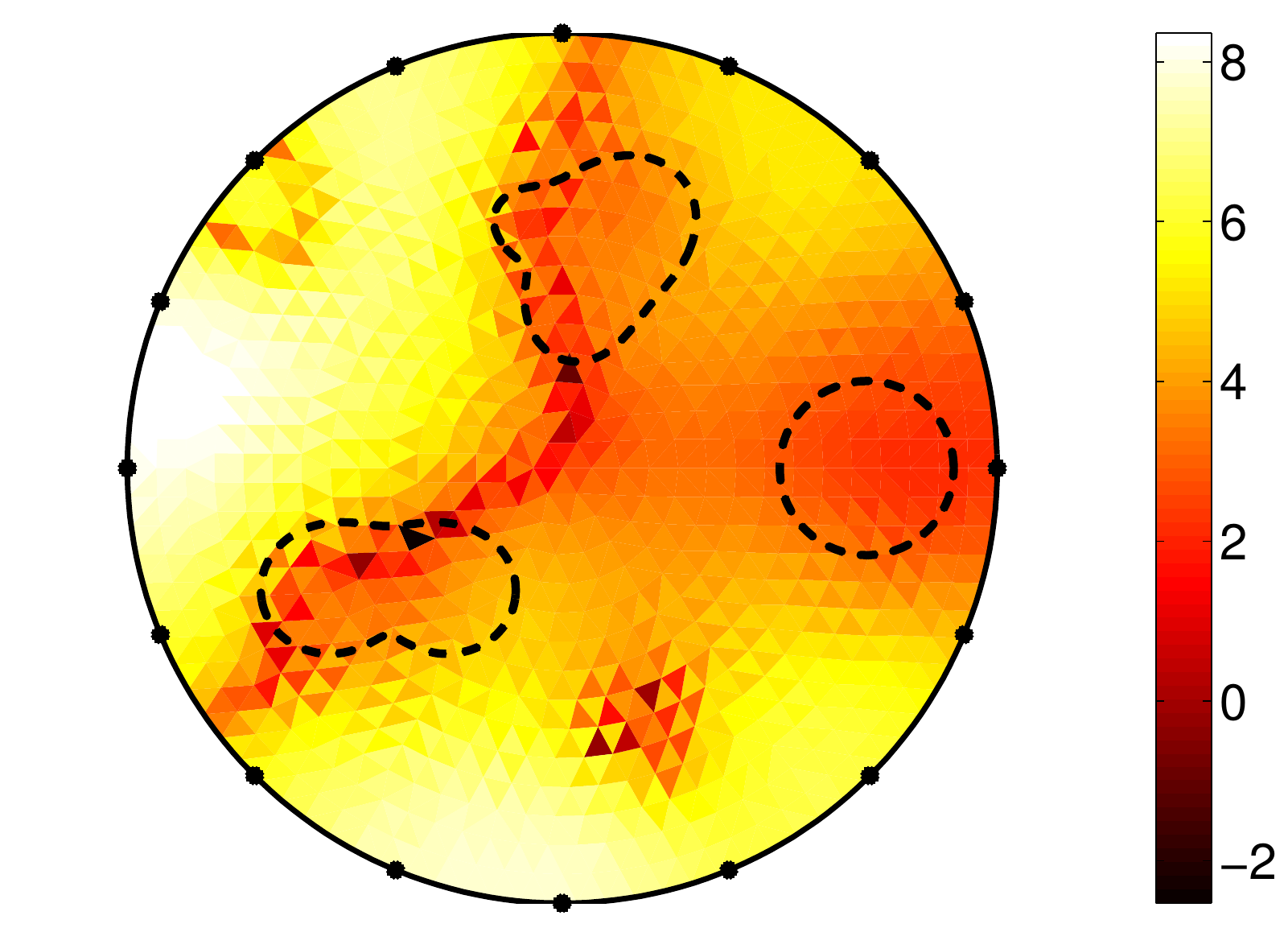}\\
\hline
\end{tabular}
\caption{ Reconstructed difference EIT images in circular domain with a data which adds 5$\%$ random noise. $\DS$: true difference image,
$\DS_{S}$: standard linearized method, cf.\ \eref{LM2},
$\DS_{B}$: naive combination of LM and S-FM, cf.\ \eref{regularization},
$\DS_{A}$: proposed combination of LM and S-FM, cf.\ \eref{recon},
$\mathbf{W1}$: S-FM alone, cf.\ \eref{SFMimage}.}
\label{recon_circle_5p}
\end{figure}

\begin{figure}
\centering
\begin{tabular}{|c|cccc|c|}
\hline
Case &  $\DS$ & $\DS_{S}$ & $\DS_{B}$ &  $\DS_{A}$ &$\mathbf{W1}$ \\
\hline
\raisebox{4ex}{\footnotesize\begin{tabular}{c}
(e)
\end{tabular}}
&
\includegraphics[keepaspectratio=true,height=1.5cm]{Fig/plot_final_20120602/deltasigma_5_woColorbar-eps-converted-to.pdf}&
\includegraphics[keepaspectratio=true,height=1.5cm]{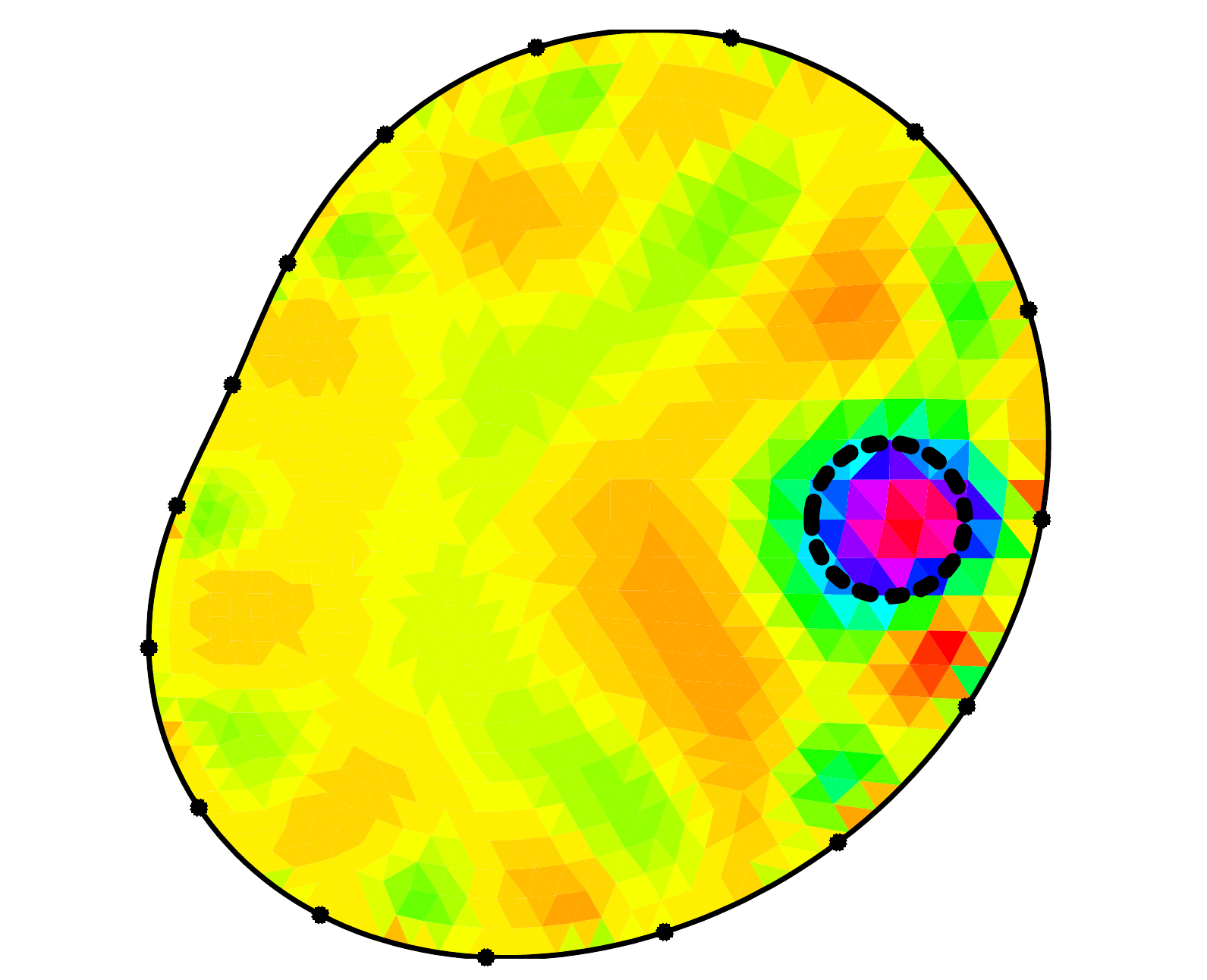}&
\includegraphics[keepaspectratio=true,height=1.5cm]{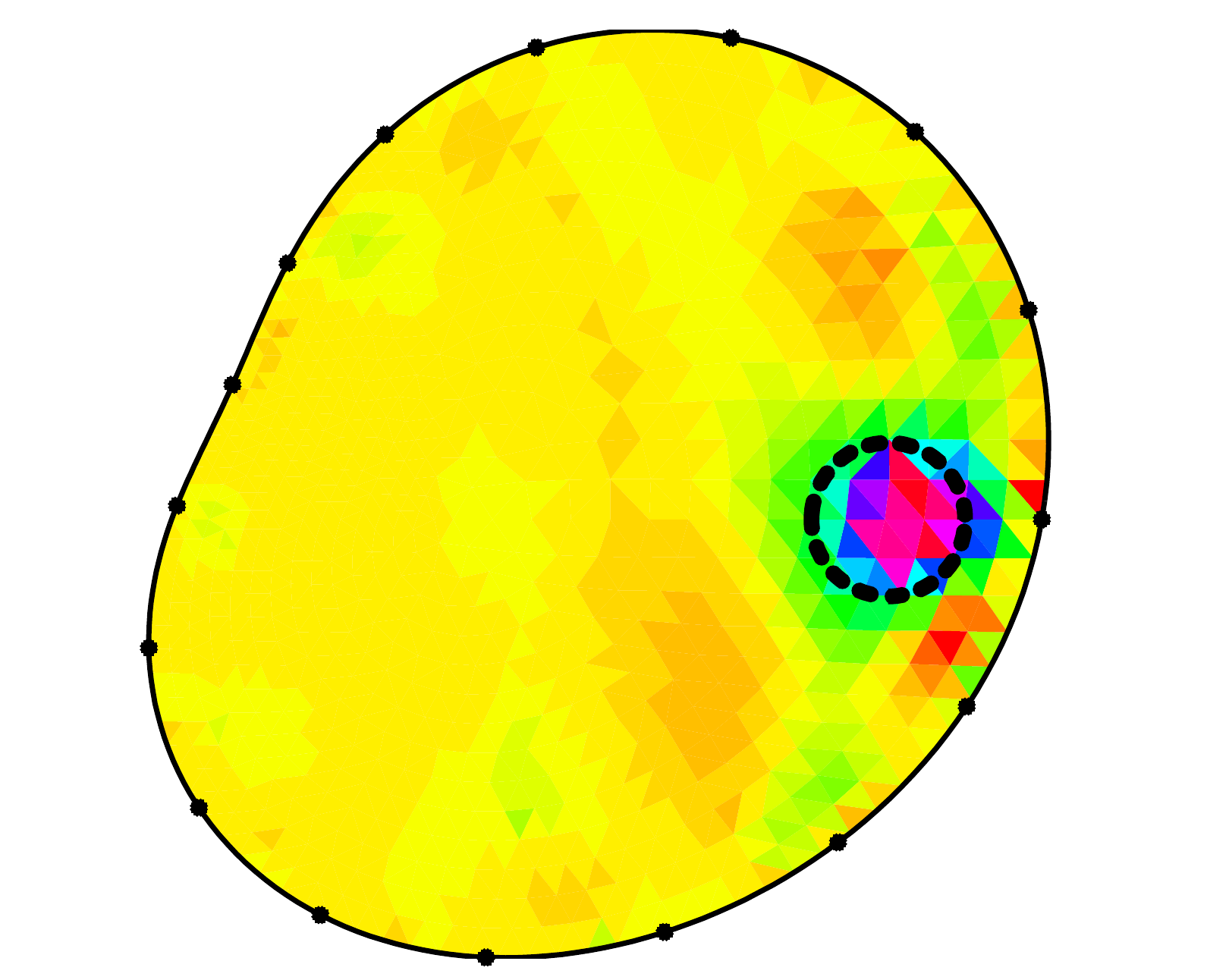}&
\includegraphics[keepaspectratio=true,height=1.5cm]{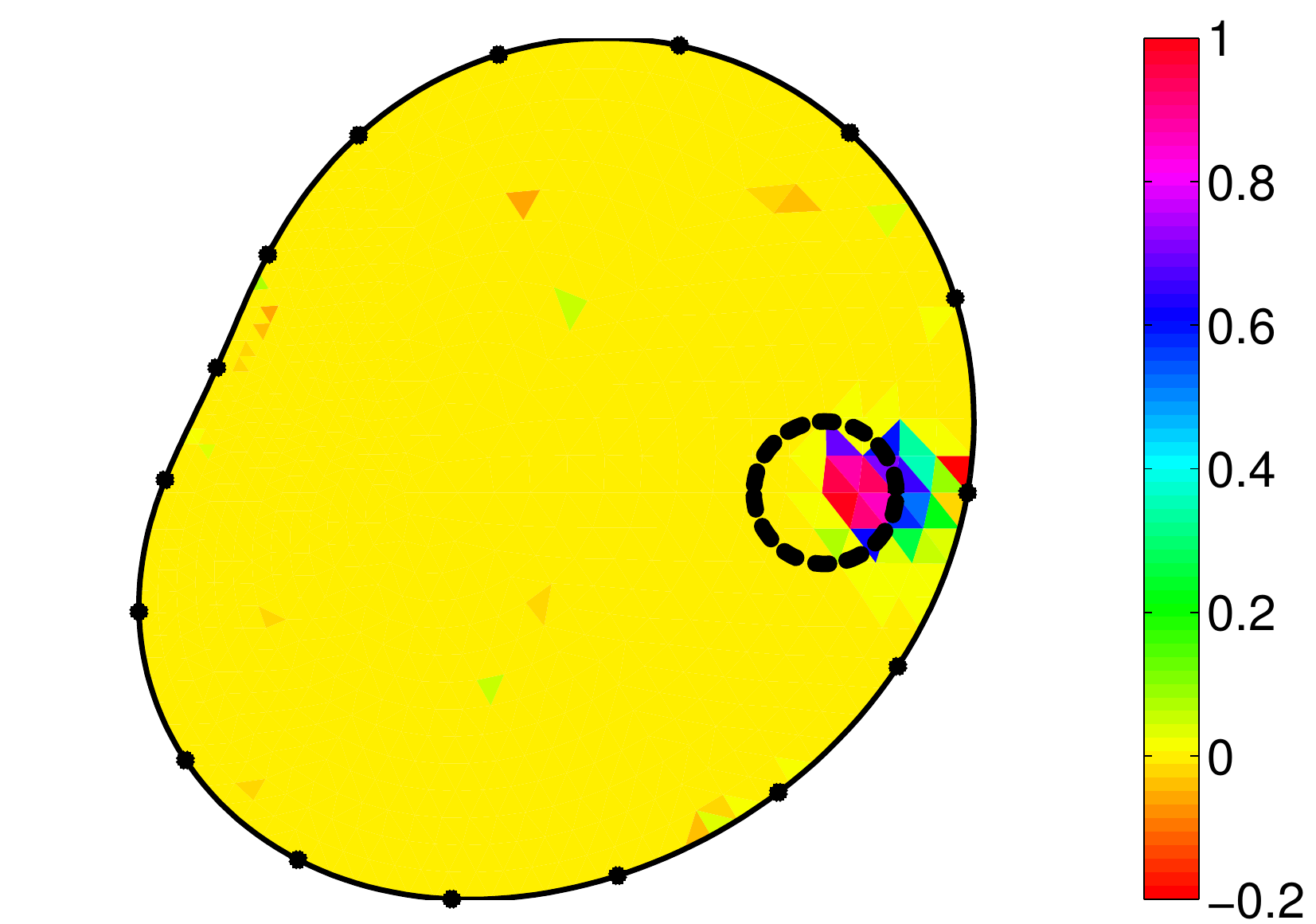}&
\includegraphics[keepaspectratio=true,height=1.45cm]{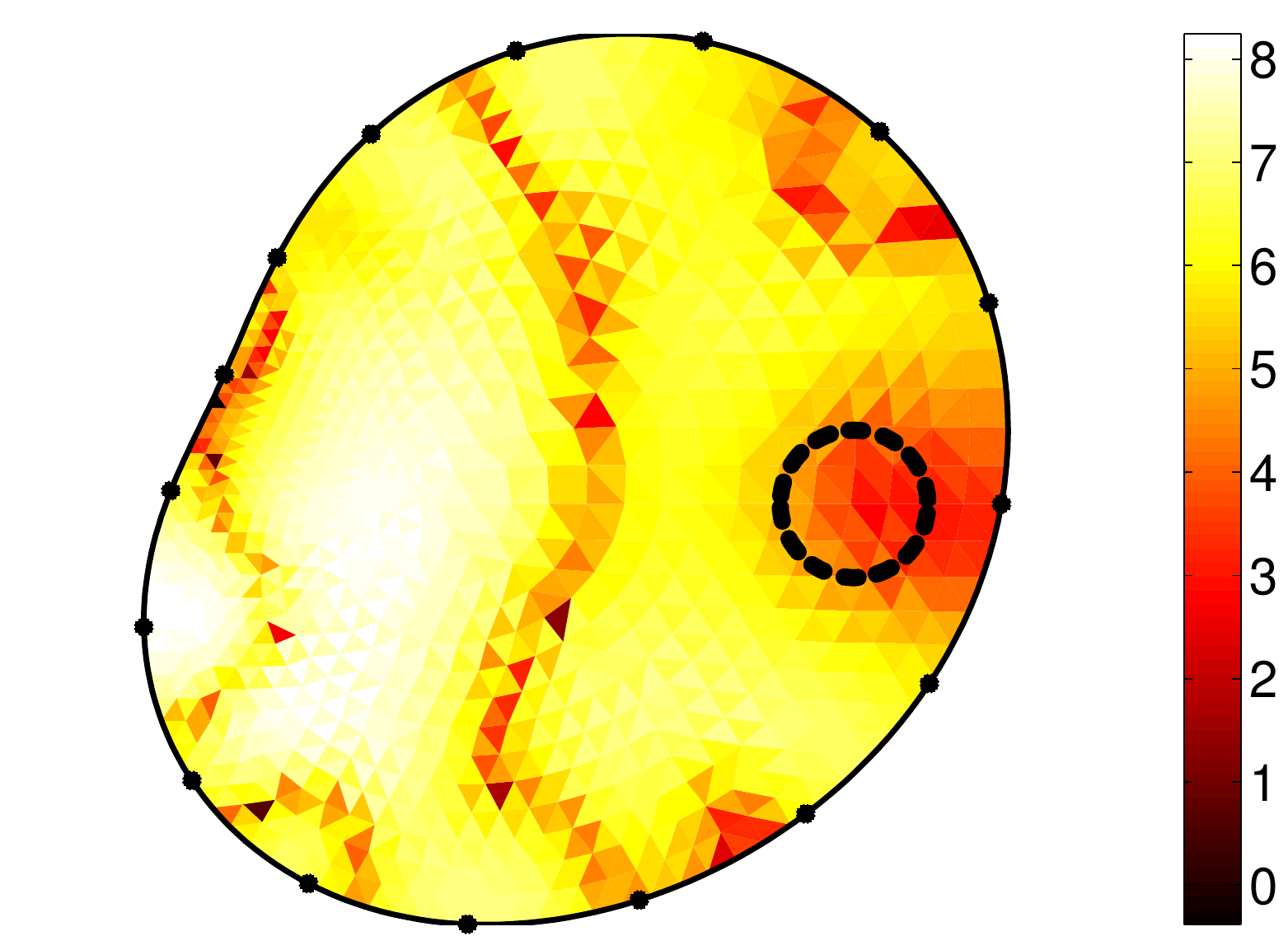}\\
\hline
\raisebox{4ex}{\footnotesize\begin{tabular}{c}
(f)
\end{tabular}}
&
\includegraphics[keepaspectratio=true,height=1.5cm]{Fig/plot_final_20120602/deltasigma_6_woColorbar-eps-converted-to.pdf}&
\includegraphics[keepaspectratio=true,height=1.5cm]{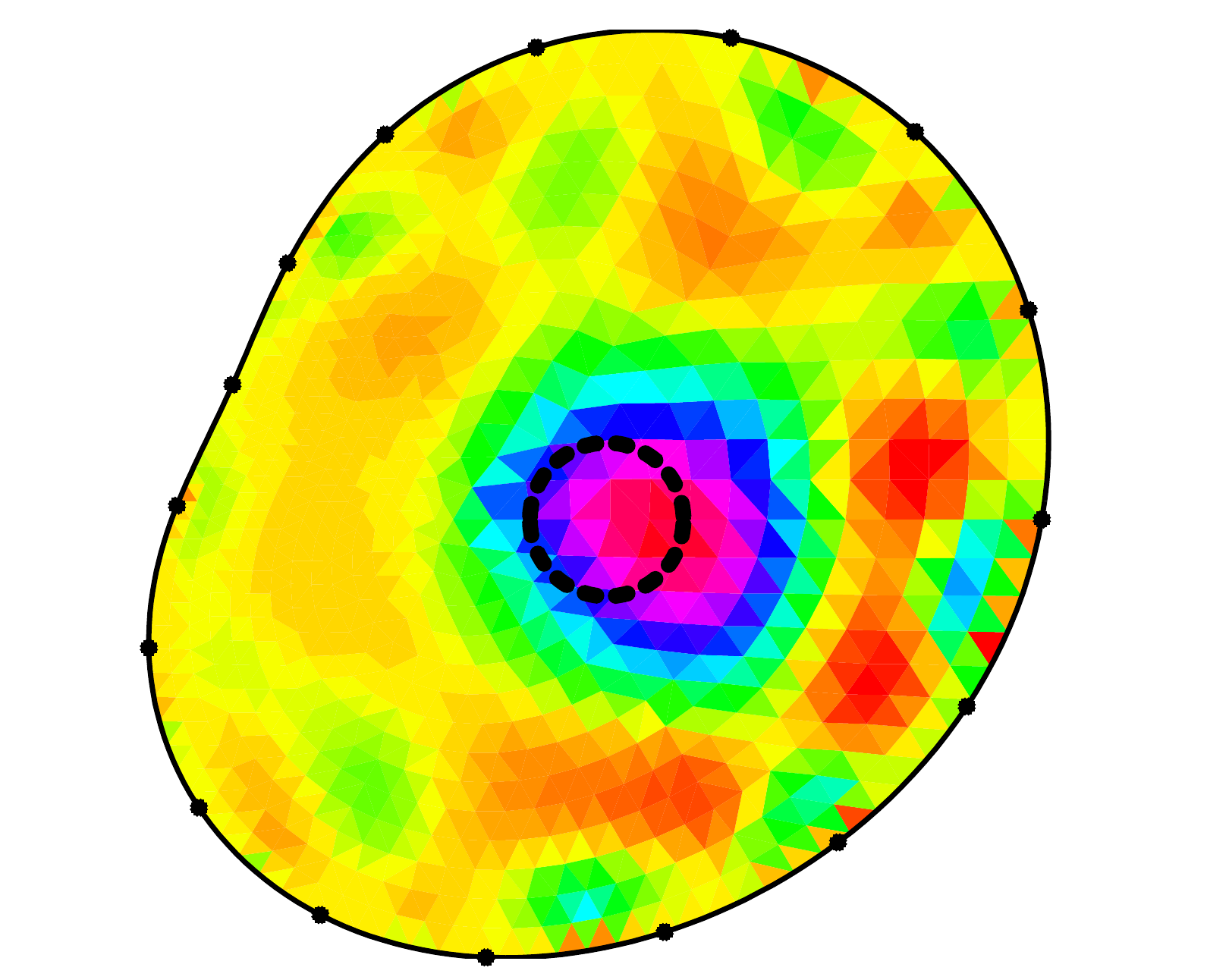}&
\includegraphics[keepaspectratio=true,height=1.5cm]{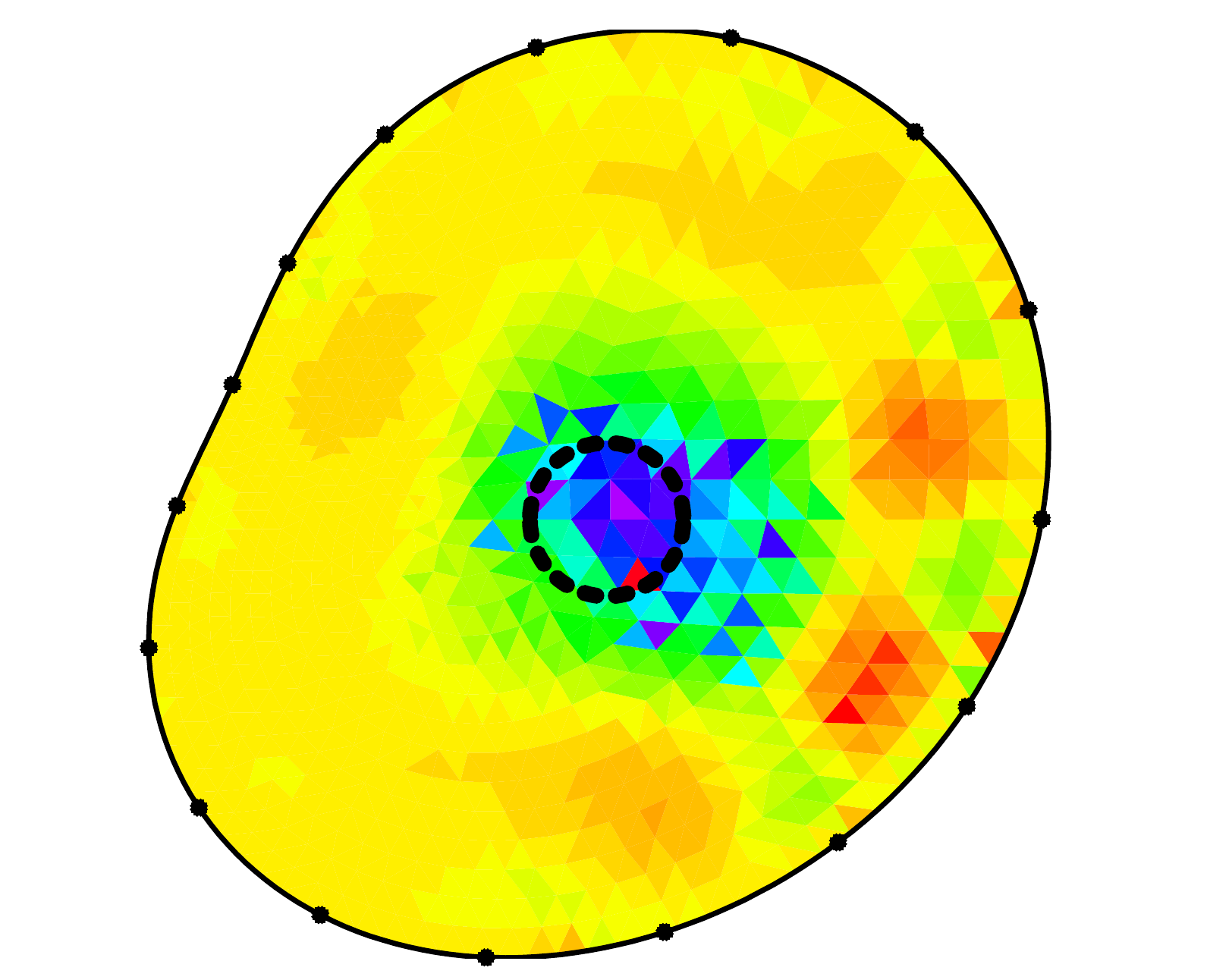}&
\includegraphics[keepaspectratio=true,height=1.5cm]{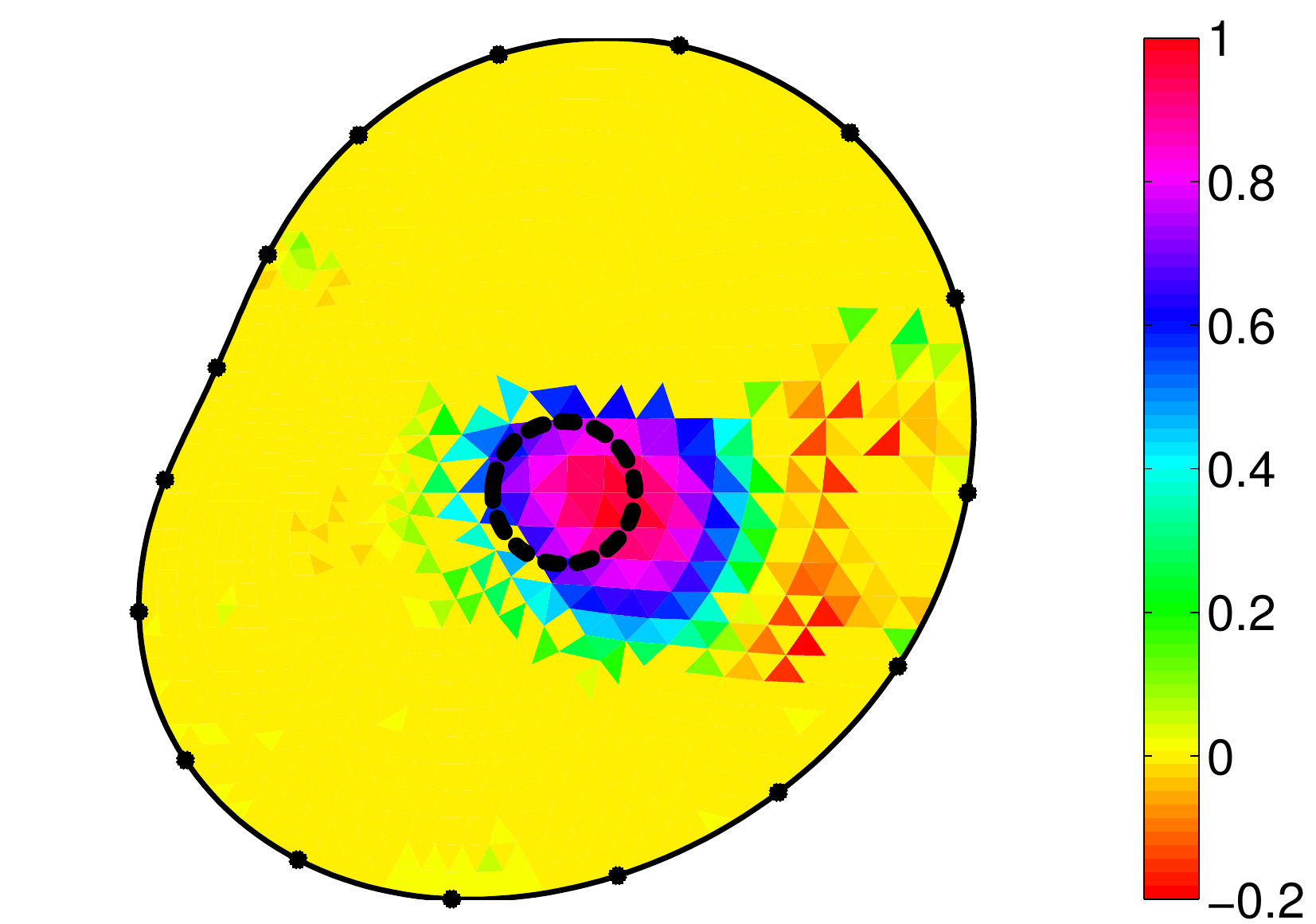}&
\includegraphics[keepaspectratio=true,height=1.5cm]{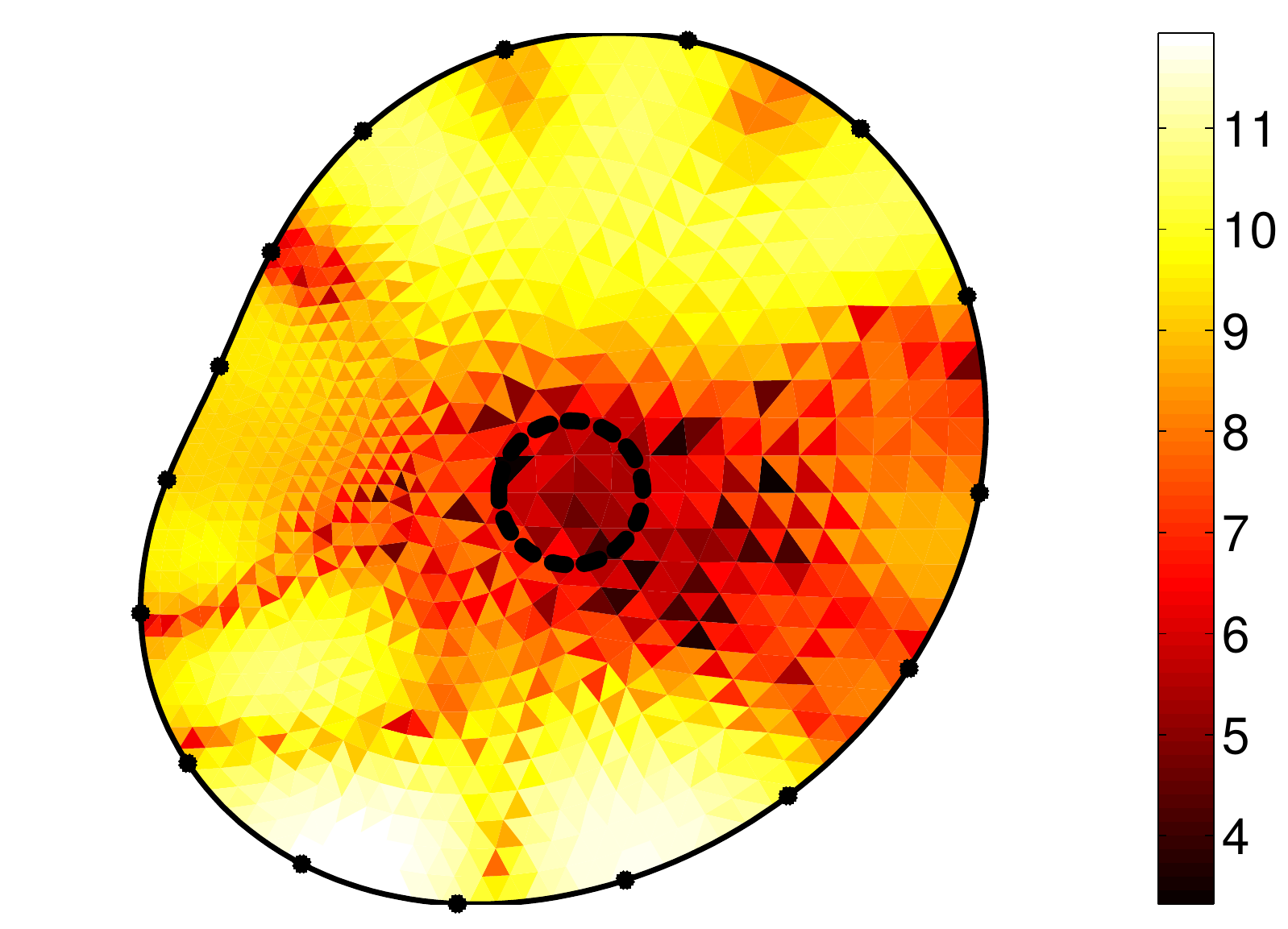}\\
\hline
\raisebox{4ex}{\footnotesize\begin{tabular}{c}
(g)
\end{tabular}}
&
\includegraphics[keepaspectratio=true,height=1.5cm]{Fig/plot_final_20120602/deltasigma_7_woColorbar-eps-converted-to.pdf}&
\includegraphics[keepaspectratio=true,height=1.5cm]{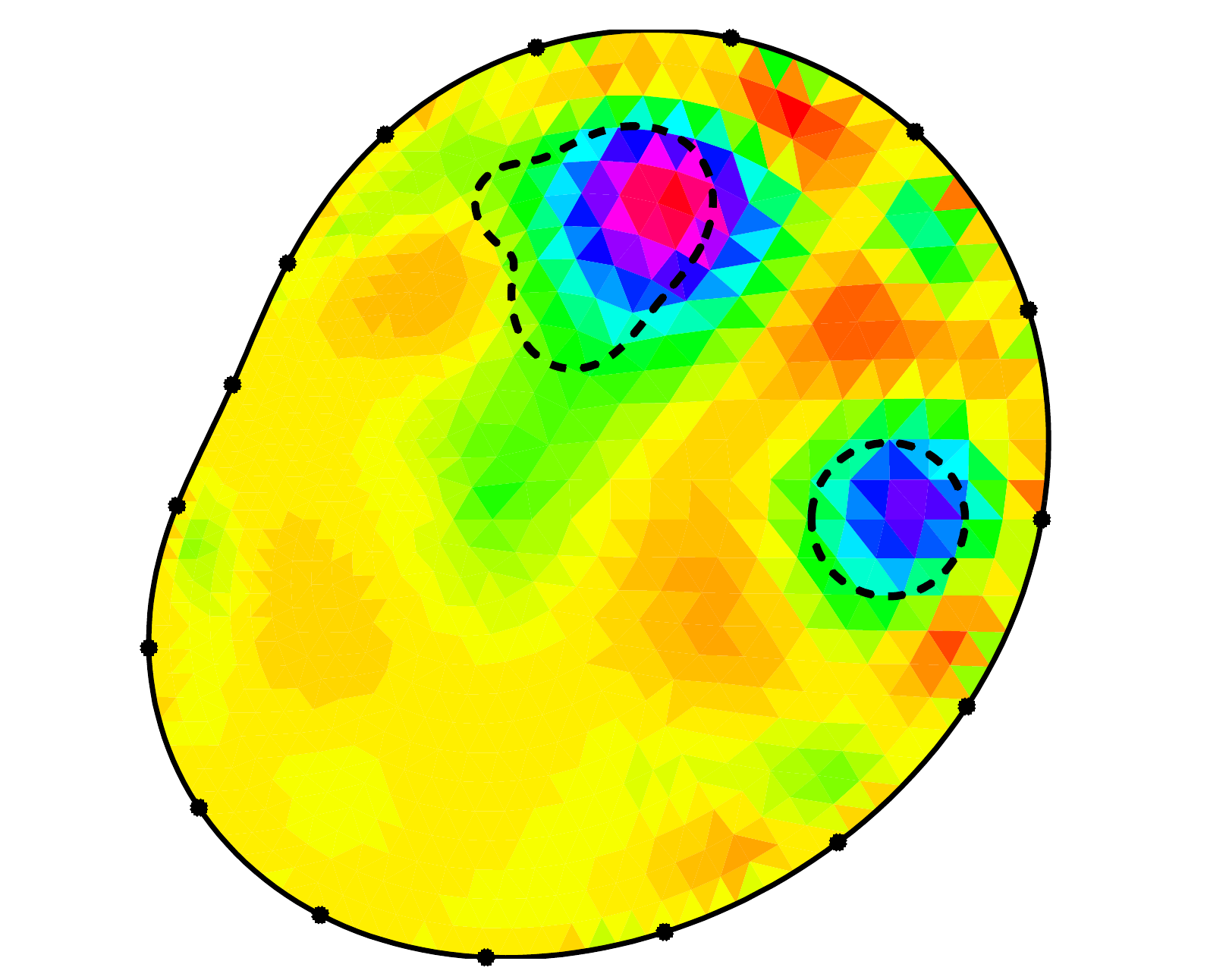}&
\includegraphics[keepaspectratio=true,height=1.5cm]{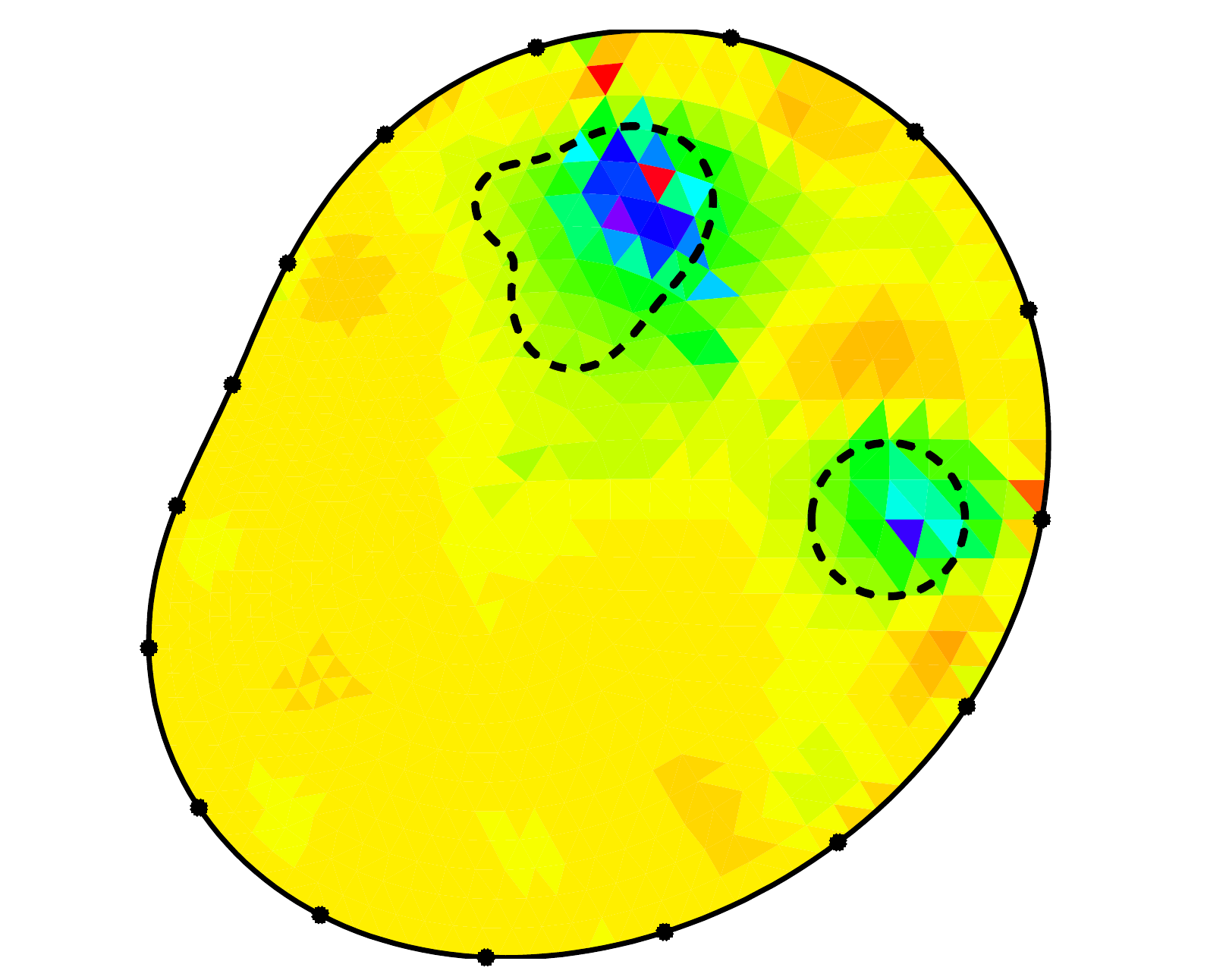}&
\includegraphics[keepaspectratio=true,height=1.5cm]{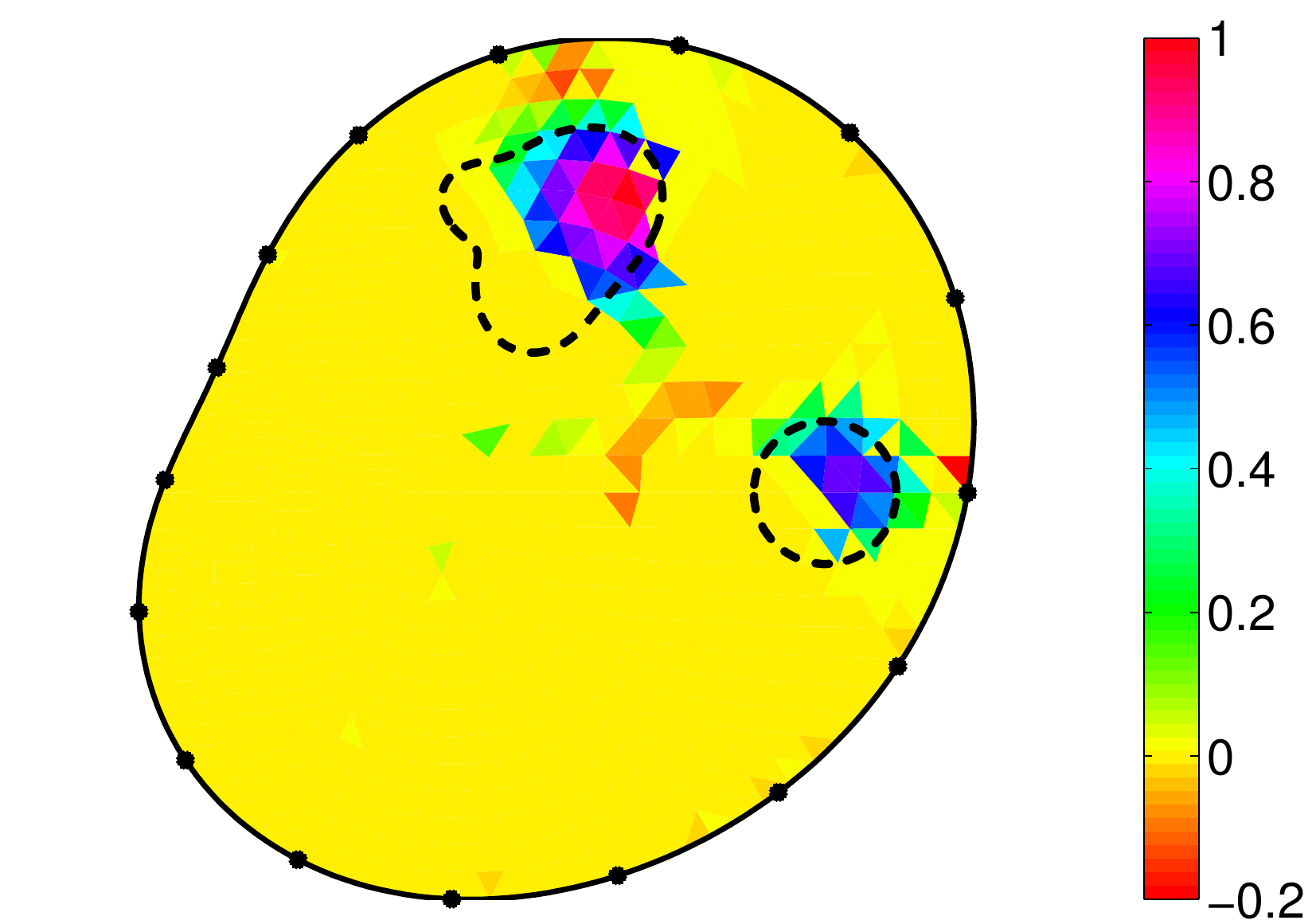}&
\includegraphics[keepaspectratio=true,height=1.5cm]{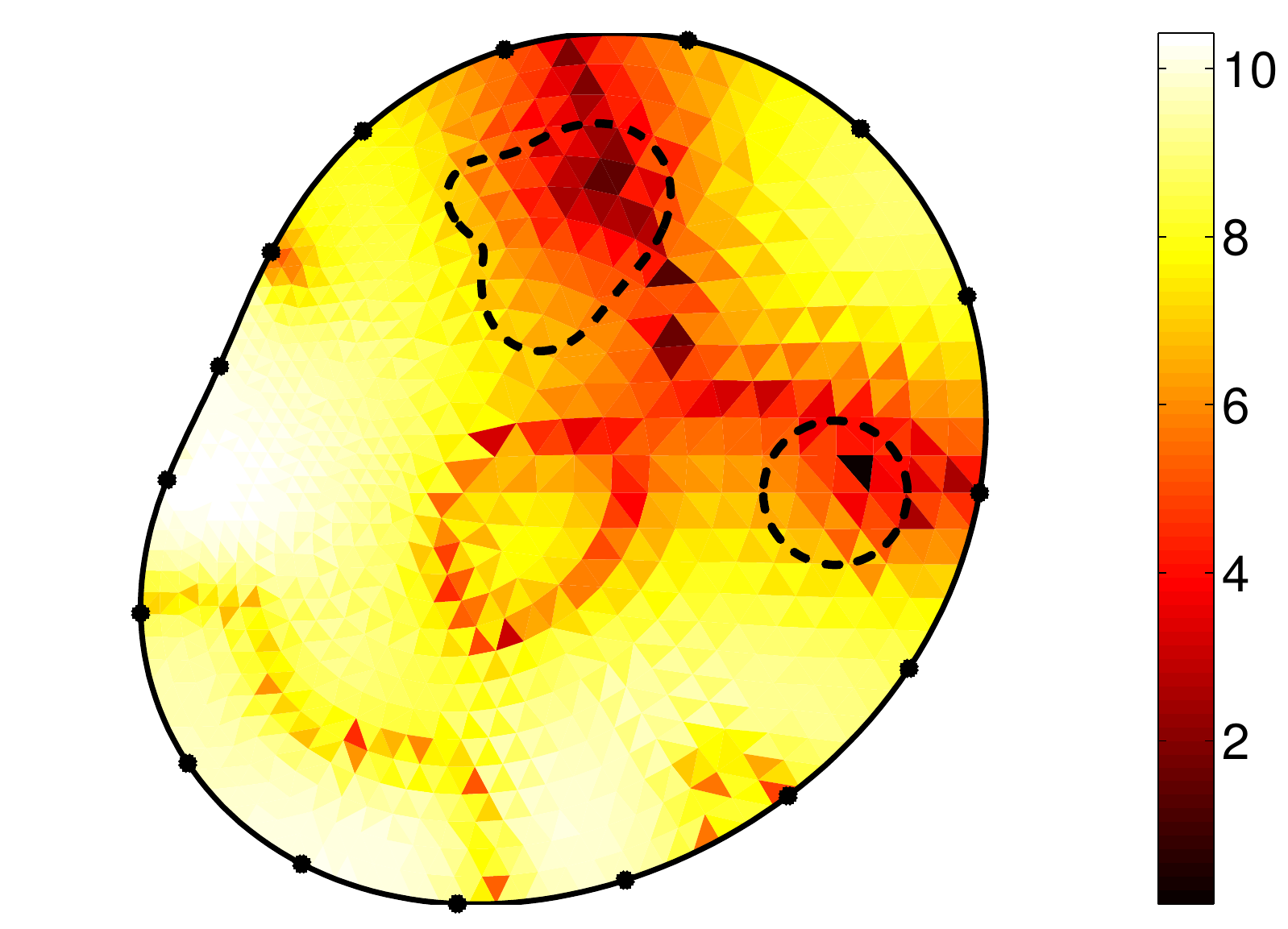}\\
\hline
\raisebox{4ex}{\footnotesize\begin{tabular}{c}
(h)
\end{tabular}}
&
\includegraphics[keepaspectratio=true,height=1.5cm]{Fig/plot_final_20120602/deltasigma_8_woColorbar-eps-converted-to.pdf}&
\includegraphics[keepaspectratio=true,height=1.5cm]{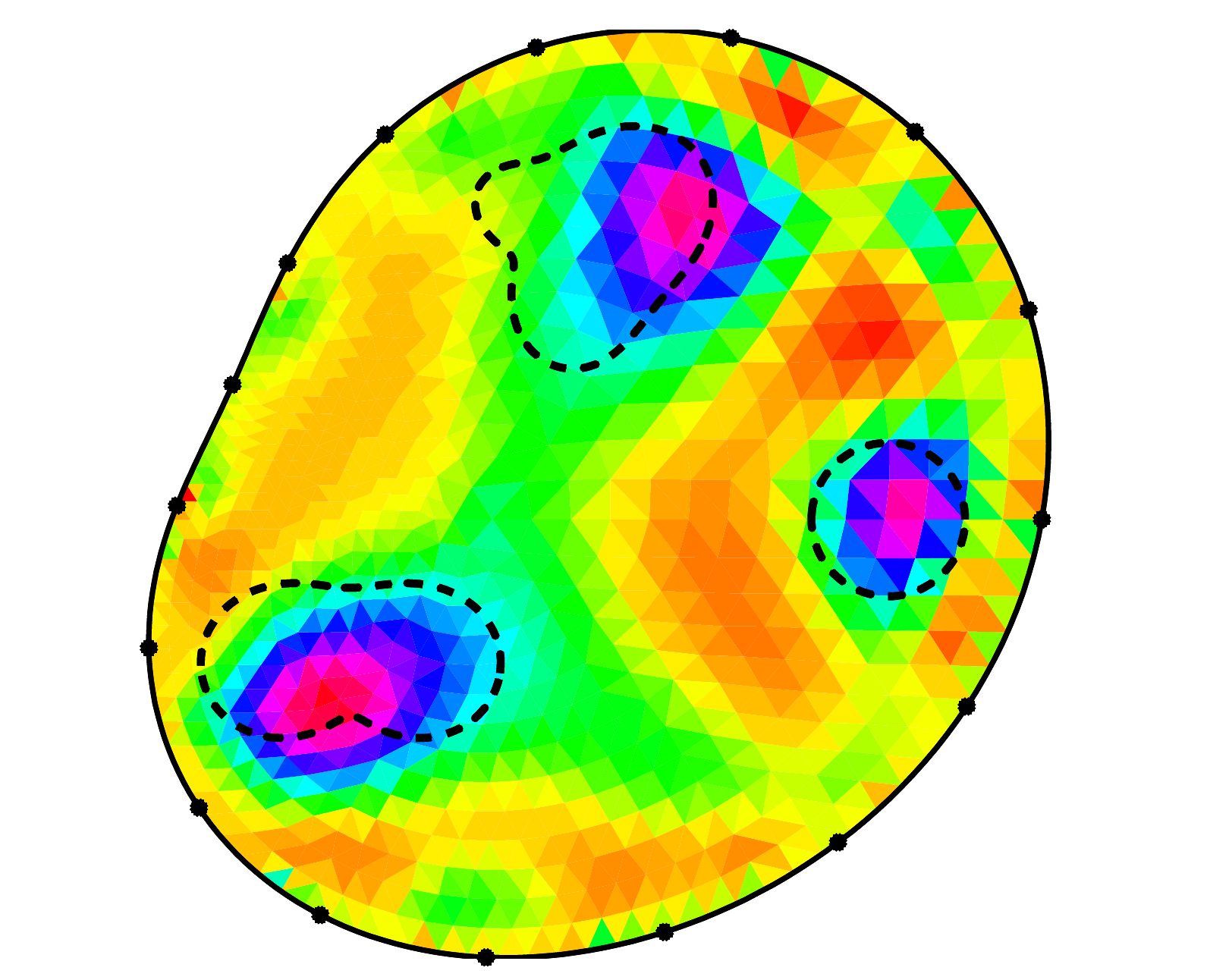}&
\includegraphics[keepaspectratio=true,height=1.5cm]{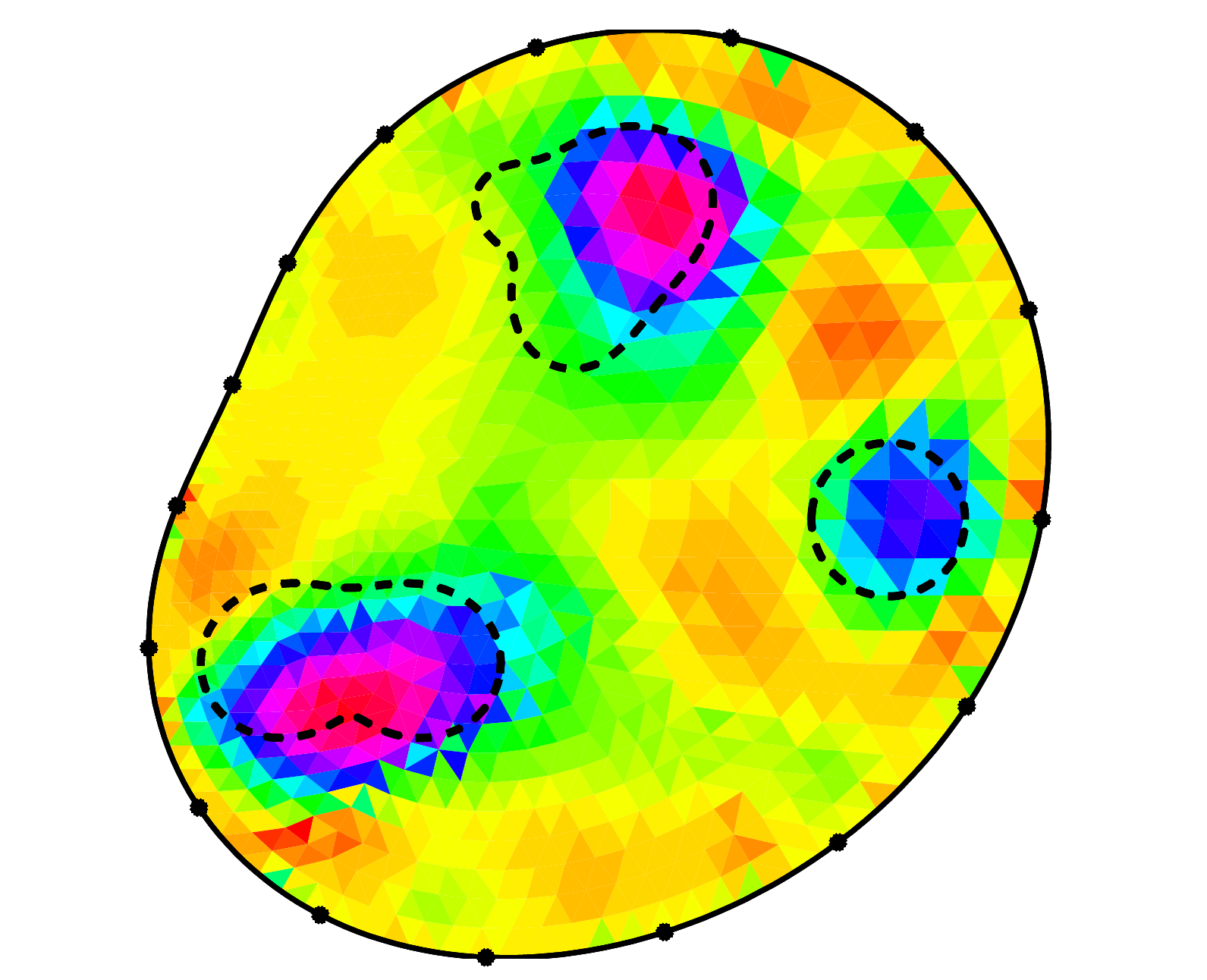}&
\includegraphics[keepaspectratio=true,height=1.5cm]{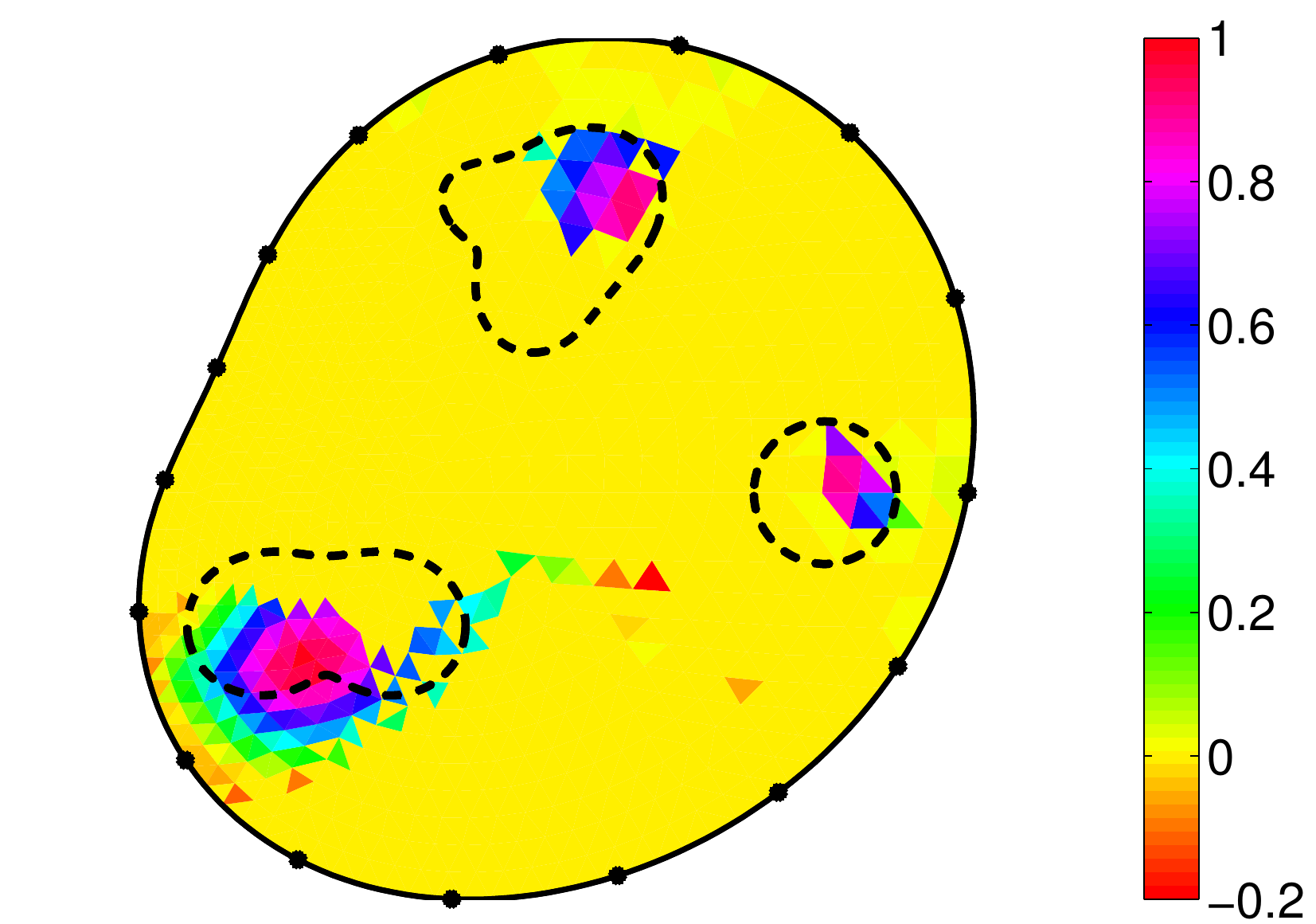}&
\includegraphics[keepaspectratio=true,height=1.5cm]{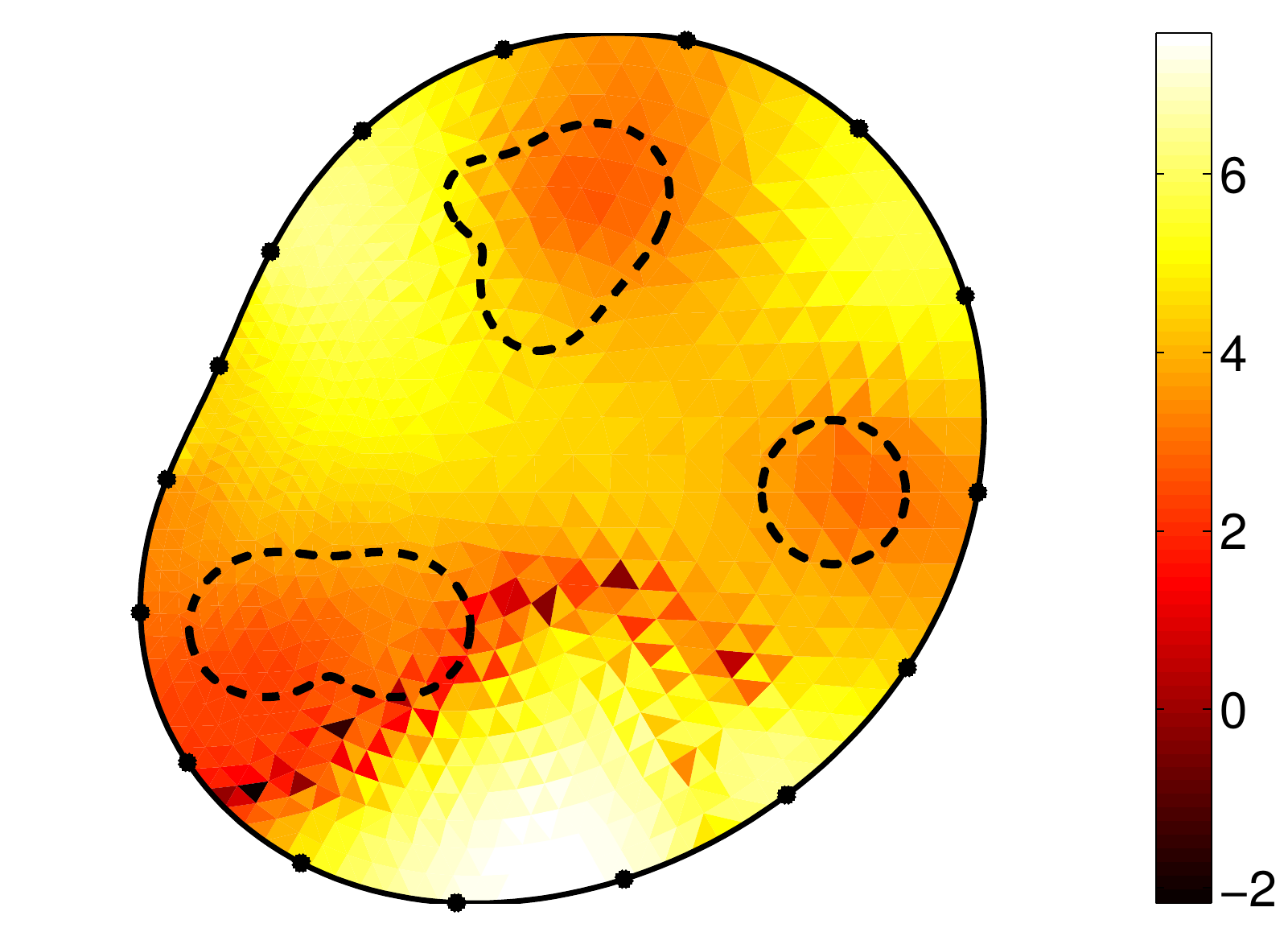}\\
\hline
\end{tabular}
\caption{ Reconstructed difference EIT images in non-circular domain with a data which adds 5$\%$ random noise. $\DS$: true difference image,
$\DS_{S}$: standard linearized method, cf.\ \eref{LM2},
$\DS_{B}$: naive combination of LM and S-FM, cf.\ \eref{regularization},
$\DS_{A}$: proposed combination of LM and S-FM, cf.\ \eref{recon},
$\mathbf{W1}$: S-FM alone, cf.\ \eref{SFMimage}.}
\label{recon_noncircle_5p}
\end{figure}

The forth column in Figure \ref{recon_circle}-\ref{recon_noncircle} shows the reconstructed images using tSVD pseudoinverse $\mathbf{B}^\dag =\hat{\mathbf{V}}_{t_1}\hat{\mathbf{\Lambda}}_{t_1}^{-1} \hat{\mathbf{U}}_{t_1}^*$ where $t_1$ is the number corresponding to $\frac{\hat\lambda_{t_1}}{\hat\lambda_1}\approx 10^{-3}$ and $\lambda_{t_0}\approx\hat{\lambda}_{t_1}$. Due to the change of singular values decay pattern, $t_1$ is close to the number of pixel $n_p$. We will denote the corresponding reconstructed image  by $\DS_{B}$ :
$$\DS_{B}=\sum_{t=1}^{t_1}\frac{1}{\hat{\lambda}_t}\langle \mathbf{b},\hat{\mathbf{u}}_t\rangle\hat\v_t.$$

The fifth column in Figure \ref{recon_circle}-\ref{recon_noncircle} shows the reconstructed images using tSVD pseudoinverse $\mathbf{A}^\dag =\tilde{\mathbf{V}}_{t_2}\tilde{\mathbf{\Lambda}}_{t_2}^{-1} \tilde{\mathbf{U}}_{t_2}^*$ :
$$\DS_{A}=\sum_{t=1}^{t_2}\frac{1}{\tilde\lambda_t}\left\langle \left[\begin{array}{c}
\DV\\
\alpha\mathbf{W}^{-1}\S^{\dag}\DV
\end{array}\right],\tilde{\mathbf{u}}_t\right\rangle\tilde\v_t
.$$
The parameter  $\alpha$ is an weighting factor to balance  between the LM and S-FM. In the special case of $\alpha=0$,  the image of $\DS_A$ is similar to the image of $\DS_S$.  Here, we choose $\alpha=1$.

The choice of $t_2$ is more involved. Since we use a TSVD, our reconstruction lies inside a $t_2$-dimensional space spanned by right singular vectors. 
Hence, it seems natural that an optimal choice of $t_2$ should be related to the number of pixels that belong to the inclusion.
To estimate the number of pixels inside the inclusion we count the pixels in the S-FM reconstruction. 
More precisely, we take the number of pixels for which the S-FM-indicator value belongs to
the upper third of the indicator values, and then, to be on the safe side, we double this number.
Hence, the truncation number $t_2$ is chosen  by
$$
t_2=\aleph\Big(\{ n ~|~w_n^{-1} ~\leq~ w_{\tiny{\min}}^{-1} -\frac{w_{\tiny{\min}}^{-1} - w_{\tiny{\max}}^{-1}}{3}\}\Big)\times 2
$$
where $\aleph(A)$ is a number of elements of a set $A$, and $w_{\tiny{\min}}=\displaystyle\min_{1\leq n\leq n_p} w_n$ and $w_{\tiny{\max}}=\displaystyle\max_{1\leq n\leq n_p} w_n$. In the Figure \ref{recon_circle}-\ref{recon_noncircle}, the set $
\{ n ~|~w_n^{-1} ~\leq~ w_{\tiny{\min}}^{-1} -\frac{w_{\tiny{\min}}^{-1} - w_{\tiny{\max}}^{-1}}{3}\}
$  occupies a region of reddish part of image ${\bf W}{\bf 1}$ (not region of yellowish part).

Figure \ref{recon_circle_1p}-\ref{recon_noncircle_5p} are structured in the same way as Figure \ref{recon_circle}-\ref{recon_noncircle} and show the effect of noise
added to the difference data  (\ref{eq:uuV1}).

\section{Conclusion}

We developed a modified discrete version of the factorization method
(FM), herein called \textit{sensitivity matrix based FM}(S-FM) that can
be used as a regularization term for the standard linearized
reconstruction method (LM). The S-FM provides a pixel-wise index
indicating the possibility of the presence of anomaly at each pixel.
Using this information for regularizing the LM
resulted in an a better localization of anomalies and alleviated Gibbs
ringing artifacts.

We should mention about the choice of $f={\bf W}^{-1}\S^{\dag}\DV$  in \eref{recon}.   Indeed, the choice is motivated by heuristical arguments and numerical evidence. The role of the augmented matrix  ${\bf W}^{-1}$ is that the right singular vectors associated with high singular values are used to reconstruct images of anomaly region ( $\delta\sigma\neq 0$).  The choice of $f$  representing the unknown ${\bf W}^{-1}\DS$ is somehow limited. 
Future studies will be concerned with rigorously justified ways of
extending the linearized method with shape reconstruction properties.

\end{document}